\documentclass[12pt]{article}
\usepackage[textwidth=16cm, textheight=21cm]{geometry}
\usepackage{times} % font
\usepackage[compress]{cite}
\usepackage{scalefnt} %\scalefont{factor} scales the size of the current font and the current value of \baselineskip
\usepackage{lipsum}
\usepackage{newtxmath} 
\usepackage{amsmath}
\usepackage{bm}
\usepackage{amsbsy}
\usepackage{subdepth} % prevents subscript from being shifted downward in presence of superscript
\usepackage{color}
\usepackage{ifpdf}
\usepackage{algorithm,algpseudocode}
\usepackage{graphicx}
\usepackage{textcomp}
\usepackage[table]{xcolor}
\usepackage{threeparttable}
\usepackage{tabularx}
\usepackage{makecell}
\usepackage{longtable}
\usepackage{colortbl,ifthen}
\usepackage{arydshln}
\usepackage{multirow}
\usepackage{bigdelim}
\usepackage{setspace}
\usepackage{mathtools,cuted}
\usepackage{fontawesome}
\usepackage{enumitem}
\setlist{nolistsep}
\usepackage{dashbox}
\usepackage{ifthen}
\usepackage{tikz}
\usetikzlibrary{shapes,arrows,automata,positioning,decorations.text,decorations.markings,bending,snakes,calc,mindmap,external}
\usepackage{bbding} %check and cross

\usepackage{marvosym} %checked and crossed boxes
\usepackage{fp}

\usepackage[amsmath]{empheq} % emphasize equations with box
\usepackage{adjustbox} % highlight texts without changing line width

\usepackage{calc}

\usepackage{floatpag}
\newcommand{\eval}[1]{\the\numexpr#1\relax}

\usepackage[defaultcolor=red]{changes} % track changes

\newcommand{\stkout}[1]{\ifmmode\text{\sout{\ensuremath{#1}}}\else\sout{#1}\fi}
\setdeletedmarkup{\stkout{#1}}

\usepackage{pdflscape}
\usepackage{afterpage}

\usepackage{tcolorbox}
\newtcolorbox{equationbox}{
	colback=white,
	colframe=black,
	boxrule=0.5pt,
	arc=0mm,
	boxsep=2pt,
	left=2pt,right=2pt,top=2pt,bottom=2pt
}

\allowdisplaybreaks[1]

\ifpdf
  \usepackage{hyperref}
\else
  \usepackage[hypertex]{hyperref}
\fi
\definecolor{dgreen}{rgb}{0,0.7,0}
\definecolor{ddgreen}{rgb}{0,0.5,0}
\definecolor{dblue}{rgb}{0,0,0.7}
\definecolor{dred}{rgb}{0.7,0,0}
\hypersetup{bookmarksnumbered=true,
            plainpages=false,
            colorlinks=true,
            linkcolor=black,%dblue,
            citecolor=black,%dgreen,
            filecolor=black,%dred,
            urlcolor=black}%blue}

\renewcommand{\eqref}[1]{(\ref{#1})}

\newcommand{\uls}[1]{\underline{\smash{#1}}}

\newcommand{\hide}[1]{} %\newcommand{\hide}[1]{#1}

\tikzstyle{startstop} = [rectangle, rounded corners, minimum width=3cm, minimum height=1cm,text centered, draw=black, thick, fill=red!30]
\tikzstyle{io} = [trapezium, trapezium left angle=70, trapezium right angle=110, minimum width=3cm, minimum height=1cm, text centered, draw=black, thick, fill=blue!30]
\tikzstyle{process} = [rectangle, minimum width=3cm, minimum height=1cm, text centered, draw=black, thick, fill=orange!30]
\tikzstyle{decision} = [diamond, aspect=3, minimum width=3cm, minimum height=1cm, text centered, draw=black, thick, fill=green!30]
\tikzstyle{arrow} = [thick,->,>=stealth]
\tikzstyle{noarrow} = [thick,-]
\tikzstyle{doublearrow} = [thick, double distance=1pt, ->, >=stealth]
\tikzstyle{doublenoarrow} = [thick, double distance=1pt]
\tikzstyle{outerwhite} = [thin, white, line width=5pt, shorten >= 4.5pt]
\tikzstyle{innerwhite} = [thin, white, line width=1pt]
\tikzset{state/.style={
           rectangle,
           rounded corners,
           draw=black, very thick,
           minimum height=2em,
           inner sep=2pt,
           },}
\tikzstyle{start} = [circle, inner sep=0.06cm, draw=black, fill=black]   
\tikzset{concept/.append style={fill={none}}}   

\tikzstyle{gain} = [draw, thick, fill=white, isosceles triangle, isosceles triangle apex angle=30, minimum height=0.5em, minimum width=0.5em, align=center]
\tikzstyle{bigblock} = [draw, thick, fill=white, rectangle, 
minimum height=4em, minimum width=6em, align=center]
\tikzstyle{smallblock} = [draw, thick, fill=white, rectangle, 
minimum height=1em, minimum width=1em, align=center]
\tikzstyle{medianblock} = [draw, thick, fill=white, rectangle, 
minimum height=3em, minimum width=3em, align=center]
\tikzstyle{bigfcnblock} = [draw, thick, fill=white, rectangle, 
minimum height=4em, minimum width=6em, double, double distance=0.5mm, align=center]
\tikzstyle{smallfcnblock} = [draw, thick, fill=white, rectangle, 
minimum height=1em, minimum width=1em, double, double distance=0.5mm, align=center]
\tikzstyle{medianfcnblock} = [draw, thick, fill=white, rectangle, 
minimum height=3em, minimum width=3em, double, double distance=0.5mm, align=center]
\tikzstyle{invisiblenode} = [draw=none, minimum size=0pt, align=center]
\tikzstyle{sum} = [draw, thick, fill=white, circle, node distance=1cm]
\makeatletter
\pgfdeclareshape{record}{
	\inheritsavedanchors[from={rectangle}]
	\inheritbackgroundpath[from={rectangle}]
	\inheritanchorborder[from={rectangle}]
	\foreach \x in {center,north east,north west,north,south,south east,south west}{
		\inheritanchor[from={rectangle}]{\x}
	}
	\foregroundpath{
		\pgfpointdiff{\northeast}{\southwest}
		\pgf@xa=\pgf@x \pgf@ya=\pgf@y
		\northeast
		\pgfpathmoveto{\pgfpoint{0}{0.33\pgf@ya}}
		\pgfpathlineto{\pgfpoint{0}{-0.33\pgf@ya}}
		\pgfpathmoveto{\pgfpoint{0.33\pgf@xa}{0}}
		\pgfpathlineto{\pgfpoint{-0.33\pgf@xa}{0}}
		\pgfpathmoveto{\pgfpointadd{\southwest}{\pgfpoint{-0.33\pgf@xa}{-0.6\pgf@ya}}}
		\pgfpathlineto{\pgfpointadd{\southwest}{\pgfpoint{-0.5\pgf@xa}{-0.6\pgf@ya}}}
		\pgfpathlineto{\pgfpointadd{\northeast}{\pgfpoint{-0.5\pgf@xa}{-0.6\pgf@ya}}}
		\pgfpathlineto{\pgfpointadd{\northeast}{\pgfpoint{-0.33\pgf@xa}{-0.6\pgf@ya}}}
	}
}
\makeatother
\tikzstyle{saturation} = [draw, record, fill=white, minimum height=3em, minimum width=3em]

\tikzstyle{spring}=[thick,decorate,decoration={zigzag,amplitude=0.2cm,pre length=0.3cm,post length=0.3cm,segment length=6}]
\tikzstyle{damper}=[thick,decoration={markings,  
  mark connection node=dmp,
  mark=at position 0.5 with 
  {
    \node (dmp) [thick,inner sep=0pt,transform shape,rotate=-90,minimum width=15pt,minimum height=3pt,draw=none] {};
    \draw [thick] ($(dmp.north east)+(2pt,0)$) -- (dmp.south east) -- (dmp.south west) -- ($(dmp.north west)+(2pt,0)$);
    \draw [thick] ($(dmp.north)+(0,-5pt)$) -- ($(dmp.north)+(0,5pt)$);
  }
}, decorate]
\tikzstyle{ground}=[fill,pattern=north east lines,draw=none,minimum width=0.75cm,minimum height=0.3cm]

\tikzset{
	pics/.cd,
	vector out/.style={
		code={
			\draw[#1] (0,0)  circle (1) (45:1) -- (225:1) (135:1) -- (315:1);
		}%end code   
	}%end style
}%end tikzset
\tikzset{
	pics/.cd,
	vector in/.style={
		code={
			\draw[#1] (0,0)  circle (1);
			\fill[#1] (0,0)  circle (.1);
		}%end code   
	}%end style
}%end tikzset

\tikzset{->-/.style={decoration={
			markings,
			mark=at position .5 with {\arrow{>}}},postaction={decorate}}}

\pgfdeclarelayer{bg}    % declare background layer
\pgfsetlayers{bg,main}  % set the order of the layers (main is the standard layer)
\tikzstyle{circlenode} = [draw, circle, minimum width=1ex]
\tikzstyle{squarenode} = [draw, regular polygon, regular polygon sides=4, minimum width=1ex]

\tikzstyle{hexagonnode} = [draw, thick, shape=regular polygon, regular polygon sides=6, minimum width=5pt]

\newcounter{line}
\newcommand\alternaterowcolor{%
	\addtocounter{line}{1}%
	\ifthenelse{\isodd{\value{line}}}{\\\rowcolor[gray]{0.95}}{\\}}

\usepackage[utf8]{inputenc}
\usepackage{subcaption}
\usepackage[labelformat=simple]{subcaption}  
\usepackage[labelfont=bf]{caption}
\usepackage{xr}
\makeatletter
\newcommand*{\addFileDependency}[1]{% argument=file name and extension
	\typeout{(#1)}
	\@addtofilelist{#1}
	\IfFileExists{#1}{}{\typeout{No file #1.}}
}
\makeatother
\newcommand*{\myexternaldocument}[1]{%
	\externaldocument[SI-]{#1}%
	\addFileDependency{#1.tex}%
	\addFileDependency{#1.aux}%
}
\myexternaldocument{supplementary}

\def\nullv{ {\bm 0} }

\def\bv{ {\bm b} }
\def\cv{ {\bm c} }

\def\ev{ {\bm e} }
\def\fv{ {\bm f} }
\def\gv{ {\bm g} }
\def\hv{ {\bm h} }

\def\pv{ {\bm p} }

\def\uv{ {\bm u} }

\def\xv{ {\bm x} }
\def\yv{ {\bm y} }

\def\vargv{ {\bm{\mathscr{g}}} }
\def\varhv{ {\bm{\mathscr{h}}} }

\def\varyv{ {\bm{\mathscr{y}}} }
\def\Bm{ {\bm B} }

\def\Gm{ {\bm G} }

\def\Km{ {\bm K} }

\def\Lambdam{ {\bm \varLambda} }

\def\zetav{ \bm{\zeta} }
\def\etav{ \bm{\eta} }

\def\lambdav{ \bm{\lambda} }
\def\muv{ \bm{\mu} }

\def\phiv{ \bm{\phi} }

\def\varphiv{ \bm{\varphi} }
\def\chiv{ \bm{\chi} }

\def\OTC{\text{DCC}} % Dual Cost-Constraint
\def\SQP{\text{SQP}}

\def\CPUsubscript{{\scriptscriptstyle\textup{C}}}
\def\OTCsubscript{{\scriptscriptstyle\textup{\OTC}}}
\def\SQPsubscript{{\scriptscriptstyle\textup{\SQP}}}

\def\FLOPsubscript{{\scriptscriptstyle\textup{FLOP}}}
\def\FcnEvalsubscript{{\scriptscriptstyle\textup{FE}}}

\def\realset{ \mathbb{R} }

\DeclareMathOperator{\principal}{pc}

\def\dimx{ d_{\mathrm{x}} }
\def\dimy{ d_{\mathrm{y}} }

\def\final{ \textup{f} }

\newcommand{\desired}[1]{ {#1}_{\mathrm{ref}} }

\def\dt{ \delta t }

\makeatletter
\DeclareRobustCommand{\optimal}{%
	\mathbin{\mathpalette\o@plus@times\relax}%
}
\newcommand{\o@plus@times}[2]{%
	\ooalign{$\m@th#1\oplus$\cr$\m@th#1\otimes$\cr}%
}
\makeatother

\def\errv{ \ev }

\def\colorsota{gray!30}
\def\coloripi{cyan!30}
\def\colorfbc{yellow!30}
\def\colorotc{green!30}
\def\colorfw{magenta!30}

\begin{document}
\baselineskip24pt % Double-space the manuscript.

\title{%
%Real-time optimal control without iteration
Optimal Control without Optimization
%Breaking the Real-Time Barrier in High-Dimensional Optimal Control
}

\author
{Tingli Hu\textsuperscript{1} and Sami Haddadin\textsuperscript{1,$\ast$}
	\\
	\normalsize{\textsuperscript{1}Mohamed bin Zayed University of Artificial Intelligence, Abu Dhabi, United Arab Emirates}
	\\
	\normalsize{\textsuperscript{$\ast$}To whom correspondence should be addressed; E-mail: \url{sami.haddadin@mbzuai.ac.ae}}
}
\date{}

\maketitle

\begin{quote}\bf
\textbf{Abstract:}
The real-time barrier in optimal control of nonlinear dynamical systems remains a longstanding limitation across science, engineering, and economics. Existing approaches rely on iterative numerical optimization and therefore cannot compute optimal control actions directly within physical time for complex, high-dimensional systems.
Here we introduce the Dual Cost-Constraint projection (DCC), a closed-form dynamical representation that enables real-time solutions for a broad class of pseudoconvex optimal control and optimization problems and demonstrates that real-time optimal control can admit a direct closed-form representation. Unlike classical formulations with Lagrange multipliers or adjoint state variables, the proposed DCC embeds constraints directly within the system dynamics. The derivation further reveals a structural equivalence between interior-point optimization iteration and nonlinear feedback control, linking constrained optimization with classical stability theory.
Through theoretical analysis and real-world-relevant numerical studies---including autonomous system control, biomechanics monitoring, and economic decision processes---we show that DCC achieves accurate optimal behavior even for highly nonlinear and high-dimensional systems. The numerical benchmark experiments suggest that DCC controls a $1000$-dimensional system at $1$~kHz using 77~\% of a modern CPU, whereas sequential quadratic programming (SQP), {a widely used} state-of-the-art solver for this class of problems, requires more than 80 processors to achieve comparable performance.
Beyond its computational advantages, DCC provides a system-level, causally deterministic interpretation of constrained optimization, revealing that optimal behavior in a broad class of optimal control problems can emerge as the stable evolution of the system itself. These results lay the foundation for extending this perspective to more general problem classes.

% Real-time optimal control of high-dimensional nonlinear dynamical systems remains a central challenge across science, engineering, and economics.
% State-of-the-art methods cannot compute optimal control actions faster than physical time. The delays can cause disasters such as aircraft crashes, irreversible injuries, and economic losses.
% %
% We introduce the \addedtwo{\emph{Dual Cost-Constraint projection} (DCC), a dynamical systems theory} that breaks through the real-time barrier for a broad spectrum of systems and problems, including (R1) autonomous system control, (R2) biomechanics monitoring, (R3) economic decision-making.
% Through theoretical analysis and real-world-relevant numerical examples (R1, R2, R3), we show that DCC achieves optimal behavior with high accuracy in actual physical time, even for highly dynamical, high-dimensional problems.
% Notably, DCC requires only $77\,\%$ CPU usage (vs. gold-standard method, $8708\,\%$) for optimally controlling a 1000-dimensional system at $1\,\textup{kHz}$ on a modern CPU.
% These advantages of DCC are crucial for every time- and safety-critical application, such as aerospace, medical, or technical infrastructure operations.
% %
% Besides, the DCC embodies a system-level, causally deterministic interpretation of constrained optimization, \addedtwo{advancing the} understanding \addedtwo{of} optimality principles in complex natural and artificial processes.\\
% \\
\end{quote}

\section{Introduction}\label{sec:introduction}

Optimal control provides a fundamental framework for understanding and governing dynamical processes across the natural and engineered sciences, yet computing optimal control actions for complex nonlinear systems in real time remains a longstanding challenge. Optimality principles govern a wide range of natural and engineered processes, including neuromechanics~\cite{ErdemirMcLHerBog2007,MarshallGlaTraAme2022}, biochemical networks~\cite{OrthThiPal2010,SchuetzZamZamHei2012,NaseriKof2020}, autonomous systems~\cite{AllenspachBodBruRin2020,JaleelSha2020,IzzoBlaFerOri2024}, economic decision-making~\cite{Steinbach2001}, energy systems and power grids~\cite{FrankReb2016,AbdiBeiSca2017,WeiWanLiMei2017,JhaInaBisSur2023}, telecommunications~\cite{GershmanSidShaBen2010}, logistics~\cite{MataiSinMit2010}, and resource allocation~\cite{BouajajaDri2017}. In many of these domains, optimal policies must be computed faster than physical time in order to maintain stability, safety, and efficiency. When optimal actions cannot be determined before the system state evolves, control decisions fall behind the physical dynamics of the system, creating what we refer to as the \emph{real-time barrier}. We consider a broad but well-defined class of pseudoconvex optimal control problems for nonlinear control-affine systems.

Breaking this barrier for complex, high-dimensional systems has proven challenging. Classical optimal control formulations--ranging from Pontryagin’s maximum principle to Hamilton-Jacobi-Bellman theory--ultimately require solving constrained optimization problems. In practice, these problems are typically addressed using iterative numerical algorithms such as active-set methods~\cite{Fletcher1987,NocedalWri2006}, interior-point methods (IPM)~\cite{PotraWri2000,BoydVan2004,NocedalWri2006}, or sequential quadratic programming (SQP)~\cite{NocedalWri2006,GillMurSau2005,BueskensWas2013}. Although these approaches have enabled major advances in optimal control and numerical optimization, their iterative nature often prevents optimal actions from being computed directly within physical time for complex dynamical systems. Consequently, many real-world implementations rely on approximations such as discretization, linearization, steady-state assumptions, or predictive schemes such as model-predictive control~\cite{MayneRawRaoSco2000,DrgonaArrCupBlu2020}.

These approximations enable useful practical solutions but introduce limitations. Linear-quadratic formulations apply primarily to simplified system classes~\cite{Mehrmann1991}, steady-state approaches neglect important dynamical behavior in systems such as metabolic networks~\cite{OrthThiPal2010}, and predictive control schemes require repeated numerical optimization over short horizons. As system dimensionality and dynamical complexity increase, the computational burden of these approaches grows rapidly, making real-time optimal control increasingly difficult. As a result, a general formulation capable of producing optimal control actions directly within the time scale of the physical system has remained elusive.

Over the past several decades, researchers have explored replacing iterative numerical optimization with compact closed-form controllers. Several promising approaches have been proposed~\cite{Yamashita1980,ZhouShi1997,XiaFen2005,FepponAllDap2020,RaveendranMahVai2023,AllibhoyCor2023}, particularly for quadratic or low-dimensional optimization problems. However, existing formulations are typically limited to specific problem classes and have not demonstrated scalability to high-dimensional nonlinear optimal control problems. Consequently, a general and scalable closed-form dynamical representation capable of addressing realistic high-dimensional systems remains an open challenge.

\begin{figure}[!h]
	\centering
	\includegraphics[trim={3cm 0 3cm 0},clip, width=16cm]{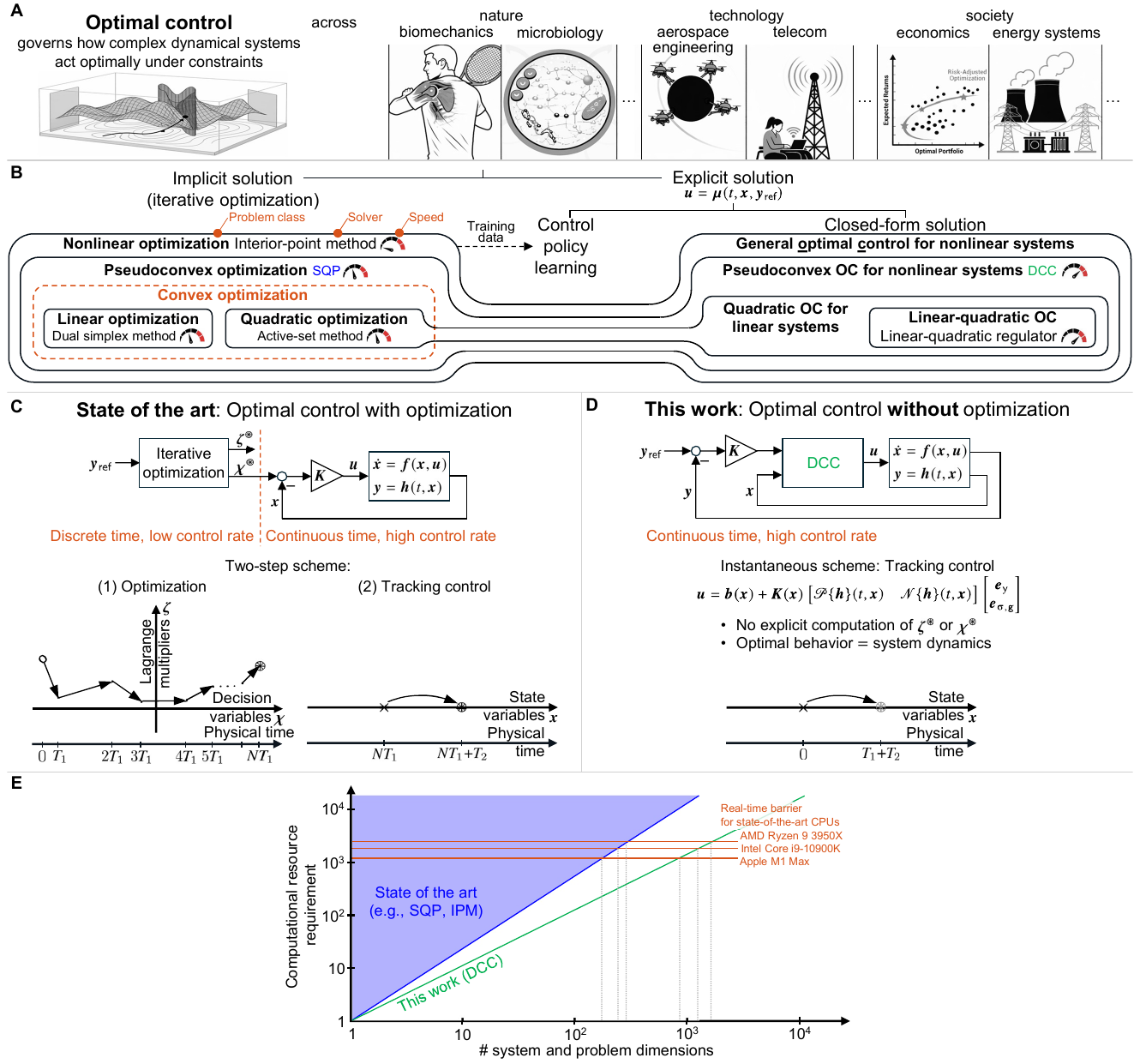}
	\caption{\label{fig:introduction}%
		\textbf{Dual Cost-Constraint projection (DCC) enables real-time optimal control without optimization.}
		(\textbf{A}) Optimal control governs how complex dynamical systems act optimally under constraints and spreads across nature, technology and society. (\textbf{B}) It is solved either implicitly or explicitly, depending on problem class (shown as Venn diagram). For (pseudo)convex problems, (\textbf{C}) state-of-the-art optimal control relies on iterative numerical optimization, e.g., Sequential Quadratic Programming (SQP), Interior-Point Method (IPM), whose (\textbf{E}) requirement on computational resource exceeds physical-time constraints (the \emph{real-time barrier}) even at low system dimensionality. (\textbf{D}) The proposed DCC is a closed-form solution that reformulates optimal control problems involving convex and certain pseudoconvex objectives as a stable dynamical system evolution, enabling (\textbf{E}) real-time optimal control for both low- and high-dimensional systems. 
	}
\end{figure}

Here we introduce the \emph{Dual Cost--Constraint projection} (DCC), a dynamical-systems formulation that enables real-time optimal control for a broad class of pseudoconvex optimal control and optimization problems across diverse domains in science and engineering{, see Figure~\ref{fig:introduction}}.
Instead of solving an optimization problem iteratively, the proposed approach embeds cost and constraint functions directly within the dynamics of the control system. Specifically, DCC constructs a \emph{dual projector} that maps cost and constraint errors into the control input space of a nonlinear dynamical system. The resulting closed-loop dynamics are designed such that their attractive equilibria coincide with the optimal solutions of the corresponding constrained optimization problem. In this sense, the proposed formulation provides a \emph{dynamical realization of constrained optimality conditions}, akin to those traditionally derived in classical optimal control theory for linear systems with quadratic cost functions~\cite{Mehrmann1991}. 

This construction establishes a structural equivalence between numerical optimization and nonlinear feedback control. The derivation reveals that interior-point optimization iterations admit an equivalent realization as nonlinear feedback dynamics. As a result, optimal behavior can be interpreted as the stable evolution of a dynamical system rather than the outcome of an iterative numerical optimization process. Thus, the proposed formulation replaces iterative optimization with a continuous-time dynamical process whose equilibria encode optimal solutions.

To evaluate the practical implications of this formulation, we compare DCC against the {widely used} state-of-the-art solver SQP on a real-time computing platform. 
{SQP is a natural benchmark because it addresses the same class of constrained nonlinear optimization problems considered in this work and exhibits an excellent convergence for smooth (pseudo)convex problems~\cite{Fletcher1987,NocedalWri2006}.
On top of that, the state-of-the-art SQP implementation combines its theoretical strength with decades of algorithmic refinement and highly optimized implementations~\cite{TheMathWorks2020,GillMurSau2005,BueskensWas2013}, including global convergence guarantees, warm-start techniques, and efficient linear-algebra routines. 
Although other optimization algorithms exist, such as IPM or simplex method, these either require longer computation times~\cite{Fletcher1987,NocedalWri2006} (see also Supplementary Information, Section~\ref{SI-sec:computational complexity}) or exploit more specialized problem structures.
As a result, SQP is widely considered a state-of-the-art reference method for (pseudo)convex optimization and therefore provides a strong baseline for evaluating DCC.
}

Across a range of benchmark problems, DCC demonstrates substantial improvements in computational efficiency while maintaining high solution quality. In particular, the cost-function error remains below $1.3\,\%$ for DCC compared with $197.2\,\%$ for SQP, task constraint violation remains below $1.2\,\%$ compared with $117.7\,\%$, and state tracking error remains below $0.2\,\%$ compared with $11.3\,\%$. For a $1000$-dimensional system operating at a control rate of $1\,\mathrm{kHz}$, DCC requires only $77\,\%$ CPU usage on a modern processor, whereas SQP requires more than $80$ processors to achieve comparable performance. CPU usage below $100\,\%$ indicates that the computation can be performed within the available real-time budget.

Beyond these benchmarks, the proposed formulation enables real-time optimal control in several safety- and time-critical domains. In particular, we demonstrate the approach in three representative examples illustrated in Figure~\ref{fig:otc_and_applications}: autonomous systems control (R1), biomechanics monitoring (R2), and optimal economic decision-making (R3). These examples highlight the practical relevance of DCC for applications in which control decisions must be generated within the time scale of the physical process.

The formulation also extends naturally to classical stationary optimization problems in fields such as microbiology~\cite{OrthThiPal2010}, energy systems~\cite{JhaInaBisSur2023}, and telecommunications~\cite{GershmanSidShaBen2010}, while enabling the study of their dynamical aspects within realistic physical models. These capabilities open new possibilities for real-time optimization in complex systems ranging from smart energy grids and communication networks to biological and socio-economic processes.

Conceptually, DCC reveals a close relationship between iterative algorithms and nonlinear dynamical systems theory in the context of constrained optimization. By embedding cost and constraint functions and their derivatives directly into system dynamics, optimal behavior can be interpreted as the stable evolution of a dynamical system rather than the outcome of an iterative numerical optimization process. This perspective demonstrates that optimal control problems may admit \emph{closed-form dynamical representation in which optimal solutions arise from system evolution itself}.

\section{Results}\label{sec:results}
These results suggest that certain classes of optimal control problems can be solved not by iterative optimization, but by appropriately constructed dynamical systems.
{%
Figure~\ref{fig:otc_and_applications} provides an overview of the proposed formulation and its application scope. Panel A illustrates the DCC formulation, while Panels B and C highlight representative application domains and benchmark comparisons with SQP.
}

\afterpage{%
	%\unsetvruler % turn off line numbering
	\newgeometry{left=0.1cm,right=0.1cm,top=0.1cm,bottom=0.1cm}
	\restoregeometry
	\begin{figure}[!h]
		\centering
		\vspace{-1cm}
		\resizebox{0.999\textwidth}{!}{
			\begin{tikzpicture}[auto,node distance=0cm]	
				\node[invisiblenode] (origin) at (0cm, 0cm) {
					\includegraphics[scale=1.0]{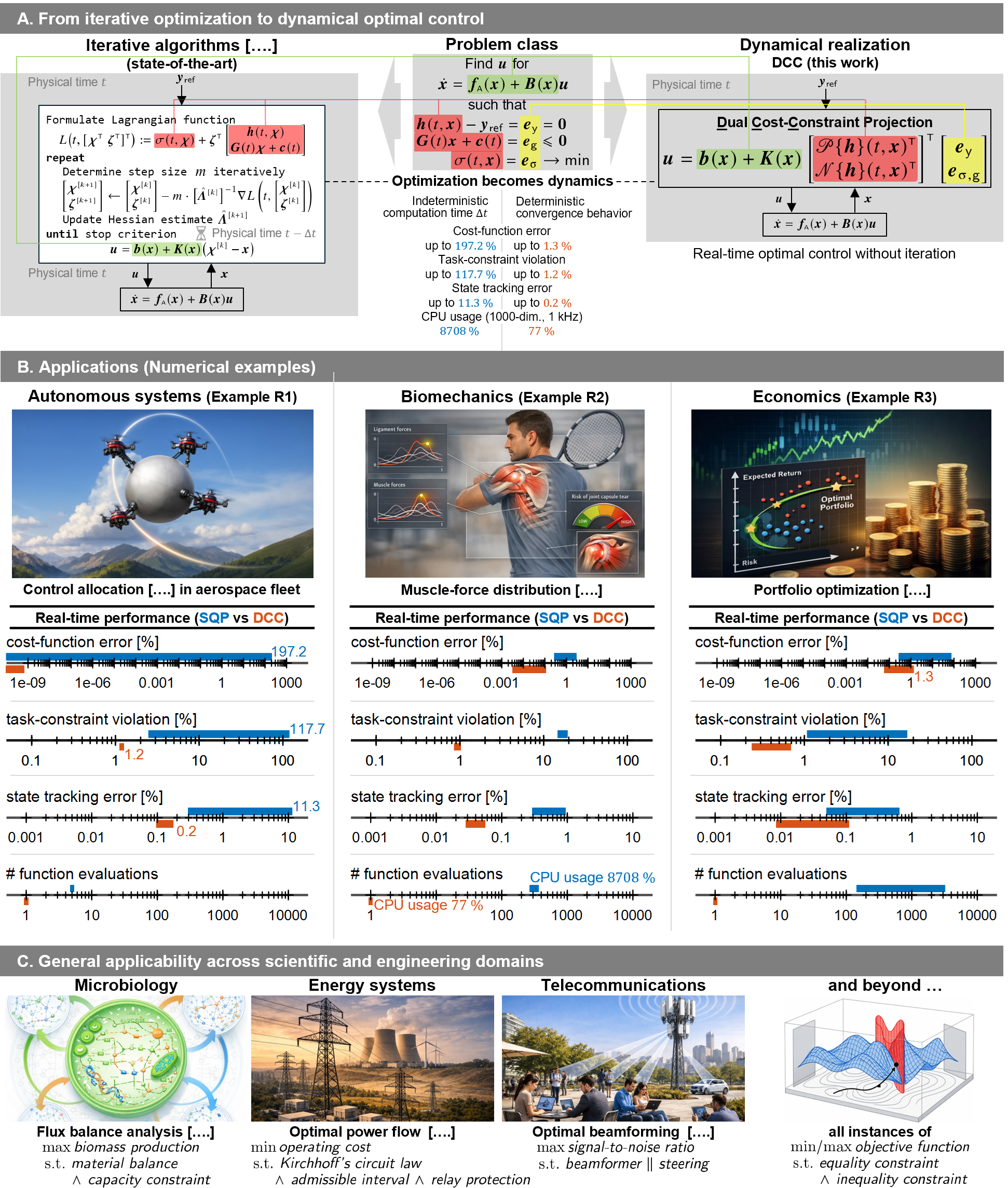}
				};
				
				% iterative algorithms
				\node[invisiblenode, above of=origin, xshift=-3.65cm, yshift=8.125cm, scale=0.5, fill=white] (B) {$\vphantom{1}\qquad$};	
				\node[invisiblenode, above of=origin, xshift=-3.325cm, yshift=8.125cm, scale=0.6, fill=none] (B) {\cite{Fletcher1987,PotraWri2000,BoydVan2004,NocedalWri2006,BueskensWas2013}};	
				
				% control allocation
				\node[invisiblenode, above of=origin, xshift=-5.125cm, yshift=0.0625cm, scale=0.5, fill=white] (B) {$\vphantom{1}\quad~$};	
				\node[invisiblenode, above of=origin, xshift=-5.16cm, yshift=0.06cm, scale=0.5, fill=none] (B) {\cite{AllenspachBodBruRin2020,IzzoBlaFerOri2024}};	
				
				% muscle-force distribution
				\node[invisiblenode, above of=origin, xshift=1.185cm, yshift=0.0625cm, scale=0.5, fill=white] (B) {$\vphantom{1}\quad~~$};	
				\node[invisiblenode, above of=origin, xshift=1.19cm, yshift=0.0525cm, scale=0.5, fill=none] (B) {\cite{ErdemirMcLHerBog2007,MarshallGlaTraAme2022}};	
				
				% portfolio
				\node[invisiblenode, above of=origin, xshift=6.12cm, yshift=0.0625cm, scale=0.5, fill=white] (B) {$\vphantom{1}\quad\quad$};	
				\node[invisiblenode, above of=origin, xshift=6.175cm, yshift=0.0525cm, scale=0.5, fill=none] (B) {\cite{Markowitz1952,SandhuGeoTan2016}};	
				% flux balance analysis
				\node[invisiblenode, above of=origin, xshift=-4.5cm, yshift=-8cm, scale=0.5, fill=white] (B) {$\vphantom{1}\quad\quad$};	
				\node[invisiblenode, above of=origin, xshift=-4.5cm, yshift=-8cm, scale=0.5, fill=none] (B) {\cite{BeckerFeiMoHan2007,OrthThiPal2010}};
				
				% optimal power flow
				\node[invisiblenode, above of=origin, xshift=-0.9cm, yshift=-8cm, scale=0.5, fill=white] (B) {$\vphantom{1}\quad\quad$};	
				\node[invisiblenode, above of=origin, xshift=-0.9cm, yshift=-8cm, scale=0.5, fill=none] (B) {\cite{WeiWanLiMei2017,JhaInaBisSur2023}};
				
				% beamforming
				\node[invisiblenode, above of=origin, xshift=2.9cm, yshift=-8cm, scale=0.5, fill=white] (B) {$\vphantom{1}\quad\quad$};	
				\node[invisiblenode, above of=origin, xshift=2.95cm, yshift=-8cm, scale=0.5, fill=none] (B) {\cite{GershmanSidShaBen2010,LauZha2015}};
			\end{tikzpicture}
		}
		%}%
	\caption{\label{fig:otc_and_applications}
		\textbf{Dual Cost--Constraint Projection enables real-time optimal control across domains.}
		\textbf{A}: The proposed Dual Cost--Constraint Projection (DCC) transforms a constrained optimization problem into a dynamical system whose attractive equilibria coincide with optimal solutions. This dynamical formulation enables optimal control actions to be generated directly in physical time.
		\textbf{B}: 
		Representative numerical examples studied in this work, including autonomous systems control (R1), muscle-force distribution in biomechanics (R2), and portfolio optimization in economic decision-making (R3). In each case, DCC is compared against the {widely used} state-of-the-art solver Sequential Quadratic Programming (SQP).
		Performance bars quantify improvements in optimality, constraint satisfaction, state tracking accuracy, and computational efficiency. DCC achieves real-time performance for high-dimensional systems where SQP exceeds the available computational budget.
		\textbf{C}: Broader application scope of optimal control and optimization problems across multiple scientific and engineering domains, including microbiology, energy systems, and telecommunications.
	}
\end{figure}
\restoregeometry 
}

\subsection{Problem class}
Almost all physical processes are in essence nonlinear systems.
In this work, we study the broad class of nonlinear time-variant control-affine systems
\begin{subequations}\label{eq:gryphon}
	\begin{align}
		\dot{\xv} &= \fv_{\scriptscriptstyle\textup{A}}(\xv) + \Bm(\xv) \uv, \label{eq:gryphon, state eqn}\\
		\yv &= \hv(t,\xv), \label{eq:gryphon, output eqn}
	\end{align}
\end{subequations}
where 
%$\xv = [\,x_{1}~~x_{2}~~\cdots~~x_{\dimx}\,]^\transp$,
%$\uv = [\,u_{1}~~u_{2}~~\cdots~~u_{\dimx}\,]^\transp$, 
%and $\yv = [\,y_{1}~~y_{2}~~\cdots~~y_{\dimy}\,]^\transp$
$\xv \in \realset^{\dimx}$,
$\uv \in \realset^{\dimx}$, 
and $\yv \in \realset^{\dimy}$
are the state, the input, and the output, respectively, of this $\dimx$-dimensional system,
and $t$ denotes time.
For instance, $\xv$ can denote electric current of a motor in an aerial vehicle~\cite{AllenspachBodBruRin2020}, activation level of muscles in a human limb~\cite{HuKueHad2020}, or capital invested on an asset market~\cite{Steinbach2001}.
Correspondingly, $\uv$ can denote voltage, neural excitation, and capital reallocation, respectively;
and $\yv$ can denote the linear and angular acceleration of the aerial vehicle, the resulting joint torques acting on the human limb, or the return on capital, respectively.
Given the reference trajectory $\desired{\yv}$ which $\yv$ should follow,
the optimality principle associated with \eqref{eq:gryphon} can be representatively and equivalently written as 
%either the pseudoconvex optimization problem
%\begin{equation}\label{eq:constrained optimization}
%	\begin{split}
%	%\chiv^{\optimal} = \arg 
%	\min_{\chiv} \sigma(t,\chiv) \qquad
%	\text{subject to}~
%	\nullv &= \hv(t,\chiv) - \desired{\yv}(t), \\[-1ex] %\label{eq:result-section, equality constraint}
%	%~\wedge~
%	\nullv &\geqslant \gv(t,\chiv) := \Gm(t) \chiv + \cv(t) %\label{eq:result-section, inequality constraint}
%	\end{split}
%\end{equation}
%for every time instant $t \geqslant 0$
%or the optimal tracking control problem
%\begin{equation}\label{eq:optimal control problem}
%\begin{split}
%	%\phiv^{\optimal} = \arg 
%	\min_{\phiv} \int_{0}^{t_{\final}} \sigma\big(t,\phiv(t)\big)\,\mathrm{d}t \qquad
%	%\quad~~
%	\text{subject to}~
%	\nullv = \hv\big(t,\,&\phiv(t)\big) - \desired{\yv}(t)~\forall t \in (0, t_{\final}], \\[-0.5ex]
%	%~\wedge~
%	\nullv \geqslant \gv\big(t,\,&\phiv(t)\big)~\forall t \in (0, t_{\final}], \\[0ex]%[0.75ex]
%	%~\wedge~
%	\xv(0) =\ &\phiv(0)
%\end{split}
%\end{equation}
%for the time period $[0, t_{\final}]$.
%
either the pseudoconvex optimization problem for every time instant $t \geqslant 0$
or the optimal tracking control problem for the time period $[0, t_{\final}]$, respectively:\par
\hspace{-0.05\textwidth}
\begin{minipage}{0.4\textwidth}
\begin{subequations}\label{eq:constrained optimization}
\begin{align}
	%\chiv^{\optimal} = \arg 
	\min_{\chiv} \sigma(t,\,&\chiv) \\
	\text{subject to}~\nonumber\\
	\nullv = \hv(t,\,&\chiv) - \desired{\yv}(t), \label{eq:constrained optimization, equality constraint} \\
	%~\wedge~
	\nullv \geqslant \gv(t,\,&\chiv) := \Gm(t) \chiv + \cv(t) \label{eq:constrained optimization, inequality constraint} \\
	~\nonumber
\end{align}
\end{subequations}
\end{minipage}
\hfill\vline\hfill
\hspace{-0.05\textwidth}
\begin{minipage}{0.56\textwidth}
\begin{subequations}\label{eq:optimal control problem}
\begin{align}
	%\phiv^{\optimal} = \arg 
	\min_{\phiv} \int_{0}^{t_{\final}} \sigma\big(t,\,&\phiv(t)\big)\,\mathrm{d}t \\
	%\quad~~
	\text{subject to}\qquad~\nonumber\\
	\nullv = \hv\big(t,\,&\phiv(t)\big) - \desired{\yv}(t)~\forall t \in (0, t_{\final}], \label{eq:optimal control problem, equality constraint} \\
	%~\wedge~
	\nullv \geqslant \gv\big(t,\,&\phiv(t)\big)~\forall t \in (0, t_{\final}], \label{eq:optimal control problem, inequality constraint} \\
	%~\wedge~
	\xv(0) =\ &\phiv(0)
\end{align}
\end{subequations}
\end{minipage}
\vspace{1ex}

\noindent
Herein, $\chiv \in \realset^{\dimx}$ denotes candidate variable and $\phiv: \realset_{\geqslant 0}\to\realset^{\dimx}$ candidate trajectory; both are to be determined.
Let $\chiv^{\optimal}(t) \in \realset^{\dimx}$ denote a local minimizer of~\eqref{eq:constrained optimization} at $t$
and $\phiv^{\optimal}: [0, t_{\final}]\to\realset^{\dimx}$ a local minimizer of \eqref{eq:optimal control problem} for the time period~$[0, t_{\final}]$,
then the state $\xv(t)$ of~\eqref{eq:gryphon} is called \emph{optimal} if $\xv(t) = \chiv^{\optimal}(t)$ or $\xv(t) = \phiv^{\optimal}(t)$.
It means that a minimal cost--either instantaneous, i.e., $\sigma(t,\xv)$, or accumulated, i.e., its time integral--is preferred as long as $\yv$ matches $\desired{\yv}$% 
, i.e., equality constraint~\eqref{eq:constrained optimization, equality constraint} or \eqref{eq:optimal control problem, equality constraint}, %(i.e., equality constraint $\nullv = \hv(t,\xv) - \desired{\yv}$) 
and $\xv$ remains within the region specified by the inequality constraint~\eqref{eq:constrained optimization, inequality constraint} or \eqref{eq:optimal control problem, inequality constraint}. %$\nullv \geqslant \gv(t,\xv) \in \realset^{\dimc}$.
Concrete examples are: 
a minimal energy consumption is desired for an aerial vehicle, as long as an assigned task can be finished while every motor operates within its payload capacity~\cite{AllenspachBodBruRin2020};
lower activation level of muscles is preferred while an intended movement can be accomplished within physiological feasibility~\cite{HuKueHad2020};
and the investment risk shall be as low as possible while the total return meets the investor's expectation according to his investment capability~\cite{Steinbach2001}.

\subsection{Dynamical realization of optimal control}
Our \emph{primary result} is the closed-form dynamical representation \emph{Dual Cost-Constraint projection} (DCC)
\begin{equation}\label{eq:OTC}
	\uv 
	= 
	\bv(\xv) 
	+ \Km(\xv) \begin{bmatrix}
		\mathscr{P}\{\hv\}(t,\xv) & \mathscr{N}\{\hv\}(t,\xv)
	\end{bmatrix} \begin{bmatrix}
		\errv_{\textup{y}} \\
		\errv_{\upsigma,\textup{g}}
	\end{bmatrix}
\end{equation}%
that solves problems~\eqref{eq:constrained optimization} and \eqref{eq:optimal control problem} for any system following the general form of \eqref{eq:gryphon}.
DCC~\eqref{eq:OTC} is composed of two components:
\begin{enumerate}
	\item the compensator $\bv(\xv)$ that cancels nonlinear dynamics~\eqref{eq:gryphon, state eqn}, and
	\item the \emph{dual projector} $\begin{bmatrix}
		\mathscr{P}\{\hv\}(t,\xv) & \mathscr{N}\{\hv\}(t,\xv)
	\end{bmatrix}$ that projects errors $\errv_{\textup{y}} := \desired{\yv} - \yv$ (related to equality constraint) and $\errv_{\upsigma,\textup{g}}$ (related to cost function and inequality constraint) to the input space,
\end{enumerate}
where $\Km(\xv) := \Bm(\xv)^{-1} K_{\textup{x}}$ is the gain scheduling, $\mathscr{P}\{\hv\}(t,\xv)$ and $\mathscr{N}\{\hv\}(t,\xv)$ are the local pseudoinverse and the associated null-space projector of $\hv$ in \eqref{eq:gryphon, output eqn}, respectively.
Due to this complementarity, \eqref{eq:OTC} nullifies $\errv_{\textup{y}}$ and minimizes $\errv_{\upsigma,\textup{g}}$.
Instead of being trivial or arbitrary, 
the dual projector $\mathscr{P}$ and $\mathscr{N}$ has a strict form that integrates the cost function $\sigma$, the equality constraint $\hv$, and the inequality constraint $\gv$ of \eqref{eq:constrained optimization} and \eqref{eq:optimal control problem}.
The details of \eqref{eq:OTC} are provided in Supplementary Equation~\eqref{SI-eq:DCC projection}. These results illustrate that optimal behavior can emerge directly from stable system evolution.

With \eqref{eq:OTC}, system~\eqref{eq:gryphon} exhibits an exponential convergence and tracking behavior in 
$\yv(t) \to \desired{\yv}(t)$ and $\xv(t) \to \chiv^{\optimal}(t) = \phiv^{\optimal}(t)$.

The reader may refer to the Supplementary Information for details of the system and problem class~\eqref{eq:gryphon}--\eqref{eq:optimal control problem}, the derivations and deductions behind~\eqref{eq:OTC} and the formal proofs of %Theorems~\ref{theorem: y to y_ref exponential}--\ref{theorem:chi optimal equals phi optimal}
the exponential convergence and tracking behavior.

\afterpage{%
	\newgeometry{left=1.5cm,right=1.5cm,top=3.85cm,bottom=3.15cm}
\begin{figure}[!h]
	\centering
	\includegraphics[width=16cm]{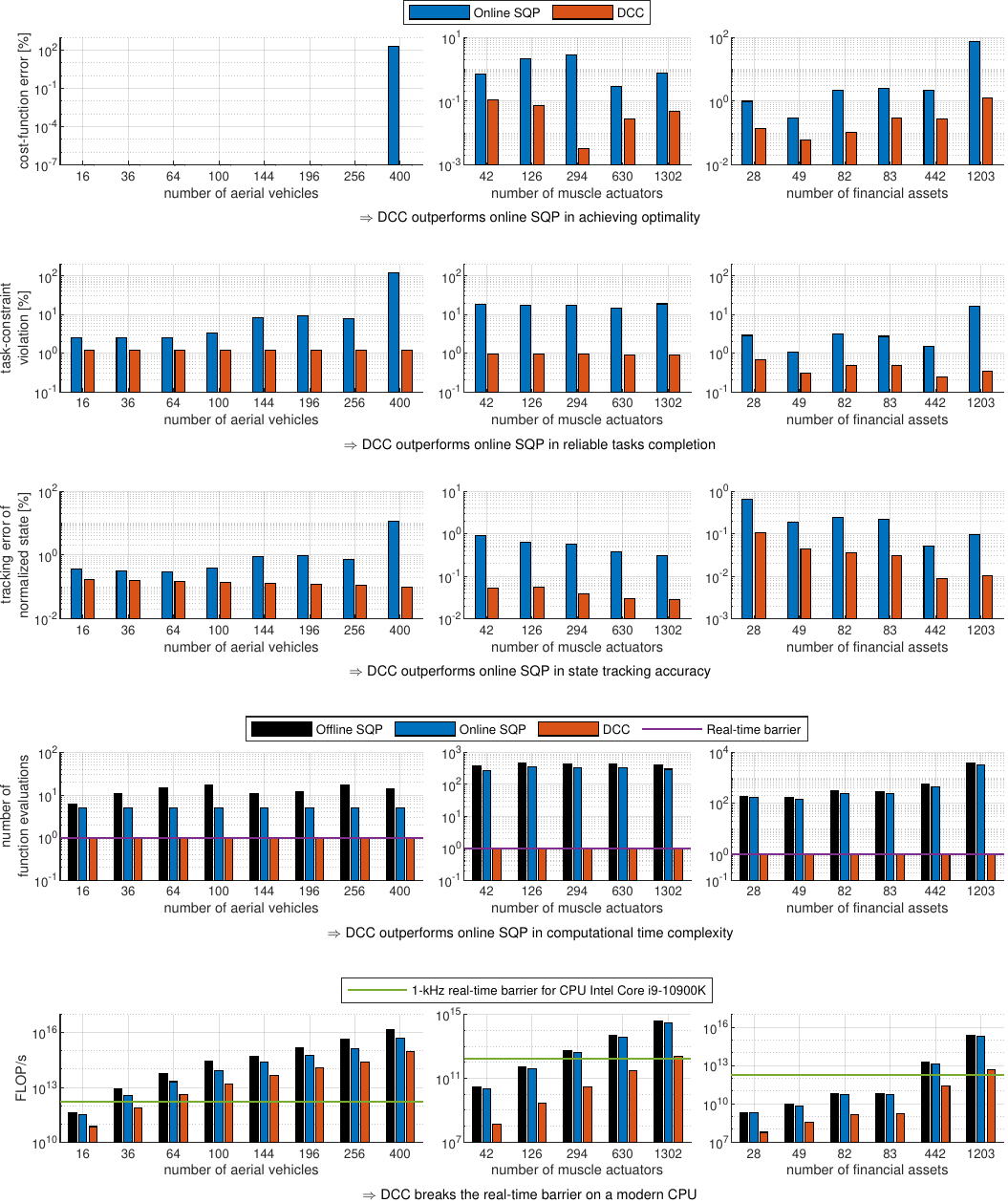}
	\caption{\label{fig:performance metrics otc vs sqp}%
		\textbf{Real-time performance metrics} of online \SQP\ (state of the art, blue) vs. DCC\ (this work, red) for three benchmark problems (R1: aerial fleet control allocation, R2: muscle-force distribution, and R3: portfolio optimization).
		Offline \SQP\ (state of the art, black) which provides the true optimal solution \smash{$\chiv^{\optimal}$} serves as reference.	
		The ordinate axis is displayed on the logarithmic scale and starts from a different value in each subplot.	
		Refer to Methods (Secion~\ref{sec:methods}) for term clarification and metric definition.
	}
\end{figure}
\restoregeometry
}

\subsection{Real-time control of high-dimensional systems}

To evaluate the proposed formulation and the theoretical convergence guarantees, we study three representative optimal control problems drawn from autonomous systems, biomechanics, and economic decision-making. 
These three scalable, realistic benchmarks (R1: aerial-fleet control allocation, R2: muscle-force distribution, R3: portfolio optimization) are analyzed with a state-space dimension $\dimx$ ranging from $28$ to $9600$, see Figure~\ref{fig:performance metrics otc vs sqp}.

Across all benchmarks and dimensions, DCC (red bars) maintains consistently low cost-function errors (first row), low task-constraint violations (second row), and low state tracking errors (third row): each stays below $1\,\%$, even at $\dimx = 9600$ (corresponding to $n = \dimx/24 = 400$ aerial vehicles).
More importantly, these metrics remain essentially invariant with respect to $\dimx$, demonstrating preserved optimality and strict task-constraint satisfaction in high-dimensional settings. These results indicate that the optimality conditions remain satisfied even as the dimensionality of the control problem increases by more than two orders of magnitude.

In contrast, online \SQP\ (blue bars) exhibits in general higher cost-function errors, higher task-constraint violations, and higher state tracking errors, with a pronounced increase beyond a moderate problem size ($144$ aerial vehicles, corresponding to a $\dimx = 3456$).
At the highest $\dimx = 9600$, the cost-function error reaches $192\,\%$ and the task-constraint violation $117\,\%$, indicating a failure to produce optimal behavior and to complete the task under real-time execution. Consequently,
\begin{itemize}
 	\item high task-constraint violations in R1 may lead to loss of flight stability or unintended load release,
 	\item insufficient tracking accuracy in R2 may prevent timely injury-risk detection, and
 	\item suboptimal decisions in R3 may lead to poor investment performance.
\end{itemize}
% \begin{itemize}
% 	\item high task-constraint violations in R1 cause aerial vehicles to crash into, e.g., a populated area,
% 	\item a bad tracking quality in R2 may fail the time-critical injury prevention, and
% 	\item the inability to behave optimally in R3 can lead to poor or even negative returns on investment.
% \end{itemize}
The emergence of these disastrous consequences will be analyzed in the respective subsections.

In summary, DCC outperforms online \SQP\ on a real-time platform by 2 orders of magnitude in 
proven optimality ($<1.3\,\%$ vs. $<197.2\,\%$), 
reliable task completion ($<1.2\,\%$ vs. $<117.7\,\%$),
state tracking accuracy ($<0.2\,\%$ vs. $<11.3\,\%$); and the corresponding metrics of DCC are less variant across $\dimx$ than those of \SQP.
This pronounced performance gap originates from DCC's ability to break the real-time barrier, whereas iterative algorithms such as \SQP\ fundamentally fail to do so. While SQP successfully solves the optimization problem offline, its computational cost exceeds the available real-time budget by more than two orders of magnitude. This disparity is elaborated in the following subsection.

The reader may refer to Section~\textit{\nameref{sec:methods}} for details of the numerical benchmarking Examples R1--R3, as well as the definition of real-time performance metrics.

\subsection{DCC breaks the real-time barrier}
Due to DCC's closed-form dynamical representation~\eqref{eq:OTC}, the functions $\sigma$, $\gv$, $\hv$ are evaluated exactly \emph{once} per time step under all circumstances; that is, $N_{\FcnEvalsubscript}^{\OTCsubscript}=1$.
In contrast, state-of-the-art optimization algorithms, such as \SQP, due to their inherent iteration process, require a substantially larger number of \uls{f}unction \uls{e}valuations, $N_{\FcnEvalsubscript}^{\SQPsubscript} \gg 1$, depending on desired accuracy (e.g., tolerance and stopping criterion), see Figure~\ref{fig:performance metrics otc vs sqp} (fourth row).
No general theoretical upper bound for $N_{\FcnEvalsubscript}^{\SQPsubscript}$ is available; only empirical values for practical convergence.
By achieving the deterministic, time-invariant and, most importantly, minimal possible number of function evaluations, DCC breaks the real-time barrier in principle.
The only question remaining is how much computing resource is required for a certain problem instance.

In practice, real-time capability is ultimately constrained by the available computing platform, such as the CPU. 
For DCC to control a $\dimx$-dimensional system in real-time at a sampling rate of $1/\dt$, the number of \uls{fl}oating-point \uls{op}erations per second [FLOP/s] must be at least
\begin{equation}\label{eq:OTC flop/s}
	F_{\OTCsubscript}
	\approx \frac{1.30\,\dimx^{3}}{\dt}
	\approx \frac{1.30}{0.67 + 0.51\,N_{\FcnEvalsubscript}^{\SQPsubscript}} \times F_{\SQPsubscript},
\end{equation}
where $F_{\OTCsubscript}$ and $F_{\SQPsubscript}$ are the minimum required FLOP-per-second for DCC and \SQP, respectively.
Their precise values
for a representative sampling interval $\dt = 1\,\textup{ms}$ and three different realistic examples are illustrated in Figure~\ref{fig:performance metrics otc vs sqp} (fifth row).
A sampling interval of $\dt = 1\,\textup{ms}$ means that optimal control actions must be computed within a strict real-time budget of $1\,\textup{ms}$.
To evaluate real-time feasibility, we consider a modern processor \emph{Intel Core i9-10900K} ($1696 \times 10^{9}$~FLOP/s). On this platform, DCC is capable of controlling systems with $\dimx \approx 1000$ state variables within the $1$~ms budget. In contrast, \SQP\ is limited to systems with $\dimx < 300$ under identical hardware constraints.

Specifically, for a system with $\dimx = 1000$, DCC requires $F_{\OTCsubscript} \approx 1.30\times 10^{12}\,\textup{FLOP/s}$, corresponding to a CPU usage of $76.65\,\%$. Hence, real-time operation is achievable on a single CPU of this type. 
For general systems and problem setups, such as Examples R2 and R3 (for reference, R1's $\sigma(t,\chiv)$ is trivial), a typical value of $N_{\FcnEvalsubscript}^{\SQPsubscript} \approx 290$ is observed (cf. Figure~\ref{fig:performance metrics otc vs sqp}, fourth row), yielding 
$F_{\SQPsubscript} \approx 1.48\times 10^{14}\,\textup{FLOP/s}$ demanded by \SQP. This value corresponds to a CPU usage of $8707.82\,\%$. This is equivalent to at least $88$ CPUs operating in parallel without considering the inter-CPU communication overhead. This exceeds the available computational budget by almost two orders of magnitude.

A further limitation of \SQP\ is that the values of $F_{\SQPsubscript}$ reported here are empirical observations from Examples R1--R3 and cannot be guaranteed for other systems or problem instances, as $N_{\FcnEvalsubscript}^{\SQPsubscript}$ behaves indeterministically.
In contrast, $F_{\OTCsubscript}$ in \eqref{eq:OTC flop/s} generalizes reliably due to deterministic and time-invariant $N_{\FcnEvalsubscript}^{\OTCsubscript} = 1$.

This pronounced superior real-time capability of DCC over the current state-of-the-art \SQP\ merits further investigation for practical implications.
To this end, the following subsections demonstrate the performance of DCC and \SQP\ in controlling and optimizing different dynamical systems across real-world-relevant application domains.
All numerical simulations are performed under identical computational budgets by emulating a computing platform with a prescribed FLOP-per-second capacity $F_{\OTCsubscript}$, scaled with $\dimx$ as given in~\eqref{eq:OTC flop/s}.
These experiments illustrate the ability of DCC to generate optimal control actions directly in real time, while enabling a fair and consistent comparison with \emph{online SQP}--it operates under the same computational budget $F_{\OTCsubscript}$ and its computation times exceeding the real-time budget manifest as control action delays. 
As a reference for evaluation, we consider further the idealized \emph{offline SQP}, which assumes unlimited computational resources and yields the delay-free optimal solution $\chiv^{\optimal}$.

\begin{figure}[!t]
	\centering
	\includegraphics[scale=0.817]{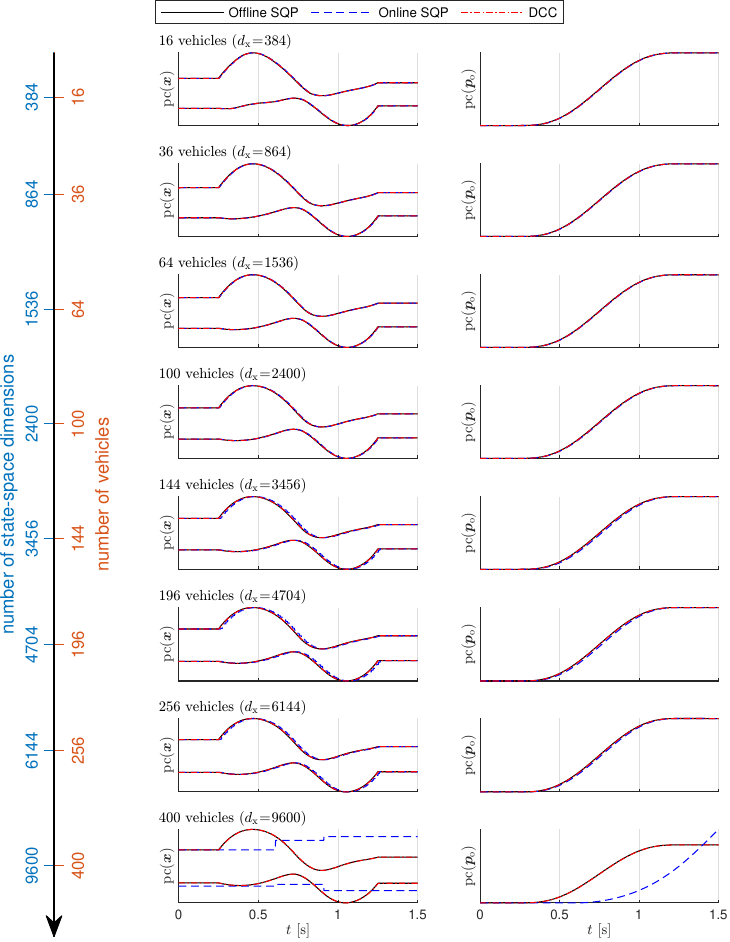}
	\caption{\label{fig:swarm tracking timeseries}%
		The actuator state trajectory $\xv(t)$ of aerial vehicles and the object's pose trajectory $\pv_{\text{o}}(t)$ resulting from \SQP~(blue) and DCC~(red) executed on an emulated real-time computer. 
		The offline computed optimal solution $\chiv^{\optimal}(t)$ (black) serves as the reference.
		Here, $\principal(\xv)$ denotes the first two \uls{p}rincipal \uls{c}omponents of $\xv$.
	}
\end{figure}

\subsection{Autonomous systems control (R1)}

The \emph{control allocation problem} of omnidirectional micro aerial vehicles aims to determine the optimal distribution of control inputs among multiple actuators in order to generate the desired wrench for flight, maneuverability, or interaction with the environment while satisfying physical constraints such as payload capacity~\cite{KirchengastSteHor2018}.

To evaluate the proposed formulation, we consider a cooperative aerial-transport scenario in which a fleet of $n$ aerial vehicles carries an object along a predefined trajectory. As the fleet size increases, the dimensionality of the control problem grows proportionally ($n \propto \dimx$), making real-time optimal control increasingly challenging. This experiment therefore provides a representative benchmark for evaluating whether DCC can generate optimal control actions within the time scale of the physical system.

Figure~\ref{fig:swarm tracking timeseries} illustrates the resulting tracking behavior. When the system is controlled using online \SQP, the delay between the computed system state $\xv^{\SQPsubscript}(t)$ and the offline optimal solution $\chiv^{\optimal}(t)$ increases with fleet size. For moderate fleet sizes ($n \leqslant 256$), this delay leads to small deviations of the object's pose trajectory $\pv_{\text{o}}^{\SQPsubscript}(t)$ from the desired trajectory. However, for larger fleets ($n = 400$), the computational delay becomes prohibitive and \SQP\ fails to control the system, resulting in unstable object motion.

In contrast, the trajectories generated by DCC remain indistinguishable from the offline optimal solution $\chiv^{\optimal}(t)$ across all tested fleet sizes. Consequently, the object's pose trajectory $\pv_{\text{o}}^{\OTCsubscript}(t)$ closely follows the desired trajectory even for large-scale systems.

These results demonstrate that while \SQP\ may remain suitable for small systems, its computational requirements prevent reliable real-time control as system dimensionality increases. Such failures can lead to severe safety risks, including aerial vehicle collisions or unintended load release. By contrast, DCC maintains stable and accurate control across all tested fleet sizes, illustrating its ability to overcome the real-time barrier in high-dimensional autonomous systems.

\subsection{Biomechanics monitoring (R2)}
\begin{figure}[!t]
	\centering
	\includegraphics[scale=0.817]{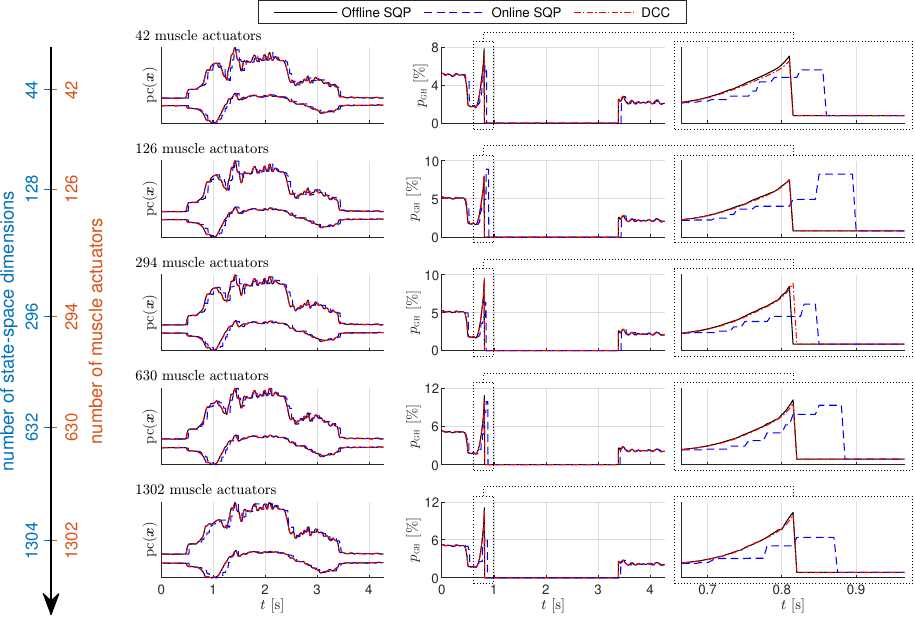}\vspace{-1ex}
	\caption{\label{fig:ulmsk tracking timeseries}%
		The muscle activation levels $\xv(t)$ and the glenohumeral joint capsule tear probability $p_{\scriptscriptstyle\textup{GH}}(t)$ estimated by \SQP~(blue) and DCC~(red) on an emulated real-time computer. 
		The offline computed optimal solution $\chiv^{\optimal}(t)$ (black) serves as the reference.
		Here, $\principal(\xv)$ denotes the first two \uls{p}rincipal \uls{c}omponents of $\xv$.
	}
\end{figure}

Due to practical limitations in measuring muscle activations and forces directly, model-based estimation is commonly used to monitor musculoskeletal conditions~\cite{ErdemirMcLHerBog2007}. Such monitoring relies on solving the \emph{muscle-force distribution problem}, which determines the muscle activations and ligament forces--within their biomechanical limits--that produce a desired multijoint limb motion while resisting possible external forces~\cite{ErdemirMcLHerBog2007,Herzog2009a,HuKueHad2020}.

To evaluate the proposed formulation in this context, we establish model-based monitors that estimate the probability $p_{\scriptscriptstyle\textup{GH}}(t)$ of a glenohumeral joint capsule tear. Each monitor integrates a shoulder-arm musculoskeletal model~\cite{HuKueHad2020} with a different number of actuators per muscle, thereby increasing the dimensionality of the optimization problem. This setup provides a representative benchmark for assessing whether DCC can compute muscle-force distributions within the time scale required for real-time biomechanical monitoring.

Figure~\ref{fig:ulmsk tracking timeseries} shows the resulting monitoring trajectories. When the estimation problem is solved using online \SQP, a noticeable delay arises between the computed probability trajectory and the offline optimal solution. As a consequence, the \SQP-based monitor fails to capture rapid variations in $p_{\scriptscriptstyle\textup{GH}}(t)$, missing short but critical peaks of injury risk.

In contrast, the DCC-based monitor exhibits negligible delay and accurately tracks all temporal changes in $p_{\scriptscriptstyle\textup{GH}}(t)$. The resulting probability estimates therefore remain closely aligned with the offline optimal solution throughout the experiment.

This improvement is essential for applications such as online monitoring during training or rehabilitation, where preventive interventions depend on timely detection of injury risk. For example, an intelligent exercise machine could adapt the applied load according to $p_{\scriptscriptstyle\textup{GH}}(t)$ in real time, thereby reducing the probability of joint injury. These results demonstrate that DCC enables reliable real-time biomechanical monitoring even for high-dimensional musculoskeletal models.

\subsection{Optimal economic decision-making (R3)}
\begin{figure}[!t]
	\centering
	\includegraphics[scale=0.817]{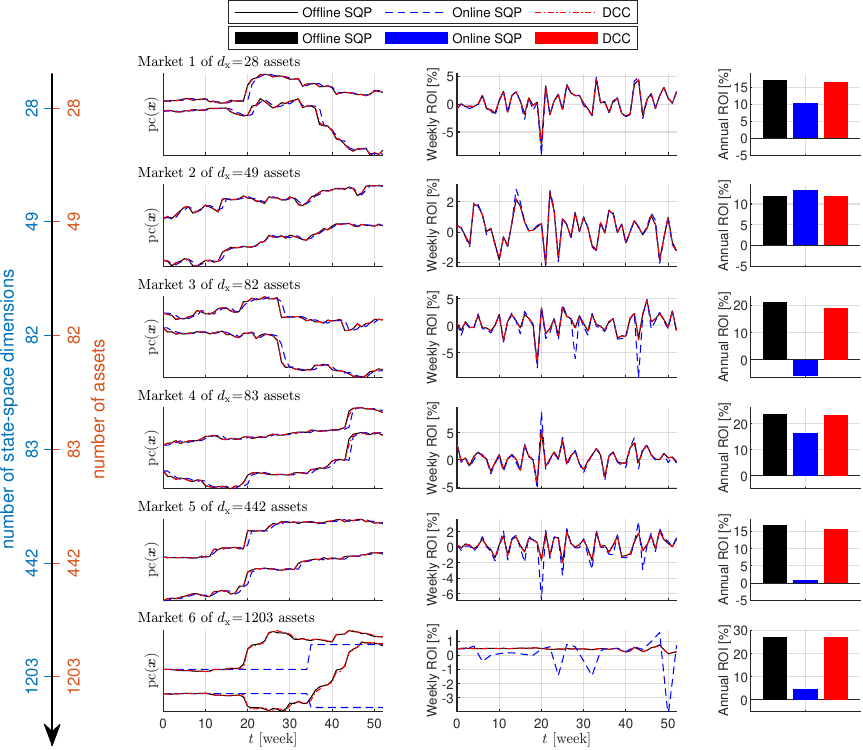}
	\caption{\label{fig:portfolio tracking timeseries}%
		The portfolio $\xv(t)$, weekly \uls{r}eturn \uls{o}n \uls{i}nvestment (ROI), and annual ROI resulting from \SQP~(blue) and DCC~(red) on an emulated real-time computer. 
		The offline computed optimal solution $\chiv^{\optimal}(t)$ (black) serves as the reference.
		Here, $\principal(\xv)$ denotes the first two \uls{p}rincipal \uls{c}omponents of $\xv$.
	}
\end{figure}

As one of the fundamental frameworks for portfolio optimization, the \emph{Markowitz mean-variance problem} determines the optimal allocation of assets by balancing the trade-off between expected return and risk, quantified by the mean and covariance of historical market data~\cite{Steinbach2001}. 

To evaluate the proposed formulation in this setting, we establish numerical market scenarios reflecting real-world financial systems~\cite{BruniCesScoTar2016}. In these scenarios, the portfolio allocation must be continuously updated based on evolving market conditions, while the dimensionality of the optimization problem increases with the number of assets. This setup therefore provides a representative benchmark for evaluating whether DCC can compute optimal portfolio allocations within the time scale required for real-time financial decision-making.

Figure~\ref{fig:portfolio tracking timeseries} illustrates the resulting investment trajectories. Across all market scenarios, the portfolio states generated by DCC, $\xv^{\OTCsubscript}(t)$, remain nearly indistinguishable from the offline optimal solution $\chiv^{\optimal}(t)$. In contrast, the trajectories obtained with online \SQP\ exhibit a noticeable time delay relative to the optimal solution. For very large markets with $\dimx = 1203$ assets, \SQP\ fails completely due to excessive computational delay.

Since Markowitz optimization relies on \emph{past} market data while investment returns depend on \emph{present} market conditions, any computational delay can lead to missed profitable opportunities. As illustrated in Fig.~\ref{fig:portfolio tracking timeseries}, such delays accumulate into substantial deviations in long-term investment performance, as observed for \SQP\ in market scenario~3.

By contrast, DCC closely tracks the Markowitz-optimal weekly return on investment (ROI), and the resulting annual ROI remains consistent with the theoretical optimum. These results demonstrate that DCC enables real-time portfolio optimization even for large-scale financial markets.

% As one of the fundamental frameworks for portfolio optimization,
% \emph{Markowitz mean-variance problem} seeks the optimal allocation of assets by balancing the trade-off between expected return and risk, as measured by mean and covariance, respectively, of the historical market data~\cite{Steinbach2001}.
% In this example, we establish numerical scenarios mirroring real-world markets~\cite{BruniCesScoTar2016}.

% Figure~\ref{fig:portfolio tracking timeseries} depicts the simulation results. 
% Across all market scenarios, there is no noticeable difference between $\xv^{\OTCsubscript}(t)$ and the optimal solution $\chiv^{\optimal}(t)$, whereas a visible time delay appears between $\xv^{\SQPsubscript}(t)$ and $\chiv^{\optimal}(t)$, and \SQP\ fails completely on the very large markets with $\dimx = 1203$ assets.
% Since Markowitz optimization inherently relies on \emph{past} market data and the actual Return on Investment (ROI) is governed by \emph{present} market conditions, any computational time delay implies missed profitable opportunities, which can accumulate into long-term investment failures--as observed with \SQP\ on market 3. In contrast, DCC tracks the Markowitz-optimal weekly ROI closely, and the  annual ROI matches closely the theoretical optimum as well.
%

\section{Discussion}\label{sec:discussion}

This work introduces Dual Cost--Constraint projection (DCC), a closed-form dynamical representation that generates optimal control actions directly within physical time. The results demonstrate that DCC maintains optimality, constraint satisfaction, and tracking accuracy even for high-dimensional nonlinear systems, thereby overcoming the real-time barrier that limits iterative optimization algorithms.

These results suggest that certain classes of optimal control problems may be solved not by iterative optimization but by appropriately constructed dynamical systems.

DCC applies to a broad class of optimal control problems involving convex and certain pseudoconvex objectives for time-variant control-affine systems with continuously differentiable dynamics. Such systems arise across many domains, including neuromechanics (muscle-force distribution~\cite{ErdemirMcLHerBog2007}, motor-unit recruitment~\cite{MarshallGlaTraAme2022}), metabolic networks (flux balance analysis~\cite{BeckerFeiMoHan2007,OrthThiPal2010}, multidimensional optimality~\cite{SchuetzZamZamHei2012}), communication systems (beamforming~\cite{GershmanSidShaBen2010,LauZha2015}), autonomous systems~\cite{AllenspachBodBruRin2020,JaleelSha2020}, electricity grids~\cite{AbdiBeiSca2017,WeiWanLiMei2017,JhaInaBisSur2023}, and financial optimization~\cite{Markowitz1952,Steinbach2001,SandhuGeoTan2016}. The breadth of this system and problem class highlights the general applicability of DCC.

The compact closed-form dynamical representation of DCC enables optimal control actions to be computed directly within each sampling interval. As demonstrated in Section~2, DCC operates within the available real-time budget even for high-dimensional nonlinear systems, whereas iterative optimization algorithms exceed the computational limits by several orders of magnitude. While existing approaches can achieve real-time performance for low-dimensional linear systems~\cite{BogertGeiEveSte2013,BodieTayKamSie2020}, the capability of DCC extends to systems with high dimensionality, nonlinear dynamics, and complex optimality principles. 

These results suggest that certain classes of optimal control problems can be solved not by iterative optimization, but by appropriately constructed dynamical systems.

Model predictive control (MPC) is currently the dominant framework for optimal decision-making in time-critical systems~\cite{MayneRawRaoSco2000}. However, MPC requires solving a constrained optimal control problem at every sampling instant using iterative numerical optimization. Replacing the iterative solver with DCC would allow MPC to operate at higher sampling frequencies and longer prediction horizons, bringing system behavior closer to the true optimum. Such improvements could increase efficiency and sustainability in industrial applications including process engineering, climate control, and manufacturing~\cite{SchwenzerAyBerAbe2021}. DCC-enhanced MPC may also enable real-time on-board decision-making in spacecraft operations~\cite{ErenPraKocRak2017,IzzoBlaFerOri2024}, where rapid responses to unpredictable environmental changes are essential.

Beyond its practical utility, DCC provides a dynamical interpretation of constrained optimization. Because the formulation is expressed directly as a causal dynamical system, it offers a new perspective on optimality principles in natural and engineered processes. This viewpoint suggests that DCC may serve as a useful model in grey-box identification, potentially enabling the discovery of unknown optimality principles in biological or physical systems such as neuromechanical coordination strategies.

A complementary interpretation reveals that DCC may also generalize a wide spectrum of control approaches to \emph{state-constrained controllers} (see Supplementary Information, Section~\ref{SI-sec:constrained system}). In particular, extending classical optimal control concepts such as the linear–quadratic regulator~\cite{Mehrmann1991} to nonlinear state-constrained systems suggests that DCC may function simultaneously as an optimal trajectory planner and feedback controller for nonlinear time-variant systems with convex cost functions.

Despite these advantages, DCC currently applies to a specific class of optimal control problems, and several extensions remain open. Hybrid dynamical systems introduce discontinuities through mechanical contacts, impacts, and friction~\cite{Cortes2008}. Control-affinity does not hold in certain systems where the number of inputs is smaller than the number of states, including many networked systems~\cite{LiuSloBar2011}, although extensions to controllable non-affine systems may be possible (see Supplementary Information, Sections~\ref{SI-sec:remark on state equation} and~\ref{SI-sec:chain of integrators}). In addition, combinatorial optimization problems common in logistics~\cite{MataiSinMit2010}, human-resource allocation~\cite{BouajajaDri2017}, and synthetic biology~\cite{NaseriKof2020} involve discrete feasible sets and are therefore nonconvex.

Exploring these directions is essential for developing a more general causally deterministic theory of optimization. Recent advances addressing nonlinear inequality constraints~\cite{RaveendranMahVai2023} and control of non-affine systems~\cite{BinazadehRah2017,ChenBaiKon2024}, together with the present work, suggest a possible pathway toward real-time optimal decision-making in increasingly complex dynamical systems.

\section{Methods}\label{sec:methods}

This section defines the real-time barrier, describes the benchmarking setup and numerical examples, and summarizes the implementation of the Dual Cost--Constraint Projection (DCC).
Detailed derivations, proofs, and extended model descriptions are provided in the Supplementary Information.

\subsection{Real-time barrier}
Real-time computation requires that the control law be evaluated before the physical system advances to the next sampling instant:
\begin{equation}\label{eq:real-time condition}
	T_{\CPUsubscript,\FLOPsubscript} \times N_{\FLOPsubscript} \leqslant \dt
\end{equation}
Here, $\dt$ is the desired sampling interval,
$N_{\FLOPsubscript}$ the number of \uls{fl}oating-point \uls{op}erations (FLOP) required by an algorithm for each time step,
and $T_{\CPUsubscript,\FLOPsubscript}$ the \uls{C}PU time per floating-point operation. Its reciprocal $T_{\CPUsubscript,\FLOPsubscript}^{-1}$ defines the available computational throughput (FLOP/s).

Real-time feasibility therefore requires that the computational demand of the algorithm does not exceed the available throughput of the computing platform.

Inserting $N_{\FLOPsubscript}^{\OTCsubscript}(\dimx)$ and $N_{\FLOPsubscript}^{\SQPsubscript}(\dimx)$ into \eqref{eq:real-time condition} results in \eqref{eq:OTC flop/s}.
Note that \eqref{eq:OTC flop/s} is an approximation for sufficiently large $\dimx \gg 1$.
The reader may refer to Supplementary Information, Section~\ref{SI-sec:computational complexity} for a detailed deduction and exact expressions.

\subsection{Benchmarking setup}

To evaluate DCC, we consider one mathematical benchmark and three application scenarios. These examples span different domains but share the same underlying structure of high-dimensional constrained optimal control.

\paragraph{Example M1}
The purely mathematical example is a lowest-dimensional representative of the problem class \eqref{eq:gryphon}--\eqref{eq:optimal control problem} with high degree of nonlinearity and suitable time-variance.

\paragraph{Example R1: aerial vehicle fleet control}
A cooperative transport task involving multiple aerial vehicles, where the key challenge is real-time optimal control allocation.

\paragraph{Example R2: musculoskeletal monitoring}
A system for monitoring the injury risk of glenohumeral joint, with the challenage of estimating coordinated muscle-force distribution in real-time.

\paragraph{Example R3: portfolio optimization}
An online investment decision-making process that selects profitable investments with minimal risk under time-critical market conditions.

\par~\par

All methods are evaluated under identical computational conditions without parallelization to ensure a fair comparison of real-time performance. DCC (with gain $K_x = 500$) is compared against state-of-the-art SQP~\cite{TheMathWorks2020}. All simulations were implemented in MATLAB R2020b on a workstation with an Intel i9-10900K CPU
and 64\,GB DDR4 RAM.
Ordinary differential equations were solved numerically using the classical Runge--Kutta method (RK4).
To maximize execution efficiency, MATLAB executable (MEX) functions were compiled for both DCC and SQP.

\subsection{Numerical examples}

\paragraph{Example R1: aerial vehicle fleet control}
Example R1 considers cooperative aerial transport using $n$ omnidirectional aerial vehicles.
The effective actuator state dimension satisfies $\dim(\xv) = 24n$.
The task is to generate the wrench required to move a payload along a predefined minimum-jerk trajectory
while satisfying actuator and contact constraints.
The quadratic cost function is used to distribute actuator effort. The results are shown in Figure~\ref{fig:swarm tracking timeseries}.

\paragraph{Example R2: musculoskeletal monitoring}
Example R2 uses a human shoulder--arm musculoskeletal model with 12 joints and 42 muscles. Each muscle is subdivided into multiple actuators, increasing the system dimension. The optimization reconstructs muscle forces consistent with observed limb motion while respecting biomechanical constraints. The resulting force distribution is used to estimate injury risk of the glenohumeral joint. The results are shown in Figure~\ref{fig:ulmsk tracking timeseries}.

\paragraph{Example R3: portfolio optimization}
Example R3 considers the classical mean--variance portfolio problem. Weekly financial market data from several markets are used to estimate time-varying
expected returns and covariance matrices. The optimization minimizes portfolio variance subject to a prescribed expected return and allocation constraints.
The results are shown in Figure~\ref{fig:portfolio tracking timeseries}.

\par~\par

Further modeling details and full parameter definitions are provided in the Supplementary Information, Section~\ref{SI-sec:extended_numerical_examples}.

\subsection{Real-time performance metrics}

To quantify real-time performance, we evaluate deviations from the offline optimal solution using the three metrics 
(i) \emph{cost-function error}, (ii) \emph{task-constraint violation}, and (iii) \emph{state tracking error}, see Supplementary Equations~\eqref{SI-eq:cost-function error}--\eqref{SI-eq:state tracking error}. 
All metrics are computed relative to the offline optimal solution. The offline optimal reference solution is obtained using SQP without computational limits. The results are shown in Figure~\ref{fig:portfolio tracking timeseries}.

\subsection{Dual cost--constraint projection}

DCC is derived by connecting numerical optimization with nonlinear control theory. Starting from formulating a Lagrangian function using barrier function and Lagrange multipliers, we {improved the interior-point optimization iteration}, which is then reformulated to DCC, a closed-loop dynamical controller with guaranteed optimal behavior.

\begin{figure}[!t]
	\centering
	\includegraphics[scale=1]{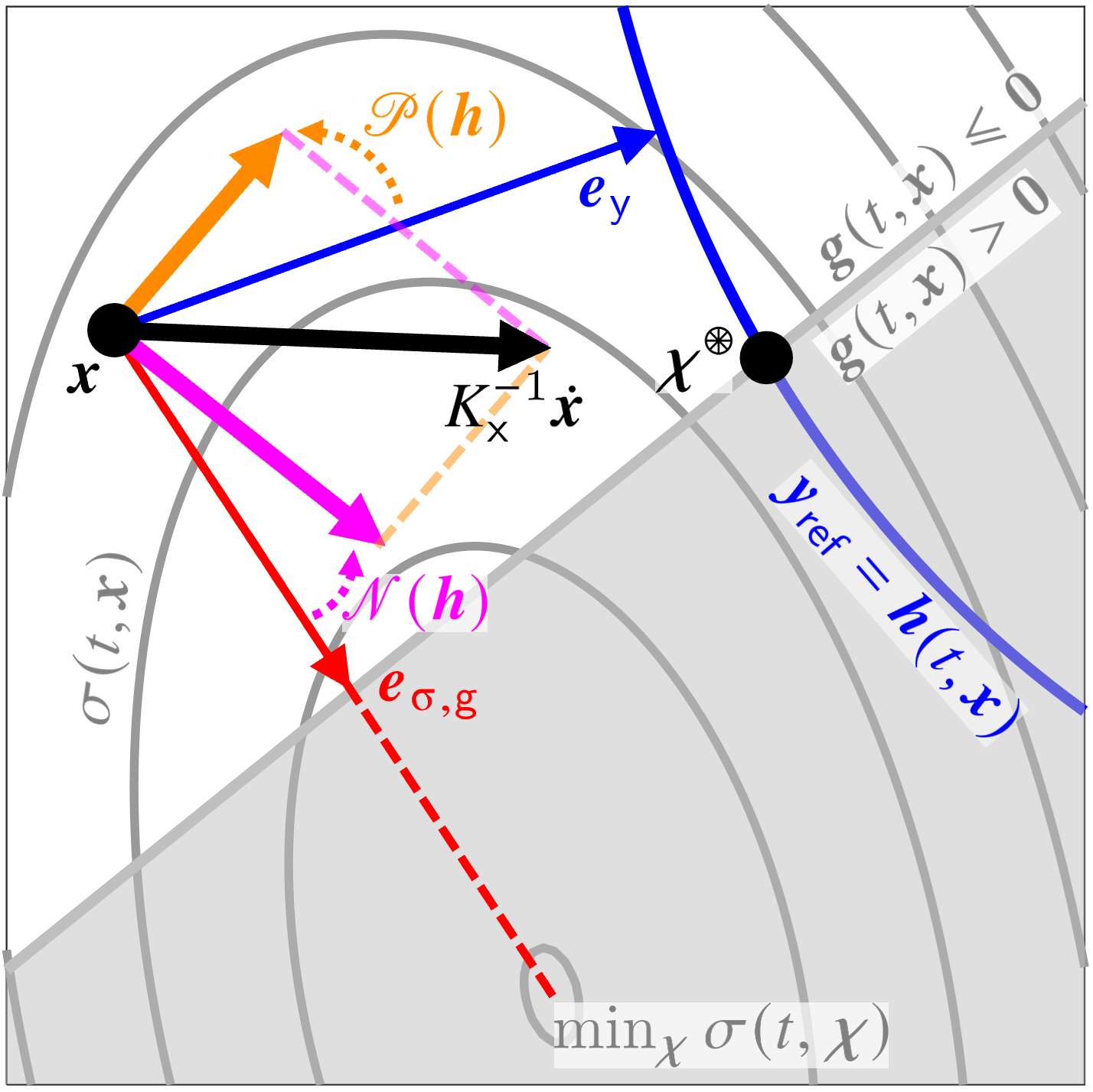}
	\caption{\label{fig:geometric interpretation}%
		\textbf{Geometric interpretation of Dual Cost--Constraint Projection.}
        The error vectors $\ev_{\textup{y}}$ and $\ev_{\upsigma,\textup{g}}$ are projected by the dual projector $\begin{bmatrix}
        \mathscr{P}\{\hv\}(t,\xv) & \mathscr{N}\{\hv\}(t,\xv)
        \end{bmatrix}$
        to an orthogonal pair (colored orange and magenta).
        The orthogonal pair creates a direction (black, $K_{\textup{x}}^{-1}\dot{\xv}$) for the state $\xv$ moving toward the optimal solution $\chiv^{\optimal}$.
	}
\end{figure}

The resulting controller contains two complementary components:
\begin{enumerate}
    \item an output-tracking component enforcing task objectives ($\ev_{\textup{y}} \rightarrow \nullv$, see Fig.~\ref{fig:geometric interpretation}),
    \item a null-space component minimizing the inequality constraint violation and the cost function $\ev_{\upsigma,\textup{g}} \rightarrow \nullv$).
\end{enumerate}

Figure~\ref{fig:geometric interpretation} illustrates a geometric interpretation of DCC.
The dual projector 
$\begin{bmatrix}
    \mathscr{P}\{\hv\}(t,\xv) & \mathscr{N}\{\hv\}(t,\xv)
\end{bmatrix}$
ensures that the two vectors (colored orange and magenta) are always orthogonal to each other.

This formulation allows optimal control actions to be generated by the closed-loop system dynamics, thereby replacing iterative optimization with a dynamical system that produces optimal behavior in physical time.

Formal derivations, lemmata, theorems, proofs--including full stability and convergence analysis and results--and computational complexity analysis supporting DCC are provided in the Supplementary Information (Sections \ref{SI-sec:lemmata_theorems_proofs}--\ref{SI-sec:computational complexity}).

%

%\backmatter
\clearpage

\section*{Supplementary information}
Mathematical deduction and derivation of DCC,
analysis of DCC's convergence behavior and computational complexity, real-time barrier analysis, remarks on barrier function and state equation, DCC's orthogonality, extended description of numerical examples, conjectures on generalizations of DCC.

\section*{Code availability}
The MATLAB codes (*.m, *.mexw64) of all numerical examples, including M1, R1, R2, R3, and the ones supporting our conjectures may be provided upon reasonable requests.
%We provide the MATLAB codes of all numerical examples, including M1, R1, R2, R3, and the ones supporting our conjectures.
%These codes are provided in the supplementary material.

\section*{Data availability}
All numerical simulation results (as MATLAB *.mat files) are accessible for download at \cite{HuHad2026}.

\section*{Acknowledgements}
The authors gratefully acknowledge the funding of this work by the Alfried Krupp von Bohlen und Halbach Foundation, the Deutsche Forschungsgemeinschaft through the Gottfried Wilhelm Leibniz Programme (award to S. Haddadin; grant no. HA7372/3-1), the European Union's Horizon 2020 research and innovation programme as part of the project SoftPro (grant no. 688857),
the Lighthouse Initiative Geriatronics by StMWi Bayern (Project X, grant no. 5140951),
and LongLeif GaPa gGmbH (Project Y, grant no. 5140953), as well as the German Federal Ministry of Education and Research (BMBF) funding as part of the project AI.D (grant no. 16ME0539K).%

\section*{Author contributions}
T. Hu and S. Haddadin developed the theory and methods. 
T. Hu developed and conducted the numerical experiments. 
Both authors interpreted the results. 
Both authors conceptualized and wrote the manuscript. Both authors read and approved the final paper.

\newpage

\clearpage
\bibliographystyle{naturemag} %\bibliographystyle{amsplain} % %
\bibliography{ref_hu}

\begin{thebibliography}{10}
\expandafter\ifx\csname url\endcsname\relax
  \def\url#1{\texttt{#1}}\fi
\expandafter\ifx\csname urlprefix\endcsname\relax\def\urlprefix{URL }\fi
\providecommand{\bibinfo}[2]{#2}
\providecommand{\eprint}[2][]{\url{#2}}

\bibitem{ErdemirMcLHerBog2007}
\bibinfo{author}{Erdemir, A.}, \bibinfo{author}{McLean, S.},
  \bibinfo{author}{Herzog, W.} \& \bibinfo{author}{van~den Bogert, A.~J.}
\newblock \bibinfo{title}{Model-based estimation of muscle force exerted during
  movements}.
\newblock \emph{\bibinfo{journal}{Clin. Biomech.}}
  \textbf{\bibinfo{volume}{22}}, \bibinfo{pages}{131--154}
  (\bibinfo{year}{2007}).

\bibitem{MarshallGlaTraAme2022}
\bibinfo{author}{Marshall, N.} \emph{et~al.}
\newblock \bibinfo{title}{Flexible neural control of motor units}.
\newblock \emph{\bibinfo{journal}{Nature Neuroscience}}
  \textbf{\bibinfo{volume}{25}}, \bibinfo{pages}{1--13} (\bibinfo{year}{2022}).

\bibitem{OrthThiPal2010}
\bibinfo{author}{Orth, J.}, \bibinfo{author}{Thiele, I.} \&
  \bibinfo{author}{Palsson, B.}
\newblock \bibinfo{title}{What is flux balance analysis?}
\newblock \emph{\bibinfo{journal}{Nature biotechnology}}
  \textbf{\bibinfo{volume}{28}}, \bibinfo{pages}{245--8}
  (\bibinfo{year}{2010}).

\bibitem{SchuetzZamZamHei2012}
\bibinfo{author}{Schuetz, R.}, \bibinfo{author}{Zamboni, N.},
  \bibinfo{author}{Zampieri, M.}, \bibinfo{author}{Heinemann, M.} \&
  \bibinfo{author}{Sauer, U.}
\newblock \bibinfo{title}{Multidimensional optimality of microbial metabolism}.
\newblock \emph{\bibinfo{journal}{Science}} \textbf{\bibinfo{volume}{336}},
  \bibinfo{pages}{601--604} (\bibinfo{year}{2012}).

\bibitem{NaseriKof2020}
\bibinfo{author}{Naseri, G.} \& \bibinfo{author}{Koffas, M.}
\newblock \bibinfo{title}{Application of combinatorial optimization strategies
  in synthetic biology}.
\newblock \emph{\bibinfo{journal}{Nature Communications}}
  \textbf{\bibinfo{volume}{11}} (\bibinfo{year}{2020}).

\bibitem{AllenspachBodBruRin2020}
\bibinfo{author}{Allenspach, M.} \emph{et~al.}
\newblock \bibinfo{title}{Design and optimal control of a tiltrotor
  micro-aerial vehicle for efficient omnidirectional flight}.
\newblock \emph{\bibinfo{journal}{The International Journal of Robotics
  Research}} \textbf{\bibinfo{volume}{39}}, \bibinfo{pages}{1305--1325}
  (\bibinfo{year}{2020}).

\bibitem{JaleelSha2020}
\bibinfo{author}{Jaleel, H.} \& \bibinfo{author}{Shamma, J.~S.}
\newblock \bibinfo{title}{Distributed optimization for robot networks: From
  real-time convex optimization to game-theoretic self-organization}.
\newblock \emph{\bibinfo{journal}{Proceedings of the IEEE}}
  \textbf{\bibinfo{volume}{108}}, \bibinfo{pages}{1953--1967}
  (\bibinfo{year}{2020}).

\bibitem{IzzoBlaFerOri2024}
\bibinfo{author}{Izzo, D.} \emph{et~al.}
\newblock \bibinfo{title}{Optimality principles in spacecraft neural guidance
  and control}.
\newblock \emph{\bibinfo{journal}{Science Robotics}}
  \textbf{\bibinfo{volume}{9}}, \bibinfo{pages}{eadi6421}
  (\bibinfo{year}{2024}).

\bibitem{Steinbach2001}
\bibinfo{author}{Steinbach, M.}
\newblock \bibinfo{title}{Markowitz revisited: Mean-variance models in
  financial portfolio analysis}.
\newblock \emph{\bibinfo{journal}{Society for Industrial and Applied
  Mathematics}} \textbf{\bibinfo{volume}{43}}, \bibinfo{pages}{31--85}
  (\bibinfo{year}{2001}).

\bibitem{FrankReb2016}
\bibinfo{author}{Frank, S.} \& \bibinfo{author}{Rebennack, S.}
\newblock \bibinfo{title}{An introduction to optimal power flow: Theory,
  formulation, and examples}.
\newblock \emph{\bibinfo{journal}{IIE Transactions}}
  \textbf{\bibinfo{volume}{48}}, \bibinfo{pages}{1172--1197}
  (\bibinfo{year}{2016}).
\newblock \urlprefix\url{https://doi.org/10.1080/0740817X.2016.1189626}.
\newblock \eprint{https://doi.org/10.1080/0740817X.2016.1189626}.

\bibitem{AbdiBeiSca2017}
\bibinfo{author}{Abdi, H.}, \bibinfo{author}{Beigvand, S.~D.} \&
  \bibinfo{author}{Scala, M.~L.}
\newblock \bibinfo{title}{A review of optimal power flow studies applied to
  smart grids and microgrids}.
\newblock \emph{\bibinfo{journal}{Renewable and Sustainable Energy Reviews}}
  \textbf{\bibinfo{volume}{71}}, \bibinfo{pages}{742--766}
  (\bibinfo{year}{2017}).
\newblock
  \urlprefix\url{https://www.sciencedirect.com/science/article/pii/S1364032116311583}.

\bibitem{WeiWanLiMei2017}
\bibinfo{author}{Wei, W.}, \bibinfo{author}{Wang, J.}, \bibinfo{author}{Li, N.}
  \& \bibinfo{author}{Mei, S.}
\newblock \bibinfo{title}{Optimal power flow of radial networks and its
  variations: A sequential convex optimization approach}.
\newblock \emph{\bibinfo{journal}{IEEE Transactions on Smart Grid}}
  \textbf{\bibinfo{volume}{8}}, \bibinfo{pages}{2974--2987}
  (\bibinfo{year}{2017}).

\bibitem{JhaInaBisSur2023}
\bibinfo{author}{Jha, R.~R.} \emph{et~al.}
\newblock \bibinfo{title}{Distribution grid optimal power flow {(D-OPF)}:
  Modeling, analysis, and benchmarking}.
\newblock \emph{\bibinfo{journal}{IEEE Transactions on Power Systems}}
  \textbf{\bibinfo{volume}{38}}, \bibinfo{pages}{3654--3668}
  (\bibinfo{year}{2023}).

\bibitem{GershmanSidShaBen2010}
\bibinfo{author}{Gershman, A.}, \bibinfo{author}{Sidiropoulos, N.},
  \bibinfo{author}{Shahbazpanahi, S.}, \bibinfo{author}{Bengtsson, M.} \&
  \bibinfo{author}{Ottersten, B.}
\newblock \bibinfo{title}{Convex optimization-based beamforming}.
\newblock \emph{\bibinfo{journal}{IEEE Signal Processing Magazine}}
  \textbf{\bibinfo{volume}{27}}, \bibinfo{pages}{62 -- 75}
  (\bibinfo{year}{2010}).

\bibitem{MataiSinMit2010}
\bibinfo{author}{Matai, R.}, \bibinfo{author}{Singh, S.} \&
  \bibinfo{author}{Mittal, M.~L.}
\newblock \bibinfo{title}{Traveling salesman problem: an overview of
  applications, formulations, and solution approaches}.
\newblock In \bibinfo{editor}{Davendra, D.} (ed.)
  \emph{\bibinfo{booktitle}{Traveling Salesman Problem}},
  chap.~\bibinfo{chapter}{1} (\bibinfo{publisher}{IntechOpen},
  \bibinfo{address}{Rijeka}, \bibinfo{year}{2010}).

\bibitem{BouajajaDri2017}
\bibinfo{author}{Bouajaja, S.} \& \bibinfo{author}{Dridi, N.}
\newblock \bibinfo{title}{A survey on human resource allocation problem and its
  applications}.
\newblock \emph{\bibinfo{journal}{Operational Research}}
  \textbf{\bibinfo{volume}{17}}, \bibinfo{pages}{339--369}
  (\bibinfo{year}{2017}).

\bibitem{Fletcher1987}
\bibinfo{author}{Fletcher, R.}
\newblock \emph{\bibinfo{title}{Practical Methods of Optimization}}
  (\bibinfo{publisher}{John Wiley \& Sons}, \bibinfo{address}{Chichester, West
  Sussex, England}, \bibinfo{year}{1987}), \bibinfo{edition}{second} edn.

\bibitem{NocedalWri2006}
\bibinfo{author}{Nocedal, J.} \& \bibinfo{author}{Wright, S.~J.}
\newblock \emph{\bibinfo{title}{Numerical optimization}}
  (\bibinfo{publisher}{Springer Science+Business Media, LLC.},
  \bibinfo{address}{New York, NY, USA}, \bibinfo{year}{2006}).

\bibitem{PotraWri2000}
\bibinfo{author}{Potra, S.~J., F. A.;~Wright}.
\newblock \bibinfo{title}{Interior-point methods}.
\newblock \emph{\bibinfo{journal}{J. Comput. Appl. Math.}}
  \textbf{\bibinfo{volume}{124}}, \bibinfo{pages}{281--302}
  (\bibinfo{year}{2000}).

\bibitem{BoydVan2004}
\bibinfo{author}{Boyd, S.} \& \bibinfo{author}{Vandenberghe, L.}
\newblock \emph{\bibinfo{title}{Convex Optimization}}
  (\bibinfo{publisher}{{Cambridge University Press}},
  \bibinfo{address}{Cambridge, UK}, \bibinfo{year}{2004}).

\bibitem{GillMurSau2005}
\bibinfo{author}{Gill, P.~E.}, \bibinfo{author}{Murray, W.} \&
  \bibinfo{author}{Saunders, M.~A.}
\newblock \bibinfo{title}{Snopt: An sqp algorithm for large-scale constrained
  optimization}.
\newblock \emph{\bibinfo{journal}{SIAM Review}} \textbf{\bibinfo{volume}{47}},
  \bibinfo{pages}{99--131} (\bibinfo{year}{2005}).
\newblock \urlprefix\url{http://www.jstor.org/stable/20453604}.

\bibitem{BueskensWas2013}
\bibinfo{author}{B{\"u}skens, C.} \& \bibinfo{author}{Wassel, D.}
\newblock \emph{\bibinfo{title}{The ESA NLP Solver WORHP}},
  \bibinfo{pages}{85--110} (\bibinfo{publisher}{Springer New York},
  \bibinfo{address}{New York, NY}, \bibinfo{year}{2013}).
\newblock \urlprefix\url{https://doi.org/10.1007/978-1-4614-4469-5_4}.

\bibitem{MayneRawRaoSco2000}
\bibinfo{author}{Mayne, D.}, \bibinfo{author}{Rawlings, J.},
  \bibinfo{author}{Rao, C.} \& \bibinfo{author}{Scokaert, P.}
\newblock \bibinfo{title}{Constrained model predictive control: Stability and
  optimality}.
\newblock \emph{\bibinfo{journal}{Automatica}} \textbf{\bibinfo{volume}{36}},
  \bibinfo{pages}{789--814} (\bibinfo{year}{2000}).
\newblock
  \urlprefix\url{https://www.sciencedirect.com/science/article/pii/S0005109899002149}.

\bibitem{DrgonaArrCupBlu2020}
\bibinfo{author}{Drgoňa, J.} \emph{et~al.}
\newblock \bibinfo{title}{All you need to know about model predictive control
  for buildings}.
\newblock \emph{\bibinfo{journal}{Annual Reviews in Control}}
  \textbf{\bibinfo{volume}{50}}, \bibinfo{pages}{190--232}
  (\bibinfo{year}{2020}).
\newblock
  \urlprefix\url{https://www.sciencedirect.com/science/article/pii/S1367578820300584}.

\bibitem{Mehrmann1991}
\bibinfo{author}{Mehrmann, V.}
\newblock \emph{\bibinfo{title}{The Autonomous Linear Quadratic Control
  Problem: Theory and Numerical Solution}} (\bibinfo{publisher}{Springer-Verlag
  Berlin}, \bibinfo{address}{Heidelberg, Germany}, \bibinfo{year}{1991}).

\bibitem{Yamashita1980}
\bibinfo{author}{Yamashita, H.}
\newblock \bibinfo{title}{A differential equation approach to nonlinear
  programming}.
\newblock \emph{\bibinfo{journal}{Mathematical Programming}}
  \textbf{\bibinfo{volume}{18}}, \bibinfo{pages}{155--168}
  (\bibinfo{year}{1980}).

\bibitem{ZhouShi1997}
\bibinfo{author}{Zhou, Z.} \& \bibinfo{author}{Shi, Y.}
\newblock \bibinfo{title}{An {ODE} method of solving nonlinear programming}.
\newblock \emph{\bibinfo{journal}{Computers \& Mathematics with Applications}}
  \textbf{\bibinfo{volume}{34}}, \bibinfo{pages}{97--102}
  (\bibinfo{year}{1997}).

\bibitem{XiaFen2005}
\bibinfo{author}{Xia, Y.} \& \bibinfo{author}{Feng, G.}
\newblock \bibinfo{title}{An improved neural network for convex quadratic
  optimization with application to real-time beamforming}.
\newblock \emph{\bibinfo{journal}{Neurocomputing}}
  \textbf{\bibinfo{volume}{64}}, \bibinfo{pages}{359--374}
  (\bibinfo{year}{2005}).

\bibitem{FepponAllDap2020}
\bibinfo{author}{Feppon, F.}, \bibinfo{author}{Allaire, G.} \&
  \bibinfo{author}{Dapogny, C.}
\newblock \bibinfo{title}{{Null space gradient flows for constrained
  optimization with applications to shape optimization}}.
\newblock \emph{\bibinfo{journal}{{ESAIM: Control, Optimisation and Calculus of
  Variations}}} \textbf{\bibinfo{volume}{26}}, \bibinfo{pages}{90}
  (\bibinfo{year}{2020}).

\bibitem{RaveendranMahVai2023}
\bibinfo{author}{Raveendran, R.}, \bibinfo{author}{Mahindrakar, A.~D.} \&
  \bibinfo{author}{Vaidya, U.}
\newblock \bibinfo{title}{Dynamical system approach for time-varying
  constrained convex optimization problems}.
\newblock \emph{\bibinfo{journal}{IEEE Transactions on Automatic Control}}
  \textbf{\bibinfo{volume}{69}}, \bibinfo{pages}{3822--3834}
  (\bibinfo{year}{2024}).

\bibitem{AllibhoyCor2023}
\bibinfo{author}{Allibhoy, A.} \& \bibinfo{author}{Cortés, J.}
\newblock \bibinfo{title}{Control-barrier-function-based design of gradient
  flows for constrained nonlinear programming}.
\newblock \emph{\bibinfo{journal}{IEEE Transactions on Automatic Control}}
  \textbf{\bibinfo{volume}{69}}, \bibinfo{pages}{3499--3514}
  (\bibinfo{year}{2024}).

\bibitem{TheMathWorks2020}
 \emph{\bibinfo{title}{Optimization Toolbox{\texttrademark} User's Guide}}
  (\bibinfo{publisher}{The MathWorks, Inc.}, \bibinfo{year}{2020}),
  \bibinfo{edition}{{R2020b}} edn.

\bibitem{HuKueHad2020}
\bibinfo{author}{Hu, T.}, \bibinfo{author}{K\"{u}hn, J.} \&
  \bibinfo{author}{Haddadin, S.}
\newblock \bibinfo{title}{Forward and inverse dynamics modeling of human
  shoulder-arm musculoskeletal system with scapulothoracic constraint}.
\newblock \emph{\bibinfo{journal}{Computer Methods in Biomechanics and
  Biomedical Engineering}}  (\bibinfo{year}{2020}).

\bibitem{Markowitz1952}
\bibinfo{author}{Markowitz, H.}
\newblock \bibinfo{title}{Portfolio selection}.
\newblock \emph{\bibinfo{journal}{The Journal of Finance}}
  \textbf{\bibinfo{volume}{7}}, \bibinfo{pages}{77--91} (\bibinfo{year}{1952}).
\newblock \urlprefix\url{http://www.jstor.org/stable/2975974}.

\bibitem{SandhuGeoTan2016}
\bibinfo{author}{Sandhu, R.}, \bibinfo{author}{Georgiou, T.} \&
  \bibinfo{author}{Tannenbaum, A.}
\newblock \bibinfo{title}{Ricci curvature: An economic indicator for market
  fragility and systemic risk}.
\newblock \emph{\bibinfo{journal}{Science Advances}}
  \textbf{\bibinfo{volume}{2}}, \bibinfo{pages}{e1501495--e1501495}
  (\bibinfo{year}{2016}).

\bibitem{BeckerFeiMoHan2007}
\bibinfo{author}{Becker, S.~A.} \emph{et~al.}
\newblock \bibinfo{title}{Quantitative prediction of cellular metabolism with
  constraint-based models: the {COBRA Toolbox}}.
\newblock \emph{\bibinfo{journal}{Nature Protocols}}
  \textbf{\bibinfo{volume}{2}}, \bibinfo{pages}{727--738}
  (\bibinfo{year}{2007}).
\newblock \urlprefix\url{https://api.semanticscholar.org/CorpusID:5687582}.

\bibitem{LauZha2015}
\bibinfo{author}{Lau, V. K.~N.} \& \bibinfo{author}{Zhang, F.}
\newblock \bibinfo{title}{Optimal beamforming for video streaming in
  multiantenna interference networks via diffusion limit}.
\newblock \emph{\bibinfo{journal}{IEEE Transactions on Information Theory}}
  \textbf{\bibinfo{volume}{61}}, \bibinfo{pages}{1819--1841}
  (\bibinfo{year}{2015}).

\bibitem{KirchengastSteHor2018}
\bibinfo{author}{Kirchengast, M.}, \bibinfo{author}{Steinberger, M.} \&
  \bibinfo{author}{Horn, M.}
\newblock \bibinfo{title}{Control allocation under actuator saturation: An
  experimental evaluation}.
\newblock \emph{\bibinfo{journal}{IFAC-PapersOnLine}}
  \textbf{\bibinfo{volume}{51}}, \bibinfo{pages}{48--54}
  (\bibinfo{year}{2018}).
\newblock \bibinfo{note}{9th IFAC Symposium on Robust Control Design ROCOND
  2018}.

\bibitem{Herzog2009a}
\bibinfo{author}{Herzog, W.}
\newblock \bibinfo{title}{Distribution problem in biomechanics}.
\newblock In \bibinfo{editor}{Binder, M.~D.}, \bibinfo{editor}{Hirokawa, N.} \&
  \bibinfo{editor}{Windhorst, U.} (eds.) \emph{\bibinfo{booktitle}{Encyclopedia
  of Neuroscience}}, \bibinfo{pages}{983--985}
  (\bibinfo{publisher}{Springer-Verlag Berlin Heidelberg},
  \bibinfo{address}{Heidelberg, Germany}, \bibinfo{year}{2009}).

\bibitem{BruniCesScoTar2016}
\bibinfo{author}{Bruni, R.}, \bibinfo{author}{Cesarone, F.},
  \bibinfo{author}{Scozzari, A.} \& \bibinfo{author}{Tardella, F.}
\newblock \bibinfo{title}{Real-world datasets for portfolio selection and
  solutions of some stochastic dominance portfolio models}.
\newblock \emph{\bibinfo{journal}{Data in Brief}} \textbf{\bibinfo{volume}{8}},
  \bibinfo{pages}{858--862} (\bibinfo{year}{2016}).

\bibitem{BogertGeiEveSte2013}
\bibinfo{author}{van~den Bogert, A.~J.}, \bibinfo{author}{Geijtenbeek, T.},
  \bibinfo{author}{Even-Zohar, O.}, \bibinfo{author}{Steenbrink, F.} \&
  \bibinfo{author}{Hardin, E.~C.}
\newblock \bibinfo{title}{A real-time system for biomechanical analysis of
  human movement and muscle function}.
\newblock \emph{\bibinfo{journal}{Medical \& biological engineering \&
  computing}} \textbf{\bibinfo{volume}{51}}, \bibinfo{pages}{1069–1077}
  (\bibinfo{year}{2013}).

\bibitem{BodieTayKamSie2020}
\bibinfo{author}{Bodie, K.}, \bibinfo{author}{Taylor, Z.},
  \bibinfo{author}{Kamel, M.} \& \bibinfo{author}{Siegwart, R.}
\newblock \bibinfo{title}{Towards efficient full pose omnidirectionality with
  overactuated mavs}.
\newblock In \bibinfo{editor}{Xiao, J.}, \bibinfo{editor}{Kr{\"o}ger, T.} \&
  \bibinfo{editor}{Khatib, O.} (eds.) \emph{\bibinfo{booktitle}{Proceedings of
  the 2018 International Symposium on Experimental Robotics}},
  \bibinfo{pages}{85--95} (\bibinfo{publisher}{Springer International
  Publishing}, \bibinfo{address}{Cham}, \bibinfo{year}{2020}).

\bibitem{SchwenzerAyBerAbe2021}
\bibinfo{author}{Schwenzer, M.}, \bibinfo{author}{Ay, M.},
  \bibinfo{author}{Bergs, T.} \& \bibinfo{author}{Abel, D.}
\newblock \bibinfo{title}{Review on model predictive control: an engineering
  perspective}.
\newblock \emph{\bibinfo{journal}{The International Journal of Advanced
  Manufacturing Technology}} \textbf{\bibinfo{volume}{117}},
  \bibinfo{pages}{1327 -- 1349} (\bibinfo{year}{2021}).
\newblock \urlprefix\url{https://api.semanticscholar.org/CorpusID:238713285}.

\bibitem{ErenPraKocRak2017}
\bibinfo{author}{Eren, U.} \emph{et~al.}
\newblock \bibinfo{title}{Model predictive control in aerospace systems:
  Current state and opportunities}.
\newblock \emph{\bibinfo{journal}{Journal of Guidance, Control, and Dynamics}}
  \textbf{\bibinfo{volume}{40}}, \bibinfo{pages}{1541--1566}
  (\bibinfo{year}{2017}).

\bibitem{Cortes2008}
\bibinfo{author}{Cortes, J.}
\newblock \bibinfo{title}{Discontinuous dynamical systems}.
\newblock \emph{\bibinfo{journal}{IEEE Control Systems Magazine}}
  \textbf{\bibinfo{volume}{28}}, \bibinfo{pages}{36--73}
  (\bibinfo{year}{2008}).

\bibitem{LiuSloBar2011}
\bibinfo{author}{Liu, Y.-Y.}, \bibinfo{author}{Slotine, J.-J.} \&
  \bibinfo{author}{Barabasi, A.-L.}
\newblock \bibinfo{title}{Controllability of complex networks}.
\newblock \emph{\bibinfo{journal}{Nature}} \textbf{\bibinfo{volume}{473}},
  \bibinfo{pages}{167--73} (\bibinfo{year}{2011}).

\bibitem{BinazadehRah2017}
\bibinfo{author}{Binazadeh, T.} \& \bibinfo{author}{Rahgoshay, M.~A.}
\newblock \bibinfo{title}{Robust output tracking of a class of non-affine
  systems}.
\newblock \emph{\bibinfo{journal}{Systems Science \& Control Engineering}}
  \textbf{\bibinfo{volume}{5}}, \bibinfo{pages}{426--433}
  (\bibinfo{year}{2017}).

\bibitem{ChenBaiKon2024}
\bibinfo{author}{Chen, Z.}, \bibinfo{author}{Bai, W.} \& \bibinfo{author}{Kong,
  L.}
\newblock \bibinfo{title}{Robust output tracking control of uncertain nonaffine
  systems with guaranteed tracking error bounds}.
\newblock \emph{\bibinfo{journal}{International Journal of Control, Automation
  and Systems}} \textbf{\bibinfo{volume}{22}}, \bibinfo{pages}{1--13}
  (\bibinfo{year}{2024}).

\bibitem{HuHad2026}
\bibinfo{author}{Hu, T.} \& \bibinfo{author}{Haddadin, S.}
\newblock \bibinfo{title}{{Simulation results of numerical examples for Dual
  Cost-Constraint projection solving optimal control}}  (\bibinfo{year}{2026}).
\newblock \urlprefix\url{https://doi.org/10.6084/m9.figshare.31577233}.

\end{thebibliography}


\begin{thebibliography}{10}
\expandafter\ifx\csname url\endcsname\relax
  \def\url#1{\texttt{#1}}\fi
\expandafter\ifx\csname urlprefix\endcsname\relax\def\urlprefix{URL }\fi
\providecommand{\bibinfo}[2]{#2}
\providecommand{\eprint}[2][]{\url{#2}}

\bibitem{NocedalWri2006}
\bibinfo{author}{Nocedal, J.} \& \bibinfo{author}{Wright, S.~J.}
\newblock \emph{\bibinfo{title}{Numerical optimization}}
  (\bibinfo{publisher}{Springer Science+Business Media, LLC.},
  \bibinfo{address}{New York, NY, USA}, \bibinfo{year}{2006}).

\bibitem{BoydVan2004}
\bibinfo{author}{Boyd, S.} \& \bibinfo{author}{Vandenberghe, L.}
\newblock \emph{\bibinfo{title}{Convex Optimization}}
  (\bibinfo{publisher}{{Cambridge University Press}},
  \bibinfo{address}{Cambridge, UK}, \bibinfo{year}{2004}).

\bibitem{LiXu2010}
\bibinfo{author}{Li, S.~J.} \& \bibinfo{author}{Xu, S.}
\newblock \bibinfo{title}{Sufficient conditions of isolated minimizers for
  constrained programming problems}.
\newblock \emph{\bibinfo{journal}{Numerical Functional Analysis and
  Optimization}} \textbf{\bibinfo{volume}{31}}, \bibinfo{pages}{715--727}
  (\bibinfo{year}{2010}).

\bibitem{TheMathWorks2020}
 \emph{\bibinfo{title}{Optimization Toolbox{\texttrademark} User's Guide}}
  (\bibinfo{publisher}{The MathWorks, Inc.}, \bibinfo{year}{2020}),
  \bibinfo{edition}{{R2020b}} edn.

\bibitem{Yamashita1980}
\bibinfo{author}{Yamashita, H.}
\newblock \bibinfo{title}{A differential equation approach to nonlinear
  programming}.
\newblock \emph{\bibinfo{journal}{Mathematical Programming}}
  \textbf{\bibinfo{volume}{18}}, \bibinfo{pages}{155--168}
  (\bibinfo{year}{1980}).

\bibitem{ZhouShi1997}
\bibinfo{author}{Zhou, Z.} \& \bibinfo{author}{Shi, Y.}
\newblock \bibinfo{title}{An {ODE} method of solving nonlinear programming}.
\newblock \emph{\bibinfo{journal}{Computers \& Mathematics with Applications}}
  \textbf{\bibinfo{volume}{34}}, \bibinfo{pages}{97--102}
  (\bibinfo{year}{1997}).

\bibitem{AllibhoyCor2023}
\bibinfo{author}{Allibhoy, A.} \& \bibinfo{author}{Cortés, J.}
\newblock \bibinfo{title}{Control-barrier-function-based design of gradient
  flows for constrained nonlinear programming}.
\newblock \emph{\bibinfo{journal}{IEEE Transactions on Automatic Control}}
  \textbf{\bibinfo{volume}{69}}, \bibinfo{pages}{3499--3514}
  (\bibinfo{year}{2024}).

\bibitem{FepponAllDap2020}
\bibinfo{author}{Feppon, F.}, \bibinfo{author}{Allaire, G.} \&
  \bibinfo{author}{Dapogny, C.}
\newblock \bibinfo{title}{{Null space gradient flows for constrained
  optimization with applications to shape optimization}}.
\newblock \emph{\bibinfo{journal}{{ESAIM: Control, Optimisation and Calculus of
  Variations}}} \textbf{\bibinfo{volume}{26}}, \bibinfo{pages}{90}
  (\bibinfo{year}{2020}).

\bibitem{Khatib1990}
\bibinfo{author}{Khatib, O.}
\newblock \bibinfo{title}{Motion/force redundancy of manipulators}.
\newblock In \emph{\bibinfo{booktitle}{Proc. of the Japan-USA Symposium on
  Flexible Automation}}, vol.~\bibinfo{volume}{1}, \bibinfo{pages}{337--342}
  (\bibinfo{address}{Kyoto, Japan}, \bibinfo{year}{1990}).

\bibitem{PotraWri2000}
\bibinfo{author}{Potra, S.~J., F. A.;~Wright}.
\newblock \bibinfo{title}{Interior-point methods}.
\newblock \emph{\bibinfo{journal}{J. Comput. Appl. Math.}}
  \textbf{\bibinfo{volume}{124}}, \bibinfo{pages}{281--302}
  (\bibinfo{year}{2000}).

\bibitem{Kelley1995}
\bibinfo{author}{Kelley, C.~T.}
\newblock \emph{\bibinfo{title}{Iterative Methods for Linear and Nonlinear
  Equations}} (\bibinfo{publisher}{Society for Industrial and Applied
  Mathematics}, \bibinfo{address}{Philadelphia, PA, USA},
  \bibinfo{year}{1995}).

\bibitem{Marinov2009}
\bibinfo{author}{Marinov, R.}
\newblock \bibinfo{title}{Convergence of the method of chords for solving
  generalized equations}.
\newblock \emph{\bibinfo{journal}{Rendiconti Del Circolo Matematico Di
  Palermo}} \textbf{\bibinfo{volume}{58}}, \bibinfo{pages}{11--27}
  (\bibinfo{year}{2009}).

\bibitem{Fletcher1987}
\bibinfo{author}{Fletcher, R.}
\newblock \emph{\bibinfo{title}{Practical Methods of Optimization}}
  (\bibinfo{publisher}{John Wiley \& Sons}, \bibinfo{address}{Chichester, West
  Sussex, England}, \bibinfo{year}{1987}), \bibinfo{edition}{second} edn.

\bibitem{LaSalle1960}
\bibinfo{author}{LaSalle, J.}
\newblock \bibinfo{title}{Some extensions of {Liapunov's} second method}.
\newblock \emph{\bibinfo{journal}{IRE Transactions on Circuit Theory}}
  \textbf{\bibinfo{volume}{7}}, \bibinfo{pages}{520--527}
  (\bibinfo{year}{1960}).

\bibitem{Gronwall1919}
\bibinfo{author}{Gronwall, T.~H.}
\newblock \bibinfo{title}{Note on the derivatives with respect to a parameter
  of the solutions of a system of differential equations}.
\newblock \emph{\bibinfo{journal}{Annals of Mathematics}}
  \textbf{\bibinfo{volume}{20}}, \bibinfo{pages}{292--296}
  (\bibinfo{year}{1919}).

\bibitem{Khalil2002}
\bibinfo{author}{Khalil, H.~K.}
\newblock \emph{\bibinfo{title}{Nonlinear systems}}
  (\bibinfo{publisher}{Prentice-Hall, Inc.}, \bibinfo{address}{Upper Saddle
  River, NJ, USA}, \bibinfo{year}{2002}), \bibinfo{edition}{third} edn.

\bibitem{CurryGarSul1983}
\bibinfo{author}{Curry, J.~H.}, \bibinfo{author}{Garnett, L.} \&
  \bibinfo{author}{Sullivan, D.}
\newblock \bibinfo{title}{{On the iteration of a rational function: computer
  experiments with Newton's method}}.
\newblock \emph{\bibinfo{journal}{Communications in Mathematical Physics}}
  \textbf{\bibinfo{volume}{91}}, \bibinfo{pages}{267 -- 277}
  (\bibinfo{year}{1983}).

\bibitem{Beardon1991}
\bibinfo{author}{Beardon, A.}
\newblock \emph{\bibinfo{title}{Iteration of Rational Functions: Complex
  Analytic Dynamical Systems}}.
\newblock Graduate Texts in Mathematics (\bibinfo{publisher}{Springer New
  York}, \bibinfo{address}{New York, NY, USA}, \bibinfo{year}{1991}).

\bibitem{EpureanuGre1998}
\bibinfo{author}{Epureanu, B.~I.} \& \bibinfo{author}{Greenside, H.~S.}
\newblock \bibinfo{title}{Fractal basins of attraction associated with a damped
  {Newton's} method}.
\newblock \emph{\bibinfo{journal}{SIAM Review}} \textbf{\bibinfo{volume}{40}},
  \bibinfo{pages}{102--109} (\bibinfo{year}{1998}).

\bibitem{Devaney2010}
\bibinfo{author}{Devaney, R.~L.}
\newblock \bibinfo{title}{{Chapter 4 - Complex Exponential Dynamics}}.
\newblock In \bibinfo{editor}{Broer, H.}, \bibinfo{editor}{Hasselblatt, B.} \&
  \bibinfo{editor}{Takens, F.} (eds.) \emph{\bibinfo{booktitle}{Handbook of
  Dynamical Systems}}, vol.~\bibinfo{volume}{3}, \bibinfo{pages}{125--223}
  (\bibinfo{publisher}{Elsevier Science}, \bibinfo{address}{Amsterdam, The
  Netherlands}, \bibinfo{year}{2010}).

\bibitem{Stone1948}
\bibinfo{author}{Stone, M.~H.}
\newblock \bibinfo{title}{The generalized {Weierstrass} approximation theorem}.
\newblock \emph{\bibinfo{journal}{Mathematics Magazine}}
  \textbf{\bibinfo{volume}{21}}, \bibinfo{pages}{167--184}
  (\bibinfo{year}{1948}).

\bibitem{BodieTayKamSie2020}
\bibinfo{author}{Bodie, K.}, \bibinfo{author}{Taylor, Z.},
  \bibinfo{author}{Kamel, M.} \& \bibinfo{author}{Siegwart, R.}
\newblock \bibinfo{title}{Towards efficient full pose omnidirectionality with
  overactuated mavs}.
\newblock In \bibinfo{editor}{Xiao, J.}, \bibinfo{editor}{Kr{\"o}ger, T.} \&
  \bibinfo{editor}{Khatib, O.} (eds.) \emph{\bibinfo{booktitle}{Proceedings of
  the 2018 International Symposium on Experimental Robotics}},
  \bibinfo{pages}{85--95} (\bibinfo{publisher}{Springer International
  Publishing}, \bibinfo{address}{Cham}, \bibinfo{year}{2020}).

\bibitem{KirchengastSteHor2018}
\bibinfo{author}{Kirchengast, M.}, \bibinfo{author}{Steinberger, M.} \&
  \bibinfo{author}{Horn, M.}
\newblock \bibinfo{title}{Control allocation under actuator saturation: An
  experimental evaluation}.
\newblock \emph{\bibinfo{journal}{IFAC-PapersOnLine}}
  \textbf{\bibinfo{volume}{51}}, \bibinfo{pages}{48--54}
  (\bibinfo{year}{2018}).
\newblock \bibinfo{note}{9th IFAC Symposium on Robust Control Design ROCOND
  2018}.

\bibitem{GhazaeiArdakaniRobJoh2015}
\bibinfo{author}{Ghazaei~Ardakani, M.~M.}, \bibinfo{author}{Robertsson, A.} \&
  \bibinfo{author}{Johansson, R.}
\newblock \bibinfo{title}{Online minimum-jerk trajectory generation}.
\newblock In \emph{\bibinfo{booktitle}{IMA Conference on Mathematics of
  Robotics}} (\bibinfo{address}{Oxford, UK}, \bibinfo{year}{2015}).

\bibitem{HuKueHad2020}
\bibinfo{author}{Hu, T.}, \bibinfo{author}{K\"{u}hn, J.} \&
  \bibinfo{author}{Haddadin, S.}
\newblock \bibinfo{title}{Forward and inverse dynamics modeling of human
  shoulder-arm musculoskeletal system with scapulothoracic constraint}.
\newblock \emph{\bibinfo{journal}{Computer Methods in Biomechanics and
  Biomedical Engineering}}  (\bibinfo{year}{2020}).

\bibitem{HuKueHad2018}
\bibinfo{author}{Hu, T.}, \bibinfo{author}{Kuehn, J.} \&
  \bibinfo{author}{Haddadin, S.}
\newblock \bibinfo{title}{Identification of human shoulder-arm kinematic and
  muscular synergies during daily-life manipulation tasks}.
\newblock In \emph{\bibinfo{booktitle}{IEEE Int. Conf. Bio. Rob.}},
  \bibinfo{pages}{1011--1018} (\bibinfo{address}{Enschede, the Netherlands},
  \bibinfo{year}{2018}).

\bibitem{ErdemirMcLHerBog2007}
\bibinfo{author}{Erdemir, A.}, \bibinfo{author}{McLean, S.},
  \bibinfo{author}{Herzog, W.} \& \bibinfo{author}{van~den Bogert, A.~J.}
\newblock \bibinfo{title}{Model-based estimation of muscle force exerted during
  movements}.
\newblock \emph{\bibinfo{journal}{Clin. Biomech.}}
  \textbf{\bibinfo{volume}{22}}, \bibinfo{pages}{131--154}
  (\bibinfo{year}{2007}).

\bibitem{Herzog2009a}
\bibinfo{author}{Herzog, W.}
\newblock \bibinfo{title}{Distribution problem in biomechanics}.
\newblock In \bibinfo{editor}{Binder, M.~D.}, \bibinfo{editor}{Hirokawa, N.} \&
  \bibinfo{editor}{Windhorst, U.} (eds.) \emph{\bibinfo{booktitle}{Encyclopedia
  of Neuroscience}}, \bibinfo{pages}{983--985}
  (\bibinfo{publisher}{Springer-Verlag Berlin Heidelberg},
  \bibinfo{address}{Heidelberg, Germany}, \bibinfo{year}{2009}).

\bibitem{MarshallGlaTraAme2022}
\bibinfo{author}{Marshall, N.} \emph{et~al.}
\newblock \bibinfo{title}{Flexible neural control of motor units}.
\newblock \emph{\bibinfo{journal}{Nature Neuroscience}}
  \textbf{\bibinfo{volume}{25}}, \bibinfo{pages}{1--13} (\bibinfo{year}{2022}).

\bibitem{AvertaBarCatHad2021}
\bibinfo{author}{Averta, G.} \emph{et~al.}
\newblock \bibinfo{title}{{U-Limb: A multi-modal, multi-center database on arm
  motion control in healthy and post-stroke conditions}}.
\newblock \emph{\bibinfo{journal}{GigaScience}} \textbf{\bibinfo{volume}{10}}
  (\bibinfo{year}{2021}).

\bibitem{ChadwickBlaKirBog2014}
\bibinfo{author}{Chadwick, E.~K.}, \bibinfo{author}{Blana, D.},
  \bibinfo{author}{Kirsch, R.~F.} \& \bibinfo{author}{van~den Bogert, A.~J.}
\newblock \bibinfo{title}{Real-time simulation of three-dimensional shoulder
  girdle and arm dynamics}.
\newblock \emph{\bibinfo{journal}{IEEE Trans. Biomed. Eng.}}
  \textbf{\bibinfo{volume}{61}}, \bibinfo{pages}{1947--1956}
  (\bibinfo{year}{2014}).

\bibitem{VidtSanMarHeg2018}
\bibinfo{author}{Vidt, M.} \emph{et~al.}
\newblock \bibinfo{title}{Modeling a rotator cuff tear: Individualized shoulder
  muscle forces influence glenohumeral joint contact force predictions}.
\newblock \emph{\bibinfo{journal}{Clinical Biomechanics}}
  \textbf{\bibinfo{volume}{60}} (\bibinfo{year}{2018}).

\bibitem{SurotoLicWibGul2022}
\bibinfo{author}{Suroto, H.} \emph{et~al.}
\newblock \bibinfo{title}{Morphology of humeral head and glenoid in normal
  shoulder of indonesian population}.
\newblock \emph{\bibinfo{journal}{Orthopedic Research and Reviews}}
  \textbf{\bibinfo{volume}{14}}, \bibinfo{pages}{459--469}
  (\bibinfo{year}{2022}).

\bibitem{HaraItoIwa1996}
\bibinfo{author}{Hara, H.}, \bibinfo{author}{Ito, N.} \&
  \bibinfo{author}{Iwasaki, K.}
\newblock \bibinfo{title}{Strength of the glenoid labrum and adjacent shoulder
  capsule}.
\newblock \emph{\bibinfo{journal}{Journal of Shoulder and Elbow Surgery}}
  \textbf{\bibinfo{volume}{5}}, \bibinfo{pages}{263--268}
  (\bibinfo{year}{1996}).
\newblock
  \urlprefix\url{https://www.sciencedirect.com/science/article/pii/S1058274696800528}.

\bibitem{Steinbach2001}
\bibinfo{author}{Steinbach, M.}
\newblock \bibinfo{title}{Markowitz revisited: Mean-variance models in
  financial portfolio analysis}.
\newblock \emph{\bibinfo{journal}{Society for Industrial and Applied
  Mathematics}} \textbf{\bibinfo{volume}{43}}, \bibinfo{pages}{31--85}
  (\bibinfo{year}{2001}).

\bibitem{BruniCesScoTar2016}
\bibinfo{author}{Bruni, R.}, \bibinfo{author}{Cesarone, F.},
  \bibinfo{author}{Scozzari, A.} \& \bibinfo{author}{Tardella, F.}
\newblock \bibinfo{title}{Real-world datasets for portfolio selection and
  solutions of some stochastic dominance portfolio models}.
\newblock \emph{\bibinfo{journal}{Data in Brief}} \textbf{\bibinfo{volume}{8}},
  \bibinfo{pages}{858--862} (\bibinfo{year}{2016}).

\bibitem{CormenLeiRivSte2009}
\bibinfo{author}{Cormen, T.~H.}, \bibinfo{author}{Leiserson, C.~E.},
  \bibinfo{author}{Rivest, R.~L.} \& \bibinfo{author}{Stein, C.}
\newblock \emph{\bibinfo{title}{Introduction to algorithms}}
  (\bibinfo{publisher}{The MIT Press}, \bibinfo{address}{Cambridge,
  Massachusetts}, \bibinfo{year}{2009}), \bibinfo{edition}{third} edn.

\bibitem{CiminiBem2017}
\bibinfo{author}{Cimini, G.} \& \bibinfo{author}{Bemporad, A.}
\newblock \bibinfo{title}{Exact complexity certification of active-set methods
  for quadratic programming}.
\newblock \emph{\bibinfo{journal}{IEEE Transactions on Automatic Control}}
  \textbf{\bibinfo{volume}{62}}, \bibinfo{pages}{6094--6109}
  (\bibinfo{year}{2017}).

\bibitem{ArnstroemAxe2022}
\bibinfo{author}{Arnström, D.} \& \bibinfo{author}{Axehill, D.}
\newblock \bibinfo{title}{A unifying complexity certification framework for
  active-set methods for convex quadratic programming}.
\newblock \emph{\bibinfo{journal}{IEEE Transactions on Automatic Control}}
  \textbf{\bibinfo{volume}{67}}, \bibinfo{pages}{2758--2770}
  (\bibinfo{year}{2022}).

\bibitem{CiminiBem2019}
\bibinfo{author}{Cimini, G.} \& \bibinfo{author}{Bemporad, A.}
\newblock \bibinfo{title}{Complexity and convergence certification of a block
  principal pivoting method for box-constrained quadratic programs}.
\newblock \emph{\bibinfo{journal}{Automatica}} \textbf{\bibinfo{volume}{100}},
  \bibinfo{pages}{29--37} (\bibinfo{year}{2019}).

\bibitem{Mehrmann1991}
\bibinfo{author}{Mehrmann, V.}
\newblock \emph{\bibinfo{title}{The Autonomous Linear Quadratic Control
  Problem: Theory and Numerical Solution}} (\bibinfo{publisher}{Springer-Verlag
  Berlin}, \bibinfo{address}{Heidelberg, Germany}, \bibinfo{year}{1991}).

\bibitem{Sommer1980}
\bibinfo{author}{Sommer, R.}
\newblock \bibinfo{title}{Control design for multivariable non-linear
  time-varying systems}.
\newblock \emph{\bibinfo{journal}{International Journal of Control}}
  \textbf{\bibinfo{volume}{31}}, \bibinfo{pages}{883--891}
  (\bibinfo{year}{1980}).
\newblock \urlprefix\url{https://doi.org/10.1080/00207178008961089}.
\newblock \eprint{https://doi.org/10.1080/00207178008961089}.

\bibitem{Zeitz1989}
\bibinfo{author}{Zeitz, M.}
\newblock \bibinfo{title}{Canonical forms for nonlinear systems}.
\newblock \emph{\bibinfo{journal}{IFAC Proceedings Volumes}}
  \textbf{\bibinfo{volume}{22}}, \bibinfo{pages}{33--38}
  (\bibinfo{year}{1989}).
\newblock
  \urlprefix\url{https://www.sciencedirect.com/science/article/pii/S147466701753606X}.
\newblock \bibinfo{note}{Nonlinear Control Systems Design, Capri, Italy, 14-16
  June 1989}.

\bibitem{BusawonDeAch2001}
\bibinfo{author}{Busawon, K.}, \bibinfo{author}{{De Leon-Morales}, J.} \&
  \bibinfo{author}{Acha-Daza, S.}
\newblock \bibinfo{title}{A controller design for a special class of nonlinear
  systems}.
\newblock \emph{\bibinfo{journal}{Applied Mathematics Letters}}
  \textbf{\bibinfo{volume}{14}}, \bibinfo{pages}{499--505}
  (\bibinfo{year}{2001}).
\newblock
  \urlprefix\url{https://www.sciencedirect.com/science/article/pii/S0893965900001841}.

\end{thebibliography}
% if required, the content of .bbl file can be included here once bbl is generated
%\input main.bbl

\end{document}

\baselineskip24pt % Double-space the manuscript.
	
\title{%
	\textbf{Supplementary Information}\\ for \\ 
	Optimal Control without Optimization
    %\added{Real-time optimal control without iteration}
	%Breaking the Real-Time Barrier in High-Dimensional Optimal Control
}

\author
{Tingli Hu\textsuperscript{1} and Sami Haddadin\textsuperscript{1,$\ast$}
	\\
	\normalsize{\textsuperscript{1}Mohamed bin Zayed University of Artificial Intelligence, Abu Dhabi, United Arab Emirates}
	\\
	\normalsize{\textsuperscript{$\ast$}To whom correspondence should be addressed; E-mail: \url{sami.haddadin@mbzuai.ac.ae}}
}
\date{}

\maketitle

%\begin{equation}
%%\mathchar'1~~
%%\mathchar'2~~
%%\mathchar'3~~
%%\mathchar'4~~
%%\mathchar'5~~
%%\mathchar'6~~
%%\mathchar'7~~
%%\mathchar'10~~
%%\mathchar'11~~
%%\mathchar'12~~
%%\mathchar'13~~
%%\mathchar'14~~
%%\mathchar'15~~
%%\mathchar'16~~
%%\mathchar'17~~
%%\mathchar'18~~
%%\mathchar'19~~
%%\mathchar'20~~
%%\mathchar'21~~
%%\mathchar'22~~
%%\mathchar'23~~
%%\mathchar'24~~
%%\mathchar'25~~
%%\mathchar'26~~
%%\mathchar'27~~
%%\mathchar'28~~
%%\mathchar'29~~
%%\mathchar'30~~
%%\mathchar'31~~
%%\mathchar'32~~
%%\mathchar'33~~
%%\mathchar'34~~
%%\mathchar'35~~
%%\mathchar'36~~
%%\mathchar'37~~
%%\mathchar'38~~
%%\mathchar'39~~
%%\mathchar'40~~
%%\mathchar'41~~
%%\mathchar'42~~
%%\mathchar'43~~
%%\mathchar'44~~
%%\mathchar'45~~
%%\mathchar'46~~
%%\mathchar'47~~
%%\mathchar'48~~
%%\mathchar'49~~
%%\mathchar'50~~
%%\mathchar'51~~
%%\mathchar'52~~
%%\mathchar'53~~
%%\mathchar'54~~
%\mathchar'55~~
%%\mathchar'56~~
%%\mathchar'57~~
%%\mathchar'58~~
%%\mathchar'59~~
%%\mathchar'60~~
%\end{equation}

\tableofcontents

\newpage
% !TEX root = ../supplementary.tex
%
\addtocounter{table}{-1}
%\begin{table}[!h]
%	\centering
%	%\caption{\label{tab:}%
	%	%}%
%	\resizebox{0.99\textwidth}{!}{%
	\renewcommand{\arraystretch}{1.01}%
	\begin{longtable}{  p{0ex} p{0.357\textwidth} p{0.55\textwidth} }
		& \textbf{Notation} & \textbf{Description} \\ \hline
		\rowcolor[gray]{0.85}\multicolumn{3}{l}{General} \\
		& $\realset$ & set of real numbers \alternaterowcolor
		& $\realset_{\geqslant 0} \subset \realset$ & set of nonnegative real numbers \alternaterowcolor
		& $\realset_{> 0} \subset \realset_{\geqslant 0}$ & set of positive real numbers \alternaterowcolor
		& $\naturalset \subset \realset$ & set of natural numbers \alternaterowcolor
		& $\nullv$ & conformable vector of zeros \alternaterowcolor
		& $\nullmat$ & conformable matrix of zeros \alternaterowcolor
		& $\Imat$ & conformable identity matrix \alternaterowcolor
		& $t \in \realset_{\geqslant 0}$ & time \alternaterowcolor
		& $\varxv \in \realset^{\dimx}$ & $\dimx$-dimensional general argument of functions \alternaterowcolor
		%
		%& $\varzv \in \realset^{\dimy}$ & $\dimy$-dimensional general argument of functions \alternaterowcolor
		%
		& $\sv \in \realset^{\dimc}$ & $\dimc$-dimensional general argument of functions \alternaterowcolor
		%
%		& \makecell[l]{$\varxv^{\astroot} \in \realset^{\dimx}$\\ $\varzv^{\astroot} \in \realset^{\dimy}$} & root of system of equations \alternaterowcolor
		& $\varxv^{\astroot} \in \realset^{\dimx}$ & root of system of equations \alternaterowcolor
		& $\varxv^{\ast} \in \realset^{\dimx}$ & equilibrium point \alternaterowcolor
		& $\hspace{0.75ex}\dot{\hspace{-0.75ex}\varxv}$ & time derivative of $\varxv(t)$, i.e., $\hspace{0.75ex}\dot{\hspace{-0.75ex}\varxv}(t) := {\mathrm{d} \varxv(t)}/{\mathrm{d}t}$ \alternaterowcolor
		& $(a, b) \subseteq \realset$ & $:=\big\{\,\varx\in\realset~\lvert~a < \varx < b\,\big\}$ \alternaterowcolor
		& $[a, b) \subset \realset$ & $:=\big\{\,\varx\in\realset~\lvert~a \leqslant \varx < b\,\big\}$ \alternaterowcolor
		& $(a, b] \subset \realset$ & $:=\big\{\,\varx\in\realset~\lvert~a < \varx \leqslant b\,\big\}$ \alternaterowcolor
		& $[a, b] \subset \realset$ & $:=\big\{\,\varx\in\realset~\lvert~a \leqslant \varx \leqslant b\,\big\}$ \alternaterowcolor
		& $\mathcal{N}(\varxv^{\ast},\delta) \subset \realset^{\dimx}$ & $:= \big\{\,\varxv\in\realset^{\dimx}~\big\lvert~\|\varxv - \varxv^{\ast}\|_{2} < \delta\,\big\}$ \alternaterowcolor
		& $\|\varxv\|$, $\|\Am\|$ & norm of vector $\varxv$ or induced norm of matrix $\Am$ \alternaterowcolor
		& $\|\varxv\|_{2}$, $\|\Am\|_{2}$ & (induced) Euclidean norm \alternaterowcolor
		& $\Am^{\#}_{\Wm}$ & $\Wm$-weighted right pseudoinverse of $\Am$ \alternaterowcolor
		& $K_{\textup{x}}$, $K_{\textup{y}} \in \realset_{> 0}$ & control gain \alternaterowcolor
		\rowcolor[gray]{0.85}\multicolumn{3}{l}{Dynamical system} \alternaterowcolor
		& $\xv \in \realset^{\dimx}$ & state \alternaterowcolor
		& $\uv \in \realset^{\dimx}$ & input \alternaterowcolor
		& $\yv \in \realset^{\dimy}$ & output \alternaterowcolor
		& $\desired{\yv} \in \realset^{\dimy}$ & reference for $\yv$ \alternaterowcolor
		& \makecell[l]{$\fv_{\textup{A}}: \realset^{\dimx} \to \realset^{\dimx}$ \\ $\Bm: \realset^{\dimx} \to \realset^{\dimx\times\dimx}$}
		 & state and input function as in $\dot{\xv} = \fv_{\textup{A}}(\xv) + \Bm(\xv) \uv$ \alternaterowcolor
		& $\hv: \realset_{\geqslant 0}\times\realset^{\dimx} \to \realset^{\dimy}$ & time-variant output function as in $\yv = \hv(t,\xv)$ \alternaterowcolor
		& $\Hm: \realset_{\geqslant 0}\times\realset^{\dimx} \to \realset^{\dimy\times\dimx}$ & Jacobian of $\hv(t,\cdot)$ \alternaterowcolor
		& $V$, $U: \realset^{\dimx} \to \realset_{\geqslant 0}$ & Lyapunov function \alternaterowcolor
		\rowcolor[gray]{0.85}\multicolumn{3}{l}{Optimization} \\
		& $\chiv \in \realset^{\dimx}$ & candidate of decision variables \alternaterowcolor
		& $\chiv^{[k]} \in \realset^{\dimx}$ & $\chiv$ at the $k$-th iteration \alternaterowcolor
		& $\chiv^{\optimal} \in \realset^{\dimx}$ & locally optimal $\chiv$ \alternaterowcolor
		& $\sigma: \realset_{\geqslant 0}\times\realset^{\dimx} \to \realset$ & time-variant cost function as in $\min_{\chiv} \sigma(t,\chiv)$ \alternaterowcolor
		& $\qv: \realset_{\geqslant 0}\times\realset^{\dimx} \to \realset^{\dimx}$ & gradient of $\sigma(t,\cdot)$ \alternaterowcolor
		& $\Qm: \realset_{\geqslant 0}\times\realset^{\dimx\times\dimx} \to \realset^{\dimx}$ & Hessian of $\sigma(t,\cdot)$ \alternaterowcolor
		& \makecell[l]{$\gv(t,\chiv) \hspace{-0.3ex}=\hspace{-0.3ex} \Gm(t) \chiv \hspace{-0.3ex}+\hspace{-0.3ex} \cv(t)$ \\ $\Gm(t) \in \realset^{\dimc\times\dimx}$ \\ $\cv(t) \in \realset^{\dimc}$}
		& \makecell[l]{time-variant inequality constraint as in \\ ~~~~$\nullv \geqslant \gv(t,\chiv)$} \alternaterowcolor
		%
		%& $\hv: \realset_{\geqslant 0}\times\realset^{\dimx} \to \realset^{\dimy}$ & time-variant constraint in optimization as in $\desired{\yv} = \hv(t,\chiv)$ \\
		%
		& \makecell[l]{$\beta: \realset\to\realset$\\ $\betav:\realset^{\dimc}\to\realset^{\dimc}$} & barrier function \alternaterowcolor
		& $\rho: \realset_{\geqslant 0}\times\realset^{\dimx} \to \realset$ & augmented $\sigma$ \alternaterowcolor%cost function as in $\min_{\chiv} \rho(t,\chiv)$ \alternaterowcolor
		& $\rv: \realset_{\geqslant 0}\times\realset^{\dimx} \to \realset^{\dimx}$ & gradient of $\rho(t,\cdot)$ \alternaterowcolor
		& $\Rm: \realset_{\geqslant 0}\times\realset^{\dimx} \to \realset^{\dimx\times\dimx}$ & Hessian of $\rho(t,\cdot)$ \alternaterowcolor
		& $\zetav \in \realset^{\dimy}$ & vector of Lagrange multipliers \alternaterowcolor
		& $L: \realset_{\geqslant 0}\times\realset^{\dimx}\times\realset^{\dimy} \to \realset$ & time-variant Lagrange function $L(t,\chiv,\zetav)$ \alternaterowcolor
		& $\lambdav:$\,\scalebox{0.9}{$\realset_{\geqslant 0}\times\realset^{\dimx}\times\realset^{\dimy} \to \realset^{1\times(\dimx+\dimy)}$} & %$:= {\partial L(t,\chiv,\zetav)}/{\partial [\chiv^\transp~~\zetav^\transp]^\transp}$ 
		Jacobian of $L(t,\cdot,\cdot)$	\alternaterowcolor
		& $\Lambdam:$\,\scalebox{0.8}{$\realset_{\geqslant 0}\times\realset^{\dimx}\times\realset^{\dimy} \to \realset^{(\dimx+\dimy)\times(\dimx+\dimy)}$} & %$:= {\partial \lambdav(t,\chiv,\zetav)}/{\partial [\chiv^\transp~~\zetav^\transp]^\transp}$ 
		Hessian of $L(t,\cdot,\cdot)$ \alternaterowcolor
		& $\Lambdam_{\textup{xx}}$:\,\scalebox{0.9}{$\realset_{\geqslant 0}\times\realset^{\dimx}\times\realset^{\dimy} \to \realset^{\dimx\times\dimx}$} & Hessian of $L(t,\cdot,\zetav)$ {\small (w.r.t. the second argument)} \alternaterowcolor
		& $\etav: \realset_{\geqslant 0}\times\realset^{\dimx}\times\realset^{\dimy} \to \realset^{\dimx}$ & \makecell[l]{iterated function as in \\ ~~~~$\chiv[k+1] = \chiv[k] + \etav(t, \chiv[k], \desired{\yv})$} \alternaterowcolor
		\rowcolor[gray]{0.85}\multicolumn{3}{l}{Optimal tracking control} \\
		& $t_{\final} \in \realset_{> 0}$ & final time instant \alternaterowcolor
		& $\phiv: \realset_{\geqslant 0} \to \realset^{\dimx}$ & candidate trajectory \alternaterowcolor
		& $\phiv^{\optimal}: \realset_{\geqslant 0} \to \realset^{\dimx}$ & locally optimal $\phiv$ \alternaterowcolor
		& $\varphiv \in \realset^{\dimx N}$ & discrete-time representation of $\phiv$ for the first $N$ time instants \alternaterowcolor
		& $\varphiv^{\optimal} \in \realset^{\dimx N}$ & discrete-time representation of $\phiv^{\optimal}$ for the first $N$ time instants \alternaterowcolor
		& $\desired{\varyv} \in \realset^{\dimy N}$ & discrete-time representation of $\desired{\yv}(t)$ for the first $N$ time instants \alternaterowcolor
		& $\varsigma: \realset^{\dimx N} \to \realset$ & sum of cost $\sigma(t,\phiv)$ for the first $N$ discrete time instants \alternaterowcolor
		& $\varqv: \realset^{\dimx N} \to \realset^{\dimx N}$ & gradient of $\varsigma$ \alternaterowcolor
		& $\varQm: \realset^{\dimx N} \to \realset^{(\dimx N)\times(\dimx N)}$ & Hessian of $\varsigma$ \alternaterowcolor
		& $\vargv: \realset^{\dimx N} \to \realset^{\dimc N}$ & discrete-time representation of $\gv(t,\cdot)$ for the first $N$ time instants \alternaterowcolor 
		& $\varGm: \realset^{\dimx N} \to \realset^{(\dimc N)\times(\dimc N)}$ & Jacobian of $\vargv$ \alternaterowcolor 
		& $\varhv: \realset^{\dimx N} \to \realset^{\dimy N}$ & discrete-time representation of $\hv(t,\cdot)$ for the first $N$ time instants \alternaterowcolor
		& $\varHm: \realset^{\dimx N} \to \realset^{(\dimy N)\times(\dimy N)}$ & Jacobian of $\varhv$ \alternaterowcolor 
		& $\varrho: \realset^{\dimx N} \to \realset$ & augmented $\varsigma$ \alternaterowcolor
		& $\varrv: \realset^{\dimx N} \to \realset^{\dimx N}$ & gradient of $\varrho$ \alternaterowcolor
		& $\varRm: \realset^{\dimx N} \to \realset^{(\dimx N)\times(\dimx N)}$ & Hessian of $\varrho$ \alternaterowcolor
		\rowcolor[gray]{0.85}\multicolumn{3}{l}{Computational complexity} \alternaterowcolor
		& $\timecomplexity(\omega)$ & time complexity of operation or procedure $\omega$ \alternaterowcolor
		& $\spacecomplexity(\omega)$ & space complexity of operation or procedure $\omega$ \alternaterowcolor
		& $\bigO$ & Bachmann-Landau big $\bigO$ notation \alternaterowcolor
		%
%		%
%		$K = \mathrm{const.} \gg 1$ & controller gain that accounts for the speed of convergence of $\tauv(t)$ to $\desired{\tauv}(t)$ \alternaterowcolor
%		%
%		$\Imat_{d \times d} \in \{1\}^{d\times d}$ & $d$-by-$d$ identity matrix \alternaterowcolor
%		$\nullmat_{\mu\times \nu} \in \{0\}^{\mu\times\nu}$ & $\mu$-by-$\nu$ zero matrix
		%
	\end{longtable}%
	\renewcommand{\arraystretch}{1}%

%\newpage
%\input{supplementary/approach}

\newpage
% !TEX root = ../supplementary.tex
%\unsetvruler % turn off line numbering
%\begin{figure}[!b]
%	\centering
%	\includegraphics[width=0.99\linewidth]{figures/intro_overview_figure.pdf}
%	%\noindent\makebox[\textwidth]{\includegraphics{figures/intro_overview_figure.pdf}}
%	\caption{\label{fig:intro figure}This work aims at solving pseudoconvex optimization problems for nonlinear dynamical systems by a closed-form tracking controller. Its efficiency and superior performance over state-of-the-art methods are proven by rigorous mathematical analysis as well as demonstrated by solving challenging high-dimensional real-world problems.}
%\end{figure}
%\section{Optimal tracking controller}\label{sec:optimal tracking controller}
%This section focuses on the details of system and problem class, the rigorous mathematical deductions and derivations leading to \emph{optimal tracking controller} (\OTC), as well as the mathematical analysis resulting in key properties and features of \OTC.
\section{System and problem class}\label{sec:lemmata_theorems_proofs}
In this work, we study the time-variant nonlinear input-affine system
\begin{subequations}\label{eq:gryphon}
	\begin{align}
		\dot{\xv} &= \fv_{\textup{A}}(\xv) + \Bm(\xv) \uv, \label{eq:gryphon, state eqn}\\
		\yv &= \hv(t,\xv), \label{eq:gryphon, output eqn}
	\end{align}
\end{subequations}
where $\uv \in \realset^{\dimx}$ denotes the control input,
$\xv \in \realset^{\dimx}$ the state, and
$\yv \in \mathcal{Y} \subset \realset^{\dimy}$ the output.
The functions $\fv_{\textup{A}}: \realset^{\dimx} \to \realset^{\dimx}$ and $\Bm: \realset^{\dimx} \to \realset^{\dimx\times\dimx}$ are such that \eqref{eq:gryphon, state eqn} satisfies existence and uniqueness theorem, i.e., there exists one and only one trajectory $\xv: \realset_{> 0} \to \realset^{\dimx}$ for any $\uv: \realset_{\geqslant 0} \to \realset^{\dimx}$ and any $\xv(0) \in \realset^{\dimx}$.
The input matrix $\Bm(\xv) \in \realset^{\dimx\times\dimx}$ is invertible $\forall \xv \in \realset^{\dimx}$.
The output function $\hv: \realset_{\geqslant 0}\times\realset^{\dimx} \to \mathcal{Y}$ is continuous in $t$ and twice continuously differentiable (differentiability class $C^{2}$) and surjective in $\xv$, for all $t \geqslant$ and all $\xv \in \realset^{\dimx}$. 
Due to the surjection $\hv$, for any reference output $\desired{\yv}$, there exist infinitely many possible solutions for the control $\uv$ such that $\lim_{t\rightarrow+\infty}\|\yv(t) - \desired{\yv}(t)\|_{2} = 0$. 
This motivates us to find the optimal solution among all these possible ones.
The optimality may be specified by
the pseudoconvex optimization problem
\begin{equation}\label{eq:constrained optimization}
	\begin{split}
		%\chiv^{\optimal} = \arg 
		\min_{\chiv} \sigma(t,\chiv)~~
		\text{subject to}~
		\nullv &= \hv(t,\chiv) - \desired{\yv}(t), \\ %\label{eq:equality constraint} 
		%~\wedge~
		\nullv &\geqslant \gv(t,\chiv) := \Gm(t) \chiv + \cv(t) %\label{eq:inequality constraint}
	\end{split}
\end{equation}
for every time instant $t \geqslant 0$
or the optimal tracking control problem
\begin{equation}\label{eq:optimal control problem}
	\begin{split}
		%\phiv^{\optimal} = \arg 
		\min_{\phiv} \int_{0}^{t_{\final}} \sigma\big(t,\phiv(t)\big)\,\mathrm{d}t~~
		%\quad~~
		\text{subject to}~
		\nullv = \hv\big(t,\,&\phiv(t)\big) - \desired{\yv}(t)~\forall t \in (0, t_{\final}], \\[-0.5ex]
		%~\wedge~
		\nullv \geqslant \gv\big(t,\,&\phiv(t)\big)~\forall t \in (0, t_{\final}], \\[0.75ex]
		%~\wedge~
		\xv(0) =\ &\phiv(0)
	\end{split}
\end{equation}
for the time period $[0, t_{\final}]$.
Herein, $\chiv \in \realset^{\dimx}$ denotes candidate variable and $\phiv: \realset_{\geqslant 0}\to\realset^{\dimx}$ candidate trajectory; both are to be determined.
The function
$\sigma: \realset_{\geqslant 0}\times\realset^{\dimx} \to \realset$ or its time-integral is the cost,
the time-variant equality constraint $\hv(t,\chiv) - \desired{\yv}(t)$ reflects the desired output tracking behavior,
and 
$\gv: \realset_{\geqslant 0}\times\realset^{\dimx} \to \realset^{\dimc}$ accounts for the time-variant linear inequality constraint.
The notation $\nullv$ denotes conformable zero vector.
In order to ensure that the constrained optimization problem~\eqref{eq:constrained optimization} is feasible 
and also to facilitate derivations and deductions, 
we impose the following assumptions:
\begin{assumptions}\label{assumption:isolated local minimizers}
	At any time instant $t \geqslant 0$,
	\begin{enumerate}
		\item the cost function $\sigma(t,\chiv)$ is twice continuously differentiable and strongly convex in $\chiv$;
		\item for~\eqref{eq:constrained optimization}, there exists at least one local minimizer;
		\item if multiple local minimizers exist, they are isolated points in $\realset^{\dimx}$.
	\end{enumerate}		
	Furthermore,
	\begin{enumerate}
		\setcounter{enumi}{3}
		\item there exists a differentiable $\chiv^{\optimal}: \realset_{\geqslant 0} \to \realset^{\dimx}$ such that for any $t \geqslant 0$, $\chiv^{\optimal}(t)$ is a local minimizer of \eqref{eq:constrained optimization}.
	\end{enumerate}		
\end{assumptions}
Herein, a \emph{local minimizer} $\chiv^{\optimal}$ of \eqref{eq:constrained optimization} refers to a point $\chiv$ that a) satisfies all constraints and b) achieves the smallest value of $\sigma(t,\cdot)$ in its neighborhood~\cite{NocedalWri2006}, i.e.,
\begin{equation*}
	\exists D > 0: \forall \chiv \in \left\{
	\chiv \in \realset^{\dimx}
	~\Bigg\lvert~\scalebox{0.9}{$
		\begin{matrix}
			\chiv \neq \chiv^{\optimal} \\
			\|\chiv - \chiv^{\optimal}\|_{2} < D \\
			\nullv = \hv(t,\chiv) - \desired{\yv}(t) \\
			\nullv \geqslant \gv(t,\chiv)
		\end{matrix}
		$}
	\right\},~
	\sigma(t,\chiv^{\optimal}) < \sigma(t,\chiv).
\end{equation*}

Assumptions~\ref{assumption:isolated local minimizers}.2 and 1.3 are satisfied under conditions given in \cite{BoydVan2004,NocedalWri2006} and \cite{LiXu2010}, respectively.
Intuitively, it seems to satisfy Assumption~\ref{assumption:isolated local minimizers}.4 if all time-variant functions $\sigma(t,\cdot)$, $\hv(t,\cdot)$, $\gv(t,\cdot)$ are differentiable in $t$.
Practically, all these assumptions are easily satisfiable and concrete examples are shown in the main document.

\begin{problem}\label{problem:find a gryphonheart}
	Let $\chiv^{\optimal}(t)$ be a local minimizer of the time-variant optimization problem~\eqref{eq:constrained optimization} at time instant $t$.
	Given the reference $\desired{\yv}(t)$ for the output of the time-variant system~\eqref{eq:gryphon}, 
	find a control 
	\begin{equation*}%\label{eq:definition,gryphonheart}
		\uv =
		\muv(t, \xv, \desired{\yv}; K_{\textup{x}})
	\end{equation*}
	with which $\xv(t)$ converges, as $t \geqslant 0$ increases, exponentially to 
	\begin{enumerate}
		\item $\chiv^{\optimal}(t)$, if $\chiv^{\optimal}(t) = \chiv^{\optimal}(0)$~$\forall t \geqslant 0$.
		\item the interior of a moving ball centered at $\chiv^{\optimal}(t)$ with a radius $\propto K_{\textup{x}}^{-1}$. %$\mathcal{N}\big(\chiv^{\optimal}(t),r\big)$ 
	\end{enumerate}
	Herein, $K_{\textup{x}} > 0$ is a design parameter.
\end{problem}

\begin{problem}\label{problem:find an optimal tracking controller}
	Let $\phiv^{\optimal}: [0, t_{\final}]\to\realset^{\dimx}$ be a local minimizer of \eqref{eq:optimal control problem} for the time period~$[0, t_{\final}]$.
	Given the reference $\desired{\yv}(t)$ for the output of the time-variant system~\eqref{eq:gryphon}, 
	find a control 
	\begin{equation*}%\label{eq:definition,gryphonheart}
		\uv =
		\muv(t, \xv, \desired{\yv}; K_{\textup{x}})
	\end{equation*}
	with which $\xv(t)$ converges, as $t \geqslant 0$ increases, exponentially to 
	\begin{enumerate}
		\item $\phiv^{\optimal}(t)$, if $\phiv^{\optimal}(t) = \phiv^{\optimal}(0)$~$\forall t \geqslant 0$.
		\item the interior of a moving ball centered at $\phiv^{\optimal}(t)$ with a radius $\propto K_{\textup{x}}^{-1}$.
	\end{enumerate}
	Herein, $K_{\textup{x}} > 0$ is a design parameter.
\end{problem}

%\begin{definition}[Optimal tracking controller]\label{definition:gryphonheart}
%	Any function $\muv: \realset_{\geqslant 0} \times \realset^{\dimx} \times \realset^{\dimy} \to \realset^{\dimx}$ that solves Problem~\ref{problem:find an optimal tracking controller} is an \emph{optimal tracking controller} (\OTC).
%\end{definition}

% Although strict local minimizers are not always isolated, all isolated local minimizers are strict \cite{NocedalWri2006}

%\begin{remark}
%	If $\hv(\cdot,\rv)$ is a linear function $\forall \rv$, there exists one and only one local minimum that is the global minimum~\cite{NocedalWri2006}. 
%	%	- Theorem 2.5
%	%	- Page 352, Exercise 12.4
%	%	- Page 373, Exercise 12.22
%	%	- Page 472, Theorem 16.2
%	%	- Page 484, "unique"
%\end{remark}
%

\newpage

\newpage
% !TEX root = ../supplementary.tex
\section{Dual cost-constraint projection -- overview}\label{sec:DCC projection}
We approach the constrained optimization problem~\eqref{eq:constrained optimization} by formulating a Lagrangian function that integrates
cost function, inequality constraint (via barrier function), and equality constraint (via Lagrange multipliers), following the interior-point algorithms~\cite{NocedalWri2006,TheMathWorks2020}.
However, unlike the state of the art where the barrier parameter and the Lagrange multipliers are updated iteratively, we are able to determine them analytically. % a) set the barrier parameter to $+\infty$ and b) eliminate Lagrange multipliers.
This results in the \emph{improved interior-point iteration} (numerical optimization), which is then reformulated into \emph{Dual Cost--Constraint projection} (DCC) through a new \emph{full-state feedback for equilibrium placement} (nonlinear control theory).
This \emph{isomorphism between numerical optimization and nonlinear control theory} concludes an exponential convergence and tracking behavior, which is formally proven using standard stability theory (Lyapunov theory, LaSalle's invariance principle, Grönwall's inequality). Figure~\ref{fig:proof flowchart} illustrates this logical workflow.
As the result, we obtain the DCC~\eqref{eq:DCC projection}. Formally, we state the following theorem:

\begin{theorem}\label{theorem:DCC projection}
Let $\chiv^{\optimal}(t)$ be a local minimizer of the time-variant optimization problem~\eqref{eq:constrained optimization} at time instant~$t$.
Let $\phiv^{\optimal}: [0, t_{\final}]\to\realset^{\dimx}$ be a local minimizer of \eqref{eq:optimal control problem} for the time period~$[0, t_{\final}]$.
Suppose Assumptions~\ref{assumption:isolated local minimizers} hold for system and problem class~\eqref{eq:gryphon}--\eqref{eq:optimal control problem}.
Suppose the initial condition $\xv(0)$ of \eqref{eq:gryphon} is sufficiently close to $\chiv^{\optimal}(0)$ (Remark~\ref{remark:sufficient closeness}). 
Given a reference $\desired{\yv}(t)$ for output $\yv$ of \eqref{eq:gryphon}, then
\begin{equationbox}
\begin{subequations}\label{eq:DCC projection}
\begin{align}%
	\uv 
	&=
	\bv(\xv) 
	+ \Km(\xv) \begin{bmatrix}
		\mathscr{P}\{\hv\}(t,\xv) & \mathscr{N}\{\hv\}(t,\xv)
	\end{bmatrix} \begin{bmatrix}
		\errv_{\textup{y}} \\
		\errv_{\upsigma,\textup{g}}
	\end{bmatrix} \\
	\intertext{with}
	\bv(\xv) &:= -\Bm(\xv)^{-1} \fv_{\scriptscriptstyle\textup{A}}(\xv), \\
	\Km(\xv) &:= \Bm(\xv)^{-1} K_{\textup{x}}, \\
	\mathscr{P}\{\hv\}(t,\xv) &:= \Hm(t,\xv)^{\#}_{\Rm(t,\xv)} = \Rm(t,\xv)^{-1} \Hm(t,\xv)^\transp \big( \Hm(t,\xv) \Rm(t,\xv)^{-1} \Hm(t,\xv)^\transp \big)^{-1}, \\
	\mathscr{N}\{\hv\}(t,\xv) &:= \Imat - \Hm(t,\xv)^{\#}_{\Rm(t,\xv)} \Hm(t,\xv), \\
	\errv_{\textup{y}} &:= \desired{\yv} - \hv(t,\xv), \\
	\intertext{and}
	\errv_{\upsigma,\textup{g}} &:= -\Rm(t,\xv)^{-1} \rv(t,\xv)
\end{align}%
\end{subequations}
\end{equationbox}
\noindent
solves Problems~\ref{problem:find a gryphonheart} and \ref{problem:find an optimal tracking controller}.
Herein, $K_{\textup{x}} > 0$ is a design parameter, 
$\Imat$ is a conformable identity matrix,
$\Hm(t,\xv) := {\partial \hv(t,\xv)}/{\partial \xv}$ is the Jacobian of $\hv(t,\xv)$, 
\smash{$\rv(t,\xv) := [{\partial \rho(t,\xv)}/{\partial \xv}]^\transp$} and $\Rm(t,\xv) := {\partial \rv(t,\xv)}/{\partial \xv}$ are the gradient and the Hessian of the \emph{augmented cost function}~\eqref{eq:augmented cost function}.
The resulting exponentially convergent tracking behavior is given in \eqref{eq:exponential convergence y}, \eqref{eq:exponential convergence x v2.5}, and \eqref{eq:exponential convergence x to phi_optimal}.
\end{theorem}
\begin{proof}
	See Section~\ref{sec:gryphonheart derivation}.
\end{proof}

\begin{figure}[!h]
\vspace{12pt}
\centering
\noindent\makebox[\textwidth]{
	\fcolorbox{white}{white}{
		\resizebox{1.25\textwidth}{!}{% \resizebox{1.525\textwidth}{!}{% 
				\scalefont{0.8500004250002125}%
				\input{figures/proof_flowchart}		
				\scalefont{1.17647}%
			}%
		}
	}
	\caption{%
		\label{fig:proof flowchart}
		The workflow to derive DCC addresses three interrelated fields, namely \emph{numerical optimization} (cyan), \emph{nonlinear control theory} (yellow), and \emph{optimal tracking control} (green). Each field contributes through a series of deductions, lemmata, and theorems to the final provably stable DCC solution, as detailed in the respective sections.}
	%\vspace{12pt}
\end{figure}

\paragraph{Optimality-emergent system dynamics}
Inserting \eqref{eq:DCC projection} into \eqref{eq:gryphon} yields the closed-loop system 
\begin{equation}\label{eq:OTC closed loop}
	\dot{\xv}
	\hspace{-0.25ex}=\hspace{-0.25ex}
	K_{\textup{x}}\overbracket[0.14ex]{
		\underbracket[0.14ex]{
			\Hm(t,\xv)^{\#}_{\Rm(t,\xv)}
		}_{\mathclap{\text{Pseudoinverse}}} \big(
		\desired{\yv} \hspace{-0.25ex}-\hspace{-0.25ex} \hv(t,\xv)
		\big)
	}^{\mathclap{\text{Component I}}}
	+
	K_{\textup{x}}\overbracket[0.14ex]{
		\underbracket[0.14ex]{
			\big(
			\Imat \hspace{-0.25ex}-\hspace{-0.25ex} \Hm(t,\xv)^{\#}_{\Rm(t,\xv)} \Hm(t,\xv)
			\big)
		}_{\mathclap{\text{Null-space projection}}} \underbracket[0.14ex]{
			\vphantom{\Hm(t,\xv)^{\#}_{\Rm(t,\xv)}} \big(\hspace{-0.25ex}-\hspace{-0.25ex} \Rm(t,\xv)^{-1} \rv(t,\xv)\big)
		}_{\mathclap{\text{Newton-Raphson step}}}
	}^{\mathclap{\text{Component II}}}.
\end{equation}
It is composed of two complementary components, where Component I ensures the output tracking behavior, while Component II guarantees the convergence towards a local minimizer of $\rho(t,\cdot)$ (by the Newton-Raphson step) without affecting the output tracking behavior. This is achieved through the null-space projector of $\hv(t,\cdot)$.

Similar complementarity structures are common in state of the art, including the differential-equation solutions~\cite{Yamashita1980,ZhouShi1997,AllibhoyCor2023} and the iterative algorithms found in~\cite{NocedalWri2006,TheMathWorks2020,FepponAllDap2020}.
However, these works did not take full advantage of the cost function's Hessian matrix.
On the contrary, we incorporate the Hessian $\Rm(t,\xv)$ in three essential spots, namely
\begin{enumerate}
	\item as weighting in the pseudoinverse in Component I,
	\item as weighting in the null-space projection in Component II,
	\item to constitute a Newton-Raphson step $-\Rm(t,\xv)^{-1} \rv(t,\xv)$ rather than the steepest descent $-\rv(t,\xv)$ in Component II.
\end{enumerate}
Without this incorporation, inequality constraints may be violated, equilibrium point $\xv^{\ast}$ of~\eqref{eq:OTC closed loop} may not be optimal, and/or the convergence may not be as rapid as exponential. Somewhat related considerations can be found in controlling kinematically redundant robot manipulators~\cite{Khatib1990}.

The details of every step from the interior-point algorithms to Theorem~\ref{theorem:DCC projection}, together with the necessary lemmeta, theorems, and their proofs, are presented in Section~\ref{sec:gryphonheart derivation}.
We also provide a systematic computational complexity analysis in Section~\ref{sec:computational complexity}, extended descriptions of numerical examples in Section~\ref{sec:extended_numerical_examples}, and possible generalization of DCC in Sections~\ref{sec:constrained system} and \ref{sec:chain of integrators}.

\newpage
% !TEX root = ../supplementary.tex
\section{%
%\deleted{Optimal tracking control}\addedtwo{Dual cost-constraint projection}
Dynamical realization of optimal control%
}\label{sec:gryphonheart derivation}
%To approach hybrid Problems~\ref{problem:find a gryphonheart} and \ref{problem:find an optimal tracking controller}, we employ three workflows, each targeting the solution of pseudoconvex optimization, feedback control, and optimal tracking control, respectively (see Main Document Figure~\ref{MD-fig:proof flowchart}).
%These workflows, through their respective lemmata and theorems, complement each other and ultimately lead to the resulting \emph{Dual Cost-Constraint projection}~(\OTC).
%%
%The detailed deduction process is presented in the following subsections. 
%
%%\begin{figure}[!h]
%%	%	\centering
%%	%	\resizebox{\textwidth}{!}{% 
%%	%		\input{figures/proof_flowchart}		
%%	%	} %
%%	\noindent\makebox[\textwidth]{
%%		\fcolorbox{white}{white}{
%%		\scalebox{0.97}{
%%		\input{figures/proof_flowchart}
%%		}
%%		}
%%	}
%%	\caption{%
%%		\label{fig:proof flowchart}
%%		The workflow of deriving the \OTC\ by addressing three interrelated aspects: pseudoconvex optimization (cyan), feedback control (yellow), and optimal tracking control (green). Each aspect contributes through a series of deductions, lemmata, and theorems, ultimately revealing $\muv$ as an \uls{o}ptimal \uls{t}racking \uls{c}ontroller.
%%	}%
%%\end{figure}

\subsection{Improved interior-point iteration}\label{sec:barrier method}
For solving optimization problems with both equality and inequality constraints such as \eqref{eq:constrained optimization},
one popular numerical method so far are the interior-point algorithms~\cite{PotraWri2000,NocedalWri2006,TheMathWorks2020}.
An interior-point algorithm iteratively updates the search direction toward a local minimizer while ensuring that the point $\chiv$ at each iteration remain in the interior of the feasible region (hence the name).
The key idea is to use \emph{barrier functions} to penalize points that violate the inequality constraint.
Following this idea, we improve the interior-point iteration and provide useful insights in deriving the \OTC\ (in Sections~\ref{sec:exponential convergence}~and~\ref{sec:optimal tracking controller}).

\begin{definition}[Barrier function, augmented cost function]\label{definition:barrier fcn}
	A \emph{barrier function} $\beta: \realset \to \realset_{\geqslant 0}$ for the open interval $(-\infty,0]$ is 
	a twice continuously differentiable function such that
	\begin{enumerate}
		\item $\beta(z) \gg \beta(0) \geqslant \beta(s) \geqslant 0$~$\forall z > 0~\wedge~s < 0$ and
		\item for all $t \geqslant 0$ and any strongly convex $\sigma(t,\cdot)$, the \emph{augmented cost function} 
		\begin{equation}\label{eq:augmented cost function}
			\rho(t,\chiv) := \sigma(t,\chiv)
			+ \frac{1}{2}\betav\big(\gv(t,\chiv)\big)^\transp \betav\big(\gv(t,\chiv)\big),
		\end{equation}
		is strongly convex in $\chiv$. Herein, $\betav(\sv)$ applies $\beta$ on each element of $\sv$.
	\end{enumerate}
\end{definition}

%The barrier functions in literature are often defined only in the interior of feasible region, such as reciprocal function \cite{Skarin2008,DenHertogRooTer1994} and logarithmic function \cite{PotraWri2000,NocedalWri2006,TheMathWorks2020}.
%While it suffices the purpose of penalizing points that are close to the boundary of the feasible region,
%our Definition~\ref{definition:barrier fcn} extends the domain of barrier functions to all real numbers.
%This extension allows an initial guess $\chiv[0]$ to be chosen from both inside and outside the feasible region. 
%%
One example of barrier function $\beta$ according to Definition~\ref{definition:barrier fcn} is
\begin{equation}\label{eq:barrier function}
	\beta(s)
	%= \beta(s; [\,p_{1}~p_{2}\,]^\transp)
	= \frac{p_{2}}{p_{1}} \ln\big(1+\exp(p_{1} s)\big),
	\quad
	p_{1} \gg 1,
	\quad
	p_{2} \gg 1,
\end{equation}
where $p_{1}$ and $p_{2}$ are \emph{barrier parameters}.
Figure~\ref{fig:barrierfcn} illustrates the graph of $\beta(s)/p_{2}$ with different values of $p_{1}$.
%\begin{equation}\label{eq:barrier function}%
%	\begin{split}
%		\beta(s)
%		&= \begin{cases}
%			0 & \text{if~}s \leqslant -E, \\
%			A\,E^{-2P}\,\left(s - E\right)^{2P} & \text{if~}s > -E,
%		\end{cases} \\
%		%
%		A &\gg 1, %\textstyle\max_{\zv \in [0,1]^{\dimx}}\big\|\Qm(t)^{-\frac{1}{2}}(t) \zv\big\|_{2}~\forall t > 0,
%		\quad
%		0 < E \ll 1, 
%		\quad
%		P \in \naturalset^{+}.
%	\end{split}
%\end{equation}
%An advantage of \eqref{eq:barrier function} over reciprocal and logarithmic functions is that it imposes no penalty on $s$ within the interval $(-\infty, -E]$. 
%Although the penalty it imposes when $x$ approaches $0^{-}$ is less than that of reciprocal and logarithmic functions, it remains sufficiently substantial.

%With the aid of a barrier function $\beta$, we combine \eqref{eq:cost function} and \eqref{eq:inequality constraint} into the \emph{augmented cost function}
%\begin{equation}\label{eq:augmented cost function}
%	\rho(t,\chiv) := \sigma(t,\chiv)
%	+ \frac{1}{2}\betav\big(\Gm(t)\chiv + \cv(t)\big)^\transp \betav\big(\Gm(t)\chiv + \cv(t)\big),
%\end{equation}
%where $\betav(\sv)$ applies $\beta$ on each element of $\sv$.
The barrier function and the augmented cost function transform \eqref{eq:constrained optimization} to the following time-variant optimization problem
\begin{equation}\label{eq:constrained optimization with barrier function}
	\begin{split}
		\min_{\chiv} \rho(t,\chiv)~~
		\text{subject to~~}
		\nullv = \hv(t,\chiv) - \desired{\yv}(t).
	\end{split}
\end{equation}
Although \eqref{eq:constrained optimization with barrier function} is not equivalent to the original \eqref{eq:constrained optimization}, 
the error occurs only for the local minimizers close to the boundary,
given Definition~\ref{definition:barrier fcn} of $\betav$ and $\rho$.\\[-2ex]

\begin{lemma}[Bound constraint as penalty]\label{corollary:bnd as penalty, v3}%\label{corollary:bnd as penalty}
	Let $\mathcal{E}^{\optimal}(t) \subset \realset^{\dimx}$ be the set of all local minimizers of \eqref{eq:constrained optimization} at instant $t$, and $\mathcal{E}^{\optimal}_{\textup{B}}(t)$ the set of all local minimizers \eqref{eq:constrained optimization with barrier function} at instant $t$.
	%
%	Let $\beta$ be a barrier function (Definition~\ref{definition:barrier fcn}).
	%
%	If $\beta(s) \gg 1$ is sufficiently large $\forall s > 0$
%	and $\beta(s) \ll 1$ is sufficiently small $\forall s \leqslant 0$ (Definition~\ref{definition:barrier fcn}),
	At any instant $t \geqslant 0$,
	for any $\chiv \in \mathcal{E}^{\optimal}(t)$, there exists $\chiv_{\textup{B}} \in \mathcal{E}^{\optimal}_{\textup{B}}(t)$ such that $\|\chiv - \chiv_{\textup{B}}\|_{2} < 10^{-n}$ with $n \gg 1$, and vice versa.
	%
%	Moreover, if there exists $E > 0$ such that $\beta(s) = 0$~$\forall s \leqslant -E$,
	Moreover, if \eqref{eq:barrier function} is chosen as the barrier function (Definition~\ref{definition:barrier fcn}), 
	then $n \rightarrow +\infty$ as $p_{1} \rightarrow +\infty$.
\end{lemma}
%
%\begin{lemma}[Bound constraint as penalty]\label{corollary:bnd as penalty, v2}
%	Suppose that all local minimizers of \eqref{eq:constrained optimization} at any time instant $t \geqslant 0$
%	are within the interval $(-\infty,\,-E]^{\dimx}$.
%	%
%	If \eqref{eq:barrier function} is chosen as the barrier function in \eqref{eq:constrained optimization with barrier function}, 
%	then \eqref{eq:constrained optimization with barrier function} is exactly the same as \eqref{eq:constrained optimization}.
%\end{lemma}

Since $p_{1}$ and $p_{2}$ can be chosen arbitrarily large, the optimization problem~\eqref{eq:constrained optimization with barrier function} with \eqref{eq:barrier function} is essentially \eqref{eq:constrained optimization}.
%We will show later in Section~\ref{sec:implementation remark} that it is manageable to set $p_{1}$ to $+\infty$ and to choose $p_{2}$ to be sufficiently large.
In fact, it is manageable to set $p_{1}$ to $+\infty$ and to choose $p_{2}$ to be sufficiently large, see Section~\ref{sec:remark on barrier fcn}.

%\begin{lemma}\label{corollary:bnd as penalty, v3}
%	The local minimizers of \eqref{eq:constrained optimization with barrier function} are the local minimizers of \eqref{eq:constrained optimization}.
%\end{lemma}

In order to solve \eqref{eq:constrained optimization with barrier function} at one arbitrary time instant, we apply the method of Lagrange multipliers~\cite{NocedalWri2006}.
The Lagrangian function is
\begin{align*}
	L(t, \chiv, \zetav)
	&= \rho(t,\chiv)
	+ \zetav^\transp \big(-\hv(t,\chiv) + \desired{\yv}(t)\big)
\end{align*}
with Lagrange multipliers $\zetav\in \realset^{\dimy}$.
The saddle points of $L(t, \chiv, \zetav)$ with respect to $\chiv$ and $\zetav$ are the roots of
\begin{subequations}\label{eq:Jacobian of Lagrangian}
	\begin{equation}\label{eq:Lagrangian's Jacobian}
		\nullv%_{\dimx+\dimy}
		= \lambdav(t,\chiv,\zetav)^\transp
		:= \begin{bmatrix}
			\frac{\partial L(t,\chiv, \zetav)}{\partial \chiv} &
			\frac{\partial L(t,\chiv, \zetav)}{\partial \zetav}
		\end{bmatrix}^\transp
		= \begin{bmatrix}
			\rv(t,\chiv) - \Hm(t,\chiv)^\transp \zetav  \\
			-\hv(t,\chiv) + \desired{\yv}(t)
		\end{bmatrix},
	\end{equation}
	where $\lambdav(t,\chiv,\zetav) \in \realset^{\dimx+\dimy}$ is the Jacobian of the Lagrangian function~$L(t,\chiv, \zetav)$,
	and
	\begin{alignat}{2}%\label{eq:definition of r, H}
		\rv(t,\varxv) &:= \left[\frac{\partial \rho(t,\varxv)}{\partial \varxv}\right]^\transp &&= \qv(t,\varxv) + \Gm(t)^\transp \xiv\big(\Gm(t)\varxv + \cv(t)\big), 
		\quad \label{eq:definition of r}\\
		&~&&\quad\qv(t,\varxv) := \left[\frac{\partial \sigma(t,\varxv)}{\partial \varxv}\right]^\transp,
		\quad
		\xiv(\sv) := \frac{\partial \betav(\sv)}{\partial \sv}^\transp \betav(\sv), \label{eq:definition of q and xi}
		\quad \\
		\Hm(t,\varxv) &:= \frac{\partial \hv(t,\varxv)}{\partial \varxv}. \label{eq:definition of H}
	\end{alignat}
\end{subequations}
Here, $\varxv \in \realset^{\dimx}$ and $\sv \in \realset^{\dimc}$ denote general arguments of functions.
The nonlinear equation~\eqref{eq:Lagrangian's Jacobian} may be solved by the Newton-Raphson method (see Section~\ref{sec:newton-raphson method}),
if the Jacobian of $\lambdav(t,\chiv,\zetav)$ (i.e, the Hessian of $L(t,\chiv, \zetav)$) 
\begin{subequations}\label{eq:Hessian of Lagrangian}
	\begin{equation}
		\begin{split}
		\Lambdam\big(t,\chiv,\zetav\big) 
		&:= \begin{bmatrix}
			\frac{\partial \lambdav(t,\chiv,\zetav)}{\partial \chiv} & \frac{\partial \lambdav(t,\chiv,\zetav)}{\partial \zetav}
		\end{bmatrix} \\
		&= \begin{bmatrix}
			\frac{\partial}{\partial \chiv}\frac{\partial L(t,\chiv, \zetav)}{\partial \chiv} & \frac{\partial}{\partial \zetav}\frac{\partial L(t,\chiv, \zetav)}{\partial \chiv} \\[1ex]
			\frac{\partial}{\partial \chiv}\frac{\partial L(t,\chiv, \zetav)}{\partial \zetav} & \frac{\partial}{\partial \zetav}\frac{\partial L(t,\chiv, \zetav)}{\partial \zetav}
		\end{bmatrix}
		= \begin{bmatrix}
			\Lambdam_{\textup{xx}}(t,\chiv,\zetav) & &-\Hm(t,\chiv)^\transp \\
			-\Hm(t,\chiv) & & \nullmat
		\end{bmatrix}
		\end{split}
	\end{equation}
	is invertible,
	where $\nullmat$ denotes the conformable zero matrix, and
	\begin{align}%\label{eq:definition of Lambda_xx, R}
		\Lambdam_{\textup{xx}}(t,\varxv,\zetav)
		&:= \frac{\partial}{\partial \varxv}\frac{\partial L(t,\varxv, \zetav)}{\partial \varxv}
		= \Rm(t,\varxv) + \frac{\partial \big(\Hm(t,\varxv)^\transp \zetav\big)}{\partial \varxv},
		\quad \label{eq:definition of Lambda_xx}\\
		\Rm(t,\varxv) &:= \frac{\partial \rv(t,\varxv)}{\partial \varxv} = \Qm(t,\varxv) + \Gm(t)^\transp \Xim\big(\Gm(t) \varxv + \cv(t)\big) \Gm(t), \quad \label{eq:definition of R}\\
		&~\quad\Qm(t,\varxv) := \frac{\partial}{\partial \varxv} \frac{\partial \sigma(t,\varxv)}{\partial \varxv},
		\quad 
		\Xim(\sv) := \frac{\partial}{\partial \sv} \left(
		\frac{\partial \betav(\sv)}{\partial \sv}^\transp \betav(\sv)
		\right). \label{eq:definition of Q and Xi}
	\end{align}
\end{subequations}
From the strong convexity of $\rho(t,\varxv)$ (Definition~\ref{definition:barrier fcn}), one concludes that $\Rm(t,\varxv)$ is \uls{s}ymmetric and \uls{p}ositive \uls{d}efinite (SPD). Consequently,
one observes
\begin{enumerate}
	\item $\Lambdam_{\textup{xx}}(t,\varxv,\zetav)$ is SPD if {$\hv(t,\varxv)$ is linear in $\varxv$, i.e.,} $\Hm(t,\varxv) = \Hm(t,\nullv)$ {is no longer a function of $\varxv$}.
	%$\hv(t,\xv)$ is linear in $\xv$ ($\Leftrightarrow \Hm(t,\xv) = \Hm(t,\nullv)$), and 
	\item {If one lets $\zetav = \nullv$, then} {$\Lambdam_{\textup{xx}}(t,\varxv,\zetav)=$}$\Lambdam_{\textup{xx}}(t,\varxv,\nullv) = \Rm(t,\varxv)$ is SPD, {as if $\hv(t,\varxv)$ were linear.}
\end{enumerate}
Therefore, for solving \eqref{eq:Lagrangian's Jacobian}, we apply a quasi-Newton method, where the Hessian $\Lambdam(t,\chiv,\zetav)$ is approximated by the invertible $\Lambdam(t,\chiv,\nullv)$.
This results in the following iteration
\begin{subequations}\label{eq:quasi-Newton sequence}
	\begin{empheq}[box=\colorbox{\coloripi}]{equation}
		\begin{bmatrix}
			\chiv[k+1] \\
			\zetav[k+1]
		\end{bmatrix}
		= \begin{bmatrix}
			\chiv[k] \\
			\zetav[k]
		\end{bmatrix} 
		- \Lambdam\big(t,\chiv[k],\nullv\big)^{-1} \lambdav(t, \chiv[k], \zetav[k])^\transp \label{eq:quasi-Newton sequence, v1}
	\end{empheq}
    \vspace{-4ex}
	\begin{align}
%		&\begin{bmatrix}
%			\chiv[k+1] \\
%			\zetav[k+1]
%		\end{bmatrix}
%		= \begin{bmatrix}
%			\chiv[k] \\
%			\zetav[k]
%		\end{bmatrix} 
%		- \Lambdam\big(t,\chiv[k],\nullv\big)^{-1} \lambdav(t, \chiv[k], \zetav[k])^\transp \label{eq:quasi-Newton sequence, v1}
%		\\[0.5ex]
		&\Leftrightarrow \nonumber\\
		&\hspace{-1ex}\chiv[k+1] = \chiv[k] 
		+ \etav_{\textup{A}}\big(t,\chiv[k],\desired{\yv}(t)\big)
		- \etav_{\textup{B}}(t,\chiv[k])
		+ \psiv_{\textup{A}}(t,\chiv[k],\zetav[k]), \label{eq:quasi-Newton sequence, chi}\\[-0.5ex]
		&\hspace{-1ex}\zetav[k+1] = \psiv_{\textup{B}}\big(t,\chiv[k],\desired{\yv}(t)\big) \label{eq:quasi-Newton sequence, zeta}
	\end{align}%
	with
	\begin{align}%
		\etav_{\textup{A}}(t, \varxv, \desired{\yv})
		&:= \Hm(t,\varxv)^{\#}_{\Rm(t,\varxv)} \big(
			\desired{\yv} - \hv(t,\varxv)
		\big), \label{eq:definition of eta_A}\\
		\etav_{\textup{B}}(t,\varxv) 
		&:= \Big(
			\Imat - \Hm(t,\varxv)^{\#}_{\Rm(t,\varxv)} \Hm(t,\varxv)
		\Big) \Rm(t,\varxv)^{-1} \rv(t,\varxv), \label{eq:definition of eta_B}\\
		\psiv_{\textup{A}}(t,\varxv,\zetav) 
		&:= \Big(
			\Imat - \Hm(t,\varxv)^{\#}_{\Rm(t,\varxv)} \Hm(t,\varxv)
		\Big) \Rm(t,\varxv)^{-1} \Hm(t,\varxv)^\transp \zetav, \\
		\psiv_{\textup{B}}(t,\varxv,\desired{\yv}) 
		&:= \Big(\Hm(t,\varxv) \Rm(t,\varxv)^{-1} \Hm(t,\varxv)^\transp\Big)^{\hspace{-0.25ex}-1\hspace{-0.25ex}} \Big(
		\desired{\yv} 
		\hspace{-0.3ex}-\hspace{-0.3ex} \hv(t,\varxv)
		\hspace{-0.3ex}+\hspace{-0.3ex} \Hm(t,\varxv) \Rm(t,\varxv)^{\hspace{-0.25ex}-1\hspace{-0.25ex}} \rv(t,\varxv)
		\Big)
	\end{align}
\end{subequations}
and conformable identity matrix $\Imat$.
%$\Rm(t,\xv)$-weighted right pseudoinverse $\Hm(t,\xv)^{\#}_{\Rm(t,\xv)}$ of $\Hm(t,\xv)$
%\begin{equation*}
%	\Hm(t,\xv)^{\#}_{\Rm(t,\xv)} 
%	= \Rm(t,\xv)^{-1} \Hm(t,\xv)^\transp \Big(
%		\Hm(t,\xv) \Rm(t,\xv)^{-1} \Hm(t,\xv)^\transp
%	\Big)^{-1}
%\end{equation*}
The operation $\Am^{\#}_{\Wm}$ is the $\Wm$-weighted right pseudoinverse of a matrix $\Am$, i.e.,
\begin{equation}\label{eq:weighted right pseudoinverse}
	\Am^{\#}_{\Wm} 
	= \Wm^{-1} \Am^\transp \big(
		\Am \Wm^{-1} \Am^\transp
	\big)^{-1}.
\end{equation}

\begin{definition}[Quotient-linear convergence~\cite{Kelley1995,NocedalWri2006,Marinov2009}]\label{definition:q-linear convergence}
	A convergent sequence $\chiv[0]$, $\chiv[1]$, $\chiv[2]$, $\ldots$ with $\chiv^{\ast} = \lim\limits_{k\rightarrow+\infty} \chiv[k]$
	is called \emph{quotient-linearly convergent} if 
	\begin{equation*}
		\exists M = \textup{const.} \in (0, 1):
		\big\|\chiv[k+1] - \chiv^{\ast}\big\|_{2}
		\leqslant 
		M\,\big\|\chiv[k] - \chiv^{\ast}\big\|_{2}
		~~
		\forall k \in \naturalset.
	\end{equation*}%\vspace{0.5ex}
	The constant $M$ is known as the rate of convergence.
\end{definition}

\begin{lemma}[Chord method \cite{Kelley1995}]\label{lemma:chord method}
	Let $[\,\chiv^{\astroot\,\transp}\,\,\zetav^{\astroot\,\transp}\,]^\transp$ denote a root of \eqref{eq:Lagrangian's Jacobian}.
	There exist $\delta_{1} > 0$ and $\delta_{2} > 0$ such that
	if \begin{equation}\label{eq:sufficient closeness of chord method}
		\left\|\begin{bmatrix}
			\chiv[0] - \chiv^{\astroot} \\
			\zetav[0] - \zetav^{\astroot}
		\end{bmatrix}\right\|_{2} < \delta_{1}
		~\wedge~
		\left\|\hat{\Lambdam} - \Lambdam(t,\chiv^{\astroot},\zetav^{\astroot})\right\|_{2} < \delta_{2}
	\end{equation}
	then the iteration
	\begin{equation}\label{eq:chord-method iteration}
		\begin{bmatrix}
			\chiv[k+1] \\
			\zetav[k+1]
		\end{bmatrix}
		= \begin{bmatrix}
			\chiv[k] \\
			\zetav[k]
		\end{bmatrix} 
		- \hat{\Lambdam}^{-1} \lambdav(t,\chiv[k], \zetav[k])
	\end{equation}
	converges quotient-linearly to $[\,\chiv^{\astroot\,\transp}\,\,\zetav^{\astroot\,\transp}\,]^\transp$ as $k \in \naturalset$ increases.
	Moreover, the rate of convergence decreases proportionally as $\delta_{1} + \delta_{2}$ decreases.
	Herein, $\|\cdot\|_{2}$ denotes the Euclidean norm (for vectors) or the matrix norm induced by the Euclidean norm.
\end{lemma}	
\begin{proof}
	See~\emph{chord method} in \cite{Kelley1995}.
	Note that $t$ is considered constant since we are considering only \emph{one} time instant and $t$ does not depend on $k$.
\end{proof}

Since \eqref{eq:quasi-Newton sequence} is a special case of \eqref{eq:chord-method iteration} with $\hat{\Lambdam} = \Lambdam(t,\chiv[k],\nullv)$, Lemma~\ref{lemma:chord method} is valid for \eqref{eq:quasi-Newton sequence}.
Furthermore, we claim that the sufficient-closeness condition~\eqref{eq:sufficient closeness of chord method} can be further simplified in case $\hat{\Lambdam} = \Lambdam(t,\chiv[k],\nullv)$.

\begin{lemma}\label{lemma:quasi Newton linear convergence}
	Let $[\,\chiv^{\astroot\,\transp}\,\,\zetav^{\astroot\,\transp}\,]^\transp$ denote a root of \eqref{eq:Lagrangian's Jacobian}.
	There exists $\delta > 0$ such that
	if \begin{equation}\label{eq:sufficient closeness of quasi Newton}
		\|\chiv[0] - \chiv^{\astroot}\|_{2} < \delta
	\end{equation}
	then the sequence $[\,\chiv[0]^\transp\,\,\zetav[0]^\transp\,]^\transp$, $[\,\chiv[1]^\transp\,\,\zetav[1]^\transp\,]^\transp$, $[\,\chiv[2]^\transp\,\,\zetav[2]^\transp\,]^\transp$, $\ldots$ generated by \eqref{eq:quasi-Newton sequence}
	converges quotient-linearly to $[\,\chiv^{\astroot\,\transp}\,\,\zetav^{\astroot\,\transp}\,]^\transp$.
	Moreover, the rate of convergence decreases as $\delta$ decreases.
\end{lemma}	
\begin{proof}
	%	Since $\zetav$ is uniquely determined by $\chiv$, cf.~\eqref{eq:quasi-Newton sequence, zeta}, 
	%	we denote $\zetav = \psiv(\chiv)$.
	Since $\rv(t,\varxv)$ is of $C^{1}$ in $\varxv$ and $\hv(t,\varxv)$ is of $C^{2}$ in $\varxv$, the functions $\Lambdam(t,\varxv,\zetav)$ and $\psiv_{\textup{B}}(t,\varxv,\desired{\yv})$ are locally Lipschitz continuous in $\varxv$.
	Furthermore, $\Lambdam(t,\varxv,\zetav)$ is linear in $\zetav$ as defined in \eqref{eq:Hessian of Lagrangian}, i.e., $\Lambdam(t,\varxv,\zetav)$ is globally Lipschitz continuous in $\zetav$.
	Hence, one obtains the following deduction:
	\begin{subequations}\label{eq:sufficient closeness of quasi Newton, zeta and lambda}
		\begin{align}
			\begin{split}
			\left\|\Lambdam(t,\chiv,\nullv) - \Lambdam(t,\chiv^{\astroot},\zetav^{\astroot}) \right\|_{2} 
			< \gamma_{\scriptscriptstyle \varLambda}(t) \left\|\begin{bmatrix}
				\chiv - \chiv^{\astroot} \\
				\nullv - \zetav^{\astroot}
			\end{bmatrix}\right\|_{2}
			= \gamma_{\scriptscriptstyle \varLambda}(t) \sqrt{
				\left\| \chiv - \chiv^{\astroot} \right\|_{2}^{2} + 	
				\left\| \nullv - \zetav^{\astroot} \right\|_{2}^{2}
			} \\[-1ex]
			< \gamma_{\scriptscriptstyle \varLambda}(t) \sqrt{
				\delta^{2} + 	
				\left\| \zetav^{\astroot} \right\|_{2}^{2}
			},
			\end{split}\\[1ex]
			\begin{split}
			\left\| \zetav[0] - \zetav^{\astroot} \right\|_{2}
			= \left\| \psiv_{\textup{B}}(t,\chiv[0],\desired{\yv}) - \psiv_{\textup{B}}(t,\chiv^{\astroot},\desired{\yv}) \right\|_{2}
			< \gamma_{\scriptscriptstyle \psi_{\textup{B}}}(t,\desired{\yv}) \left\| \chiv - \chiv^{\astroot} \right\|_{2} \\
			= \gamma_{\scriptscriptstyle \psi_{\textup{B}}}(t,\desired{\yv}) \delta,
			\end{split}
		\end{align}
	\end{subequations}
	where $\gamma_{\scriptscriptstyle \varLambda}(t)$ is the Lipschitz constant of $\Lambdam(t,\cdot,\cdot)$ in the considered neighborhood of $\chiv^{\astroot}$,
	and \smash{$\gamma_{\scriptscriptstyle \psi_{\textup{B}}}(t,\desired{\yv})$} the Lipschitz constant of $\psiv_{\textup{B}}(t,\cdot,\desired{\yv})$ in the considered neighborhood of \smash{$\chiv^{\astroot}$}.	
	As mentioned before, $t$ is considered constant, and so is $\desired{\yv}$, since we are considering only \emph{one} time instant and $t$ does not depend on $k$.
	It is obvious that \eqref{eq:sufficient closeness of quasi Newton} and \eqref{eq:sufficient closeness of quasi Newton, zeta and lambda} together are \eqref{eq:sufficient closeness of chord method}.
	Therefore, the conclusion of Lemma~\ref{lemma:chord method} is applicable here.
\end{proof}

\begin{lemma}\label{lemma:chi_astroot = chi_optimal}
	%	Suppose that $\Lambdam(\chiv^{\astroot},\zetav^{\astroot}$ is SPD for any root $[\,\chiv^{\astroot\,\transp}\,\,\zetav^{\astroot\,\transp}\,]^\transp$ denote an arbitrary root of \eqref{eq:Lagrangian's Jacobian}.
	The roots ($\chiv$ component) of \eqref{eq:Lagrangian's Jacobian} are the local minimizers of \eqref{eq:constrained optimization}.
\end{lemma}	
\begin{proof}
	All roots ($\chiv$ component) of \eqref{eq:Lagrangian's Jacobian} are the local minimizers of \eqref{eq:constrained optimization with barrier function} (the Karush-Kuhn-Tucker conditions \cite{Fletcher1987,BoydVan2004,NocedalWri2006}), % Theorems 12.1 (KKT), 12.5, 12.6 in NocedalWri2006,
	and the latter are the local minimizers of \eqref{eq:constrained optimization} (Lemma~\ref{corollary:bnd as penalty, v3}).
\end{proof}

This means, all $\chiv^{\astroot}$ above can be replaced by $\chiv^{\optimal}$.
For the sake of readability, we use the single notation $\chiv^{\optimal}$ for all these three notions (roots of \eqref{eq:Lagrangian's Jacobian}, local minimizers of \eqref{eq:constrained optimization with barrier function}, and local minimizers of \eqref{eq:constrained optimization}) from this point forth.

Since $\zetav[k]$ converges as $k$ increases (Lemma~\ref{lemma:quasi Newton linear convergence}), we apply $\zetav[k+1] = \zetav[k]$ to \eqref{eq:quasi-Newton sequence}. By inserting \eqref{eq:quasi-Newton sequence, zeta} into \eqref{eq:quasi-Newton sequence, chi}, one obtains the \emph{interior-point iteration}
\begin{align}
	\chiv[k+1] &= \chiv[k] + \etav(t, \chiv[k], \desired{\yv}), \label{eq:barrier-method sequence} \\
	\etav(t, \varxv, \desired{\yv}) 
	&:= \etav_{\textup{A}}(t, \varxv, \desired{\yv})
	- \etav_{\textup{B}}(t, \varxv), \label{eq:definition of eta}
\end{align}
where $\etav_{\textup{A}}$ and $\etav_{\textup{B}}$ are defined in \eqref{eq:definition of eta_A} and \eqref{eq:definition of eta_B}, respectively.
Rigorously, \eqref{eq:barrier-method sequence} shall be written as
\begin{empheq}[box=\colorbox{\coloripi}]{equation}
%\begin{equation}
	\label{eq:barrier-method sequence, v2}
	\chiv^{[k+1]}(t) = \chiv^{[k]}(t) + \etav\big(t, \chiv^{[k]}(t), \desired{\yv}(t)\big)
%\end{equation}
\end{empheq}
to reflect the time-variant property of \eqref{eq:constrained optimization}.

\begin{theorem}[Interior-point iteration]\label{theorem:barrier-method iteration}
	For any local minimizer $\chiv^{\optimal}$ of \eqref{eq:constrained optimization},
	there exists $\delta > 0$ such that
	if $\|\chiv^{[0]}(t) - \chiv^{\optimal}\|_{2} < \delta$
	then $\chiv^{[k]}(t)$ of~\eqref{eq:barrier-method sequence, v2}
	converges quotient-linearly (Definition~\ref{definition:q-linear convergence}) to $\chiv^{\optimal}(t)$ as $k \in \naturalset$ increases, i.e.,
	\begin{empheq}[box=\colorbox{\coloripi}]{equation}
%	\begin{equation*}
	\label{eq:barrier-method iteration, convergence}
		\exists M \in (0, 1):
		\big\|\chiv^{[k+1]}(t) - \chiv^{\optimal}(t)\big\|_{2}
		\leqslant 
		M\,\big\|\chiv^{[k]}(t) - \chiv^{\optimal}(t)\big\|_{2}
		~~
		\forall k \in \naturalset.
%	\end{equation*}
	\end{empheq}	
\end{theorem}
\begin{proof}
	As elaborated above, there is no difference between \eqref{eq:quasi-Newton sequence} and \eqref{eq:barrier-method sequence}. Hence, one concludes the local convergence of $\chiv[k]$ to a root $\chiv^{\astroot}$ of \eqref{eq:Lagrangian's Jacobian} as in Lemma~\ref{lemma:quasi Newton linear convergence},
	where $\chiv^{\astroot}$ is as well a local minimizer $\chiv^{\optimal}$ of \eqref{eq:constrained optimization} (Lemma~\ref{lemma:chi_astroot = chi_optimal}).
	%Inserting \eqref{eq:definition of eta_A} and \eqref{eq:definition of eta_B} into \eqref{eq:definition of eta} and then into \eqref{eq:barrier-method sequence, v2} results in \eqref{eq:barrier-method sequence, v3}.
	%Therefore, the local quotient-linear convergence of $\chiv[k]$ to a local minimizer $\chiv^{\optimal}$ of \eqref{eq:constrained optimization} is concluded.
\end{proof}

\begin{definition}[Equilibrium points at instant $t$]\label{definition:equilibrium chi_ast}
	A point $\chiv^{\ast} \in \realset^{\dimx}$ is called an \emph{equilibrium point} of \eqref{eq:barrier-method sequence} \emph{at instant $t$},
	if and only if
	\begin{equation*}
		\nullv = \etav(t, \chiv^{\ast}, \desired{\yv})
		~\Leftrightarrow~
		\chiv^{\ast}
		= \chiv^{\ast}
		+ \etav(t, \chiv^{\ast}, \desired{\yv}).
	\end{equation*}
	In other words, any 
	\begin{equation}\label{eq:set of equilibrium points, barrier-method sequence}
		\chiv^{\ast} \in
		\mathcal{E}_{\chi}(t,\desired{\yv}) := \big\{\,
		\varxv\in\realset^{\dimx} ~\big\lvert~ \nullv = \etav(t, \varxv, \desired{\yv})
		\,\big\}
	\end{equation}
	is an \emph{equilibrium point} of \eqref{eq:barrier-method sequence} \emph{at instant $t$}.
\end{definition}

It is worth noting that \begin{itemize}
	\item the set $\mathcal{E}_{\chi}(t,\desired{\yv})$ of equilibrium points is time-variant and depends on $\desired{\yv}$, and
	\item $\mathcal{E}_{\chi}(t,\desired{\yv}) \neq \emptyset$ as long as Assumptions~\ref{assumption:isolated local minimizers}.2 and \ref{assumption:isolated local minimizers}.3 are satisfied.
\end{itemize}

\begin{lemma}[Local minimizers]\label{corollary:local minimizer is a  asymptotically stable equilibrium, barrier-method sequence}
	A local minimizer of \eqref{eq:constrained optimization}
	is an asymptotically stable equilibrium point of \eqref{eq:barrier-method sequence}.
\end{lemma}
\begin{proof}
	It can be concluded immediately from Lemma~\ref{lemma:chi_astroot = chi_optimal} and Theorem~\ref{theorem:barrier-method iteration}.
\end{proof}

\subsection{Full-state feedback for equilibrium placement}\label{sec:feedback linearization}
Instead of considering the system~\eqref{eq:gryphon} immediately, 
we firstly consider a simple, yet, very special system that is
\begin{subequations}\label{eq:special gryphon}
	\begin{align}
		\dot{\xv} &= \fv_{\textup{A}}^\prime(t,\xv)
		+ \uv^\prime, \\
		\yv &= \hv(t,\xv), \label{eq:special gryphon, output fcn}
	\end{align}
\end{subequations}
where $\uv^\prime \in \realset^{\dimx}$.
The function $\fv_{\textup{A}}^\prime: \realset_{\geqslant 0}\times\realset^{\dimx} \to \realset^{\dimx}$ is of a certain form, which will be revealed later.
It is easy to conclude that if 
\begin{equation}\label{eq:gryphon, fbl-2}
	\uv = \Bm(\xv)^{-1} \Big(-\fv_{\textup{A}}(\xv) 
	+ \fv_{\textup{A}}^\prime(t,\xv)
	+ \uv^\prime\Big),
\end{equation}
is chosen, then \eqref{eq:gryphon} becomes \eqref{eq:special gryphon}.

To control the output $\yv$ of \eqref{eq:special gryphon} to follow the reference $\desired{\yv}$, we deploy the standard feedback linearization.
\begin{lemma}\label{lemma:h(x) converges to zero, v2}
	Suppose %$\desired{\yv}: \realset_{\geqslant 0} \to \realset^{\dimy}$ 
	$\desired{\yv}(t)$ and $\hv(t,\xv)$ are
	differentiable in $t$.
	If 
	\begin{equation}\label{eq:special gryphon, fbl v2}
		%		\begin{split}
			\uv^\prime = 
			\Hm(t,\xv)^{\#}_{\Qm^\prime} \big( K_{\textup{y}}\big(\desired{\yv} - \hv(t,\xv)\big) - \Hm(t,\xv)\fv_{\textup{A}}^\prime(t,\xv) \big),
			\quad 
			\Qm^\prime\text{ is SPD},
			\quad
			K_{\textup{y}} > 0
			\quad
			%		\end{split}
	\end{equation}
	is chosen for \eqref{eq:special gryphon},
	then $\yv(t)$ converges exponentially{, as $t > 0$ increases, with time constant $K_{\textup{y}}^{-1}$ to a ball neighborhood of $\desired{\yv}(t)$ with a radius $\propto K_{\textup{y}}^{-1}$.}
\end{lemma}
\begin{proof}
	Applying Lie derivative to \eqref{eq:special gryphon, output fcn} results in
	\begin{equation}\label{eq:Lie derivative of special gryphon}
		\dot{\yv} 
		= \frac{\partial \hv(t,\xv)}{\partial t} + \Hm(t,\xv) \dot{\xv}
		= \frac{\partial \hv(t,\xv)}{\partial t} + \Hm(t,\xv) \fv_{\textup{A}}^\prime(t,\xv) + \Hm(t,\xv) \uv^\prime.
	\end{equation}
	{%
	Inserting \eqref{eq:special gryphon, fbl v2} into \eqref{eq:Lie derivative of special gryphon} yields
	\begin{align}
		\dot{\yv} 
		&= \frac{\partial \hv(t,\xv)}{\partial t} 
		+ \Hm(t,\xv) \fv_{\textup{A}}^\prime(t,\xv) 
		+ \overbracket[0.14ex]{
			\Hm(t,\xv) \Hm(t,\xv)^{\#}_{\Qm^\prime}
		}^{=\,\Imat} \big( K_{\textup{y}}\big(\desired{\yv} - \yv\big) - \Hm(t,\xv)\fv_{\textup{A}}^\prime(t,\xv) \big) \nonumber\\
		&= \frac{\partial \hv(t,\xv)}{\partial t} 
		+ K_{\textup{y}}\big(\desired{\yv} - \yv\big). \label{eq:FBL special gryphon}
	\end{align}
	%where $\gammav(t,\varxv) := {\partial\hv(t,\varxv)}/{\partial t}$.
	Let \eqref{eq:FBL special gryphon} be rewritten in terms of error dynamics; this results in
	\begin{align}
		\widetilde{\yv}(t) &:= \yv(t) - \desired{\yv}(t), \\
		\wv(t) &:=
		\desired{\dot{\yv}}(t) - \frac{\partial \hv\big(\tau,\xv(t)\big)}{\partial \tau}\bigg|_{\tau=t}, \\
		\dot{\widetilde{\yv}}(t)
		&= - K_{\textup{y}} \widetilde{\yv}(t) - \wv(t),
	\end{align}
	%
%	\begin{align}
%		\dot{\widetilde{\yv}}
%		+ K_{\textup{y}} \widetilde{\yv} 
%		&= - \left(\desired{\dot{\yv}} - \frac{\partial \hv(t,\xv)}{\partial t}\right)
%	\end{align}
	%
	whose solution is
	\begin{empheq}[box=\colorbox{\colorfbc}]{equation}
	%\begin{equation}
		\label{eq:exponential convergence y}
		\widetilde{\yv}(t) = \exp(\hspace{-0.25ex}-\hspace{-0.25ex}K_{\textup{y}}t)\,\widetilde{\yv}(0)
		- \int_{0}^{t}\hspace{-1ex} \exp\big(\hspace{-0.25ex}-\hspace{-0.25ex}K_{\textup{y}}(t\hspace{-0.25ex}-\hspace{-0.25ex}\tau)\big)\,\wv(\tau)
		\,\mathrm{d}\tau.
	%\end{equation}
	\end{empheq}
	}
\end{proof}

Inserting \eqref{eq:special gryphon, fbl v2} into \eqref{eq:special gryphon} results in the closed-loop system
\begin{equation}\label{eq:special gryphon, fbl, closed loop}
%\begin{split}
	\dot{\xv} 
	= 
	\Hm(t,\xv)^{\#}_{\Qm^\prime} K_{\textup{y}}\big(\desired{\yv} - \hv(t,\xv)\big)
	+ \Big( \Imat - \Hm(t,\xv)^{\#}_{\Qm^\prime} \Hm(t,\xv) \Big)	\fv_{\textup{A}}^\prime(t,\xv).
%\end{split}
\end{equation}
By comparing \eqref{eq:special gryphon, fbl, closed loop} and \eqref{eq:barrier-method sequence},
it is easy to determine that
with
\begin{subequations}\label{eq:choice for elegant gryphon candidate}
	\begin{align}
		\Qm^\prime &= \Rm(t,\xv), \\
		\fv_{\textup{A}}^\prime(t,\xv) &= - K_{\textup{x}} \Rm(t,\xv)^{-1} \rv(t,\xv),
		\quad
		K_{\textup{x}} = K_{\textup{y}},
	\end{align}
\end{subequations}
\eqref{eq:special gryphon, fbl, closed loop} becomes
\begin{empheq}[box=\colorbox{\colorfbc}]{equation}
%\begin{equation}
	\label{eq:elegant gryphon candidate}
	%\vspace{0.25ex}
	\dot{\xv} 
	= K_{\textup{x}} \etav(t, \xv, \desired{\yv}),
%\end{equation}
\end{empheq}
which shares certain similarities with \eqref{eq:barrier-method sequence} in equilibrium points and convergence behavior, which are outlined as follows.

\begin{definition}[Equilibrium points at instant $t$]\label{definition:equilibrium x_ast}
	A point $\xv^{\ast} \in \realset^{\dimx}$ is called an \emph{equilibrium point} of \eqref{eq:elegant gryphon candidate} \emph{at instant $t$}, if and only if
	\begin{equation*}
		\nullv = K_{\textup{x}} \etav(t, \xv^{\ast}, \desired{\yv}).
	\end{equation*}
	In other words, any
	\begin{equation}\label{eq:set of equilibrium points, elegant gryphon candidate}
		\xv^{\ast} \in
		\mathcal{E}_{\textup{x}}(t,\desired{\yv}) := \big\{\,
		\varxv \in \realset^{\dimx} ~\big\lvert~ \nullv = K_{\textup{x}} \etav(t, \varxv, \desired{\yv})
		\,\big\}
	\end{equation}
	is an \emph{equilibrium point} of \eqref{eq:elegant gryphon candidate} \emph{at instant $t$}.
\end{definition}

\begin{lemma}\label{corollary:equibrium points are same}
	The differential equation \eqref{eq:elegant gryphon candidate} and the difference equation \eqref{eq:barrier-method sequence} have the same equilibrium points, i.e.,
	\begin{equation*}
		\mathcal{E}_{\textup{x}}(t,\desired{\yv}) = \mathcal{E}_{\chiv}(t,\desired{\yv})~~
		\forall t \geqslant 0.
	\end{equation*}
\end{lemma}	
\begin{proof}
	It can be concluded immediately by comparing
	\eqref{eq:set of equilibrium points, barrier-method sequence} and \eqref{eq:set of equilibrium points, elegant gryphon candidate}.
\end{proof}

\begin{lemma}[Local minimizers]\label{corollary:local minimizer is a  asymptotically stable equilibrium, elegant gryphon candidate}
	A local minimizer of \eqref{eq:constrained optimization}
	is an asymptotically stable equilibrium point of \eqref{eq:elegant gryphon candidate}.
\end{lemma}
\begin{proof}
	It can be concluded immediately from Lemmata~\ref{corollary:local minimizer is a  asymptotically stable equilibrium, barrier-method sequence} and \ref{corollary:equibrium points are same}.
\end{proof}

Evidently, \eqref{eq:choice for elegant gryphon candidate} places the closed-loop system's equilibrium points $\xv^{\ast}$ at the desired coordinates $\chiv^{\optimal}$. 
Therefore, \eqref{eq:choice for elegant gryphon candidate} serves as an \emph{equilibrium placement}.
{It must be noted that feedback linearization is used here as a facilitative tool in the deduction and it does not appear in the resulting equation~\eqref{eq:elegant gryphon candidate}; in fact, \eqref{eq:elegant gryphon candidate} is a nonlinear system.}

\subsection{Exponential convergence}\label{sec:exponential convergence}
\begin{lemma}\label{lemma:x to x_ast}%[$\xv \rightarrow \xv^{\ast}$]
	Consider the system~\eqref{eq:elegant gryphon candidate} with $K_{\textup{x}} > 0$.
	Suppose $\desired{\yv}(t) = \desired{\yv}(0)$~$\forall t \geqslant 0$.
	Suppose for all $t \geqslant 0$ and all $\xv \in \realset^{\dimx}$,
	$\sigma(t,\xv) = \sigma(0,\xv)$,
	$\hv(t,\xv) = \hv(0,\xv)$, 
	and $\gv(t,\xv) = \gv(0,\xv)$.
	%$\Gm(t) = \Gm(0)$, and $\cv(t) = \cv(0)$.
	Let 
	\begin{equation}\label{eq:ball about x_ast}		
		\mathcal{N}(\xv^{\ast},\delta) := \big\{\,\xv \in \realset^{\dimx}~\big\lvert~\|\xv - \xv^{\ast}\|_{2} < \delta\,\big\}.
	\end{equation}
	There exists $\delta: \realset^{\dimx} \to \realset_{>0}$ and $\xv^{\ast} \in \mathcal{E}_{\textup{x}}$ such that
	if %the initial condition 
	$\xv(0) \in \mathcal{N}\big(\xv^{\ast},\delta(\xv^{\ast})\big)$
	then $\xv(t)$ converges exponentially to $\xv^{\ast}$ as $t > 0$ increases.
	%	then $\lim\limits_{t\rightarrow+\infty} \big\|\xv(t) - \xv^{\ast}\big\|_{2} = 0$.
\end{lemma}
\begin{proof}
	%	Since $K_{\textup{y}} = K_{\textup{x}}$,
	%	Equation~\eqref{eq:elegant gryphon candidate} becomes
	%	\begin{equation*}
		%		\dot{\xv} = K_{\textup{x}} \etav(\xv,\nullv,\Qm).
		%	\end{equation*}
	We choose \adjustbox{bgcolor=\colorfbc}{$\xv^{\ast}  = \chiv^{\optimal}$} (this is possible due to Lemma~\ref{corollary:local minimizer is a  asymptotically stable equilibrium, elegant gryphon candidate})
	and $\delta$ to be such that $\mathcal{N}\big(\xv^{\ast},\delta(\xv^{\ast})\big) \cap \mathcal{E}_{\textup{x}}(0,\desired{\yv}) = \{\xv^{\ast}\}$ (such a $\delta(\xv^{\ast})$ exists due to Assumptions~\ref{assumption:isolated local minimizers}).
	Therefore, the time-invariant scalar function
	\begin{equation}\label{eq:Lyapunov function}
		V(\xv) = \frac{1}{2} (\xv - \xv^{\ast})^\transp K_{\textup{x}}^{-1} (\xv - \xv^{\ast})
		\begin{cases}
			= 0 &\text{if~} \xv = \xv^{\ast}, \\
			> 0 &\text{if~} \xv \in \mathcal{N}\big(\xv^{\ast},\delta(\xv^{\ast})\big)\setminus\{\xv^{\ast}\}		
		\end{cases}
	\end{equation}
	is a Lyapunov candidate function on $\mathcal{N}\big(\xv^{\ast},\delta(\xv^{\ast})\big)$.
	Its time derivative is
	\begin{equation}\label{eq:time derivative of Lyapunov function}
		\frac{\mathrm{d}}{\mathrm{d}t}V(\xv)
		= \frac{\partial V(\xv)}{\partial \xv} \dot{\xv}
		= (\xv - \xv^{\ast})^\transp \etav(t,\xv,\desired{\yv}).
	\end{equation}
	Note that $\xv^{\ast} = \textup{const.}$ due to the time-invariant $\sigma$, $\hv$, $\gv$ and constant $\desired{\yv}$.
	
	In order to analyze the bound of ${\mathrm{d}V(\xv)}/{\mathrm{d}t}$, we further let $\delta$ be such that 
	if $\|\chiv[0] - \chiv^{\optimal}\|_{2} < \delta(\chiv^{\optimal})$
	then $\chiv[k]$ of \eqref{eq:barrier-method sequence}
	converges quotient-linearly to a local minimizer $\chiv^{\optimal}$ of \eqref{eq:constrained optimization} as $k \in \naturalset$ increases (Theorem~\ref{theorem:barrier-method iteration}). 
	From this, with \eqref{eq:barrier-method sequence} and Definition~\ref{definition:q-linear convergence}, one concludes that
	$\exists M \in (0, 1)$:
	\begin{align}
		&\forall \chiv \in \mathcal{N}\big(\chiv^{\optimal},\delta(\chiv^{\optimal})\big),~ \nonumber\\
		&\quad\big(\chiv + \etav(t,\chiv,\desired{\yv}) - \chiv^{\optimal}\big)^\transp \big(\chiv + \etav(t,\chiv,\desired{\yv}) - \chiv^{\optimal}\big)
		\leqslant 
		\big(\chiv - \chiv^{\optimal}\big)^\transp M^{2} \big(\chiv - \chiv^{\optimal}\big) \nonumber \\
		&\quad\quad\Leftrightarrow
		(\chiv - \chiv^{\optimal})^\transp \etav(t,\chiv,\desired{\yv}) 
		\leqslant
		-\frac{1}{2}\big(1-M^{2}\big) \big\|\chiv - \chiv^{\optimal}\big\|_{2}^{2} - \frac{1}{2} \big\|\etav(t,\chiv,\desired{\yv})\big\|_{2}^{2}. \label{eq:qlin converg vs asymp converg, step 1}
	\end{align}
	Since $\xv^{\ast}  = \chiv^{\optimal}$ has been chosen,
	inequality~\eqref{eq:qlin converg vs asymp converg, step 1} can be written as	
	\begin{equation}\label{eq:qlin converg vs asymp converg, step 2}
		(\xv - \xv^{\ast})^\transp \etav(t,\xv,\desired{\yv}) 
		\leqslant 
		-\frac{1}{2}\big(1-M^{2}\big) \big\|\xv - \xv^{\ast}\big\|_{2}^{2} - \frac{1}{2} \big\|\etav(t,\xv,\desired{\yv})\big\|_{2}^{2}.
	\end{equation}
	Inserting \eqref{eq:qlin converg vs asymp converg, step 2} into \eqref{eq:time derivative of Lyapunov function} results in
	\begin{equation}\label{eq:qlin converg vs asymp converg, step 3}
		\frac{\mathrm{d}}{\mathrm{d}t}V(\xv)
		\leqslant 
		-\frac{1}{2}\big(1-M^{2}\big) \big\|\xv - \xv^{\ast}\big\|_{2}^{2} - \frac{1}{2} \big\|\etav(t,\xv,\desired{\yv})\big\|_{2}^{2}
	\end{equation}
	and consequently, ${\mathrm{d}V(\xv)}/{\mathrm{d}t} \leqslant 0$.
	This concludes the Lyapunov stability.
	Furthermore, ${\mathrm{d}} V(\xv)/{\mathrm{d}t}$ reaches zero if and only if $\xv$ reaches $\xv^{\ast}$, the only equilibrium point in $\mathcal{N}\big(\xv^{\ast},\delta(\xv^{\ast})\big)$.
	Therefore, the asymptotic convergence of $\xv(t)$ is concluded (LaSalle's invariance principle~\cite{LaSalle1960}).
	
	Since \eqref{eq:Lyapunov function} can be written as $2 K_{\textup{x}} V(\xv) = \|\xv - \xv^{\ast}\|_{2}^{2}$, inequality~\eqref{eq:qlin converg vs asymp converg, step 3} becomes
	\begin{equation*}%\label{eq:qlin converg vs asymp converg, step 4}
		\frac{\mathrm{d}}{\mathrm{d}t}V(\xv)
		\leqslant 
		-\frac{1}{2}\big(1-M^{2}\big) 2 K_{\textup{x}} V(\xv) - \frac{1}{2} \big\|\etav(t,\xv,\desired{\yv})\big\|_{2}^{2}
		\leqslant 
		-\frac{1}{2}\big(1-M^{2}\big) 2 K_{\textup{x}} V(\xv).
	\end{equation*}
	Solving this differential inequality using Gr\"{o}nwall's inequality~\cite{Gronwall1919} yields
	\begin{align*}%\label{eq:exponential convergence x}
		&V\big(\xv(t)\big) \leqslant V\big(\xv(0)\big) \exp\big(\hspace{-0.5ex}-\hspace{-0.5ex}(1\hspace{-0.5ex}-\hspace{-0.5ex}M^{2}) K_{\textup{x}} t\big) \\
		&~\Leftrightarrow~
		\big\|\xv(t) - \xv^{\ast}\big\|_{2} \leqslant \big\|\xv(0) - \xv^{\ast}\big\|_{2}\,\exp\big(\hspace{-0.25ex}-\hspace{-0.25ex}0.5 (1\hspace{-0.5ex}-\hspace{-0.5ex}M^{2}) K_{\textup{x}} t\big).
	\end{align*}
\end{proof}

Since $\xv^{\ast} = \chiv^{\optimal}$ is chosen as an arbitrary local minimizer $\chiv^{\optimal}$ explicitly in the proof of Lemma~\ref{lemma:x to x_ast}, 
one can replace $\xv^{\ast}$ with $\chiv^{\optimal}$ in Lemma~\ref{lemma:x to x_ast} and its proof.
Also, it is worth noting that, in general, $M$ in \eqref{eq:qlin converg vs asymp converg, step 1} depends on $\chiv^{\optimal}$, especially when there are multiple local minimizers.
This results in the following lemma:
\begin{lemma}\label{lemma:x to x_ast, v2}
	Consider the system~\eqref{eq:elegant gryphon candidate} with $K_{\textup{x}} > 0$.
	Suppose that the optimization problem~\eqref{eq:constrained optimization} is time-invariant, i.e., for all $t \geqslant 0$ and all $\xv \in \realset^{\dimx}$,
	$\sigma(t,\xv) = \sigma(0,\xv)$,
	$\hv(t,\xv) = \hv(0,\xv)$,
	$\desired{\yv}(t) = \desired{\yv}(0)$,
	and $\gv(t,\xv) = \gv(0,\xv)$.
	%$\Gm(t) = \Gm(0)$, and $\cv(t) = \cv(0)$.
	There exists a pair $\delta: \realset^{\dimx} \to \realset_{>0}$ and $M: \realset^{\dimx} \to (0, M_{\sup}]$ with $M_{\sup} < 1$
	such that,
	for any local minimizer $\chiv^{\optimal}$,
	if $\xv(0) \in \mathcal{N}\big(\chiv^{\optimal},\delta(\chiv^{\optimal})\big)$
	then %$\xv(t)$ converges exponentially to $\chiv^{\optimal}$ as $t > 0$ increases.
	%then $\lim\limits_{t\rightarrow+\infty} \|\xv(t) - \chiv^{\optimal}\|_{2} = 0$.
	\begin{empheq}[box=\colorbox{\colorfbc}]{equation}
%	\begin{equation}
		\label{eq:exponential convergence x}
		\big\|\xv(t) - \chiv^{\optimal}\big\|_{2} \leqslant \big\|\xv(0) - \chiv^{\optimal}\big\|_{2}\,\exp\Big(\hspace{-0.25ex}-\hspace{-0.25ex}0.5\, \big(1\hspace{-0.5ex}-\hspace{-0.5ex}M(\chiv^{\optimal})^{2}\big) K_{\textup{x}} t\Big).
%	\end{equation}
	\end{empheq}
\end{lemma}
%\begin{proof}
%	It is concluded immediately from Lemmata~\ref{corollary:local minimizer is a  asymptotically stable equilibrium, elegant gryphon candidate} and \ref{lemma:x to x_ast}.
%\end{proof}

Besides, from \eqref{eq:qlin converg vs asymp converg, step 1} and \eqref{eq:qlin converg vs asymp converg, step 2}, one concludes:
\begin{lemma}\label{lemma: inequality for Lyapunov}
	For the pair $\delta$ and $M$ in Lemma~\ref{lemma:x to x_ast, v2} 
	and any local minimizer $\chiv^{\optimal}(t)$ of \eqref{eq:constrained optimization} at instant $t$,
	if $\big\|\xv - \chiv^{\optimal}(t)\big\|_{2} < \delta\big(\chiv^{\optimal}(t)\big)$ then
	\begin{equation*}%
		\big(\xv - \chiv^{\optimal}(t)\big)^\transp \etav\big(t,\xv,\desired{\yv}(t)\big) 
		\leqslant
		-\frac{1}{2}\big(1-M\big(\chiv^{\optimal}(t)\big)^{2}\big) \big\|\xv - \chiv^{\optimal}(t)\big\|_{2}^{2} - \frac{1}{2} \big\|\etav\big(t,\xv,\desired{\yv}(t)\big)\big\|_{2}^{2}
	\end{equation*}
	$\forall t \geqslant 0$.
\end{lemma}\vspace{0.1ex}

Due to the exponential convergence stated in Lemma~\ref{lemma:x to x_ast, v2}, one may infer that similar convergence behavior holds even in the presence of time-variant $\sigma$, $\hv$, $\desired{\yv}$, and $\gv$, provided that $\chiv^{\optimal}(t)$ remains band-limited in the frequency domain {(equivalently, \smash{$\big\|\dot{\chiv}^{\optimal}(t)\big\|_{2} \leqslant D_{\chi}$} is uniformly bounded)}.
To formalize this inference, we present the following lemma.

\begin{lemma}\label{lemma:x to x_ast, v3}%[$\xv \rightarrow \xv^{\ast}$]
	Consider the system~\eqref{eq:elegant gryphon candidate} with $K_{\textup{x}} > 0$.
	Consider an arbitrary local minimizer $\chiv^{\optimal}: \realset_{\geqslant 0} \to \realset^{\dimx}$ of \eqref{eq:constrained optimization}.
	Suppose that
	\begin{equation}\label{eq:slowly varying chi_optimal}
		\big\|\dot{\chiv}^{\optimal}(t)\big\|_{2} \leqslant D_{\chi}~~
		\forall t \geqslant 0.
	\end{equation}
	There exists $\delta: \realset^{\dimx} \to \realset_{>0}$
	and $M: \realset^{\dimx} \to (0, M_{\sup}]$ with $M_{\sup} < 1$
	such that 
	if 
	%	\begin{subequations}
		\begin{align}
			\label{eq: x(0) is sufficiently close to chi_optimal(0)}
			\big\|\xv(0) - \chiv^{\optimal}(0)\big\|_{2} &< \delta\big(\chiv^{\optimal}(0)\big), \\
			\label{eq: x(t) is sufficiently close to chi_optimal(t)}
			\big\|\xv(0) - \chiv^{\optimal}(0)\big\|_{2} \exp (-\kappa t)
			+ \kappa^{-1} D_{\chi}
			&<
			\delta\big(\chiv^{\optimal}(t)\big)~~
			\forall t > 0
		\end{align}
		%	\end{subequations}
	then 
	\begin{empheq}[box=\colorbox{\colorfbc}]{equation}
	%\begin{equation}
		\label{eq:exponential convergence x v2.5} 
		\big\|\xv(t) \hspace{-0.25ex}-\hspace{-0.25ex} \chiv^{\optimal}(t)\big\|_{2}
		\leqslant
		\big\|\xv(0) \hspace{-0.25ex}-\hspace{-0.25ex} \chiv^{\optimal}(0)\big\|_{2} \exp (-\kappa t)
		+ \int_{0}^{t} \exp\big(\hspace{-0.25ex}-\hspace{-0.25ex}\kappa (t\hspace{-0.25ex}-\hspace{-0.25ex}\tau)\big) \big\|\dot{\chiv}^{\optimal}(\tau)\big\|_{2}\,\mathrm{d}\tau,
	%\end{equation}
	\end{empheq}
	where 
	\begin{equation}\label{eq:definition of kappa}
		%\kappa = K_{\textup{x}} \inf_{t\geqslant 0} \frac{1}{2} \Big(1-M\big(\chiv^{\optimal}(t)\big)^{2}\Big),
		%
		\kappa = \frac{1}{2} \big(1-M_{\sup}^{2}\big) K_{\textup{x}}.
	\end{equation}
\end{lemma}
\begin{proof}
	Let \eqref{eq:elegant gryphon candidate} be rewritten in terms of error dynamics;
	this results in 
	\begin{align}
		\widetilde{\xv} &= \xv - \chiv^{\optimal}, \label{eq:error of x} \\
		\dot{\widetilde{\xv}}
		&= K_{\textup{x}} \etav\big(t,
		\widetilde{\xv} + \chiv^{\optimal},
		\desired{\yv}
		\big)
		- \dot{\chiv}^{\optimal}. \label{eq:error dynamics of x}
	\end{align}
	This system can be considered as the nominal system
	\begin{equation}\label{eq:nominal system}
		\dot{\widetilde{\xv}}
		= K_{\textup{x}} \etav\big(t,
		\widetilde{\xv} + \chiv^{\optimal},
		\desired{\yv}
		\big)
	\end{equation}
	perturbed by the perturbation $\dot{\chiv}^{\optimal}$~\cite{Khalil2002}.
	For the nominal system~\eqref{eq:nominal system},
	we choose \eqref{eq:Lyapunov function} as the Lyapunov function, which, with respect to \eqref{eq:error of x}, becomes
	\begin{equation}\label{eq:Lyapunov function v2}
		U(\widetilde{\xv}) 
		= \frac{1}{2} \widetilde{\xv}^\transp K_{\textup{x}}^{-1} \widetilde{\xv}.
		%= \frac{1}{2 K_{\textup{x}}} \big\|\widetilde{\xv}\big\|_{2}^{2}
	\end{equation}
	
	In order to proceed, we introduce an Ansatz 
	\begin{equation}\label{eq:Ansatz x follows chi within a ball}
		\big\|\widetilde{\xv}(t)\big\|_{2} < \delta\big(\chiv^{\optimal}(t)\big)~~
		\forall t > 0,
	\end{equation}
	meaning that $\xv(t) = \widetilde{\xv}(t) + \chiv^{\optimal}(t)$ is always within the ball $\mathcal{N}\big(\chiv^{\optimal}(t),\delta({\small\chiv^{\optimal}(t)})\big)$ around $\chiv^{\optimal}(t)$. 
	At the end of this proof, we will show that this Ansatz is satisfied.
	
	With \eqref{eq: x(0) is sufficiently close to chi_optimal(0)} and Ansatz~\eqref{eq:Ansatz x follows chi within a ball},
	one can apply Lemma~\ref{lemma: inequality for Lyapunov} and obtains
	\begin{equation}
		\widetilde{\xv}^\transp \etav\big(t,
		\widetilde{\xv} + \chiv^{\optimal},
		\desired{\yv}
		\big)
		\leqslant
		-\frac{1}{2}\big(1-M(\chiv^{\optimal})^{2}\big) \big\|\widetilde{\xv}\big\|_{2}^{2} - \frac{1}{2} \big\|\etav\big(t,
		\widetilde{\xv} + \chiv^{\optimal},
		\desired{\yv}
		\big)\big\|_{2}^{2}.
	\end{equation}
	Consequently, 
	the time derivative
	\begin{equation*}%\label{eq:time derivative of Lyapunov function v2}
		\frac{\mathrm{d}}{\mathrm{d}t} U(\widetilde{\xv}) 
		= \frac{\partial U(\widetilde{\xv})}{\partial\widetilde{\xv}} \dot{\widetilde{\xv}}
		= \widetilde{\xv}^\transp \etav\big(t,
		\widetilde{\xv} + \chiv^{\optimal},
		\desired{\yv}
		\big)
		\leqslant
		0
	\end{equation*}
	is nonpositive.
	Analogously to the proof of Lemma~\ref{lemma:x to x_ast},
	one concludes in sequence
	the Lyapunov stability, asymptotic stability (LaSalle's invariance principle~\cite{LaSalle1960}), and exponential stability (Gr\"{o}nwall's inequality~\cite{Gronwall1919}) of the nominal system \eqref{eq:nominal system}'s equilibrium point $\widetilde{\xv}^{\ast} = \nullv$.
	
	Since the eqiulibrium points of the nominal system~\eqref{eq:nominal system} are exponentially stable and the corresponding Lyapunov function~\eqref{eq:Lyapunov function v2} satisfies 
	\begin{subequations}\label{eq:x to x_ast, v3, Lyapunov inequalities}
	\begin{align}
		b_{1} \big\|\widetilde{\xv}\big\|_{2}^{2}
		\leqslant
		U(\widetilde{\xv}) 
		\leqslant
		b_{2} \big\|\widetilde{\xv}\big\|_{2}^{2}, 
		&\quad b_{1} = b_{2} = \frac{1}{2} K_{\textup{x}}^{-1}, \\
		\frac{\mathrm{d}}{\mathrm{d}t} U(\widetilde{\xv}) 
		\leqslant
		-b_{3} \big\|\widetilde{\xv}\big\|_{2}^{2}, 
		&\quad
		%	b_{3} 
		%	= \inf_{t\geqslant 0} \frac{1}{2} \Big(1-M\big(\chiv^{\optimal}(t)\big)^{2}\Big), 
		b_{3} = \frac{1}{2} \big(1-M_{\sup}^{2}\big)\\
		\left\|\frac{\partial U(\widetilde{\xv})}{\partial\widetilde{\xv}}\right\|_{2}
		= K_{\textup{x}}^{-1} \big\|\widetilde{\xv}^\transp \big\|_{2}
		\leqslant b_{4} \big\|\widetilde{\xv}^\transp \big\|_{2},
		&\quad 
		b_{4} = K_{\textup{x}}^{-1}, \\
		\left\|\frac{\partial U(\widetilde{\xv})}{\partial \big[\,t~~\desired{\yv}^{\,\transp}~~\chiv^{\optimal\transp}\,\big]^\transp}\right\|_{2} = 0 \leqslant b_{5} \big\|\widetilde{\xv}^\transp \big\|_{2},
		&\quad
		b_{5} = 0,
	\end{align}
	\end{subequations}
	with constants $b_{1}$, $b_{2}$, $b_{3}$, $b_{4}$, $b_{5}$,
	the stability theorem for \emph{slowly varying systems} \cite{Khalil2002}~(see also Theorem~\ref{theorem:slowly varying system}) is applicable to the perturbed system~\eqref{eq:error dynamics of x} with slowly varying perturbation $\chiv^{\optimal}$~\eqref{eq:slowly varying chi_optimal} (cf.~Remark~\ref{remark:application of theorm for slowly varying systems}).
	This results in
	\begin{subequations}\label{eq:exponential convergence x v2+3}
		\begin{align}
			\big\|\widetilde{\xv}(t)\big\|_{2}
			&\leqslant
			\big\|\widetilde{\xv}(0)\big\|_{2} \exp (-\kappa t)
			+ \frac{b_{4}}{2 b_{1}} \int_{0}^{t} \exp\big(\hspace{-0.25ex}-\hspace{-0.25ex}\kappa (t\hspace{-0.25ex}-\hspace{-0.25ex}\tau)\big) \big\|\dot{\chiv}^{\optimal}(\tau)\big\|_{2}\,\mathrm{d}\tau \label{eq:exponential convergence x v2} \\
			&\leqslant
			\|\widetilde{\xv}(0)\|_{2} \exp (-\kappa t)
			+ \frac{b_{4}}{2 b_{1}} \frac{D_{\chi}}{\kappa}
			= \|\widetilde{\xv}(0)\|_{2} \exp (-\kappa t)
			+ \kappa^{-1} D_{\chi}, \label{eq:exponential convergence x v3}
		\end{align}
	\end{subequations}
	where $\kappa \in \realset_{>0}$ is defined in \eqref{eq:definition of kappa}.
	Note that \eqref{eq:exponential convergence x v2} is identical to \eqref{eq:exponential convergence x v2.5}.
	Furthermore, the inequality~\eqref{eq:exponential convergence x v3}, together with \eqref{eq: x(t) is sufficiently close to chi_optimal(t)}, confirms Ansatz~\eqref{eq:Ansatz x follows chi within a ball}.
\end{proof}

As expected, when $\chiv^{\optimal} = \textup{const.}$ ($\Leftrightarrow \dot{\chiv}^{\optimal} = \nullv$), Lemma~\ref{lemma:x to x_ast, v3} degenerates to Lemma~\ref{lemma:x to x_ast, v2}, matching the desired behavior in Problem~\ref{problem:find a gryphonheart}.1. 
Furthermore, the upper bound~\eqref{eq:exponential convergence x v3} of $\big\|\widetilde{\xv}(t)\big\|_{2}$ converges exponentially to $\kappa^{-1} D_{\chi} = K_{\textup{x}}^{-1} \big(1-M_{\sup}^{2}\big)^{-1} 2 D_{\chi}$ as $t$ increases.
Since the constants $M_{\sup}$ and $D_{\chi}$ are related to the problem instance only, one concludes that the behavior in Problem~\ref{problem:find a gryphonheart}.2 is addressed.

{%
%	\vspace{-5ex}
\begin{remark}%
	``Slowly'' is a relative notion: any finite $D_{\chi} \geqslant 0$, even very large, satisfies the definition. 
	Since all systems within classical physics evolve with finite rates, this \emph{slowly-varying-system} assumption holds universally. %for meaningful physical systems
\end{remark}
}%

One may notice that Condition~\eqref{eq: x(t) is sufficiently close to chi_optimal(t)} does not seem to be trivial, in contrast to the other two conditions~\eqref{eq:slowly varying chi_optimal} and \eqref{eq: x(0) is sufficiently close to chi_optimal(0)}.
Let us start considering the worst case where $M_{\sup}$ is very close to~1. It is still manageable to make $\kappa \gg \max\{1, D_{\chi}\}$ by choosing a sufficiently large $K_{\textup{x}} > 0$.
As a result, the right-hand side of \eqref{eq:exponential convergence x v2+3} converges virtually instantaneously.
Thus, one may express \eqref{eq:exponential convergence x v2+3} as
\begin{equation*}
	\forall t \gg \kappa^{-1} \approx 0^{+},~
	%		\big\|\xv(t) - \chiv^{\optimal}(t)\big\|_{2} =
	\big\|\widetilde{\xv}(t)\big\|_{2}
	\lesssim
	%	\frac{b_{4}}{2 b_{1} \kappa} \big\|\dot{\chiv}^{\optimal}(t)\big\|_{2}
	%	= 
	\kappa^{-1} \big\|\dot{\chiv}^{\optimal}(t)\big\|_{2}
	\leqslant 
	\kappa^{-1} D_{\chi}
	\approx 0^{+},
\end{equation*}
indicating a virtually immediate closeness of $\xv(t)$ to $\chiv^{\optimal}(t)$.
Therefore, the inequality~\eqref{eq: x(t) is sufficiently close to chi_optimal(t)} is satisfied from a practical point of view.
Furthermore, since $\xv(t)$ is virtually immediately close to $\chiv^{\optimal}(t)$, 
the supremum $M_{\sup}$ of convergence rate
%the rate $M\big(\chiv^{\optimal}(t)\big) \in (0,1)$ of convergence (see Definition~\ref{definition:q-linear convergence}) 
is small (see Lemma~\ref{lemma:quasi Newton linear convergence}).
A small $M_{\sup}$ leads to a larger $\kappa$, see \eqref{eq:definition of kappa}. 
In summary, we conclude the following remark.

\begin{remark}\label{remark: sufficiently large Kx ensures the closeness for all t > 0}
A larger $K_{\textup{x}}$ practically increases the likelihood of satisfying \eqref{eq: x(t) is sufficiently close to chi_optimal(t)}
and a sufficiently large $K_{\textup{x}} > 0$ practically makes \eqref{eq: x(t) is sufficiently close to chi_optimal(t)} satisfied.
\end{remark}
Consequently, we have the following theorem.

\begin{theorem}\label{theorem:optimal effort distributor}
	Consider the time-variant system~\eqref{eq:gryphon} with an initial condition $\xv(0)$ for state and a reference $\desired{\yv}(t)$ for output.
	Consider the time-variant optimization problem~\eqref{eq:constrained optimization} with a differentiable trajectory $\chiv^{\optimal}: \realset_{\geqslant 0} \to \realset^{\dimx}$ of local minimizer.
	Consider some constant $K_{\textup{x}} \in \realset_{>0}$.
	Suppose that Assumptions~\ref{assumption:isolated local minimizers} hold.
	If \eqref{eq:slowly varying chi_optimal}, \eqref{eq: x(0) is sufficiently close to chi_optimal(0)}, and \eqref{eq: x(t) is sufficiently close to chi_optimal(t)} are satisfied,
	then 
	\begin{empheq}[box=\colorbox{\colorfbc}]{equation}
%	\begin{equation}
		\label{eq:gryphonheart}
		\uv 
		= \muv\big(t, \xv, \desired{\yv}; K_{\textup{x}}\big)
		= \Bm(\xv)^{-1} \big(
		-\fv_{\textup{A}}(\xv) 
		+ K_{\textup{x}} \etav(t, \xv, \desired{\yv})
		\big)
		%
		%	\quad
		%	K_{\textup{x}} > 0
%	\end{equation}
	\end{empheq}
	with $\etav$ defined in \eqref{eq:definition of eta}\eqref{eq:definition of eta_A}\eqref{eq:definition of eta_B}
%	\begin{equation}
%		\begin{split}
%			\uv 
%			= \muv\big(t, \xv, \desired{\yv}; K_{\textup{x}}\big)
%			%
%			= \Bm(\xv)^{-1} \big(
%			-\fv_{\textup{A}}(\xv) 
%			&+ K_{\textup{x}} \Hm(t,\xv)^{\#}_{\Rm(t,\xv)} \big(
%			\desired{\yv}(t) - \hv(t,\xv)
%			\big) \\
%			&+ K_{\textup{x}} \Big(
%			\Imat - \Hm(t,\xv)^{\#}_{\Rm(t,\xv)} \Hm(t,\xv)
%			\Big) \Rm(t,\xv)^{-1} \rv(t,\xv)
%			\big)
%		\end{split}
%	\end{equation}
	solves Problem~\ref{problem:find a gryphonheart}.
	%is an \OTC~(Definition~\ref{definition:gryphonheart}).
	The resulting closed-loop trajectory $\xv(t)$ follows $\chiv^{\optimal}(t)$ in an exponentially convergent manner described by~\eqref{eq:exponential convergence x v2.5}.
\end{theorem}
\begin{proof}
	Equations~\eqref{eq:gryphon, fbl-2}, \eqref{eq:special gryphon, fbl v2}, and \eqref{eq:choice for elegant gryphon candidate} constitute \eqref{eq:gryphonheart}.
	Lemma~\eqref{lemma:x to x_ast, v3} delineates the convergence behavior of the closed-loop system \eqref{eq:gryphon}\eqref{eq:gryphonheart}.
\end{proof}

Although Theorem~\ref{theorem:optimal effort distributor} solves Problem~\ref{problem:find a gryphonheart}, it requires the initial value $\xv(0)$ of system's state to be sufficiently close to a local minimizer $\chiv^{\optimal}(0)$ of \eqref{eq:constrained optimization} at the initial instant $t = 0$, cf. \eqref{eq: x(0) is sufficiently close to chi_optimal(0)}.

\begin{remark}\label{remark:sufficient closeness}
	In all lemmata and theorems above, 
	the notation $\mathcal{N}(\chiv^{\optimal},\delta) = \big\{\,\xv \in \realset^{\dimx}~\big\lvert~\|\xv - \chiv^{\optimal}\|_{2} < \delta\,\big\}$ shall be read as ``the interior of a ball of radius $\delta$ and center $\chiv^{\optimal}$''
	rather than ``an open neighborhood of $\chiv^{\optimal}$''.
	Hence, $\xv(0) \in \mathcal{N}(\chiv^{\optimal},\delta)$ or, equivalently, $\big\|\xv(0) - \chiv^{\optimal}\big\|_{2} < \delta$
	shall be read as ``$\xv(0)$ being \emph{sufficiently close} to $\chiv^{\optimal}$'' rather than ``$\xv(0)$ being very close to $\chiv^{\optimal}$''.
\end{remark}

It needs to be emphasized that $\xv^{\ast}$ and $\delta$ are chosen as pairs in the proof of Lemma~\ref{lemma:x to x_ast} to specify
\begin{itemize}
	\item the local minimizer (denoted $\chiv^{\optimal}$) that is being concerned (although it can be any local minimizer), and
	\item the region in which the quotient-linear convergence of \eqref{eq:barrier-method sequence} is guaranteed (Theorem~\ref{theorem:barrier-method iteration}).	
\end{itemize}
The chosen pairs of $\xv^{\ast}$ and $\delta$ display a sufficient condition for the exponential convergence of $\xv(t)$ as stated in Lemmata~\ref{lemma:x to x_ast}, \ref{lemma:x to x_ast, v2}, and \ref{lemma:x to x_ast, v3},
This condition does not exclude other regions in $\realset^{\dimx}$ that ensure this kind of convergence.

However, the above statement does not address whether there exist points that lie outside the basin of attraction for any equilibrium point of \eqref{eq:elegant gryphon candidate}.
For exploring this further,
insights from researches on Newton fractals \cite{CurryGarSul1983,Beardon1991,EpureanuGre1998,Devaney2010} become relevant, considering that \eqref{eq:elegant gryphon candidate} originates from a quasi-Newton method~\eqref{eq:barrier-method sequence}.
State-of-the-art research in this field has primarily focused on one-dimensional complex polynomials \cite{CurryGarSul1983,Beardon1991,EpureanuGre1998,Devaney2010}. These works have revealed that the basin boundaries exhibit a fractal structure, dividing the complex plane into three types of regions: i) regions where the Newton-Raphson iteration converges to a root, ii) regions where the iteration converges to a periodic cycle, and iii) regions where no convergence occurs.
While a one-dimensional complex polynomial can be expressed as a two-dimensional real polynomial, generalizing these findings to multidimensional polynomials is not trivial.
Nevertheless, we hypothesize that certain properties of fractal basin boundaries in the two-dimensional case may extend to more general multidimensional polynomials.
Building on this hypothesis and the Stone-Weierstrass approximation theorem~\cite{Stone1948}, which states that any continuous function of multiple variables can be approximated by a polynomial,
we propose that analogous fractal structures and convergence behaviors may also be present in the case with multidimensional real-valued continuously differentiable functions as in \eqref{eq:Lagrangian's Jacobian}.
Accordingly, we posit the following proposition:%
\begin{proposition}\label{proposition:newton fractal}
	The real domain $\realset^{\dimx} = \mathcal{X}_{1} \cup \mathcal{X}_{2} \cup \mathcal{X}_{3}$ is composed of three complementary domains $\mathcal{X}_{1}$, $\mathcal{X}_{2}$, and $\mathcal{X}_{3}$ such that
	\begin{enumerate}
		\item If $\chiv[0] \in \mathcal{X}_{1}$, then $\chiv[k]$ of the interior-point iteration~\eqref{eq:barrier-method sequence} converges to one of the local minimizers of \eqref{eq:constrained optimization}.
		\item If $\chiv[0] \in \mathcal{X}_{2}$, then $\chiv[k]$ converges to a periodic cycle.
		The points that form these cycles (summarized in $\mathcal{X}_{2}^{\prime} \subseteq \mathcal{X}_{2}$), due to the discrete natural of the sequence $\chiv[k]$, are isolated points in $\realset^{\dimx}$.
		\item If $\chiv[0] \in \mathcal{X}_{3}$, then $\chiv[k]$ does not converge at all. The domain $\mathcal{X}_{3}$ (Julia set) is a closed set and is nowhere dense in $\realset^{\dimx}$.
	\end{enumerate}
\end{proposition}

Proposition~\ref{proposition:newton fractal} accounts for the behavior of the discrete-time system~\eqref{eq:barrier-method sequence}.
Even though being related to \eqref{eq:barrier-method sequence}, the continuous-time system~\eqref{eq:elegant gryphon candidate} is able to overcome the undesired Cases 2 and~3:
\begin{itemize}
	\item Due to the continuous natural of \eqref{eq:elegant gryphon candidate}, $\dot{\xv}(t)$ always drives the point $\xv(t)$ into the neighborhood of $\xv(t)$. This means, if $\xv(t) \in \mathcal{X}_{2}^{\prime}$, it is able to leave $\mathcal{X}_{2}^{\prime}$. Remember that $\mathcal{X}_{2}^{\prime}$ is composed of isolated points.
	\item Since $\xv(t)$ converges to some hyperplane
	%$\{\,\chiv \in \realset^{\dimx}~\lvert~\nullv = \big(\hv(t,\chiv) - \desired{\yv}(t)\big)\,\}$ 
	from \emph{any} initial condition (see Lemma~\ref{lemma:h(x) converges to zero, v2}),	
	Case 3 does not exist for \eqref{eq:elegant gryphon candidate}.
\end{itemize}
These observations conclude the following theorem:
\begin{theorem}\label{lemma:any x is sufficiently close}
	If Proposition~\ref{proposition:newton fractal} is true, then \adjustbox{bgcolor=\colorfw}{any $\xv(0) \in \realset^{\dimx}$ is sufficiently close} (as explained in Remark~\ref{remark:sufficient closeness}) to at least one local minimizer of \eqref{eq:constrained optimization} at $t=0$, meeting the requirement for $\xv(0)$ specified as \eqref{eq: x(0) is sufficiently close to chi_optimal(0)} in Theorem~\ref{theorem:optimal effort distributor}.
\end{theorem}

We leave rigorously proofing Proposition~\ref{proposition:newton fractal} (and consequently, Theorem~\ref{lemma:any x is sufficiently close}) for future work.

\subsection{Isomorphism between optimization and optimal tracking}\label{sec:optimal tracking controller}
It is obvious that the optimal tracking control problem~\eqref{eq:optimal control problem} exhibits a high similarity with the time-variant pseudoconvex optimization problem~\eqref{eq:constrained optimization}.
The result $\phiv^{\optimal}: [0, t_{\final}] \to \realset^{\dimx}$ (a local minimizer) of \eqref{eq:optimal control problem} can be used to constitute the control 
%$\uv: \realset_{\geqslant 0} \to \realset^{\dimx}$ 
\begin{equation*}
	\uv = \Bm(\phiv^{\optimal})^{-1} \big(\dot{\phiv}^{\optimal} - \fv_{\textup{A}}(\phiv^{\optimal})\big)
\end{equation*}
for the time-variant system~\eqref{eq:gryphon}
such that its state $\xv(t)$ behaves optimally, i.e., $\xv(t) = \phiv^{\optimal}(t)~\forall t \geqslant 0$.

Numerical methods may be applied to obtain an approximated result which is close to a local minimizer. 
In doing so, one converts the cost functional to a cost function in discrete time; that is
\begin{equation}\label{eq:optimal control cost, rectangle rule}
	\int_{0}^{t_{\text{f}}} \sigma\big(t,\phiv(t)\big)\,\mathrm{d}t
	\approx 
	\sigma\big(0, \phiv(0)\big) \dt
	+ \sum_{n=1}^{N} \sigma\big(n\dt, \phiv(n\dt)\big) \dt
\end{equation}
using the rectangle rule,
where $\dt \in \{\,\delta \in \realset_{>0}~\lvert~t_{\final}/\delta \in \naturalset^{+}~\wedge~t_{\final}/\delta \gg 1\,\}$ is a small positive number.
Since $\dt$ is a positive constant and $\phiv(0) = \xv(0)$ (cf.~\eqref{eq:optimal control problem}) is determined, the cost function \eqref{eq:optimal control cost, rectangle rule} can be replaced by
%$\varsigma: \realset^{\dimx N} \to \realset_{\geqslant 0}$
\begin{align*}
	\varsigma(\varphiv) :=\sum_{n=1}^{N} \sigma\big(n\dt, \phiv(n\dt)\big).
\end{align*}
in the context of minimization.
Herein,
\begin{align*}
	\varphiv := \big[\,\phiv(1\dt)^\transp~~\phiv(2\dt)^\transp~~\cdots~~\phiv(N\dt)^\transp\,\big]^\transp \in \realset^{\dimx N}
\end{align*}
is the discrete representation of $\phiv: \realset_{> 0}\to\realset^{\dimx}$.
Ergo, the optimal tracking control problem~\eqref{eq:optimal control problem} can be approximately written as the nonlinear optimization problem
\begin{subequations}\label{eq:optimal control convex optimization}
%	\begin{empheq}[box=\colorbox{\colorotc}]{equation}
	\begin{equation}
		\min_{\varphiv} \varsigma(\varphiv)
	\end{equation}
%	\end{empheq}
	subject to
%	\begin{empheq}[box=\colorbox{\colorotc}]{align}
	\begin{align}
		\nullv &= \varhv(\varphiv) - \desired{\varyv}
		\in \realset^{\dimy N}, 
		\quad
		\varhv(\varphiv) := \begin{bmatrix}
			\hv\big(1\dt,\phiv(1\dt)\big) \\
			\hv\big(2\dt,\phiv(2\dt)\big) \\
			\vdots \\
			\hv\big(N\dt,\phiv(N\dt)\big)
		\end{bmatrix}, \quad
		\desired{\varyv} := \begin{bmatrix}
			\desired{\yv}(1\dt) \\
			\desired{\yv}(2\dt) \\
			\vdots \\
			\desired{\yv}(N\dt)
		\end{bmatrix}, \\
		\nullv &\geqslant \vargv(\varphiv) := \begin{bmatrix}
			\gv\big(1\dt,\phiv(1\dt)\big) \\
			\gv\big(2\dt,\phiv(2\dt)\big) \\
			\vdots \\
			\gv\big(N\dt,\phiv(N\dt)\big)
		\end{bmatrix}
		\in \realset^{\dimc N}.
	\end{align}
%	\end{empheq}
\end{subequations}
Since $\dt$ can be chosen arbitrarily close to zero, the problems~\eqref{eq:optimal control problem} and \eqref{eq:optimal control convex optimization} are isomorphic to each other.

It is reasonable to believe that the trajectory $\chiv^{\optimal}: \realset_{\geqslant 0} \to \realset^{\dimx}$ with $\chiv^{\optimal}(t)$ being a local minimizer of \eqref{eq:constrained optimization} at time instant $t$ is a local minimizer $\phiv^{\optimal}$ of the optimal control problem~\eqref{eq:optimal control problem}, i.e., $\chiv^{\optimal}(t) = \phiv^{\optimal}(t)$~$\forall t \in (0, t_{\final}]$.
Since \eqref{eq:optimal control problem} is not always analytically solvable, we analyze its discrete-time approximation \eqref{eq:optimal control convex optimization} and present its equivalency to \eqref{eq:constrained optimization}.

\begin{lemma}\label{lemma:equivalency optimal control and convex programming}
	Let $\mathcal{E}^{\optimal}(t) \subset \realset^{\dimx}$ be the set of all local minimizers of the optimization problem~\eqref{eq:constrained optimization} at instant $t$.
	\adjustbox{bgcolor=\colorotc}{For} any local minimizer 
	\begin{equation*}
		\varphiv^{\optimal} = \begin{bmatrix}
			\phiv^{\optimal}(1\dt)^\transp &
			\phiv^{\optimal}(2\dt)^\transp &
			\cdots &
			\phiv^{\optimal}(N\dt)^\transp
		\end{bmatrix}^\transp
	\end{equation*}
	of the optimization problem~\eqref{eq:optimal control convex optimization}
	%any $\dt \in \{\,\delta \in \realset_{>0}~\lvert~t_{\final}/\delta \in \naturalset^{+}~\wedge~t_{\final}/\delta \gg 1\,\}$, 
	and \adjustbox{bgcolor=\colorotc}{any $n \in \naturalset^{+} \cap [1, N]$},
	there exists $\chiv^{\optimal}(n\dt) \in \mathcal{E}^{\optimal}(n\dt)$ such that \adjustbox{bgcolor=\colorotc}{$\chiv^{\optimal}(n\dt) = \phiv^{\optimal}(n\dt)$}.
\end{lemma}
\begin{proof}
For later use, we derive the partial derivatives of $\varsigma$, $\varhv$, and $\vargv$ as follows.
\begin{align*}
	\varqv(\varphiv) 
	:= \left[\frac{\partial \varsigma(\varphiv)}{\partial \varphiv}\right]^\transp
	\hspace{-0.5ex}= \begin{bmatrix}
		\qv\big(1\dt, \phiv(1\dt)\big)^\transp &
		\qv\big(2\dt, \phiv(2\dt)\big)^\transp &
		\cdots &
		\qv\big(N\dt, \phiv(N\dt)\big)^\transp
	\end{bmatrix}^\transp, \quad
	%
	%\qv(t,\xv) := \left[\frac{\partial \sigma(t,\xv)}{\partial \xv}\right]^\transp,
\end{align*}
\begin{align*}
	\varQm(\varphiv) 
	:= \frac{\partial \varqv(\varphiv) }{\partial \varphiv} 
	= \begin{bmatrix}
		\Qm\big(1\dt,\phiv(1\dt)\big) & \nullmat & \vphantom{\ddots}\cdots & \nullmat \\
		\nullmat & \Qm\big(2\dt,\phiv(2\dt)\big) & \ddots & \nullmat \\
		\vdots & \nullmat & \ddots & \nullmat \\
		\nullmat & \nullmat & \vphantom{\ddots}\cdots & \Qm\big(N\dt,\phiv(N\dt)\big)
	\end{bmatrix}, \quad
	%\Qm(t,\xv) := \frac{\partial \qv(t,\xv)}{\partial \xv}.
\end{align*}
\begin{align}\label{eq:optimal-ctrl, Jacobian of equality constraint}
	\varHm(\varphiv) 
	:= \frac{\partial \varhv(\varphiv)}{\partial \varphiv}
	= \begin{bmatrix}
		\Hm\big(1\dt,\phiv(1\dt)\big) & \nullmat & \vphantom{\ddots}\cdots & \nullmat \\
		\nullmat & \Hm\big(2\dt,\phiv(2\dt)\big) & \ddots & \nullmat \\
		\vdots & \nullmat & \ddots & \nullmat \\
		\nullmat & \nullmat & \vphantom{\ddots}\cdots & \Hm\big(N\dt,\phiv(N\dt)\big)
	\end{bmatrix}, \quad
	%
	%\Hm(t,\xv) := \frac{\partial \hv(t,\xv)}{\partial \xv}.
\end{align}
\begin{align*}
	\varGm 
	:= \frac{\partial \vargv(\varphiv)}{\partial \varphiv}
	= \begin{bmatrix}
		\Gm(1\dt) & \nullmat & \vphantom{\ddots}\cdots & \nullmat \\
		\nullmat & \Gm(2\dt) & \ddots & \nullmat \\
		\vdots & \nullmat & \ddots & \nullmat \\
		\nullmat & \nullmat & \vphantom{\ddots}\cdots & \Gm(N\dt)
	\end{bmatrix},
\end{align*}
where $\qv(t,\xv)$, $\Hm(t,\xv)$, and $\Qm(t,\xv)$ are given in \eqref{eq:definition of q and xi}, \eqref{eq:definition of H}, and \eqref{eq:definition of Q and Xi}, respectively.

Let the optimization problem~\eqref{eq:optimal control convex optimization} be solved by the interior-point iteration~\eqref{eq:barrier-method sequence}; it means that 
\begin{align}
	\varphiv^{[k]} &= \begin{bmatrix}
		\phiv^{[k]}(1\dt)^\transp &
		\phiv^{[k]}(2\dt)^\transp &
		\cdots &
		\phiv^{[k]}(N\dt)^\transp
	\end{bmatrix}^\transp, \nonumber\\
	\begin{split}
	\varphiv^{[k+1]} 
	&= \varphiv^{[k]} 
	+ \varHm(\varphiv^{[k]})^{\#}_{\varRm(\varphiv^{[k]})} \big(\desired{\varyv} - \varhv(\varphiv^{[k]})\big) \\
	&\qquad\quad- \Big(\Imat - \varHm(\varphiv^{[k]})^{\#}_{\varRm(\varphiv^{[k]})}\varHm(\varphiv^{[k]})\Big) \varRm(\varphiv^{[k]})^{-1} \varrv(\varphiv^{[k]})
	\end{split} \label{eq:optimal control, barrier-method sequence}
\end{align}
converges to a local minimizer $\varphiv^{\optimal}$ as $k \in \naturalset$ increases (Theorem 1).
Herein, $\varrv(\varphiv) := \big[{\partial \varrho(\varphiv)}/{\partial \varphiv}\big]^\transp $ and $\varRm(\varphiv) := {\partial \varrv(\varphiv)}/{\partial \varphiv}$ are the gradient and the Hessian, respectively, of
the augmented cost function
\begin{align*}
	\varrho(\varphiv) = \varsigma(\varphiv) + \frac{1}{2} \betav\big(\vargv(\varphiv)\big)^\transp \betav\big(\vargv(\varphiv)\big),
\end{align*}
cf.~Definition~\ref{definition:barrier fcn}.
Easily, one can deduct
\begin{align*}
	\varrv(\varphiv) 
	%:= \left[\frac{\partial \varrho(\varphiv)}{\partial \varphiv}\right]^\transp 
	&= \varqv(\varphiv) + \varGm^\transp \xiv\big(\vargv(\varphiv)\big), 
	\quad
	\xiv(\sv) := \frac{\partial \betav(\sv)}{\partial \sv}^\transp \betav(\sv) = \left[\frac{\partial \beta(s_{i})}{\partial s_{i}} \beta(s_{i})\right], \\
	\varRm(\varphiv) 
	%:= \frac{\partial \varrv(\varphiv)}{\partial \varphiv}
	&= \varQm(\varphiv) + \varGm^\transp \Xim\big(\vargv(\varphiv)\big) \varGm, 
	\quad
	\Xim(\sv) := \frac{\partial \xiv(\sv)}{\partial \sv}.
\end{align*}
Herein, $\xiv: \realset^{\dims} \to \realset^{\dims}$ is an element-wise function,
and $\Xim: \realset^{\dims} \to \realset^{\dims\times\dims}$ returns a diagonal and positive semidefinite matrix.
This leads to
\begin{align*}
	\xiv\big(\vargv(\varphiv)\big) = \begin{bmatrix}
		\xiv\Big(\gv\big(1\dt, \phiv(1\dt)\big)\Big) \\
		\xiv\Big(\gv\big(2\dt, \phiv(2\dt)\big)\Big) \\
		\vdots \\
		\xiv\Big(\gv\big(N\dt, \phiv(N\dt)\big)\Big)
	\end{bmatrix},
\end{align*}
\begin{align*}
	\Xim\big(\vargv(\varphiv)\big) 
	=
	\begin{bmatrix}
		\Xim\big(\gv(1\dt, \phiv(1\dt))\big) & \nullmat & \vphantom{\ddots}\cdots & \nullmat \\
		\nullmat & \Xim\big(\gv(2\dt, \phiv(2\dt))\big) & \ddots & \nullmat \\
		\vdots & \nullmat & \ddots & \nullmat \\
		\nullmat & \nullmat & \vphantom{\ddots}\cdots & \Xim\big(\gv(N\dt, \phiv(N\dt))\big)
	\end{bmatrix}.
\end{align*}
Consequently, one obtains
\begin{align*}
	\varrv(\varphiv) 
	= \begin{bmatrix}
		\rv\big(1\dt,\phiv(1\dt)\big) \\
		\rv\big(2\dt,\phiv(2\dt)\big) \\
		\vdots \\
		\rv\big(N\dt,\phiv(N\dt)\big)
	\end{bmatrix}, \quad
	%
	%\rv(t,\xv) := \qv(t, \xv) + \Gm(t)^\transp \xiv\big(\gv(t, \xv)\big)
\end{align*}
\begin{align}\label{eq:optimal-ctrl, Hessian of augmented cost}
	\varRm(\varphiv) 
	=
	\begin{bmatrix}
		\Rm\big(1\dt,\phiv(1\dt)\big) & \nullmat & \vphantom{\ddots}\cdots & \nullmat \\
		\nullmat & \Rm\big(2\dt,\phiv(2\dt)\big) & \ddots & \nullmat \\
		\vdots & \nullmat & \ddots & \nullmat \\
		\nullmat & \nullmat & \vphantom{\ddots}\cdots & \Rm\big(N\dt,\phiv(N\dt)\big)
	\end{bmatrix}, \quad
	%
	%\Rm(t,\xv) := \Qm(t,\xv) + \Gm(t)^\transp \Xim\big(\gv(t,\xv)\big) \Gm(t),
\end{align}
and that $\varRm(\varphiv)$ is a block diagonal matrix where each block matrix $\Rm(t,\xv)$ is SPD.
Herein, $\rv(t,\xv)$ and $\Rm(t,\xv)$ are given in \eqref{eq:definition of r} and \eqref{eq:definition of R}, respectively.

Since the weighted right pseudoinverse $\varHm(\varphiv)^{\#}_{\varRm(\varphiv)}$ is the key component of \eqref{eq:optimal control, barrier-method sequence},
we explore firstly its properties. 
By inserting
\eqref{eq:optimal-ctrl, Jacobian of equality constraint} and \eqref{eq:optimal-ctrl, Hessian of augmented cost} into $\varHm(\varphiv)^{\#}_{\varRm(\varphiv)}$ that is defined in \eqref{eq:weighted right pseudoinverse}, one obtains
\begin{align*}
	&\varHm(\varphiv)^{\#}_{\varRm(\varphiv)}
	= \varRm(\varphiv)^{-1} \varHm(\varphiv)^\transp \Big(\varHm(\varphiv) \varRm(\varphiv)^{-1} \varHm(\varphiv)^\transp\Big)^{-1} \\
	&= \resizebox{0.95\textwidth}{!}{$
	\begin{bmatrix}
		\Rm_{\langle 1\rangle} & \nullmat & \vphantom{\ddots}\cdots & \nullmat \\
		\nullmat & \Rm_{\langle 2\rangle} & \ddots & \nullmat \\
		\vdots & \nullmat & \ddots & \nullmat \\
		\nullmat & \nullmat & \vphantom{\ddots}\cdots & \Rm_{\langle N\rangle}
	\end{bmatrix}^{-1} 
	\begin{bmatrix}
		\Hm_{\langle 1\rangle} & \nullmat & \vphantom{\ddots}\cdots & \nullmat \\
		\nullmat & \Hm_{\langle 2\rangle} & \ddots & \nullmat \\
		\vdots & \nullmat & \ddots & \nullmat \\
		\nullmat & \nullmat & \vphantom{\ddots}\cdots & \Hm_{\langle N\rangle}
	\end{bmatrix}^\transp
	\left(
	\begin{bmatrix}
		\Hm_{\langle 1\rangle} & \nullmat & \vphantom{\ddots}\cdots & \nullmat \\
		\nullmat & \Hm_{\langle 2\rangle} & \ddots & \nullmat \\
		\vdots & \nullmat & \ddots & \nullmat \\
		\nullmat & \nullmat & \vphantom{\ddots}\cdots & \Hm_{\langle N\rangle}
	\end{bmatrix}
	\begin{bmatrix}
		\Rm_{\langle 1\rangle} & \nullmat & \vphantom{\ddots}\cdots & \nullmat \\
		\nullmat & \Rm_{\langle 2\rangle} & \ddots & \nullmat \\
		\vdots & \nullmat & \ddots & \nullmat \\
		\nullmat & \nullmat & \vphantom{\ddots}\cdots & \Rm_{\langle N\rangle}
	\end{bmatrix}^{-1} 
	\begin{bmatrix}
		\Hm_{\langle 1\rangle} & \nullmat & \vphantom{\ddots}\cdots & \nullmat \\
		\nullmat & \Hm_{\langle 2\rangle} & \ddots & \nullmat \\
		\vdots & \nullmat & \ddots & \nullmat \\
		\nullmat & \nullmat & \vphantom{\ddots}\cdots & \Hm_{\langle N\rangle}
	\end{bmatrix}^\transp
	\right)^{-1}
	$} \\
	&= \resizebox{0.95\textwidth}{!}{$
		\begin{bmatrix}
		\Rm_{\langle 1\rangle}^{-1} \Hm_{\langle 1\rangle}^{~\transp} \big(\Hm_{\langle 1\rangle}^{\vphantom{\transp}} \Rm_{\langle 1\rangle}^{-1} \Hm_{\langle 1\rangle}^{~\transp}\big)^{-1} & \nullmat & \vphantom{\ddots}\cdots & \nullmat \\
		\nullmat & \Rm_{\langle 2\rangle}^{-1} \Hm_{\langle 2\rangle}^{~\transp} \big(\Hm_{\langle 2\rangle}^{\vphantom{\transp}} \Rm_{\langle 2\rangle}^{-1} \Hm_{\langle 2\rangle}^{~\transp}\big)^{-1} & \ddots & \nullmat \\
		\vdots & \nullmat & \ddots & \nullmat \\
		\nullmat & \nullmat & \vphantom{\ddots}\cdots & \Rm_{\langle N\rangle}^{-1} \Hm_{\langle N\rangle}^{~\transp} \big(\Hm_{\langle N\rangle}^{\vphantom{\transp}} \Rm_{\langle N\rangle}^{-1} \Hm_{\langle N\rangle}^{~\transp}\big)^{-1}
	\end{bmatrix}
	= \begin{bmatrix}
		\big[\Hm_{\langle 1\rangle}\big]^{\#}_{\Rm_{\langle 1\rangle}} & \nullmat & \vphantom{\ddots}\cdots & \nullmat \\
		\nullmat & \big[\Hm_{\langle 2\rangle}\big]^{\#}_{\Rm_{\langle 2\rangle}} & \ddots & \nullmat \\
		\vdots & \nullmat & \ddots & \nullmat \\
		\nullmat & \nullmat & \vphantom{\ddots}\cdots & \big[\Hm_{\langle N\rangle}\big]^{\#}_{\Rm_{\langle N\rangle}}
	\end{bmatrix},
	$}
\end{align*}
where $\Hm_{\langle n\rangle} := \Hm\big(n\dt,\phiv(n\dt)\big)$ and $\Rm_{\langle n\rangle} := \Rm\big(n\dt,\phiv(n\dt)\big)$.
It shows evidently that the pseudoinverse of $\varHm(\varphiv)$ for the entire time interval is equivalent to the pseudoinverse of $\Hm\big(t,\phiv(t)\big)$ at each time instant and there is no interference among time instants.
In fact, this observation is valid for the entire equation~\eqref{eq:optimal control, barrier-method sequence} (deduction is omitted).
As a result, one concludes that for all $n \in \naturalset^{+} \cap [1,N]$
%\begin{align*}
%	\phiv^{[k+1]}(n\dt) 
%	&= \phiv^{[k]}(n\dt) - \Hm\big(n\dt,\phiv^{[k]}(n\dt)\big)^{\#}_{\Rm\big(n\dt,\phiv^{[k]}(n\dt)\big)} \hv\big(n\dt,\phiv^{[k]}(n\dt)\big) \\
%	&\quad- \Big(
%	\Imat - \Hm\big(n\dt,\phiv^{[k]}(n\dt)\big)^{\#}_{\Rm\big(t,\phiv^{[k]}(n\dt)\big)} \Hm\big(n\dt,\phiv^{[k]}(n\dt)\big)
%	\Big) \Rm\big(n\dt,\phiv^{[k]}(n\dt)\big)^{-1} \rv\big(n\dt,\phiv^{[k]}(n\dt)\big)
%\end{align*}
%\begin{align*}
%	\phiv^{[k+1]}(t_{n}) 
%	&= \phiv^{[k]}(t_{n}) - \Hm\big(t_{n},\phiv^{[k]}(t_{n})\big)^{\#}_{\Rm\big(t_{n},\phiv^{[k]}(t_{n})\big)} \hv\big(t_{n},\phiv^{[k]}(t_{n})\big) \\
%	&\quad- \Big(
%	\Imat - \Hm\big(t_{n},\phiv^{[k]}(t_{n})\big)^{\#}_{\Rm\big(t,\phiv^{[k]}(t_{n})\big)} \Hm\big(t_{n},\phiv^{[k]}(t_{n})\big)
%	\Big) \Rm\big(t_{n},\phiv^{[k]}(t_{n})\big)^{-1} \rv\big(t_{n},\phiv^{[k]}(t_{n})\big),
%\end{align*}
%with $t_{n} = n\dt$,
\begin{align*}
	\phiv^{[k+1]}(n\dt) 
	&= \phiv^{[k]}(n\dt) + \big[\Hm_{\langle n\rangle}\big]^{\#}_{\Rm_{\langle n\rangle}} \Big(\desired{\yv}(n\dt) - \hv\big(n\dt,\phiv^{[k]}(n\dt)\big)\Big) \\
	&\qquad\qquad\quad\,- \Big(
		\Imat - \big[\Hm_{\langle n\rangle}\big]^{\#}_{\Rm_{\langle n\rangle}} \Hm_{\langle n\rangle}
	\Big) \big[\Rm_{\langle n\rangle}\big]^{-1} \rv\big(n\dt,\phiv^{[k]}(n\dt)\big) \\
	&= \phiv^{[k]}(n\dt) + \etav\big(n\dt, \phiv^{[k]}(n\dt), \desired{\yv}(n\dt)\big)
\end{align*}
converges to $\phiv^{\optimal}(n\dt)$ as $k \in \naturalset$ increases.
This behavior is the same as that of $\chiv^{[k]}(t)$, see~\eqref{eq:barrier-method sequence, v2}.
Therefore, one concludes that $\forall n \in \naturalset^{+} \cap [1, N]$, if $\phiv^{[0]}(n\dt) = \chiv^{[0]}(n\dt)$, then $\phiv^{\optimal}(n\dt) = \chiv^{\optimal}(n\dt)$.
\end{proof}

Lemma~\ref{lemma:equivalency optimal control and convex programming} shows the equivalency between $\phiv^{\optimal}$ and $\chiv^{\optimal}$ in discrete time, and, consequently, the numerical equivalency between \eqref{eq:constrained optimization} and \eqref{eq:optimal control problem}.
Since $\dt$ can be chosen arbitrarily close to zero, we conclude
$\phiv^{\optimal}(t) = \chiv^{\optimal}(t)~\forall t \in (0, t_{\final}]$.
Therefore, Lemma~\ref{lemma:x to x_ast, v3}, Remark~\ref{remark: sufficiently large Kx ensures the closeness for all t > 0}, and Theorem~\ref{theorem:optimal effort distributor} can be analogously applied here.
%the optimal control problem \eqref{eq:optimal control problem}.
Formally, we formulate the following theorem and remark.

\begin{theorem}\label{theorem:x to phi_ast}
	Consider the time-variant system~\eqref{eq:gryphon} with an initial condition $\xv(0)$ for state and a reference trajectory $\desired{\yv}$ for output.
	Consider the controller~\eqref{eq:gryphonheart} with $K_{\textup{x}} = \textup{const.} > 0$.
	Consider an arbitrary local minimizer 
	\begin{equation*}
		\varphiv^{\optimal} = \begin{bmatrix}
			\phiv^{\optimal}(1\dt)^\transp &
			\phiv^{\optimal}(2\dt)^\transp &
			\cdots &
			\phiv^{\optimal}(N\dt)^\transp
		\end{bmatrix}^\transp
	\end{equation*}
	of \eqref{eq:optimal control convex optimization}.
	Suppose that
	\begin{align}
		\phiv^{\optimal}(0) &= \xv(0), \label{eq: x(0) is phi_optimal(0)} \\
		\big\|\dot{\phiv}^{\optimal}(t)\big\|_{2} &\leqslant D_{\phi}~~
		\forall t \geqslant 0. \label{eq:slowly varying phi_optimal}
	\end{align}
	There exists $\delta: \realset^{\dimx} \to \realset_{>0}$
	and $M: \realset^{\dimx} \to (0, M_{\sup}]$ with $M_{\sup} < 1$
	such that 
	if 
	\begin{equation}\label{eq: x(t) is sufficiently close to phi_optimal(t)}
%		\big\|\xv(0) - \phiv^{\optimal}(0)\big\|_{2} \exp (-\kappa t)
%		+ \kappa^{-1} D_{\phi}
%		<
%		\delta\big(\phiv^{\optimal}(t)\big)~~
%		\forall t > 0 
		%
		\kappa^{-1} D_{\phi}
		<
		\delta\big(\phiv^{\optimal}(t)\big)~~
		\forall t > 0 
	\end{equation}
	then 
	\begin{empheq}[box=\colorbox{\colorotc}]{equation}
%	\begin{equation}
	\label{eq:exponential convergence x to phi_optimal}
%		\big\|\xv(t) \hspace{-0.25ex}-\hspace{-0.25ex} \phiv^{\optimal}(t)\big\|_{2}
%		\leqslant
%		\big\|\xv(0) \hspace{-0.25ex}-\hspace{-0.25ex} \phiv^{\optimal}(0)\big\|_{2} \exp (-\kappa t)
%		+ \int_{0}^{t} \exp\big(\hspace{-0.25ex}-\hspace{-0.25ex}\kappa (t\hspace{-0.25ex}-\hspace{-0.25ex}\tau)\big) \big\|\dot{\phiv}^{\optimal}(\tau)\big\|_{2}\,\mathrm{d}\tau,
		%
		\big\|\xv(t) \hspace{-0.25ex}-\hspace{-0.25ex} \phiv^{\optimal}(t)\big\|_{2}
		\leqslant
		\int_{0}^{t} \exp\big(\hspace{-0.25ex}-\hspace{-0.25ex}\kappa (t\hspace{-0.25ex}-\hspace{-0.25ex}\tau)\big) \big\|\dot{\phiv}^{\optimal}(\tau)\big\|_{2}\,\mathrm{d}\tau,
%	\end{equation}
	\end{empheq}
	where 
	\begin{equation*}%\label{eq:definition of kappa}
		\kappa = \frac{1}{2} \big(1-M_{\sup}^{2}\big) K_{\textup{x}}.
	\end{equation*}
\end{theorem}
\begin{proof}
It can be concluded immediately from Lemma~\ref{lemma:x to x_ast, v3}, Theorem~\ref{theorem:optimal effort distributor}, and Lemma~\ref{lemma:equivalency optimal control and convex programming}.
It must be noted that \eqref{eq: x(t) is sufficiently close to chi_optimal(t)} and \eqref{eq:exponential convergence x v2.5}  degenerate to \eqref{eq: x(t) is sufficiently close to phi_optimal(t)} and \eqref{eq:exponential convergence x to phi_optimal}, respectively, due to~\eqref{eq: x(0) is phi_optimal(0)}.
\end{proof}

\begin{remark}\label{remark: sufficiently large Kx ensures the closeness for all t > 0, optimal-ctrl}
	A larger $K_{\textup{x}}$ practically increases the likelihood of satisfying \eqref{eq: x(t) is sufficiently close to phi_optimal(t)}
	and a sufficiently large $K_{\textup{x}} > 0$ practically makes \eqref{eq: x(t) is sufficiently close to phi_optimal(t)} satisfied.
\end{remark}

Remember that \eqref{eq: x(0) is phi_optimal(0)} is an equality constraint in \eqref{eq:optimal control problem}, and it can be easily guaranteed.
{Therefore, Theorem~\ref{theorem:x to phi_ast} addresses Problem~\ref{problem:find an optimal tracking controller} completely.}
%Since Theorem~\ref{theorem:x to phi_ast} addresses Problem~\ref{problem:find an optimal tracking controller} completely, we conclude that $\muv$ in \eqref{eq:gryphonheart} is an \OTC~(Definition~\ref{definition:gryphonheart}).

%\addedtwo{The interpretation of the resulting \OTC~\eqref{eq:gryphonheart} is given in the Main Document.}

% !TEX root = ../supplementary.tex
\subsection{Orthogonality}

\begin{theorem}[DCC's orthogonality]
	There exists a weighted inner product space in which the two components $\etav_{\textup{A}}(t,\xv,\desired{\yv})$ and $\etav_{\textup{B}}(t,\xv)$, cf. \eqref{eq:definition of eta_A} and \eqref{eq:definition of eta_B}, of \eqref{eq:elegant gryphon candidate} are orthogonal to each other for all $t \geqslant 0$ and all $\xv \in \realset^{\dimx}$.
\end{theorem}
\begin{proof}
Let the weighted inner product
\begin{equation}\label{eq:inner product, step 1}
	w 
	= \etav_{\textup{A}}(t,\xv,\desired{\yv})^\transp 
	\Wm
	\etav_{\textup{B}}(t,\xv)
\end{equation}
be considered with some symmetric and positive definite matrix $\Wm$.
Inserting \eqref{eq:definition of eta_A} and \eqref{eq:definition of eta_B} into \eqref{eq:inner product, step 1} results in
\begin{align}
	w &= \big(
	\desired{\yv} - \hv(t,\xv)
	\big)^\transp \left(\Hm(t,\xv)^{\#}_{\Rm(t,\xv)}\right)^\transp 
	\Wm
	\left(
	\Imat - \Hm(t,\xv)^{\#}_{\Rm(t,\xv)} \Hm(t,\xv)
	\right) \Rm(t,\xv)^{-1} \rv(t,\xv) \nonumber\\
	&= \big(
	\desired{\yv} - \hv
	\big)^\transp \left(
		(\Hm \Rm^{-1} \Hm^\transp)^{-\transp}\Hm\Rm^{-\transp}
		\Wm
		\big(\Imat - \Rm^{-1}\Hm^\transp(\Hm \Rm^{-1} \Hm^\transp)^{-1}\Hm\big)	
	\right)
	\Rm^{-1}\rv, \label{eq:inner product, step 2}
\end{align}
where the dependency is omitted for readability.
Let \begin{equation*}
	\Wm = \Rm(t,\xv)
\end{equation*}
be chosen, then \eqref{eq:inner product, step 2} becomes
\begin{align*}
	w &= \big(
	\desired{\yv} - \hv
	\big)^\transp \Big(
	(\Hm \Rm^{-1} \Hm^\transp)^{-\transp}\Hm
	\overbracket[0.14ex]{\Rm^{-\transp}\Rm}^{=\,\Imat}
	\big(\Imat - \Rm^{-1}\Hm^\transp(\Hm \Rm^{-1} \Hm^\transp)^{-1}\Hm\big)	
	\Big)
	\Rm^{-1}\rv \\
	&= \big(
	\desired{\yv} - \hv
	\big)^\transp \Big(
	(\Hm \Rm^{-1} \Hm^\transp)^{-\transp}\Hm
	\big(\Imat - \Rm^{-1}\Hm^\transp(\Hm \Rm^{-1} \Hm^\transp)^{-1}\Hm\big)	
	\Big)
	\Rm^{-1}\rv \\
	&= \big(
	\desired{\yv} - \hv
	\big)^\transp \Big(
	\big(
		(\Hm \Rm^{-1} \Hm^\transp)^{-\transp}\Hm 
		- (\Hm \Rm^{-1} \Hm^\transp)^{-\transp} 
		\overbracket[0.14ex]{\Hm\Rm^{-1}\Hm^\transp(\Hm \Rm^{-1} \Hm^\transp)^{-1}}^{=\,\Imat}
		\Hm
	\big)	
	\Big)
	\Rm^{-1}\rv \\
	&= \big(
	\desired{\yv} - \hv
	\big)^\transp \Big(
	\overbracket[0.14ex]{
		\big(
		(\Hm \Rm^{-1} \Hm^\transp)^{-\transp}\Hm 
		- (\Hm \Rm^{-1} \Hm^\transp)^{-\transp}	\Hm
		\big)
	}^{=\,\nullmat}
	\Big)
	\Rm^{-1}\rv \\
	&= 0.
\end{align*}
Note that $\Rm(t,\xv)$ is symmetric. % and so is $\Hm(t,\xv) \Rm(t,\xv)^{-1} \Hm(t,\xv)^\transp$.
\end{proof}

%	This zero-valued inner product means that the eigendecomposition of the chosen $\Am$ in \eqref{eq:choice of A} provides an affine transformation which transforms $\vv_{\mathrm{A}}$ and $\vv_{\mathrm{B}}$ to a pair of vectors that are orthogonal to each other.
%	Since affine transformations preserve parallelism, one can conclude that a parallelism between $\vv_{\mathrm{A}}$ and $\vv_{\mathrm{B}}$ does not exist.
% !TEX root = ../supplementary.tex
%%%%%%%%%%%%%%%%%%%%%%%%%%%%%%%%%%%%%%%%%%%%%%%%%%%%%%%%%%%%%%%%%%%%%%%%%%%%%%%%%%%%%%%%%%
\subsection{Remark on barrier function}\label{sec:remark on barrier fcn}
The two barrier parameters $p_{1}$ and $p_{2}$ of the barrier function~\eqref{eq:barrier function} need to be infinitely and sufficiently large, respectively, according to Lemma~\ref{corollary:bnd as penalty, v3}. 
This subsection presents a way to design $p_{1}$ and $p_{2}$ that satisfy this requirement.

The barrier function $\beta$ is not explicitly required by the \OTC~\eqref{eq:gryphonheart}, rather implicitly through $\rv(t,\xv)$ and $\Rm(t,\xv)$, cf.~\eqref{eq:definition of eta}; they are %which are \eqref{eq:augmented cost function}'s gradient and Hessian, i.e.,
\begin{subequations}\label{eq:definition of r and R}
	\begin{align}
		\rv(t,\xv) 
		&:= \left[\frac{\partial \rho(t,\xv)}{\partial \xv}\right]^\transp 
		= \qv(t,\xv) + \Gm(t)^\transp \xiv\big(\gv(t,\xv)\big), \label{eq:definition of r}\\
		\Rm(t,\xv) 
		&:= \frac{\partial \rv(t,\xv)}{\partial \xv} 
		= \Qm(t,\xv) + \Gm(t)^\transp \Xim\big(\gv(t,\xv)\big) \Gm(t), \quad \label{eq:definition of R}	
	\end{align}
\end{subequations}
respectively, 
where 
\begin{align*}
	\qv(t,\xv) &:= \left[\frac{\partial \sigma(t,\xv)}{\partial \xv}\right]^\transp, \\
	\Qm(t,\xv) &:= \frac{\partial}{\partial \xv} \frac{\partial \sigma(t,\xv)}{\partial \xv}
\end{align*}
%are the gradient and Hessian of the cost function $\sigma(t,\xv)$,
and
%\begin{align}
%	\xiv(\sv) 
%	&:= \frac{\partial \betav(\sv)}{\partial \sv}^\transp \betav(\sv)
%	= \begin{bmatrix}
	%		\displaystyle\frac{\partial \beta(s_{i})}{\partial s_{i}}^\transp \beta(s_{i})
	%	\end{bmatrix}, \label{eq:definition of xi} \\
%	%
%	\Xim(\sv) 
%	&:= \frac{\partial}{\partial \sv} \left(
%	\frac{\partial \betav(\sv)}{\partial \sv}^\transp \betav(\sv)
%	\right)
%	= \diag\left\{
%		\frac{\partial}{\partial s_{i}}\left(\frac{\partial \beta(s_{i})}{\partial s_{i}}^\transp \beta(s_{i})\right)
%	\right\}. \label{eq:definition of Q and Xi}
%\end{align}
\begin{subequations}\label{eq:definition of xi and Xi}
	\begin{alignat}{2}
		\xiv(\sv) &:= \begin{bmatrix}
			\xi(s_{1}) & \xi(s_{2}) & \cdots & \xi(s_{\dimc})
		\end{bmatrix}^\transp,~\quad\quad
		\xi(s) 
		&&:= \frac{\partial \beta(s)}{\partial s}^\transp \beta(s), \label{eq:definition of xi} \\
		\Xim(\sv) 
		&:= \begin{bmatrix}
			\varXi(s_{1}) & 0 & \cdots & 0 \\
			0 & \varXi(s_{2}) & \cdots & 0 \\
			\vdots & \ddots & \ddots & \vdots \\
			0 & 0 & \cdots & \varXi(s_{\dimc})
		\end{bmatrix},~\quad
		\varXi(s) 
		&&:= \frac{\partial}{\partial s}\left(\frac{\partial \beta(s)}{\partial s}^\transp \beta(s)\right) \label{eq:definition of Xi}
	\end{alignat}
\end{subequations}
are related to the derivatives of $\beta$.
By analytically deriving $\xi$ and $\varXi$, one is able to let $p_{1} \rightarrow +\infty$; this results in
\begin{subequations}\label{eq:barrier metric functions xi and Xi}
	\begin{align}
		%	\lim\limits_{p_{1} \rightarrow +\infty} \xi(s; [\,p_{1}~p_{2}\,]^\transp) 
		\lim\limits_{p_{1} \rightarrow +\infty} \xi(s) 
		&= \begin{cases}
			p_{2}^2\,s & \text{if~}s > 0, \\
			0 & \text{if~}s \leqslant 0,
		\end{cases} \\
		%
		%	\lim\limits_{p_{1} \rightarrow +\infty} \varXi(s; [\,p_{1}~p_{2}\,]^\transp)
		\lim\limits_{p_{1} \rightarrow +\infty} \varXi(s)
		&= \begin{cases}
			p_{2}^2 & \text{if~}s > 0, \\
			\frac{1}{4} \big(\ln(2)+1\big) p_{2}^2 & \text{if~}s = 0, \\
			0 & \text{if~}s < 0,
		\end{cases} \label{eq:barrier metric functions Xi}
	\end{align}
\end{subequations}
see Figure~\ref{fig:barrierfcn} for graphs.
\begin{figure}[!h]
	\centering
	%	\includegraphics[width=0.99\linewidth]{figures/barrierfcn}
	%	\noindent\makebox[\textwidth]{\includegraphics{figures/barrierfcn.pdf}}%
	\noindent\makebox[\textwidth]{\input{figures/barrierfcn}}
	\vspace{-2ex}
	\caption{\label{fig:barrierfcn}Barrier function $\beta$ and its derivative-related functions $\xi$ and $\varXi$ with various values for parameter $p_{1}$.}
\end{figure}
{%
	Inserting \eqref{eq:barrier metric functions xi and Xi} into \eqref{eq:definition of r and R} results in%
}
\begin{subequations}\label{eq:implementation of r and R} 
	\begin{align}
		\rv(t,\xv) &= \begin{cases}
			\qv(t,\xv) & \text{if~}\gv(t,\xv) \leqslant \nullv, \\
			\qv(t,\xv) + p_{2}^{2} \breve{\Gm}(t)^\transp \breve{\gv}(t,\xv) & \text{otherwise},
		\end{cases} \label{eq:implementation of r} \\
		\Rm(t,\xv) &= \begin{cases}
			\Qm(t,\xv) & \text{if~}\gv(t,\xv) \leqslant \nullv, \\
			\Qm(t,\xv) + p_{2}^{2} \breve{\Gm}(t)^\transp \breve{\Gm}(t) & \text{otherwise},
		\end{cases} \label{eq:implementation of R}
	\end{align}
\end{subequations}
where 
$\breve{\gv}(t,\xv) = \breve{\Gm}(t) \xv + \breve{\cv}(t) \in \realset^{\dimcbar}$
is the vector of positive elements in $\gv(t,\xv)$.
{%
	Equation~\eqref{eq:implementation of r and R} shows that the value of $p_{2}$ matters, only when $\xv$ is in the exterior of $\{\varxv~\lvert~\gv(t,\varxv) \leqslant \nullv\}$. 
	A sufficiently large $p_{2}$ makes the Newton-Raphson step $-\Rm(t,\xv)^{-1} \rv(t,\xv)$ be dominated by the component directing to the interior, i.e.,
}
\begin{align}
	- \big(\breve{\Gm}(t)^\transp \breve{\Gm}(t)\big)^{-1} \breve{\Gm}(t)^\transp \breve{\gv}(t,\xv)
	&\approx
	- \Rm(t,\xv)^{-1} \rv(t,\xv) \label{eq:sufficiently large p2 squared, design equation} \\
	&\stackrel{\mathclap{\scalebox{0.6}{\eqref{eq:implementation of r and R}}}}{=\joinrel=}
	-\Big(
	\Qm(t,\xv) + p_{2}^{2} \breve{\Gm}(t)^\transp \breve{\Gm}(t)
	\Big)^{\hspace{-0.3ex}-1} \Big(
	\qv(t,\xv) + p_{2}^{2} \breve{\Gm}(t)^\transp \breve{\gv}(t,\xv)
	\Big). \nonumber
\end{align}
The left-hand side of \eqref{eq:sufficiently large p2 squared, design equation} is a Newton-Raphson step in finding the root of $\nullv = \breve{\gv}(t,\xv)$. %(see SI, Section~\ref{SI-sec:newton-raphson method}).
This dominance can be ensured by
\begin{equation*}%\label{eq:R much greater than Q}
	\big\|\Rm(t,\xv)\big\| 
	= \big\| \Qm(t,\xv) + p_{2}^{2} \breve{\Gm}(t)^\transp \breve{\Gm}(t) \big\|
	\gg 
	\big\|\Qm(t,\xv)\big\|
	~\Leftrightarrow~
	K_{\RoverQ} := \frac{\|\Rm(t,\xv)\|}{\|\Qm(t,\xv)\|}  \gg 1.
\end{equation*}
Here, $\|\cdot\|$ denotes an arbitrary matrix norm.
Using the triangle inequality 
\begin{align*}
	\big\|\Rm(t,\xv)\big\| 
	%	&\stackrel{\mathclap{\text{\eqref{eq:definition of R}}}}{=\joinrel=} \big\|\Qm(t,\xv) + \Gm(t)^\transp \Xim\big(\gv(t,\xv)\big) \Gm(t)\big\| \\
	%	&\stackrel{\mathclap{\text{\eqref{eq:barrier metric functions Xi}}}}{=\joinrel=} \big\|\Qm(t,\xv) + \Gm(t)^\transp p_{2}^2\,\diag\big\{\breve{\varXi}\big(\gv(t,\xv)\big)\big\} \Gm(t)\big\| \\
	%	&\leqslant \big\|\Qm(t,\xv)\big\| +	p_{2}^2\,\big\|\Gm(t)^\transp \diag\big\{\breve{\varXi}\big(\gv(t,\xv)\big)\big\} \Gm(t)\big\|
	%
	= \big\|\Qm(t,\xv) + p_{2}^2 \breve{\Gm}(t)^\transp \breve{\Gm}(t)\big\|
	\leqslant \big\|\Qm(t,\xv)\big\| +	p_{2}^2\,\big\|\breve{\Gm}(t)^\transp \breve{\Gm}(t)\big\|
\end{align*}
leads to
\begin{align*}
	K_{\RoverQ}
	%= \frac{\|\Rm(t,\xv)\|}{\|\Qm(t,\xv)\|}
	%	\leqslant 1 + \frac{p_{2}^2\,\big\|\Gm(t)^\transp \diag\big\{\breve{\varXi}\big(\gv(t,\xv)\big)\big\} \Gm(t)\big\|}{\|\Qm(t,\xv)\|},
	\leqslant 1 + \frac{p_{2}^2\,\big\|\breve{\Gm}(t)^\transp \breve{\Gm}(t)\big\|}{\|\Qm(t,\xv)\|},
\end{align*}
which is equivalently
\begin{equation}\label{eq:sufficiently large p2 squared}
	p_{2}^2
	%	\geqslant \frac{(K_{\RoverQ} - 1) \|\Qm(t,\xv)\|}{\big\|\Gm(t)^\transp \diag\big\{\breve{\varXi}\big(\gv(t,\xv)\big)\big\} \Gm(t)\big\|}.
	\geqslant \frac{(K_{\RoverQ} - 1) \|\Qm(t,\xv)\|}{\big\|\breve{\Gm}(t)^\transp \breve{\Gm}(t)\big\|}.
\end{equation}
Any $p_{2}$ that satisfies \eqref{eq:sufficiently large p2 squared} with a sufficiently large design parameter $K_{\RoverQ} \gg 1$ (typically, $K_{\RoverQ} > 1000$) can be considered sufficiently large in the concerned numerical context.
%\deleted{Moreover, for a minimum computational overhead, this work chooses the maximum norm (maximum absolute row sum) for $\|\cdot\|$ and the equal sign in \eqref{eq:sufficiently large p2 squared} for finding $p_{2}$.}
%

{%
	As indicated by \eqref{eq:implementation of r and R}, 
}%
the notion of barrier function $\beta$ was introduced as a facilitative tool for deriving the proposed \OTC~theory and is not required in the final \OTC~formulation.

%\deleted{%
	%For a large $p_{2}$, it is likely that the SPD matrix $\Rm(t,\xv)$~\eqref{eq:implementation of R} becomes ill-conditioned. %, when at least one element of $\gv(t,\xv)$ is positive.
	%This may cause non-negligible numerical errors in computing its inverse. %(in $\etav$, cf.~\eqref{eq:gryphonheart}, \eqref{eq:definition of eta}, and \eqref{eq:definition of eta_B}).
	%To overcome this, the following numerical procedure is applied:}
%\begin{itemize}
%	\item \deleted{If $\Rm(t,\xv)$ is diagonal, a diagonal matrix is constructed with diagonal elements that are reciprocals of the corresponding elements of $\Rm(t,\xv)$.}
%	\item \deleted{Otherwise, the Cholesky decomposition is applied on $\Rm(t,\xv)$ followed by forward and back substitutions.}
%\end{itemize}
%%This procedure suffices for the numerical examples presented in this work.
%%In general, there are very rare cases where the Cholesky decomposition fails.
%%In these cases, one may use the Sherman-Morrison formula to compute the inverse of $\Rm(t,\xv)$, taking advantage of a well posed $\Qm(t,\xv)$ and a diagonal $\Xim\big(\gv(t,\xv)\big)$.

{%
	\subsection{Remark on state equation}\label{sec:remark on state equation}
	It is worth noting that $\fv_{\scriptscriptstyle\textup{A}}$ and $\Bm$ in \eqref{eq:gryphon} can be generalized as
	\begin{equation}\label{eq:non-affine gryphon}
		\dot{\xv}
		= \fv(\xv, \uv),
	\end{equation}
	where $\dimu:=\dim(\uv)$ does not have to equal $\dimx$.
	The proposed \OTC\ theory is applicable
	as long as 
	%\begin{equation*}
	%	\uv = \psiv(\zv,\xv^{(r)}) ~\big|~ \xv^{(r)} = \fv\big(\zv, \psiv(\zv,\xv^{(r)})\big)
	%\end{equation*}
	\begin{equation}\label{eq:inverse of generalized gryphon}
		\left\{
		\psiv: \realset^{\dimx}\times\realset^{\dimx} \to \realset^{\dimu} ~\big|~ \ellv = \fv\big(\xv, \psiv(\xv,\ellv)\big)~\forall 
		\begin{bmatrix}
			\xv^\transp & \ellv^\transp
		\end{bmatrix}^\transp 
		\in \realset^{\dimx+\dimu}
		\right\}
		\neq \emptyset.
	\end{equation}
	Consequently, the expression $\Bm(\xv)^{-1}\big(-\fv_{\scriptscriptstyle\textup{A}}(\xv) + \cdots\big)$ in \eqref{eq:gryphon, fbl-2} and consequently in \eqref{eq:gryphonheart} is then replace by $\psiv\big(\xv, \cdots\big)$, i.e., \eqref{eq:gryphonheart} becomes
	\begin{equation*}
		\uv = \psiv\big(\xv,~ K_{\textup{x}} \etav(t, \xv, \desired{\yv})\big).
	\end{equation*}
	Note that only the existence~\eqref{eq:inverse of generalized gryphon} of $\psiv$ is required, not the uniqueness;
	and this existence is related to the controllability of \eqref{eq:non-affine gryphon}.
	In-depth investigation is left for future work.
}%

\newpage
\section{Numerical examples -- extended description}\label{sec:extended_numerical_examples}
% !TEX root = ../supplementary.tex
\subsection{Example M1}\label{sec:numerical example M1}
%
%To demonstrate \OTC's behavior and performance, we select\added{, besides the three realistic Examples R1--R3 in the main document,}
%one comprehensive, purely \uls{m}athematical example (Example M1).
%
In terms of the system and problem class~\eqref{eq:gryphon}--\eqref{eq:optimal control problem}, we construct this comprehensive example as follows:
\begin{subequations}\label{eq:example m1}
	\begin{align}
		&\xv = \begin{bmatrix}
			x_{1} & x_{2}
		\end{bmatrix}^\transp, 
		\quad
		\fv_{\scriptscriptstyle\textup{A}}(\xv) = \frac{1}{\xv^\transp \xv + 0.001},
		\quad
		\Bm(\xv) = \big(\xv^\transp \xv + 1\big) \Imat,
		\quad \\
		&h(t,\xv) = -6\,t+{\big(5\,x_{1}+25\,x_{2}^2-7\big)}^2+{\big(25\,x_{1}^2+5\,x_{2}-11\big)}^2 - 1,
		\quad
		\desired{y} = 0,
		\quad \\
		&\gv(t,\xv) = \begin{bmatrix} 1 & 0\\[1ex] -1 & 0\\[1ex] 0 & 1\\[1ex] 0 & -1\\[1ex] \frac{209}{100\,\sin(\pi t/10)-143} & 1\\[1ex] \frac{\sin(\pi t/10)}{2}+\frac{73}{75} & 1\\[1ex] \frac{77}{8\,\left(5\,\sin(\pi t/10)+7\right)} & -1\\[1ex] \frac{15\,\sin(\pi t/10)}{38}-\frac{67}{152} & -1 \end{bmatrix} \xv 
		+ \begin{bmatrix} -\frac{19}{20}\\[1ex] -\frac{19}{20}\\[1ex] -\frac{19}{20}\\[1ex] -\frac{19}{20}\\[1ex] -\frac{209\,\left(3\,\sin(\pi t/10)-10\right)}{10\,\left(100\,\sin(\pi t/10)-143\right)}\\[1ex] -\frac{3\,\sin(\pi t/10)}{10}-1\\[1ex] -\frac{77\,\left(\sin(\pi t/10)+5\right)}{40\,\left(5\,\sin(\pi t/10)+7\right)}\\[1ex] \frac{3\,\sin(\pi t/10)}{10}-1 \end{bmatrix}, 
		\quad \\
		&\sigma(t,\xv) = \sqrt{
			\big(\xv - \pv(t)\big)^\transp \Pm(t) \big(\xv - \pv(t)\big) + 0.001
		} + 0.001\, \big(\xv - \pv(t)\big)^\transp \big(\xv - \pv(t)\big),
		\quad \\
		&\quad\Pm(t) = \begin{bmatrix} \frac{9 t + 9 t \cos\big(4 \pi^{2}\,\sin(\pi t/100)\big)+40}{18 t + 40} & \frac{9 t \sin\big(4 \pi^{2}\,\sin(\pi t/100)\big)}{18 t + 40}\\[1ex] \frac{9 t \sin\big(4 \pi^{2}\,\sin(\pi t/100)\big)}{18 t + 40} & \frac{9 t {\sin\big(2\,\pi ^2\,\sin(\pi t/100)\big)}^2+20}{9 t + 20} \end{bmatrix},
		\quad \\%[1ex]
		&\quad\pv(t) = \begin{bmatrix} -\frac{2\,\sin(\pi t/5)}{3}-\frac{10\,\sin(3 \pi t/10)}{27}\\[1ex] \frac{2\,\cos(\pi t/5)}{3}-\frac{10\,\cos(3 \pi t/10)}{27} \end{bmatrix}.
	\end{align}
\end{subequations}
This example is constructed to have high degrees of nonlinearity in the state, output, and cost functions, along with time-variance in these functions and the inequality constraints. These features make it a suitable lowest-dimensional representative of the system and problem class~\eqref{eq:gryphon}--\eqref{eq:optimal control problem}.
The time evolution of the cost function $\sigma(t,\xv)$, equality constraint $\desired{y} = h(t,\xv)$, and inequality constraint $\nullv \geqslant \gv(t,\xv)$ are illustrated in Figure~\ref{fig:example 0 result}.
\afterpage{%
	\newgeometry{left=1.5cm,right=1.5cm,top=2cm,bottom=6cm}
	\begin{figure}[!t]
		\centering
		\includegraphics[width=15cm]{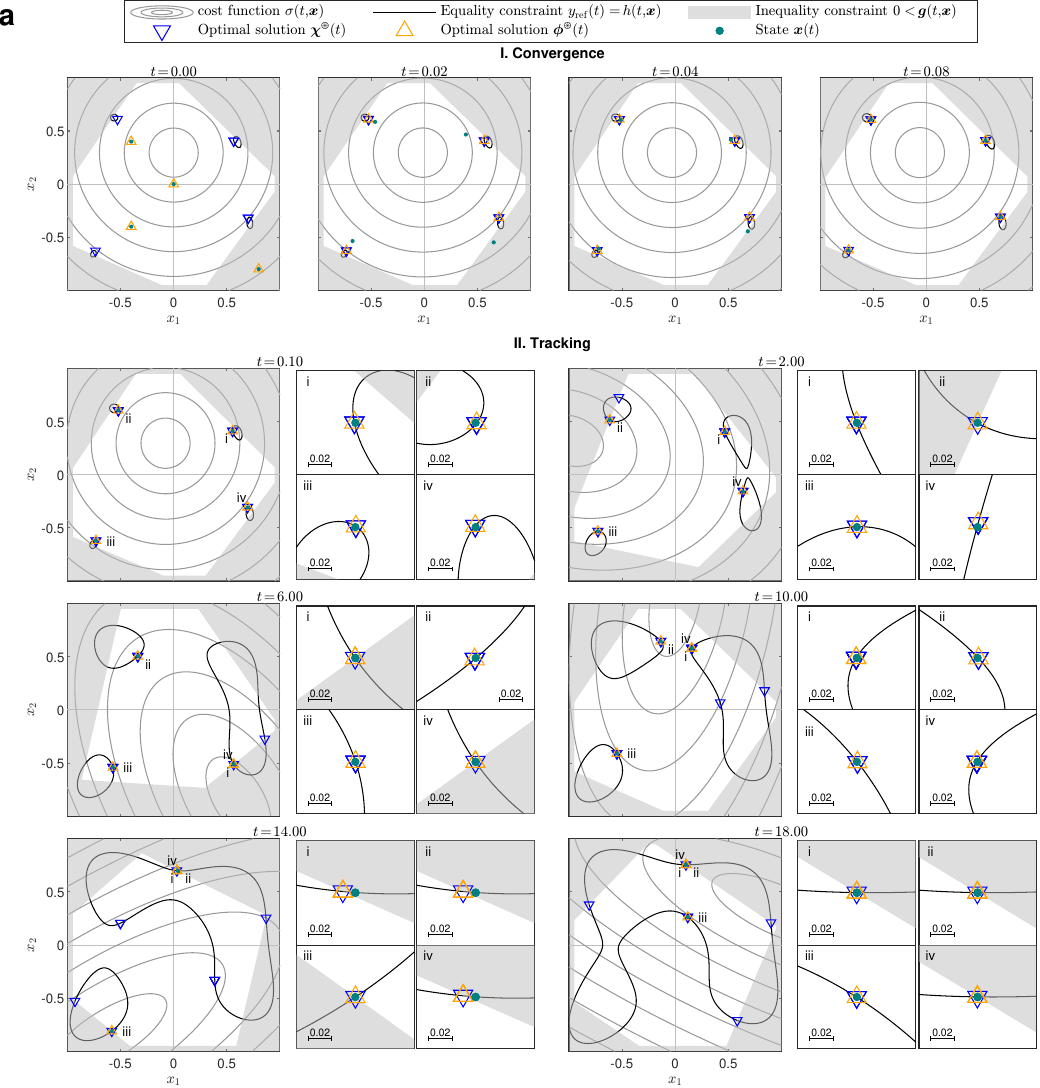}
		\includegraphics[width=15cm]{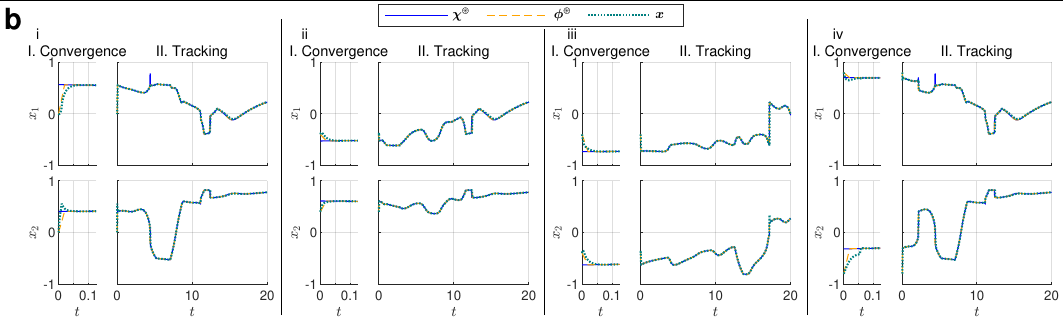}
		%\thisfloatpagestyle{empty}
		%	\vspace{-2ex}
		\caption{\label{fig:example 0 result}
			Example M1: A two-dimensional system~\eqref{eq:example m1} controlled by the \OTC\ with $K_{\textup{x}} = 100$.
			The trajectories $\xv(t)$ (marked as i, ii, iii, iv) starting from different initial conditions $\xv(0)$ follow the corresponding $\chiv^{\optimal}(t)$ and $\phiv^{\optimal}(t)$ over time $t$. 
			\textbf{a,} Trajectories of $\xv(t)$ (teal), $\chiv^{\optimal}(t)$ (blue), and $\phiv^{\optimal}(t)$ (orange) in the two-dimensional space along with cost function and constraints. Each subplot corresponds to a snapshot at time instant $t$ as indicated.
			An animation is available in the Supplementary Video (\url{https://youtu.be/tzvtecj1I5I}).
			\textbf{b,} Exponential convergence and tracking of $\xv(t)$ to the closest $\chiv^{\optimal}(t)$ and $\phiv^{\optimal}(t)$. Each column represents the state's time evolution starting from a different initial condition.
			A zoomed-in perspective to the initial transient phase is provided on the left in each column.
		}
	\end{figure}
	\restoregeometry
}

Figure~\ref{fig:example 0 result} shows that the state $\xv(t)$ of %\deleted{a two-dimensional, yet, comprehensive system, cf.~}
\eqref{eq:example m1} controlled by \eqref{eq:gryphonheart}, follows $\chiv^{\optimal}(t)$ and $\phiv^{\optimal}(t)$ closely for all $t > 0$.
These trajectories of $\xv(t)$, starting from distinct initial points $\xv(0)$, demonstrate an exponential convergence to the corresponding $\chiv^{\optimal}(t)$ and $\phiv^{\optimal}(t)$ with an initial transient phase of $< 0.1$.
Despite the occasional abrupt jumps in trajectories $\chiv^{\optimal}(t)$ and $\phiv^{\optimal}(t)$, \OTC\ consistently guarantees a rapid first-order convergence behavior of $\xv$ towards a local minimizer.
This fast and robust convergence behavior underscores the effectiveness of \OTC\ in swiftly achieving desired optimal outcomes in time-variant contexts.
% !TEX root = ../supplementary.tex
\subsection{Example R1 - aerial vehicle fleet model}

\paragraph*{Problem statement}
In Example R1, the state vector of $n$ VoliroX omnidirectional aerial vehicles~\cite{BodieTayKamSie2020} jointly carrying a hollow sphere (object) along a predefined path is defined as 
\begin{equation}
\xv = [\,\xv_{\scriptscriptstyle\langle 1\rangle}^{\,\transp}~~\xv_{\scriptscriptstyle\langle 2\rangle}^{\,\transp}~~\cdots~~\xv_{\scriptscriptstyle\langle n\rangle}^{\,\transp}\,]^\transp,
%x = [x_1^\top \; x_2^\top \; \dots \; x_n^\top]^\top
\end{equation}
where $\xv_{\scriptscriptstyle\langle i\rangle} \in \realset^{24}$ represents the actuator states of one vehicle~\cite{BodieTayKamSie2020}, implying $\dimx = 24\,n$. The system dynamics are modeled as
\begin{equation}
	\dot{\xv} = -100\,\xv + 100\,\uv, 
\end{equation}
approximating the low-level actuator dynamics of the vehicles. The cost function 
\begin{equation}
	\sigma(t,\xv) = \frac{1}{2}\xv^\transp \Imat \xv,
\end{equation}
assumes actuators to be considered equally~\cite{KirchengastSteHor2018,BodieTayKamSie2020}.

The payload trajectory follows a known minimum-jerk reference path from start to end pose. Forward dynamics simulation is used to compute the resulting object motion. In each considered scenario, the sphere is of different size and weight, and the number of vehicles is increased $n \in \{16, 36, 64, 100, 144, 196, 256, 400\}$.

\paragraph*{Details}
Through the known inverse-dynamics models of the object and vehicles, the wrench $\hv(t,\xv)$ that is required to move the \uls{o}bject along the reference trajectory $\pv_{\textup{o,ref}}(t)$ is found to be
	\begin{equation}
		\desired{\yv}(t) := \hv_{\textup{o}}\big(\pv_{\textup{o,ref}}(t), \dot{\pv}_{\textup{o,ref}}(t), \ddot{\pv}_{\textup{o,ref}}(t)\big),
	\end{equation}
and the wrench that a moving \uls{v}ehicle (along the trajectory $\pv_{\textup{v}\langle i\rangle}(t)$ for the $i$-th vehicle) is still able to exert to its environment~\cite{BodieTayKamSie2020} as
	\begin{equation}\label{eq:single-drone dynamics}
		\yv_{\langle i\rangle} := \hv_{\textup{v}}\big(\pv_{\textup{v}\langle i\rangle}(t), \dot{\pv}_{\textup{v}\langle i\rangle}(t), \ddot{\pv}_{\textup{v}\langle i\rangle}(t), \xv_{\langle i\rangle}\big).
	\end{equation}
All wrenches and poses are expressed in the same inertial frame of reference.	We further assume that there is no relative movement between any vehicle and the object, i.e., there exists a fixed relation
	\begin{equation}\label{eq:drone-object kinematics}
		\pv_{\textup{v}\langle i\rangle} = \psiv_{\textup{v}\langle i\rangle}(\pv_{\textup{o}})
	\end{equation}
between any $\pv_{\textup{v}\langle i\rangle}$ and $\pv_{\textup{o}}$. Provided that the reference $\pv_{\textup{o,ref}}(t)$ of $\pv_{\textup{o}}(t)$ is given explicitly as a function of $t$ (see below), inserting it into \eqref{eq:drone-object kinematics} and then into \eqref{eq:single-drone dynamics} yields
	\begin{equation}
		\yv_{\langle i\rangle} 
		:= \hv_{\textup{v}}\hspace{-0.5ex}\left(
		\psiv_{\textup{v}\langle i\rangle}\big(\pv_{\textup{o,ref}}(t)\big), 
		\frac{\mathrm{d}}{\mathrm{d}t} \psiv_{\textup{v}\langle i\rangle}\big(\pv_{\textup{o,ref}}(t)\big), 
		\frac{\mathrm{d}^{2}}{\mathrm{d}t^{2}} \psiv_{\textup{v}\langle i\rangle}\big(\pv_{\textup{o,ref}}(t)\big),
		\xv_{\langle i\rangle}(t)
		\right).
	\end{equation}
Its sum then yields the time variant output function
	\begin{equation}
		\begin{split}
			\yv 
			%= \sum_{i=1}^{n} \yv_{\langle i\rangle}
			\hspace{-0.25ex}=\hspace{-0.25ex} \hv(t,\xv) 
			\hspace{-0.25ex}:=\hspace{-0.25ex}
			\sum_{i=1}^{n} \hv_{\textup{v}}\hspace{-0.5ex}\left(
			\psiv_{\textup{v}\langle i\rangle}\hspace{-0.25ex}\big(\pv_{\textup{o,ref}}(t)\big), 
			\frac{\mathrm{d}}{\mathrm{d}t} \psiv_{\textup{v}\langle i\rangle}\hspace{-0.25ex}\big(\pv_{\textup{o,ref}}(t)\big), 
			\frac{\mathrm{d}^{2}}{\mathrm{d}t^{2}} \psiv_{\textup{v}\langle i\rangle}\hspace{-0.25ex}\big(\pv_{\textup{o,ref}}(t)\big),
			\xv_{\langle i\rangle}(t)
			\right).
		\end{split}
	\end{equation}

The object's pose reference $\pv_{\textup{o,ref}}(t)$ is designed to be a minimum-jerk trajectory~\cite{GhazaeiArdakaniRobJoh2015} for $t \in [0\,\textup{s},~1.5\,\textup{s}]$ subject to the boundary conditions
\begin{align*}
	\pv_{\textup{o,ref}}(t \leqslant 0.25\,\textup{s}) &= [\,
	0\,\textup{m}~~
	0\,\textup{m}~~
	0\,\textup{m}~~
	0\,\textup{rad}~~
	0\,\textup{rad}~~
	0\,\textup{rad}
	\,]^\transp, \\
	\pv_{\textup{o,ref}}(t \geqslant 1.25\,\textup{s}) &= [\,
	2\,\textup{m}~~
	2\,\textup{m}~~
	2\,\textup{m}~~
	0\,\textup{rad}~~
	0\,\textup{rad}~~
	0\,\textup{rad}
	\,]^\transp.
\end{align*}	
The trajectory is created numerically at $2$~kHz. 

The chosen inequality constraint $\gv(t,\xv)$ requires that the vehicles actuator state $\xv$ must be such that \smash{\textsuperscript{a)}}every vehicle is allowed to push -- but not pull -- the object with a normal force $\geqslant 1\,\textup{N}$. This means that vehicles and object remain in contact and no relative sliding between them occurs. Further, \smash{\textsuperscript{b)}}every actuator must operate within its payload capacity. This leads to a total of $\dimc = n + 2\dimx$ inequality constraints. The resulting expression is a function of $\xv$ and $\pv_{\textup{v}\langle i\rangle\textup{,ref}}$, i.e., of $\pv_{\textup{o,ref}}$, cf.~\eqref{eq:drone-object kinematics}, and ultimately, explicitely of $t$.

More specific implementation details about R1 are that
\begin{itemize}
	\item The contact points between the object and the vehicles are equidistant on the spherical object surface.
	\item The number of vehicles $n$ is chosen corresponding to the object radius such that there is adequate spacing between each two adjacent vehicles.
	\item The object mass ($33.90\,\textup{kg}$, $76.29\,\textup{kg}$, $135.65\,\textup{kg}$, $211.97\,\textup{kg}$, $305.26\,\textup{kg}$, $415.51\,\textup{kg}$, $542.73\,\textup{kg}$, $848.06\,\textup{kg}$) is chosen realistically in accordance with the its radius ($1.0\,\textup{m}$, $1.5\,\textup{m}$, $2.0\,\textup{m}$, $2.5\,\textup{m}$, $3.0\,\textup{m}$, $3.5\,\textup{m}$, $4.0\,\textup{m}$, $5.0\,\textup{m}$).
\end{itemize}

Once $\xv: \realset_{\geqslant 0} \to \realset^{\dimx}$ is computed by DCC or \SQP, the forward-dynamics simulation is run to obtain the trajectory of the object pose
\begin{align*}
	\pv_{\textup{o}} \in \Big\{
	\pv: \realset_{\geqslant 0} \to \realset^{6} 
	~\big|~
	\hv_{\textup{o,v}}\big(\pv(t), \dot{\pv}(t), \ddot{\pv}(t), \xv(t)\big)~\forall t \in [0\,\textup{s},~1.5\,\textup{s}]
	\Big\}.
\end{align*}
using the dynamics model
\begin{equation*}%\label{eq:drone-object dynamics compact}
	\begin{split}
		&\hv_{\textup{o,v}}\big(\pv(t), \dot{\pv}(t), \ddot{\pv}(t), \xv(t)\big) \\
		&:= \hv_{\textup{o}}\big(\pv(t), \dot{\pv}(t), \ddot{\pv}(t)\big)
		+ \sum_{i=1}^{n} \hv_{\textup{v}}\left(
		\psiv_{\textup{v}\langle i\rangle}\big(\pv(t)\big), 
		\frac{\mathrm{d}}{\mathrm{d}t} \psiv_{\textup{v}\langle i\rangle}\big(\pv(t)\big), 
		\frac{\mathrm{d}^{2}}{\mathrm{d}t^{2}} \psiv_{\textup{v}\langle i\rangle}\big(\pv(t)\big),
		\xv_{\langle i\rangle}(t)
		\right)
	\end{split}
\end{equation*}
of the object-vehicle complex. The results are shown in the Main Document, Figure~\ref{MD-fig:swarm tracking timeseries}.

\subsection{Example R2 - musculoskeletal model}
\paragraph*{Problem statement} 
The shoulder-arm model includes $12$ joints, $42$ muscles, and $2$ ligament elements~\cite{HuKueHad2020}. Each muscle is subdivided into multiple actuator units that may be considered as motor units. The functions $\fv_{\scriptscriptstyle\textup{A}}(\xv)$ and $\Bm(\xv)$ for state dynamics are obtained from the linearized muscle-activation dynamics \cite{HuKueHad2020} using a first-order approximation at $\uv = \nullv$.

The resulting optimization determines muscle forces consistent with the observed joint motion while satisfying activation bounds $\gv(t,\xv)$ that is $\xv \in [0,1]^{\dimx}$ ($\Rightarrow \dimc = 2\dimx$). 
Further, the nonlinear function $\hv(t,\xv)$ of $\xv$ produces the joint-space generalized force $\yv$, where the joint-space generalized coordinates $\pv(t)$ and the change rate $\dot{\pv}(t)$ are given explicitly as functions of $t$. 
Finally, the cost function is defined as 
\begin{equation}
	\sigma(t,\xv) = \frac{1}{2}\xv^\transp \Qm \xv
\end{equation}
with motor-task-specific muscle synergy matrix $\Qm$~\cite{HuKueHad2018}.
The glenohumeral joint reaction force is computed from the musculoskeletal model. An injury-risk probability is estimated using a cumulative distribution model based on experimentally observed joint capsule failure forces.

\paragraph*{Details} 
Each muscle is composed of $2^{n}-1$ actuators with different sizes of \uls{p}hysiological \uls{c}ross-\uls{s}ection \uls{a}rea (PCSA), while the combined force producing capability and characteristics of each muscle remain unchanged. This results in altogether $\dimx = 2 + 42\times(2^{n}-1)$ force producing elements. This extension combines in principle the muscle-force distribution~\cite{ErdemirMcLHerBog2007,Herzog2009a,HuKueHad2020} and the motor-unit recruitment~\cite{MarshallGlaTraAme2022}.

We calculate the reference signal $\desired{\yv}(t)$ for the joint-space generalized force from recorded human shoulder-arm movement $\pv(t)$ in \cite{HuKueHad2018,AvertaBarCatHad2021} for the duration $t_{\final} \approx 4.27\,\textup{s}$ and interpolated at $2$~kHz through the joint-space inverse-dynamics model~\cite{HuKueHad2020}.

\paragraph*{Capsule tear} 
Using intermediate quantities in the aforementioned $\hv(t,\xv)$ and joint-space inverse dynamics, we compute the \uls{g}leno\uls{h}umeral joint reaction force $\forcev_{\hspace{-0.51ex}\scriptscriptstyle\textup{GH}} \in \realset^{3}$,
which is exerted by connective tissues to maintain the joint structure~\cite{ChadwickBlaKirBog2014,VidtSanMarHeg2018}, improper force magnitude or direction may cause joint dislocation or tissue tear. As an exemplary monitoring function, we assess the probability $p_{\scriptscriptstyle\textup{GH}}$ of glenohumeral joint capsule tear as follows:

\begin{equation}
	p_{\scriptscriptstyle\textup{GH}} = \begin{cases}
		0 & \text{if~}\forcev_{\hspace{-0.51ex}\scriptscriptstyle\textup{GH}}\text{~points inside the \uls{g}lenoid \uls{f}ossa~\cite{ChadwickBlaKirBog2014,SurotoLicWibGul2022}}, \\
		\psi_{\hspace{-0.31ex}\scriptscriptstyle\textup{GH}}\big(\big\|\forcev_{\hspace{-0.51ex}\scriptscriptstyle\textup{GH}} - (\forcev_{\hspace{-0.51ex}\scriptscriptstyle\textup{GH}}^{\,\transp}\,\ev_{\hspace{-0.21ex}\scriptscriptstyle\textup{GH}}) \ev_{\hspace{-0.21ex}\scriptscriptstyle\textup{GH}}\big\|_{2}\big) & \text{if~}\forcev_{\hspace{-0.51ex}\scriptscriptstyle\textup{GH}}\text{~points outside GF}~\wedge~\forcev_{\hspace{-0.51ex}\scriptscriptstyle\textup{GH}}^{\,\transp}\,\ev_{\hspace{-0.21ex}\scriptscriptstyle\textup{GH}} > 0, \\
		\psi_{\hspace{-0.31ex}\scriptscriptstyle\textup{GH}}\big(\|\forcev_{\hspace{-0.51ex}\scriptscriptstyle\textup{GH}}\|_{2}\big) & \text{otherwise}.
	\end{cases}
\end{equation}

Herein, $\ev_{\hspace{-0.21ex}\scriptscriptstyle\textup{GH}} \in \realset^{3}$ is the unit vector normal to \uls{g}lenoid \uls{f}ossa (GF) and directing towards humeral head.
Therefore, a positive $\forcev_{\hspace{-0.51ex}\scriptscriptstyle\textup{GH}}^\transp\,\ev_{\hspace{-0.21ex}\scriptscriptstyle\textup{GH}}$ means that GF is able to produce the normal component $(\forcev_{\hspace{-0.51ex}\scriptscriptstyle\textup{GH}}^\transp\,\ev_{\hspace{-0.21ex}\scriptscriptstyle\textup{GH}}) \ev_{\hspace{-0.21ex}\scriptscriptstyle\textup{GH}}$ of $\forcev_{\hspace{-0.51ex}\scriptscriptstyle\textup{GH}}$.
The function $\psi_{\hspace{-0.31ex}\scriptscriptstyle\textup{GH}}: \realset \to [0,1]$ is the cumulative distribution function
\begin{equation*}
	\psi_{\hspace{-0.31ex}\scriptscriptstyle\textup{GH}}(F) = \frac{1}{2} + \frac{1}{2} \erf\left(\frac{1}{\sqrt{2}} \frac{F - 699.44\,\text{N}}{182.14\,\text{N}}\right)
\end{equation*}
of the normal distribution with mean $699.44\,\text{N}$ and standard deviation $182.14\,\text{N}$,
which is the force magnitude by which the glenohumeral joint capsule is torn~\cite{HaraItoIwa1996}.
The achieved results are depicted in the Main Document, Figure~\ref{MD-fig:ulmsk tracking timeseries}.

\subsection{Example R3 - portfolio model}

\paragraph*{Problem statement}
Financial market data consist of weekly returns for several asset markets. Expected returns and covariance matrices are estimated using a moving average window. We investigate \emph{the mean-variance problem with prescribed reward}~\cite{Steinbach2001}. 
We formulate the system dynamics as
\begin{equation}
	\dot{\xv} = -100\,\xv + 100\,\uv, 
\end{equation}
with a cost function 
\begin{equation}
	\sigma(t,\xv) = \frac{1}{2} \xv^\top \Qm(t) \xv
\end{equation}
that is subject to expected return and allocation constraints (see below). Portfolio performance is evaluated using weekly and annual return on investment.

\paragraph*{Details}
We deploy real-world data $\thetav(t)$ of weekly returns on 6 markets with $\dimx \in \{28, 49, 82,\allowbreak 83, 442, 1203\}$ assets over $> 10$ years until 2015/2016~\cite{BruniCesScoTar2016}.
We apply moving-average filter (window size 52 weeks) on $\thetav(t)$ to obtain time-varying expected returns $\pv(t) \in \realset^{\dimx}$ of and covariance matrix $\Qm(t)$ among $\dimx$ assets.  We select the data of the latest year where the established time-variant optimization problem~\eqref{eq:constrained optimization} is feasible and interpolated them to 2000 samples/year.

The system's equality and inequality constraints as well as the reference trajectory are defined as
\begin{align*}
	&\hv(t,\xv) = \begin{bmatrix}
		\pv(t)^\transp \xv \\
		\sum_{i=1}^{\dimx} x_{i}
	\end{bmatrix},
	\quad 
	\desired{\yv}(t) = \begin{cases}
		[\,0~~0\,]^\transp & \text{if~}\nexists i: p_{i}(t) > 0\\
		[\,p_{\textup{am}}(t)~~1\,]^\transp & \text{otherwise}
	\end{cases},
	\quad \\
	&\gv(t,\xv) = \begin{bmatrix}
		\xv - 1 \\
		-\xv
	\end{bmatrix}.
\end{align*}
where $p_{\textup{am}}(t)$ is the arithmetic mean of the positive elements in $\pv(t)$.

The actual \uls{r}eturn \uls{o}n \uls{i}nvestment (ROI) $R$ following the portfolio state $\xv$ is
\begin{equation*}
	R_{\scriptscriptstyle\textup{weekly}}(t) = \begin{cases}
		~0 & \text{if~} \|\xv(t)\|_{1} = 0, \\[1ex]
		\displaystyle\frac{\thetav(t)^\transp \xv(t)}{\|\xv(t)\|_{1}} & \text{otherwise},
	\end{cases}
\end{equation*}
which is the weekly ROI, since the input data are given on a weekly basis.
Consequently, the annual ROI is
\begin{equation*}
	R_{\scriptscriptstyle\textup{annual}}
	= -1 + \prod_{k=0}^{52}\Big( 1 + R_{\scriptscriptstyle\textup{weekly}}\big(k \Delta T\big) \Big)
\end{equation*}
with $\Delta T = 1\,\textup{week}$.
The results are shown in the Main Document, Figure~\ref{MD-fig:portfolio tracking timeseries}.

\begin{table}[!h]
	\centering
	\caption{\label{tab:real-world examples}Computational behavior of all real-world-relevant examples.}
	\begin{minipage}{\textwidth}
		\scalefont{0.75000000187}%
		% !TEX root = ../main.tex
%\begin{threeparttable}
%\setlength{\tabcolsep}{2pt}
%\setlength{\extrarowheight}{1pt}
%\begin{tabularx}{1.15\textwidth}{|r||>{\arraybackslash}X|>{\arraybackslash}X|>{\arraybackslash}X|>{\centering\arraybackslash}X|>{\centering\arraybackslash}X|>{\centering\arraybackslash}X||c|c||r|r||l|}
\begin{tabular}{|r||r|r|r|c|c|c||c|c||l|}
	\hline
	& \multicolumn{6}{c||}{$\min_{\chiv} \sigma(t,\chiv)~~\textup{s.t.}~\nullv=\hv(t,\chiv)-\desired{\yv}(t)\,\wedge\,\nullv \geqslant \gv(t,\chiv)$} & \multicolumn{2}{c||}{Solvers} & \\
	$n$ & $\dimx$ & $\dimy$ & $\dimc$ & \makecell{Hessian\\ $\Qminv:=\frac{\partial}{\partial \chiv}\frac{\partial \sigma(t,\chiv)}{\partial \chiv}$} & \makecell{Equality\\ constraint} & \makecell{Inequality\\ constraint} & \OED & \SQP & Remark \\
	\hline\hline
	\multicolumn{9}{|l||}{Example R1: control allocation among $n$ drones moving a spherical object} & Object's radius \& mass \\ \hline
	16 &  \eval{24*16} & \multirow{8}{*}{6} &  \eval{(1+2*24)*16} &  \multirow{8}{*}{\makecell{identity\\ matrices}} & \multirow{8}{*}{\makecell{time-variant\\ linear}} & \multirow{8}{*}{\makecell{time-variant\\ linear}} & \Checkmark & \Checkmark & 1.0\,m, 33.90\,kg \\
	36 &  \eval{24*36} &                    & \eval{(1+2*24)*36}  &                           & & & \Checkmark & \Checkmark & 1.5\,m, 76.29\,kg \\
    64 &  \eval{24*64} &                    & \eval{(1+2*24)*64}  &                           & & & \Checkmark & \Checkmark & 2.0\,m, 135.65\,kg \\
   100 & \eval{24*100} &                    & \eval{(1+2*24)*100} &                           & & & \Checkmark & \Checkmark & 2.5\,m, 211.97\,kg \\
   144 & \eval{24*144} &                    & \eval{(1+2*24)*144} &                           & & & \Checkmark & \Checkmark & 3.0\,m, 305.26\,kg \\
   196 & \eval{24*196} &                    & \eval{(1+2*24)*196} &                           & & & \Checkmark & \Checkmark & 3.5\,m, 415.51\,kg \\
   256 & \eval{24*256} &                    & \eval{(1+2*24)*256} &                           & & & \Checkmark & \Checkmark & 4.0\,m, 542.73\,kg \\
   400 & \eval{24*400} &                    & \eval{(1+2*24)*400} &                           & & & \Checkmark & \Checkmark$\strut^{\dagger}$ & 5.0\,m, 848.06\,kg \\
   \hline\hline
   \multicolumn{10}{|l|}{Example R2: muscle-force distribution among $\approx \dimx$ muscle actuators for performing a limb movement} \\ \hline
%   1 &  \eval{2+42*(2-1)} & \multirow{20}{*}{12} &  \eval{2*(2+42*(2-1))} &  \multirow{5}{*}{\makecell{identity\\ matrices}} & \multirow{20}{*}{\makecell{time-variant\\ nonlinear}} & \multirow{20}{*}{\makecell{constant\\ box}} & \Checkmark & \Checkmark & \\
%   2 &  \eval{2+42*(4-1)} &                     &  \eval{2*(2+42*(4-1))} &                           & & & \Checkmark & \Checkmark & \\
%   3 &  \eval{2+42*(8-1)} &                     &  \eval{2*(2+42*(8-1))} &                           & & & \Checkmark & \Checkmark & \\
%   4 & \eval{2+42*(16-1)} &                     & \eval{2*(2+42*(16-1))} &                           & & & \Checkmark & \Checkmark & \\
%   5 & \eval{2+42*(32-1)} &                     & \eval{2*(2+42*(32-1))} &                           & & & \Checkmark & \Checkmark & \\
%   \cline{1-2}\cline{4-5}\cline{8-9}
%   %
%   1 &  \eval{2+42*(2-1)} &                     &  \eval{2*(2+42*(2-1))} & \multirow{5}{*}{\makecell{constant\\ diagonal\\ matrices}} & & & \Checkmark & \Checkmark & \\
%   2 &  \eval{2+42*(4-1)} &                     &  \eval{2*(2+42*(4-1))} &                           & & & \Checkmark & \Checkmark & \\
%   3 &  \eval{2+42*(8-1)} &                     &  \eval{2*(2+42*(8-1))} &                           & & & \Checkmark & \Checkmark & \\
%   4 & \eval{2+42*(16-1)} &                     & \eval{2*(2+42*(16-1))} &                           & & & \Checkmark & \Checkmark & \\
%   5 & \eval{2+42*(32-1)} &                     & \eval{2*(2+42*(32-1))} &                           & & & \Checkmark & \Checkmark & \\
%   \cline{1-2}\cline{4-5}\cline{8-9}
   %
   1 &  \eval{2+42*(2-1)} & \multirow{5}{*}{12} &  \eval{2*(2+42*(2-1))} & \multirow{5}{*}{\makecell{constant\\ fully populated\\ matrices}} & \multirow{5}{*}{\makecell{time-variant\\ nonlinear}} & \multirow{5}{*}{\makecell{constant\\ box}} & \Checkmark & \Checkmark & \\
   2 &  \eval{2+42*(4-1)} &                     &  \eval{2*(2+42*(4-1))} &                           & & & \Checkmark & \Checkmark & \\
   3 &  \eval{2+42*(8-1)} &                     &  \eval{2*(2+42*(8-1))} &                           & & & \Checkmark & \Checkmark & \\
   4 & \eval{2+42*(16-1)} &                     & \eval{2*(2+42*(16-1))} &                           & & & \Checkmark & \Checkmark & \\
   5 & \eval{2+42*(32-1)} &                     & \eval{2*(2+42*(32-1))} &                           & & & \Checkmark & \Checkmark & \\
%   \cline{1-2}\cline{4-5}\cline{8-9}
   %
%   1 &  \eval{2+42*(2-1)} &                     &  \eval{2*(2+42*(2-1))} & \multirow{5}{*}{\makecell{constant\\ fully populated\\ matrices\smash{$^{\S}$}}} & & & \Checkmark & \Checkmark & \\
%   2 &  \eval{2+42*(4-1)} &                     &  \eval{2*(2+42*(4-1))} &                           & & & \Checkmark & \Checkmark & \\
%   3 &  \eval{2+42*(8-1)} &                     &  \eval{2*(2+42*(8-1))} &                           & & & \Checkmark & \Checkmark & \\
%   4 & \eval{2+42*(16-1)} &                     & \eval{2*(2+42*(16-1))} &                           & & & \Checkmark & \Checkmark & \\
%   5 & \eval{2+42*(32-1)} &                     & \eval{2*(2+42*(32-1))} &                           & & & \Checkmark & \Checkmark & \\
   \hline\hline
   \multicolumn{10}{|l|}{Example R3: portfolio optimization for $\dimx$ assets} \\ \hline
   1 & 28 & \multirow{6}{*}{2} & \eval{2*28} & \multirow{5}{*}{\makecell{~\\time-variant\\ fully populated\\ matrices}}& \multirow{6}{*}{\makecell{time-variant\\ linear}} & \multirow{6}{*}{\makecell{constant\\ box}} & \Checkmark & \Checkmark & \\
   2 & 49 &                    & \eval{2*49} &                           & & & \Checkmark & \Checkmark & \\	
   3 & 82 &                    & \eval{2*82} &                           & & & \Checkmark & \Checkmark & \\
   4 & 83 &                    & \eval{2*83} &                           & & & \Checkmark & \Checkmark & \\
   5 & 442 &                   & \eval{2*442} &                          & & & \Checkmark & \Checkmark & \\
   6 & 1203 &                  & \eval{2*1203} &                         & & & \Checkmark$\strut^{\ddagger}$ & \Checkmark$\strut^{\ddagger}$ & \\		
   \hline
\end{tabular}%
\begin{tablenotes}
	\item[~] All input data were created at or interpolated to 2000 samples per second (Examples R1 and R2) or per year (Example R3), except
	those with $^\dagger$500 samples/second and 
	$^\ddagger$1000 samples/year %, since $\Qm(t) \in \realset^{1203\times 1203}$ with 2000 samples/year exceeds the maximum size of an array in MEX.
	due to practical limitations.
	%\item[$^\S$] These two types of constant fully populated matrices reflect different muscle-coactivation pattern.
\end{tablenotes}
%\end{threeparttable}%
		\scalefont{1.33333333}%
	\end{minipage}
\end{table}
~\newpage~

\newpage
% !TEX root = ../supplementary.tex
%
\section{Computational complexity}\label{sec:computational complexity}
In this section, we delve into a comprehensive analysis and comparison of computational complexity of \OTC\ alongside the well-established and popular state-of-the-art methods \uls{i}nterior-\uls{p}oint \uls{a}lgorithm (\IPA)~\cite{BoydVan2004,NocedalWri2006,TheMathWorks2020} and \uls{s}equential \uls{q}uadratic \uls{p}rogramming (\SQP) \cite{Fletcher1987,NocedalWri2006,TheMathWorks2020}.
%By examining the computational complexities, we show that the \OTC\ not only exhibits the lowest complexity but also offers reliable performance for real-time computing.
{In the analysis, we adopt} the Bachmann-Landau big $\bigO$ notation {to denote the asymptotic upper
bound of complexity}~\cite{CormenLeiRivSte2009}{, and} $\timecomplexity(\Omega)$ and $\spacecomplexity(\Omega)$ to denote the time and space complexity, respectively, of performing operation or procedure $\Omega$.
%We suppose that $\Qm^{-1}$ is given (no need to calculate the inverse of $\Qm$), and we do not count the memory space required for the input data ($\Qm$, $\desired{\yv}$, $\Gm$, and $\cv$).

\subsection{Dual cost--constraint projection}
\paragraph{Time complexity}
The time complexity of \OTC~\eqref{eq:gryphonheart} can be derived by counting the number of arithmetic operations, {and is dominated by 
\begin{itemize}
	\item inverting $\Bm(\xv) \in \realset^{\dimx\times\dimx}$ using lower-upper decomposition followed by forward and back substitution,
	\item constituting $\rv(t,\xv) \in \realset^{\dimx}$ and $\Rm(t,\xv) \in \realset^{\dimx\times\dimx}$,
	\item inverting $\Lambdam(t,\xv,\nullv) \in \realset^{(\dimx+\dimy)\times(\dimx+\dimy)}$ using Cholesky decomposition followed by forward and back substitution.
\end{itemize}
}%
This results in
{%
\begin{equation}\label{eq:time complexity OTC}
\begin{split}
	\timecomplexity(\OTC)
	&= \left(
	\frac{1}{3}\dimx^{3} 
	+ \frac{1}{6}(\dimx+\dimy)^{3} + (\dimx+\dimy)^{2} 
	+ \dimx + 2\dimcbar\dimx
	+ \dimx^{2} + \dimcbar\dimx^{2}
	\right) \big(\timecomplexity(\times) + \timecomplexity(+)\big)
	\\[-0.5ex]
	&\quad+ \left(
	\frac{1}{2}\dimx^{2}
	+ \frac{1}{2}(\dimx+\dimy)^{2} + 2(\dimx+\dimy)
	\right) \timecomplexity(\div) %\tau_{\CPUsubscript,\div}
	+ \big(
	\dimx + \dimy
	\big) \timecomplexity(\sqrt{\cdot}\,) \\
	&\quad+ \timecomplexity(\qv) + \timecomplexity(\Qm) + \timecomplexity(\hv) + \timecomplexity(\Hm) + \timecomplexity(\Bm) + \timecomplexity(\fv_{\scriptscriptstyle\textup{A}}),
\end{split}
\end{equation}%
where $\dimcbar \in [0, \dimc]$ is the time-variant number of positive elements in $\gv(t,\xv)$.
%$\dimcbar := \dim\big(\breve{\gv}(t,\xv)\big) \leqslant \dimc$
}

\paragraph{Space complexity}
The required memory space is dominated by the matrices $\Rm(t,\xv) \in \realset^{\dimx\times\dimx}$ and $\Hm(t,\xv) \in \realset^{\dimy\times\dimx}$.
This means
\begin{equation}\label{eq:space complexity OTC}
	\begin{split}
		\spacecomplexity(\OTC)
		= \bigO(\dimx^{2} + \dimy\dimx).
	\end{split}
\end{equation}

\subsection{Interior-point algorithm}
Although the state-of-the-art \IPA~\cite{NocedalWri2006,TheMathWorks2020} differs from the one we constructed in Section~\ref{sec:barrier method}, it is still a sequence of Newton-Raphson iterations {for finding the saddle point of a Lagrangian function} with a $(\dimx+\dimy+2\dimc)$-by-$(\dimx+\dimy+2\dimc)$ {symmetric Hessian matrix}. 
After sufficiently many iterations, the resulting $\chiv^{\optimal}$ must be converted to $\uv$ in order to drive \eqref{eq:gryphon}. The simplest choice would be
\begin{equation}\label{eq:chi_optimal to u}
	\uv = \Bm(\xv)^{-1} \big(
	-\fv_{\textup{A}}(\xv) 
	+ K_{\chi} (\chiv^{\optimal} - \xv)
	\big),
	\quad
	K_{\chi} > 0.
\end{equation}
\paragraph{Time complexity}
By assuming the symmetric {Hessian} to be also positive definite, one concludes {%
\begin{equation}\label{eq:time complexity IPA per iteration}
	\begin{split}\hspace{-1ex}
		T_{\CPUsubscript,\textup{iter}}^{\IPAsubscript}(\dimx,\dimy,\dimc)
		&= \left(
		\frac{1}{6}(\dimx \hspace{-0.3ex}+\hspace{-0.3ex} \dimy \hspace{-0.3ex}+\hspace{-0.3ex} 2\dimc)^{3} + (\dimx \hspace{-0.3ex}+\hspace{-0.3ex} \dimy \hspace{-0.3ex}+\hspace{-0.3ex} 2\dimc)^{2}
		\right) \big(\timecomplexity(\times) + \timecomplexity(+)\big)
		\\[-0.5ex]
		&\quad+ \left(
		\frac{1}{2}(\dimx \hspace{-0.3ex}+\hspace{-0.3ex} \dimy \hspace{-0.3ex}+\hspace{-0.3ex} 2\dimc)^{2} + 2(\dimx \hspace{-0.3ex}+\hspace{-0.3ex} \dimy \hspace{-0.3ex}+\hspace{-0.3ex} 2\dimc)
		\right) \timecomplexity(\div) \\
		&\quad+ \big(
		\dimx \hspace{-0.3ex}+\hspace{-0.3ex} \dimy \hspace{-0.3ex}+\hspace{-0.3ex} 2\dimc
		\big) \timecomplexity(\sqrt{\cdot}\,)
		+ \timecomplexity(\qv) + \timecomplexity(\Qm) + \timecomplexity(\hv) + \timecomplexity(\Hm)
	\end{split}
\end{equation}%
}arithmetic operations per \uls{iter}ation.
{%
Under the best-case assumption, i.e.,
\begin{itemize}
	\item \eqref{eq:constrained optimization} is a convex optimization problem and 
	\item the initial guess $\chiv^{[0]}(t)$ at every time instant $t \geqslant 0$ is sufficiently close to a local minimizer and satisfies all inequality and equality constraints,
\end{itemize}%
}%
the \IPA\ requires \begin{equation}\label{eq:iteration complexity IPA}
	N_{\FcnEvalsubscript}^{\IPAsubscript}(t,\dimc) = \bigO\Big(\sqrt{1+\dimc} \log_{2}\big( (1+\dimc) \|\chiv^{[0]} - \chiv^{\optimal}(t)\|_{2} / \varepsilon \big)\Big)
\end{equation}
iterations to achieve \smash{$\|\chiv^{[N_{\text{\tiny FE}}^{\text{\tiny IPA}}]} - \chiv^{\optimal}(t)\|_{2} \leqslant \varepsilon$} \cite{PotraWri2000,BoydVan2004}.
While the target precision $\varepsilon$ can be set constant, the local minimizer $\chiv^{\optimal}(t)$ is time-variant due to the time-variant problem~\eqref{eq:constrained optimization}.
In general, one cannot choose the initial guess to be immediately close to $\chiv^{\optimal}(t)$ without introducing further algorithmic procedures; therefore, \smash{$\|\chiv^{[0]} - \chiv^{\optimal}(t)\|_{2} / \varepsilon \gg 1$} is a large number.
As a conclusion, the \IPA-based controller's time complexity is
{%
\begin{equation}\label{eq:time complexity IPA}
\begin{split}
	\timecomplexity(\IPA) 
	&= \frac{1}{3}\dimx^{3} \big(\timecomplexity(\times) + \timecomplexity(+)\big)
	+ \frac{1}{2}\dimx^{2} \timecomplexity(\div)
	+ \timecomplexity(\Bm) + \timecomplexity(\fv_{\scriptscriptstyle\textup{A}}) \\[-1ex]
	&\quad+ N_{\FcnEvalsubscript}^{\IPAsubscript}(t,\dimc)\, T_{\CPUsubscript,\textup{iter}}^{\IPAsubscript}(\dimx,\dimy,\dimc).
\end{split}
\end{equation}
}

\paragraph{Space complexity}
The required memory space is dominated by the {Hessian matrix}. This means
\begin{equation}\label{eq:space complexity IPA}
	\spacecomplexity(\IPA) = \bigO\big(\dimx^{2}+(\dimy+2\dimc)\dimx\big).
\end{equation}
{Note that the Hessian matrix is symmetric and it has zero-valued elements.}

\subsection{Sequential quadratic programming}
As an iterative method, \SQP~\cite{Fletcher1987,NocedalWri2006,TheMathWorks2020} approximates a general nonlinear optimization problem with a sequence of \uls{q}uadratic \uls{p}rogramming (QP) subproblems and solves these QP subproblems 
%(quadratic cost function and linear constraints) 
at each iteration. 
In each iteration, the state-of-the-art \SQP\ implementation \cite{TheMathWorks2020} deploys an active-set method,
%(also known as projection method), 
which undertakes three steps:
\begin{itemize}
	\item Step 1: Establish $\dimy + \dimcbar$ \emph{active constraints} from all $\dimy$ equality constraints and the $\dimcbar \in [0, \dimc]$ inequality constraints on whose boundaries the current $\chiv$ is. This turns the QP subproblem with equality and inequality constraints to an equality-constrained QP problem.
	\item Step 2: Solve the equality-constrained QP problem analytically using the method of Lagrange multipliers, {i.e., find a saddle point of a Lagrangian function} with a $(\dimx+\dimy+\dimcbar)$-by-$(\dimx+\dimy+\dimcbar)$ symmetric {Hessian} matrix. The resulting $\chiv$ and the associated Lagrange multiplier $\zetav$, relatively to those in the previous iteration, provides a search direction.
	\item Step 3: Given the search direction, a line-search procedure iteratively determines the step length that provides a sufficient decrease in a merit function. The merit function is constituted from the cost function and \emph{all} constraints of the general nonlinear optimization problem.
\end{itemize}
After {a sufficient number $N_{\textup{iter}}^{\SQPsubscript}$ of \SQP} \uls{iter}ations, the resulting $\chiv^{\optimal}$ must be converted to $\uv$ via, e.g., \eqref{eq:chi_optimal to u}.
\paragraph{Time complexity}
Hence, we conclude the time complexity of {one single} \SQP\ {iteration}
{%
\begin{equation}\label{eq:time complexity SQP per iteration}
\begin{split}\hspace{-1ex}
	T_{\CPUsubscript,\textup{iter}}^{\SQPsubscript}(\dimx,\dimy,\dimcbar,N_{\textup{ls}}^{\SQPsubscript})
	&= \left(
	\frac{1}{6}(\dimx \hspace{-0.3ex}+\hspace{-0.3ex} \dimy \hspace{-0.3ex}+\hspace{-0.3ex} \dimcbar)^{3} + (\dimx \hspace{-0.3ex}+\hspace{-0.3ex} \dimy \hspace{-0.3ex}+\hspace{-0.3ex} \dimcbar)^{2}
	\right) \big(\timecomplexity(\times) + \timecomplexity(+)\big)
	\\[-0.5ex]
	&\quad+ \left(
	\frac{1}{2}(\dimx \hspace{-0.3ex}+\hspace{-0.3ex} \dimy \hspace{-0.3ex}+\hspace{-0.3ex} \dimcbar)^{2} + 2(\dimx \hspace{-0.3ex}+\hspace{-0.3ex} \dimy \hspace{-0.3ex}+\hspace{-0.3ex} \dimcbar)
	\right) \timecomplexity(\div) \\
	&\quad+ \big(
	\dimx \hspace{-0.3ex}+\hspace{-0.3ex} \dimy \hspace{-0.3ex}+\hspace{-0.3ex} \dimcbar
	\big) \timecomplexity(\sqrt{\cdot}\,)
	+ \timecomplexity(\qv) + \timecomplexity(\Qm) + \timecomplexity(\hv) + \timecomplexity(\Hm)\\
	&\quad+ N_{\textup{ls}}^{\SQPsubscript} \times \big(\timecomplexity(\sigma) + \timecomplexity(\gv) + \timecomplexity(\hv)\big),
\end{split}
\end{equation}%
where $N_{\textup{ls}}^{\SQPsubscript} \geqslant 0$ is the number of \uls{l}ine-\uls{s}earch iterations.
}%

Since the active constraints (both size and content) vary over iterations, it is difficult to determine {formulae for upper bounds of $N_{\textup{iter}}^{\SQPsubscript}$ and $N_{\textup{ls}}^{\SQPsubscript}$}~\cite{CiminiBem2017,ArnstroemAxe2022}. 
While several works \cite{CiminiBem2019,ArnstroemAxe2022} focus on improving the update strategy of active constraints,
we analyze in this work the \SQP's iteration complexity under the following considerations:
\begin{enumerate}
	\item \SQP\ is a quasi-Newton method whose convergence behavior is described by Lemma~\ref{lemma:chord method} (mathematically proven \cite{Fletcher1987,NocedalWri2006}).
	\item If the initial guess satisfies all constraints, then \SQP\ requires fewer iterations than \IPA~(empirical knowledge \cite{NocedalWri2006}).
\end{enumerate}
Therefore, we conclude
{$N_{\textup{ls}}^{\SQPsubscript} = 1$ for all \SQP\ iterations}
and %
{$N_{\textup{iter}}^{\SQPsubscript} = {N_{\FcnEvalsubscript}^{\SQPsubscript}}/{N_{\textup{ls}}^{\SQPsubscript}} \leqslant N_{\FcnEvalsubscript}^{\IPAsubscript}$} (cf. \eqref{eq:iteration complexity IPA}), respectively, as the best-case upper-bound estimate.
In summary, {
\begin{equation}\label{eq:time complexity SQP}
	\timecomplexity(\SQP) 
	= N_{\FcnEvalsubscript}^{\SQPsubscript} \times T_{\CPUsubscript,\textup{iter}}^{\SQPsubscript}(\dimx,\dimy,\dimcbar,1)
	\leqslant
	N_{\FcnEvalsubscript}^{\IPAsubscript}(t,\dimc) \times T_{\CPUsubscript,\textup{iter}}^{\SQPsubscript}(\dimx,\dimy,\dimcbar,1),
\end{equation}
}%
if the initial guess is sufficiently close to a local minimizer and satisfies all inequality and equality constraints.
{The upper bound of} $\timecomplexity(\SQP)$ is time-variant, inheriting from {$N_{\FcnEvalsubscript}^{\IPAsubscript}(t,\dimc)$}.

\paragraph{Space complexity}
The required memory space is dominated by the {Hessian} matrix in Step 2. This means
%\begin{equation}\label{eq:space complexity, SQP}
%	\begin{split}
%		\spacecomplexity(\SQP) 
%		&= \bigO\big((\dimx+\dimy+\dimc)^{2}\big).
%	\end{split}
%\end{equation}
\begin{equation}\label{eq:space complexity SQP}
	\spacecomplexity(\IPA) =
	\bigO\big(\dimx^{2} + (\dimy\hspace{-0.3ex}+\hspace{-0.3ex}\dimc)\dimx\big).
\end{equation}
{Note that the Hessian matrix is symmetric and it has zero-valued elements.}

\subsection{Summary}
The comparison \eqref{eq:time complexity SQP} vs \eqref{eq:time complexity IPA} shows that \SQP\ requires less number of \uls{fl}oating-point \uls{op}erations (FLOP) per iteration than \IPA, whereas the comparison \eqref{eq:time complexity OTC} vs \eqref{eq:time complexity SQP} concludes that \OTC\ requires less number of FLOP per time step than \SQP, if $N_{\FcnEvalsubscript}^{\SQPsubscript} > 1$.
Further, the memory space that \OTC\ requires is less than \SQP.

The number $N_{\FLOPsubscript}^{\OTCsubscript}$ of FLOP per time step required by \OTC\ can be estimated as
\begin{equation}\label{eq:time complexity OTC, simplified}
	N_{\FLOPsubscript}^{\OTCsubscript} 
	= \frac{\timecomplexity(\OTC)}{\timecomplexity(\times)}
	\approx \frac{2}{3}\dimx^{3} + 2 \dimcbar\dimx^{2} + \frac{1}{3}\big(\dimx+\dimy\big)^{3}
\end{equation}
for sufficiently large $\dimx \gg 1$, assuming $\timecomplexity(\times) = \timecomplexity(+)$ and that $\timecomplexity(\qv) + \timecomplexity(\Qm) + \timecomplexity(\hv) + \timecomplexity(\Hm) + \timecomplexity(\Bm) + \timecomplexity(\fv_{\scriptscriptstyle\textup{A}})$ does not grow faster than $\bigO(\dimx^{3})$.
Analogously,
\begin{equation}\label{eq:time complexity SQP, simplified}
	N_{\FLOPsubscript}^{\SQPsubscript} 
	= \frac{\timecomplexity(\SQP)}{\timecomplexity(\times)}
	\approx \frac{2}{3}\dimx^{3} + N_{\FcnEvalsubscript}^{\SQPsubscript} \times \frac{1}{3}\big(\dimx \hspace{-0.3ex}+\hspace{-0.3ex} \dimy \hspace{-0.3ex}+\hspace{-0.3ex} \dimcbar\big)^{3}.
\end{equation}
For concrete examples, the empirical values for \smash{$\breve{\dimc}/\dimx$} % ca. 0.15
and $\dimy/\dimx$ may be known and (approximately) fixed, allowing one to further simplify \eqref{eq:time complexity OTC, simplified} and \eqref{eq:time complexity SQP, simplified} to functions of $\dimx$.

\newpage
% !TEX root = ../supplementary.tex
%
\section{Real-time barrier analysis}\label{sec:realtime barrier analysis}

\subsection{Real-time computer emulation}
To ensure a fair comparison under identical computational budgets, we emulate a computing platform whose throughput is set to ${1.30\,\dimx^{3}}/{\dt}$ [FLOP/s], corresponding to the minimum requirement of DCC (see the Main Document, Equation~\eqref{MD-eq:OTC flop/s}). 
This emulation is achieved via suitable postprocessing of simulation results obtained on arbitrary computing platforms.

The postprocessing requires three quantities:
\begin{enumerate}
	\item CPU time $T_{{\scriptscriptstyle\textup{C},\OTCsubscript}}(t) = T_{{\scriptscriptstyle\textup{C},\OTCsubscript}}$ for a time step required by offline simulated DCC,
	\item output $\chiv^{[k^{\prime}]}(t)$ of a prematurely terminated SQP and
	\item its corresponding CPU time $T_{{\scriptscriptstyle\textup{C},\SQPsubscript}}(t)$.
\end{enumerate}
Premature termination means that the SQP iteration is stopped as soon as the number of function evaluations reaches $k^{\prime}$, which is the smallest value such that $\big\|\chiv^{[k^{\prime}]}(t) - \chiv^{\optimal}(t)\big\|_{2} \leqslant\allowbreak \big\|\xv^{\OTCsubscript}(t) - \chiv^{\optimal}(t)\big\|_{2}$~$\forall t > 0$ and is determined empirically for each problem instance.

With these quantities, the equivalent time on a real-time computing platform is
\begin{equation}\label{eq:online SQP time delay}
	t + \frac{T_{{\scriptscriptstyle\textup{C},\SQPsubscript}}(t)}{T_{{\scriptscriptstyle\textup{C},\OTCsubscript}}} \dt
	= t + \gamma T_{{\scriptscriptstyle\textup{C},\SQPsubscript}}(t)
\end{equation}
with $\gamma := \dt/T_{{\scriptscriptstyle\textup{C},\OTCsubscript}}$ and control sampling interval $\dt$.
Note that the required CPU time $T_{{\scriptscriptstyle\textup{C},\OTCsubscript}}$ of DCC serves as the baseline reference, effectively emulating a computing platform whose throughput matches the DCC's minimum requirement.	

The time~\eqref{eq:online SQP time delay} represents the delay to which the SQP result becomes available:
\begin{equation}\label{eq:online SQP delayed control action}
	\xv^{\SQPsubscript}\big(t + \gamma T_{{\scriptscriptstyle\textup{C},\SQPsubscript}}(t)\big)
	= \chiv^{[k^{\prime}]}(t).
\end{equation}
While \SQP\ processes the latest input data at $t + \gamma T_{{\scriptscriptstyle\textup{C},\SQPsubscript}}(t)$, the value of $\xv^{\SQPsubscript}$ holds constant until the subsequent \SQP\ result becomes available: 
\begin{equation}\label{eq:online SQP skip time steps} 
	\xv^{\SQPsubscript}(\tau) 
	= \chiv^{[k^{\prime}]}(t)
	\quad
	\forall \tau \in \Big[
	t + \gamma T_{{\scriptscriptstyle\textup{C},\SQPsubscript}}(t),~ 
	t + \gamma T_{{\scriptscriptstyle\textup{C},\SQPsubscript}}(t) + \gamma T_{{\scriptscriptstyle\textup{C},\SQPsubscript}}\big(t + \gamma T_{{\scriptscriptstyle\textup{C},\SQPsubscript}}(t)\big)
	\Big).
\end{equation}

For DCC there exists no difference between offline and online $\xv^{\OTCsubscript}(t)$.

A practical algorithm for sorting $\chiv^{[k^{\prime}]}(t)$ into $\xv^{\SQPsubscript}(t)$ as described above is provided in Algorithm~\ref{alg:emulation of realistic x_SQP}. The resulting time series $\xv^{\SQP}(t)$ is referred to as \emph{online SQP}. The delayed control actions~\eqref{eq:online SQP delayed control action} and skipping multiple consecutive time steps~\eqref{eq:online SQP skip time steps} are the primary source of tracking errors and constraint violations.

Note that this emulation does not require any access to high-performance computing resources as it can be directly used on any computational hardware, including low-performance devices.
Therefore, this emulation serves the fair and consistent evaluation of real-time performance between high-dimensional DCC and SQP, independently of a particular technology.

\begin{algorithm}[!t]
	\caption{
		Emulation of the trajectory $\xv^{\scriptscriptstyle\text{\SQP}}(t)$ from $t = 0$ to $t = t_{\final}$ with consideration of response delay due to the high computation time $T_{{\scriptscriptstyle\textup{C},\SQPsubscript}}$ of \SQP.
	}
	\label{alg:emulation of realistic x_SQP}
	\begin{algorithmic}[1]
		\Require $T_{{\scriptscriptstyle\textup{C},\SQPsubscript}}: \realset_{\geqslant 0} \to \realset_{> 0}$,~~
		$T_{{\scriptscriptstyle\textup{C},\OTCsubscript}}: \realset_{\geqslant 0} \to \realset_{> 0}$,~~
		$\chiv^{[k^{\prime}]}: \realset_{\geqslant 0} \to \realset^{\dimx}$
		\Ensure $\xv^{\scriptscriptstyle\text{\SQP}}: \realset_{\geqslant 0} \to \realset^{\dimx}$
		\State $N \gets \lfloor t_{\final}/\dt \rfloor$;
		\State $n \gets 1$;
		\While{$n \leqslant N$} \Comment{$n_{1}$: the time step when the result $\chiv^{[k^{\prime}]}(n\dt)$ for time step $n$ is}
		\State $n_{1} \gets n + \lceil T_{{\scriptscriptstyle\textup{C},\SQPsubscript}}(n\dt) / T_{{\scriptscriptstyle\textup{C},\OTCsubscript}}(n\dt) \rceil$; \Comment{delivered}
		\If{$n_{1} > N$}
		\State $\xv^{\scriptscriptstyle\text{\SQP}}(n\dt) \gets \chiv^{[k^{\prime}]}(n\dt)$;
		\Else \Comment{$n_{2}$: the time step when $\chiv^{[k^{\prime}]}(n_{1}\dt)$ is delivered}
		\State $n_{2} \gets n_{1} + \lceil T_{{\scriptscriptstyle\textup{C},\SQPsubscript}}(n_{1}\dt) / T_{{\scriptscriptstyle\textup{C},\OTCsubscript}}(n_{1}\dt) \rceil$; 
		\ForAll{$\nu \in [n_{1}, n_{2}]$} 
		\State $\xv^{\scriptscriptstyle\text{\SQP}}(\nu\dt) \gets \chiv^{[k^{\prime}]}(n\dt)$; \Comment{$\chiv^{[k^{\prime}]}(n\dt)$ is the latest available result}
		\EndFor
		\EndIf
		\State $n \gets n_{1}$;
		\EndWhile
		\State \Return $\xv^{\scriptscriptstyle\text{\SQP}}$;
	\end{algorithmic}
\end{algorithm}

\subsection{Real-time performance metrics}
In order to quantitatively assess the fundamental real-time performance of DCC and online \SQP, we define the \emph{real-time performance metrics}
\begin{itemize}
	\item Cost-function error:
	\begin{equation}\label{eq:cost-function error}
		\left|\,
		\frac{\int_{0}^{t_{\final}}\sigma\big(t,\xv(t)\big)}{\int_{0}^{t_{\final}}\sigma\big(t,\chiv^{\optimal}(t)\big)} - 1\,
		\right|\times 100\%,
	\end{equation}
	\item Task-constraint violation:
	\begin{equation}\label{eq:task-constraint violation}
		\frac{\int_{0}^{t_{\final}}\big\|\desired{\yv}(t) - \hv\big(t,\xv(t)\big)\big\|_{2}\mathrm{d}t}{\int_{0}^{t_{\final}}\big\|\desired{\yv}(t)\big\|_{2}\mathrm{d}t} \times 100~\%,~\text{and}
	\end{equation}
	\item State tracking error:
	\begin{equation}\label{eq:state tracking error}
		\frac{1}{t_{\final}}\int_{0}^{t_{\final}}\frac{1}{\dimx}\big\|\xv(t) - \chiv^{\optimal}(t)\big\|_{1}\mathrm{d}t.
	\end{equation}
\end{itemize}
Herein,	
$\xv$ is evaluated for $\xv^{\OTCsubscript}$ or $\xv^{\SQPsubscript}$. The true optimal solution $\chiv^{\optimal}$ is obtained through \emph{offline \SQP} which has access to unlimited computational resources. Numerical results of these metrics are illustrated in the Main Document, Figure~\ref{MD-fig:performance metrics otc vs sqp} (first--third rows).

%\newpage
%\input{supplementary/numerical_implementation}

%\newpage
%\input{supplementary/extended_numerical_results}

%\newpage
%\input{supplementary/delay_caused_tracking_error}

\newpage
% !TEX root = ../supplementary.tex

\section{Outlook I: State-constrained systems}\label{sec:constrained system}
In absence of cost function $\sigma$ the system and problem class~\eqref{eq:gryphon}--\eqref{eq:optimal control problem} can be considered as a state-constrained system
\begin{equation}\label{eq:state-constrained system}
	\dot{\xv} = \fv_{\scriptscriptstyle\textup{A}}(\xv) + \Bm(\xv) \uv, 
	\quad
	\xv(t) \in \mathcal{X}(t) := \big\{\,
	\varxv \in \realset^{\dimx}
	~\big\lvert~
	\hv(t,\varxv) = \nullv
	~\wedge~
	\gv(t,\varxv) \leqslant \nullv
	\,\big\}. %~\forall t \geqslant 0
\end{equation}
Even though $\sigma$ is not present, the augmented cost function $\rho$~\eqref{eq:augmented cost function} can still be constructed solely from $\gv(t,\xv)$, allowing the proposed \OTC\ to be applied on \eqref{eq:state-constrained system}. This results in the controller
\begin{subequations}\label{eq:gryphonheart for constrained system}
	\begin{equation}
		\uv = \omegav(t,\xv,\vv) 
		:= \Bm(\xv)^{-1} \Big(
		-\fv_{\scriptscriptstyle\textup{A}}(\xv) 
		+ K_{\textup{x}} \omegav_{\scriptscriptstyle\textup{A}}(t, \xv)
		+ \Omegam_{\scriptscriptstyle\textup{B}}(t, \xv) \vv
		\Big)
	\end{equation}
	with
	\begin{align}
		&\omegav_{\scriptscriptstyle\textup{A}}(t, \xv) :=
		\begin{cases}
			\!\begin{aligned}
				-\Hm(t,\xv)^{\#}_{\Imat} \hv(t,\xv)
			\end{aligned} & \text{if~}\gv(t,\xv)\leqslant\nullv \\[4ex]
			\hspace{-0.5ex}
			\!\begin{aligned}
				&-\Hm(t,\xv)^{\#}_{\breve{\Gm}(t)^\transp \breve{\Gm}(t)} \hv(t,\xv) \\[-0.5ex]
				&\quad-\left(\Imat \hspace{-0.5ex}-\hspace{-0.5ex} \Hm(t,\xv)^{\#}_{\breve{\Gm}(t)^\transp \breve{\Gm}(t)} \Hm(t,\xv)\right) \big(\breve{\Gm}(t)^\transp \breve{\Gm}(t)\big)^{\hspace{-0.5ex}-1} \breve{\Gm}(t)^\transp \breve{\gv}(t,\xv)
			\end{aligned}
			& \text{otherwise},
		\end{cases} \\
		&\Omegam_{\scriptscriptstyle\textup{B}}(t, \xv) :=\begin{cases}
			\Imat \hspace{-0.5ex}-\hspace{-0.5ex} \Hm(t,\xv)^{\#}_{\Imat} \Hm(t,\xv) & \text{if~}\gv(t,\xv)\leqslant\nullv, \\[2ex]
			\Imat - 
			\begin{bmatrix}
				\breve{\Gm}(t) \\[-1ex]
				\Hm(t,\xv)
			\end{bmatrix}^{\#}_{\Imat} 
			\begin{bmatrix}
				\breve{\Gm}(t) \\[-1ex]
				\Hm(t,\xv)
			\end{bmatrix} & \text{otherwise},
		\end{cases}
	\end{align}
\end{subequations}
and a new input $\vv \in \realset^{\dimx}$.
Herein, $\breve{\gv}(t,\xv) = \breve{\Gm}(t) \xv + \breve{\cv}(t)$ is the vector of positive elements in $\gv(t,\xv)$, %i.e., the components of inequality constraints that are active at $t$ and $\xv$,
$\omegav_{\scriptscriptstyle\textup{A}}$ is a variant of $\etav$~\eqref{eq:definition of eta} in absence of cost function $\sigma$ and when $\desired{\yv} = \nullv$, 
and $\Omegam_{\scriptscriptstyle\textup{B}}(t,\xv)$ projects $\vv$ into a local null space of $\hv(t,\cdot)$ and $\breve{\gv}(t,\cdot)$.
Therefore, $\xv(t)$ of \eqref{eq:state-constrained system} being controlled by \eqref{eq:gryphonheart for constrained system} converges exponentially to $\mathcal{X}(t)$ as $t \geqslant 0$ increases, regardless of $\vv$'s value, i.e., the control $\uv$ provided by~\eqref{eq:gryphonheart for constrained system} is admissible.

Inserting \eqref{eq:gryphonheart for constrained system} into \eqref{eq:state-constrained system} reduces the state-constrained control problem to the simpler unconstrained one; that is
\begin{equation}\label{eq:gryphonheart for constrained system, inner-closed-loop}
	\dot{\xv} = K_{\textup{x}} \omegav_{\scriptscriptstyle\textup{A}}(t, \xv) + \Omegam_{\scriptscriptstyle\textup{B}}(t, \xv)\vv.
\end{equation}
One can use standard control methods~\cite{Khalil2002} such as feedback linearization, sliding mode control, or, when applicable, \uls{l}inear-\uls{q}uadratic \uls{r}egulator (LQR)~\cite{Mehrmann1991} to find $\vv$. 
This means, \OTC\ provides a way for generalizing a wide spectrum of control methods to accommodate state constraints.
%In SI, Section~\ref{SI-sec:constrained system}, we present a state-constrained LQR with both mathematical derivation and numerical example.
As an example, we demonstrate the generalization of LQR to a \emph{state-constrained LQR} later in this section.
We leave generalization of other control methods for future work.

Moreover, \eqref{eq:state-constrained system} is commonly accompanied with an output function
\begin{equation*}
	\zv = \fv_{\scriptscriptstyle\textup{C}}(t,\xv)
\end{equation*}
and it has been a challenge to find an optimal trajectory for constrained systems such that the output $\zv$ eventually reaches a goal \smash{$\zv_{\textup{g}}$}.
{Suppose there exists optimal control \smash{$\vv = \upsilonv^{\optimal}(t,\xv)$} (e.g., LQR) for \eqref{eq:gryphonheart for constrained system, inner-closed-loop}, and consequently,}
\begin{equation}\label{eq:state-constrained optimal controller}
	\uv = \omegav(t,\xv,\vv) 
	:= \Bm(\xv)^{-1} \Big(
	-\fv_{\scriptscriptstyle\textup{A}}(\xv) 
	+ K_{\textup{x}} \omegav_{\scriptscriptstyle\textup{A}}(t, \xv)
	+ \Omegam_{\scriptscriptstyle\textup{B}}(t, \xv) \upsilonv^{\optimal}(t,\xv)
	\Big)
\end{equation}
for \eqref{eq:state-constrained system}.
If $\zv_{\textup{g}}$ is somehow on the path that is generated by \eqref{eq:state-constrained optimal controller},
i.e.,
\begin{equation*}
	\left\{\tau \geqslant 0~\Bigg\lvert~
	\begin{matrix}
		\zv_{\textup{g}} = \fv_{\scriptscriptstyle\textup{C}}\big(\tau,\xv(\tau)\big) \\
		%\wedge~
		\dot{\xv}(t) = K_{\textup{x}} \omegav_{\scriptscriptstyle\textup{A}}\big(t, \xv(t)\big) + \Omegam_{\scriptscriptstyle\textup{B}}\big(t, \xv(t)\big)\upsilonv^{\optimal}\big(t,\xv(t)\big) \\
		%\wedge~
		\xv(0) \in \mathcal{X}(0)
	\end{matrix}
	\right\}
	\neq \emptyset,
\end{equation*}
then \eqref{eq:state-constrained optimal controller} is the optimal trajectory planner.
This insight points to a promising direction in solving optimal trajectory planning problems using closed-form feedback controllers.
For instance, the controller~\eqref{eq:state-constrained optimal controller} can be applied to each agent in a multi-agent system setting, enabling the entire swarm to chart a whole state space by labeling each point with the reachability of $\zv_{\textup{g}}$ and the associated cost starting from that point.
The resulting \emph{map} of the labeled state space can serve as valuable prior knowledge for the higher-level optimal controller targeting $\zv_{\textup{g}}$.
A numerical example illustrating this idea of multi-agent cartography is shown later in this section. %in SI, Section~\ref{SI-sec:constrained system}.

\subsection{State-constrained linear-quadratic regulator}%\label{sec:constrained system}
Let us consider the following state-constrained \uls{l}inear-\uls{q}uadratic (LQ) problem
\begin{subequations}\label{eq:state-constrained LQ problem}
\begin{equation}\label{eq:state-constrained LQ problem, cost}
	\min_{\xv,\uv} \int_{0}^{\infty} \xv(t)^\transp \Qm_{\textup{xx}} \xv(t) + \uv(t)^\transp \Qm_{\textup{uu}} \uv(t)\,\mathrm{d}t
\end{equation}
subject to
\begin{align}
	\dot{\xv} &= \Am \xv + \Bm \uv, \\
	\nullv &= \Hm \xv. %\\
	%\nullv \geqslant \nullmat \xv + \nullv
\end{align}
\end{subequations}
In order to approach~\eqref{eq:state-constrained LQ problem}, 
we apply the admissible transformation
\begin{equation}\label{eq:admissible transform, SCLQ}
	\uv = \omegav(\xv,\vv) 
	= \Bm^{-1} \Big(
	-\Am \xv  
	- K_{\textup{x}} \Hm^{\#}_{\Imat} \Hm \xv
	+ \big(\Imat \hspace{-0.5ex}-\hspace{-0.5ex} \Hm^{\#}_{\Imat} \Hm\big) \vv
	\Big),
\end{equation}
which is a special case of \eqref{eq:gryphonheart for constrained system}
%\added{the controller in Main Document Section~\textit{\nameref{MD-sec:constrained system}}}
with $\gv(t,\xv) = \nullmat \xv - 1 \leqslant \nullv$~$\forall \xv$.
The transformation~\eqref{eq:admissible transform, SCLQ} reduces \eqref{eq:state-constrained LQ problem} to the unconstrained LQ problem
\begin{subequations}
%\begin{equation}
%	\xv^\transp \Qm_{\textup{xx}} \xv + \uv^\transp \Qm_{\textup{uu}} \uv
%	= \xv^\transp \big(\Qm_{\textup{xx}} + \Omegam_{\textup{A}}^{\hspace{0.5ex}\transp} \Qm_{\textup{uu}} \Omegam_{\textup{A}}\big) \xv
%	+ \vv^\transp \Omegam_{\scriptscriptstyle\textup{B}}^{\hspace{0.5ex}\transp} \Qm_{\textup{uu}} \Omegam_{\scriptscriptstyle\textup{B}} \vv
%	+ 2 \xv^\transp \Omegam_{\textup{A}}^{\hspace{0.5ex}\transp} \Qm_{\textup{uu}} \Omegam_{\scriptscriptstyle\textup{B}} \vv
%\end{equation}
\begin{equation}\label{eq:state-constrained LQ problem, transformed cost}
%\scalebox{0.95}{$
	\min_{\xv,\vv} \int_{0}^{\infty}\hspace{-0.5ex}
	\xv(t)^\transp \big(\Qm_{\textup{xx}} \hspace{-0.5ex}+\hspace{-0.5ex} \Omegam_{\textup{A}}^{\hspace{0.5ex}\transp} \Qm_{\textup{uu}} \Omegam_{\textup{A}}\big) \xv(t)
	+ \vv(t)^\transp \Omegam_{\scriptscriptstyle\textup{B}}^{\hspace{0.5ex}\transp} \Qm_{\textup{uu}} \Omegam_{\scriptscriptstyle\textup{B}} \vv(t)
	+ 2 \xv(t)^\transp \Omegam_{\textup{A}}^{\hspace{0.5ex}\transp} \Qm_{\textup{uu}} \Omegam_{\scriptscriptstyle\textup{B}} \vv(t)
	\,\mathrm{d}t
%$}
\end{equation}
subject to
\begin{equation}
	\dot{\xv} = \Omegam_{\textup{A}} \xv + \Omegam_{\scriptscriptstyle\textup{B}} \vv,
	\quad
	\Omegam_{\textup{A}} := -K_{\textup{x}} \Hm^{\#}_{\Imat} \Hm, 
	\quad
	\Omegam_{\scriptscriptstyle\textup{B}} := \Imat - \Hm^{\#}_{\Imat} \Hm,
\end{equation}
\end{subequations}
which can be solved by the LQR~\cite{Mehrmann1991}.
We denote the solution as
\begin{equation}\label{eq:LQR solution, v}
	\vv = \upsilonv_{\scriptscriptstyle\text{LQR}}(\xv).
\end{equation}
Inserting \eqref{eq:LQR solution, v} into \eqref{eq:admissible transform, SCLQ} yields the \emph{\uls{s}tate-\uls{c}onstrained \uls{l}inear-\uls{q}uadratic \uls{r}egulator} (SCLQR)
\begin{equation}\label{eq:state-constrained LQR}
	\uv 
	= \Bm^{-1} \Big(
		-\Am \xv  
		- K_{\textup{x}} \Hm^{\#}_{\Imat} \Hm \xv
		+ \big(\Imat \hspace{-0.5ex}-\hspace{-0.5ex} \Hm^{\#}_{\Imat} \Hm\big) \upsilonv_{\scriptscriptstyle\text{LQR}}(\xv)
	\Big).
\end{equation}
A numerical example demonstrating this SCLQR's effectiveness will be shown later in this section.

It must be noted that incorporating linear inequality constraint in \eqref{eq:state-constrained LQ problem} results in an unconstrained optimal control problem with a \emph{non-quadratic} cost functional and a \emph{nonlinear} system.
Addressing this class of problem requires advanced optimal control methods, which are beyond the scope of this work.

\subsection{Numerical example and result}
\begin{figure}[!h]
	\centering
	%	\includegraphics[width=0.99\textwidth]{figures/consys_x_3d2d.pdf}
	%	\includegraphics[width=0.99\textwidth]{figures/consys_x_timeseries.pdf}	
%	\noindent\makebox[\textwidth]{\includegraphics[width=18cm]{figures/consys_x_3d2d.pdf}}
%	\noindent\makebox[\textwidth]{\includegraphics[width=18cm]{figures/consys_x_timeseries.pdf}}%
	\noindent\makebox[\textwidth]{\includegraphics[width=17.25cm]{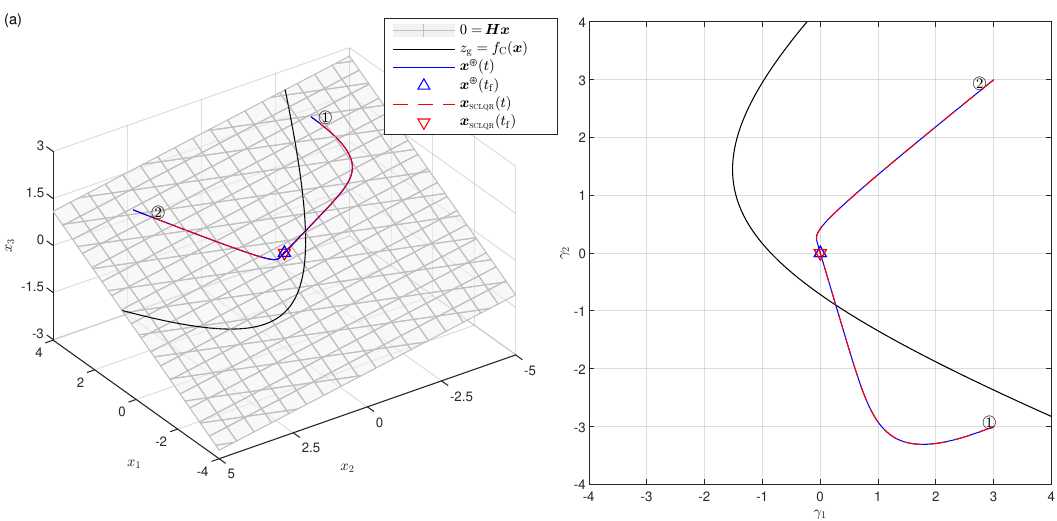}}
	\noindent\makebox[\textwidth]{\includegraphics[width=17.25cm]{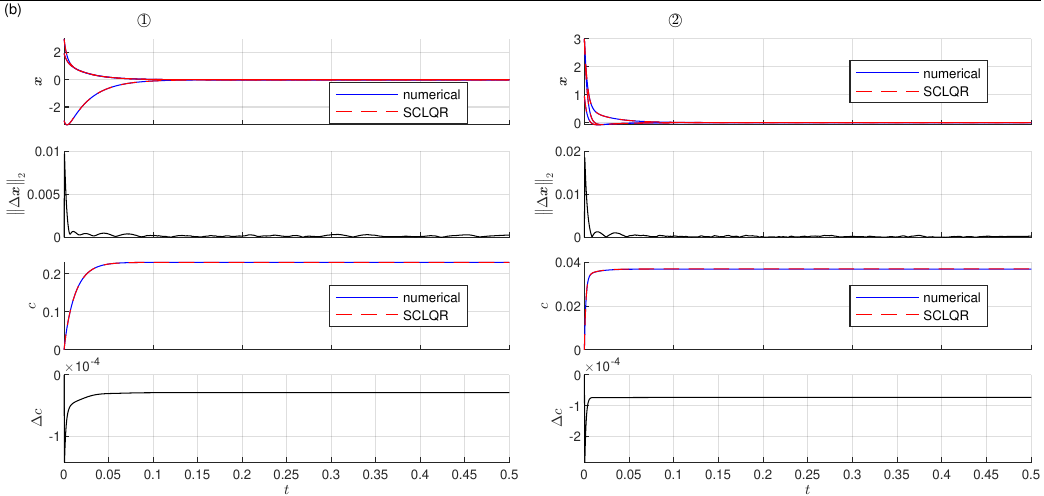}}%
	%	\thisfloatpagestyle{empty}
	\caption{\label{fig:example constrained system result}
		Exemplary three-dimensional state-constrained linear-quadratic problem~\eqref{eq:state-constrained LQ problem} solved numerically and by SCLQR~\eqref{eq:state-constrained LQR}.
		The difference between the resulting state trajectories, colored with blue and red, respectively, is negligible, being tested with two initial conditions $\xv(0)$.
		(a) State trajectories in the 3-D coordinate system (left) and in the 2-D parametric coordinate system (right).
		(b) Time evolution of relevant quantities. Refer to texts for detailed description. %The $t$-axis is on the base-10 logarithmic scale.
	}
\end{figure}

In order to illustrate the ideas mentioned above, we construct an example for \eqref{eq:state-constrained LQ problem}:
\begin{align*}
	&\xv = \begin{bmatrix}
		x_{1} \\ x_{2} \\ x_{3}
	\end{bmatrix}, 
	\quad
	\uv = \begin{bmatrix}
		u_{1} \\ u_{2} \\ u_{3}
	\end{bmatrix}, 
	\quad \\
	&\Qm_{\textup{xx}} = \begin{bmatrix}
		1.04 & -0.01695 & 0.2303\\ -0.01695 & 0.7284 & 0.2473\\ 0.2303 & 0.2473 & 1.898
	\end{bmatrix},
	\quad
	\Qm_{\textup{uu}} = \begin{bmatrix}
		7.331 & 0.1877 & 2.067\\ 0.1877 & 2.328 & -0.6628\\ 2.067 & -0.6628 & 3.956
	\end{bmatrix} \times 10^{-5},
	\quad \\
	&\Am = \begin{bmatrix}
		0.1086 & -0.2032 & -0.02073\\ 0.1763 & 0.6136 & 0.5626\\ -0.5076 & -0.3963 & -0.1000
	\end{bmatrix},
	\quad
	\Bm= \begin{bmatrix}
		0.6665 & 0.9614 & -0.8088\\ 0.3533 & 0.1329 & -0.875\\ -0.8267 & -0.5846 & -0.9484
	\end{bmatrix},
	\quad \\
	&\Hm = \begin{bmatrix}
		-0.3317 & 0.1136 & 0.6919
	\end{bmatrix},
	\quad \\
	&
	f_{\scriptscriptstyle\textup{C}}(t,\xv) = {\left(x_{1}+\frac{46}{25}\right)}^2-{\left(x_{2}-2\right)}^2-\left(3\,x_{3}-\frac{12}{25}\right)\,\left(x_{3}-\frac{4}{25}\right)+4, 
	\quad
	z_{\textup{g}} = 0.
\end{align*}
The matrix $\Hm$ defines $\mathcal{X}$ to be a two-dimensional plane in the three-dimensional space. This 2-D plane can be parametrically represented
as
\begin{align*}
	\xv = \psiv(\gammav) 
	= \begin{bmatrix}
		\gamma _{1}\\ \gamma _{2}\\ 0.4795\,\gamma _{1}-0.1642\,\gamma _{2}
	\end{bmatrix}
\end{align*}
with parametric coordinates $\gammav = [\,\gamma_{1}~~\gamma_{2}\,]^\transp$.
We use this coordinate transformation to project the constraint into a 2-D space, providing a more informative and readable illustration in Figure~\ref{fig:example constrained system result}(a).

Figure~\ref{fig:example constrained system result} shows the resulting trajectories of state $\xv(t)$ starting from two distinct initial conditions $\xv(0)$ and driven by a control $\uv$ that is \begin{itemize}
	\item determined \uls{num}erically using iterative algorithm for constrained optimization~\cite{TheMathWorks2020}.
	\item provided by the SCLQR~\eqref{eq:state-constrained LQR} with $K_{\textup{x}} = 500$.
\end{itemize}
Merely a slight difference between these two trajectories $\xv_{\textup{num}}$ and $\xv_{\scriptscriptstyle\text{SCLQR}}$ is noticeable, see Figure~\ref{fig:example constrained system result}. In fact, the difference 
\begin{equation*}%\label{eq:difference u_NQR vs u_OED}
	\big\|\Delta\xv(t)\big\|_{2} := \big\|\xv_{\textup{num}}(t) - \xv_{\scriptscriptstyle\text{SCLQR}}(t)\big\|_{2}
\end{equation*}
is on the order of magnitude of $10^{-3}$, see Figure~\ref{fig:example constrained system result}(b).
Also, the values of running cost
\begin{equation*}%\label{eq:difference c_NQR vs c_OED}
	c_{\Box}(t) := \int_{0}^{t}\hspace{-0.5ex}
	\big(\xv_{\Box}(\tau)^\transp \Qm_{\textup{xx}} \xv_{\Box}(\tau) + \uv_{\Box}(\tau)^\transp \Qm_{\textup{uu}} \uv_{\Box}(\tau)\big)
	\,\mathrm{d}\tau
\end{equation*}
are almost the same, with a difference 
\begin{equation*}%\label{eq:difference c_NQR vs c_OED}
	\Delta c(t) := c_{\textup{num}}(t) - c_{\scriptscriptstyle\text{SCLQR}}(t)
\end{equation*}
on the order of magnitude of $10^{-4}$, see Figure~\ref{fig:example constrained system result}(b).
These very small differences $\big\|\Delta\xv(t)\big\|_{2}$ and $\Delta c(t)$ are presumably numerical in nature, probably originating from
\begin{itemize}
	\item discrete-time integral (numerical solver and SCLQR),
	\item the convex but not strongly convex cost functional~\eqref{eq:state-constrained LQ problem, transformed cost} due to rank-deficient $\Omegam_{\scriptscriptstyle\textup{B}}$ (SCLQR),
	\item finite horizon $t \in [0, 0.5]$ (numerical solver) instead of infinite horizon $t \in [0,+\infty)$ (SCLQR), and
	\item premature termination of iteration due to numerical resolution (numerical solver).
\end{itemize}

Note, this example can be viewed as the simplest multi-agent system with two agents only. 
The trajectory of agent~{\small $\raisebox{.5pt}{\textcircled{\raisebox{-.9pt} {1}}}$} passes through the curve that defines all solutions in the state space that correspond to the goal $z_{\textup{g}}$, see Figure~\ref{fig:example constrained system result}(a).

These numerical results showcase the effectiveness of the proposed SCLQR that combines \OED\ and LQR, suggesting \OED's capability in generalizing standard control methods to incorporate state constraints. 
Also, they demonstrate the potential of \OED\ for  simultaneous optimal trajectory planning and tracking control, which will be further investigated in the future.

\newpage
% !TEX root = ../supplementary.tex
%
\section{Outlook II: DCC for nonlinear systems with a chain of integrators}\label{sec:chain of integrators}
\begin{figure}[!h]
	\centering
	\resizebox{\textwidth}{!}{%
		\scalefont{0.8500004250002125}%
		% !TEX root = ../../supplementary.tex
\begin{tikzpicture}[auto,node distance=0cm]	

\node[invisiblenode] at (0,0) (origin) {};

% state fcn
\node[medianfcnblock, right of=origin, xshift=0cm, yshift=0cm] (n-th-component) {
$\fv_{\scriptscriptstyle\textup{A},n}(\xv) + \Bm_{n}(\xv)\uv$
};

% chain of integrators
\node[medianblock, right of=n-th-component, xshift=2.99cm, yshift=0cm] (integrator1) {%
$\int(\cdot)\mathrm{d}t$%
};
\node[medianblock, right of=integrator1, xshift=2.4cm, yshift=0cm] (integrator2) {%
$\int(\cdot)\mathrm{d}t$%
};
\node[invisiblenode, right of=integrator2, xshift=2.1cm, yshift=0cm] (integrator3) {%
$\cdots$%
};
\node[medianblock, right of=integrator3, xshift=1.8cm, yshift=0cm] (integrator4) {%
$\int(\cdot)\mathrm{d}t$%
};

% output fcn
\node[medianfcnblock, right of=integrator4, xshift=2cm, yshift=0cm] (outfcn) {%
$\hv_{\uptheta}(t,\thetav)$%
};

\draw[doublearrow] ($(n-th-component)+(-2.25cm,0cm)$) --node[]{$\uv$} (n-th-component);
\draw[doublearrow] (n-th-component) --node[]{$\thetav^{(n)}$} (integrator1);
\draw[doublearrow] (integrator1) --node[]{$\thetav^{(n-1)}$} (integrator2);
\draw[doublearrow] (integrator2) --node[]{$\thetav^{(n-2)}$} (integrator3);
\draw[doublearrow] (integrator3) --node[]{$\dot{\thetav}$} (integrator4);
\draw[doublearrow] (integrator4) --node[]{$\thetav$} (outfcn);
\draw[doublearrow] (outfcn) --node[]{$\yv$} ($(outfcn.east)+(0.75cm,0cm)$);

\draw[doublearrow] (integrator1.east) -| ($(integrator1.east)+(0.6cm,-1cm)$) -| ($(n-th-component.south)+(1cm,0cm)$);
\draw[doublearrow] (integrator2.east) -| ($(integrator2.east)+(0.5cm,-1.4cm)$) -| ($(n-th-component.south)+(0.5cm,0cm)$);
\draw[doublearrow] (integrator3.east) -| ($(integrator3.east)+(0.3cm,-1.8cm)$) -| ($(n-th-component.south)+(-0.5cm,0cm)$);
\draw[doublearrow] (integrator4.east) -| ($(integrator4.east)+(0.2cm,-2.2cm)$) -| ($(n-th-component.south)+(-1cm,0cm)$);
\node[invisiblenode, below of=n-th-component, xshift=0cm, yshift=-0.9cm] () {%
$\cdots$%
};

\draw[innerwhite] ($(integrator1.east)+(0.2cm,0cm)$) -- ($(integrator2.west)+(-0.2cm,0cm)$);
\draw[innerwhite] ($(integrator2.east)+(0.2cm,0cm)$) -- ($(integrator3.west)+(-0.2cm,0cm)$);
\draw[innerwhite] ($(integrator3.east)+(0.2cm,0cm)$) -- ($(integrator4.west)+(-0.2cm,0cm)$);
\draw[innerwhite] ($(integrator4.east)+(0.2cm,0cm)$) -- ($(outfcn.west)+(-0.2cm,0cm)$);

\end{tikzpicture}%
		\scalefont{1.17647}%
	}%
	\caption{\label{fig:nlsys chain integrators}
		Block diagram of a nonlinear system with a chain of $n$ integrators.
	}
\end{figure}
\vspace{2ex}
In this section, we consider the following nonlinear system with a chain of $n$ integrators:
\begin{subequations}\label{eq:chain of integrators}
\vspace{-2ex}
	\begin{align}
		\overbracket[0.14ex]{%
			\begin{bmatrix}
				\dot{\thetav} \\
				\ddot{\thetav} \\
				\vdots \\
				\thetav^{(n-1)} \\ \hdashline[1pt/1pt]
				\thetav^{(n)}
			\end{bmatrix}%
		}^{\dot{\xv}}
		&= \overbracket[0.14ex]{
			\begin{bmatrix}
				\begin{bmatrix}
					\nullmat & \Imat & \cdots & \nullmat & \nullmat \\
					\nullmat & \nullmat & \ddots & \nullmat & \nullmat \\
					\vdots & \ddots & \ddots & \ddots & \vdots \\
					\nullmat & \nullmat & \cdots & \nullmat & \Imat
				\end{bmatrix}  
				\begin{bmatrix}
					\thetav \\
					\dot{\thetav} \\
					\vdots \\
					\thetav^{(n-2)} \\
					\thetav^{(n-1)}
				\end{bmatrix} \\ \hdashline[1pt/1pt]
				\fv_{\scriptscriptstyle\textup{A},n}\left(\big[\,
					\thetav^{\transp} ~~
					\dot{\thetav}^{\transp} ~~
					\cdots ~~
					\thetav^{(n-2)\transp} ~~
					\thetav^{(n-1)\transp}
				\,\big]^{\transp}\right)
			\end{bmatrix}
		}^{\fv_{\scriptscriptstyle\text{A}}(\xv)}
		+
		\overbracket[0.14ex]{
			\begin{bmatrix}
				\nullmat \\
				\nullmat \\
				\vdots \\
				\nullmat \\ \hdashline[1pt/1pt]
				\Bm_{n}(\xv)
			\end{bmatrix}
		}^{\Bm(\xv)} \uv, \\
		\yv &= \hv(t,\xv) = \hv_{\uptheta}(t,\thetav),
	\end{align}
\end{subequations}
where
$\xv = \big[\,
\thetav^{\transp} ~~
\dot{\thetav}^{\transp} ~~
\cdots ~~
\thetav^{(n-2)\transp} ~~
\thetav^{(n-1)\transp}
\,\big]^{\transp}$
with $\thetav^{(n)} := {\mathrm{d}^{n}\thetav(t)}/{\mathrm{d}t^{n}}$
is the state,
$\Bm_{n}(\xv)$ is invertible $\forall\xv$ (implying $\dim(\uv) = \dim(\thetav)$).
System~\eqref{eq:chain of integrators} is related to a nonlinear canonical form~\cite{Sommer1980,Zeitz1989,BusawonDeAch2001}.
Its block diagram is depicted in Figure~\ref{fig:nlsys chain integrators}.

For \eqref{eq:chain of integrators}, we specify the optimality principle as follows:
\begin{equation}\label{eq:constrained optimization theta}
	\begin{split}
		%\varthetav^{\optimal} = \arg 
		\min_{\varthetav} \sigma_{\uptheta}(t,\varthetav)~~
		\text{subject to}~
		\nullv = \hv_{\uptheta}(t,\varthetav) - \desired{\yv}(t)%, \\ %\label{eq:result-section, equality constraint}
		~\wedge~
		\nullv \geqslant \gv_{\uptheta}(t,\varthetav) %\label{eq:result-section, inequality constraint}
	\end{split}
\end{equation}
and $\varthetav^{\optimal}$ denotes one of its local minimizers. 
The optimization problem~\eqref{eq:constrained optimization theta} is a special case of \eqref{eq:constrained optimization} with 
\begin{equation}\label{eq:constrained optimization chi->theta}
	\begin{split}\hspace{-1cm}
%		&\chiv = \big[
%		\begin{matrix}
%			\varthetav{}^\transp &
%			\dot{\varthetav}{}^\transp &
%			\cdots &
%			\varthetav^{(n-2)}{}^\transp &
%			\varthetav^{(n-1)}{}^\transp
%		\end{matrix}\big]^\transp, \\
		%
		&\chiv = \big[\,
		\varthetav^{\transp} ~~
		\dot{\varthetav}^{\transp} ~~
		\cdots ~~
		\varthetav^{(n-2)\transp} ~~
		\varthetav^{(n-1)\transp}
		\,\big]^{\transp}, \\
		&\sigma(t,\chiv) = \sigma_{\uptheta}(t,\varthetav), \\
		&\hv(t,\chiv) = \hv_{\uptheta}(t,\varthetav), \\
		&\gv(t,\chiv) = \gv_{\uptheta}(t,\varthetav)
	\end{split}
\end{equation}
and $%\begin{equation}
\chiv^{\optimal} = \big[\,
	\varthetav^{\optimal\transp} ~~
	\nullv^\transp ~~
	\cdots ~~
	\nullv^\transp ~~
	\nullv^\transp
\,\big]^\transp
$%\end{equation}
is one of its %\eqref{eq:constrained optimization}\eqref{eq:constrained optimization chi->theta}'s 
local minimizers.

Inspired by \OTC~\eqref{eq:gryphonheart} for the system and problem class~\eqref{eq:gryphon}--\eqref{eq:optimal control problem}, we propose the controller
\begin{align}
	\uv
	&= \muv(t,\xv,\desired{\yv}; K_{\textup{x}}) \nonumber\\
	\begin{split}\label{eq:OTC for chain of integrators}
		&= \Bm_{n}(\xv)^{-1} \bigg(
		-\fv_{\scriptscriptstyle\textup{A},n}(\xv)
		+ K_{0} \Hm_{\uptheta}(t,\thetav)^{\#}_{\Rm_{\uptheta}(t,\thetav)} \big(
		\desired{\yv} - \hv_{\uptheta}(t,\thetav)
		\big) \\[-2ex]
		&\qquad\qquad\qquad- K_{0} \Big(
		\Imat - \Hm_{\uptheta}(t,\thetav)^{\#}_{\Rm_{\uptheta}(t,\thetav)} \Hm_{\uptheta}(t,\thetav)
		\Big) \Rm_{\uptheta}(t,\thetav)^{-1} \rv_{\uptheta}(t,\thetav)
		- \sum_{i=1}^{n-1} K_{i} \thetav^{(i)}
		\bigg)
	\end{split}%
\end{align}%
for the system and problem class~\eqref{eq:chain of integrators}--\eqref{eq:constrained optimization chi->theta}.
Herein, $\Hm_{\uptheta}(t,\thetav) := {\partial \hv_{\uptheta}(t,\thetav)}/{\partial \thetav}$ is the Jacobian of $\hv_{\uptheta}(t,\thetav)$;
the functions \smash{$\rv_{\uptheta}(t,\thetav) := [{\partial \rho_{\uptheta}(t,\thetav)}/{\partial \thetav}]^\transp$} and
$\Rm_{\uptheta}(t,\thetav) := {\partial \rv_{\uptheta}(t,\thetav)}/{\partial \thetav}$ are the gradient and Hessian of 
$\rho_{\uptheta}(t,\thetav) := \sigma_{\uptheta}(t,\thetav) + \frac{1}{2}\betav\big(\gv_{\uptheta}(t,\thetav)\big)^\transp \betav\big(\gv_{\uptheta}(t,\thetav)\big)$, respectively;
the controller gains $K_{0}$, $K_{1}$, $\cdots$, $K_{n}$ are functions of $K_{\textup{x}}$ (placement of dominant pole $-K_{\textup{x}}$, see below).
%Remark: Neither $\Bm(\xv)$'s invertibility nor $\dim(\xv) = \dim(\uv)$ are required.

On the ground of deduction and numerical simulation (in the following subsections), we claim that 
%\eqref{eq:OTC for chain of integrators} is an \OTC~(Definition~\ref{definition:gryphonheart}) for \eqref{eq:chain of integrators}--\eqref{eq:constrained optimization chi->theta}, i.e., 
the state $\xv(t)$ of the closed-loop system~\eqref{eq:chain of integrators}\eqref{eq:OTC for chain of integrators} converges, as $t \geqslant 0$ increases, exponentially to the interior of a moving ball centered at $\chiv^{\optimal}(t)$ with a radius $\propto K_{\textup{x}}^{-1}$. We leave the proof for future work.

\subsection{Deduction}
For controlling \eqref{eq:chain of integrators}, we apply a cascaded control: the inner controller
\begin{equation}\label{eq:chain of integrators, inner controller}
	\uv 
	= \Bm_{n}(\xv)^{-1} \left(
	-\fv_{\scriptscriptstyle\textup{A}}(\xv)
	+ K_{0} (\vv - \thetav)
	- \sum_{i=1}^{n-1} K_{i} \thetav^{(i)}
	\right) 
\end{equation}
shapes the linear inner closed-loop system
\begin{equation}\label{eq:chain of integrators, inner closed loop dynamics}
	\begin{bmatrix}
		\dot{\thetav} \\
		\ddot{\thetav} \\
		\vdots \\
		\thetav^{(n-1)} \\
		\thetav^{(n)}
	\end{bmatrix}
	= \begin{bmatrix}
		\nullmat & \Imat & \cdots & \nullmat & \nullmat \\
		\nullmat & \nullmat & \ddots & \nullmat & \nullmat \\
		\vdots & \ddots & \ddots & \ddots & \vdots \\
		\nullmat & \nullmat & \cdots & \nullmat & \Imat \\
		~-K_{0} \Imat~ & ~-K_{1} \Imat~ & \cdots & ~-K_{n-2} \Imat~ & -K_{n-1} \Imat~
	\end{bmatrix} 
	\begin{bmatrix}
		\thetav \\
		\dot{\thetav} \\
		\vdots \\
		\thetav^{(n-2)} \\
		\thetav^{(n-1)}
	\end{bmatrix}
	+ \begin{bmatrix}
		\nullmat \\
		\nullmat \\
		\vdots \\
		\nullmat \\
		K_{0} \Imat
	\end{bmatrix} \vv
\end{equation}
with a new input $\vv$.
%Herein, all identity ($\Imat$) and zero matrices ($\nullmat$) are $\dimtheta$-by-$\dimtheta$.
%The controller gains $K_{0}$, $K_{1}$, $K_{2}$, $\cdots$, $K_{n-1}$ are designed such that \eqref{eq:chain of integrators, inner closed loop dynamics}'s transfer function is
%\begin{equation}
%	\frac{\mathscr{L}\{x_{i}(t)\}(s)}{\mathscr{L}\{v_{i}(t)\}(s)}
%	= \frac{1}{\big(1 + K_{\textup{x}}^{-1} s\big) \big(1 + (q K_{\textup{x}})^{-1} s\big)^{n-1}},
%\end{equation}
%where $\mathscr{L}\{x_{i}(t)\}(s)$ denotes Laplace transformation of $x$ from time domain ($t$) to frequency domain ($s$), and
%the higher the constant $q/n \gg 1$ is, the more dominant the pole $-K_{\textup{x}}$ is among all poles.
%With a sufficiently large $q$, the response behavior of $\xv(t)$ to $\vv(t)$ in \eqref{eq:high-order system, inner closed loop dynamics} can be approximated by
%\begin{equation}\label{eq:chain of integrators, inner closed loop dynamics approx}
%	\dot{\hat{\thetav}} = - K_{\textup{x}} \hat{\thetav} + K_{\textup{x}} \vv,
%\end{equation}
%on which, \OTC~\eqref{eq:gryphonheart} is applicable; that is
The controller gains $K_{0}$, $K_{1}$, $K_{2}$, $\cdots$, $K_{n-1}$ are designed such that $-K_{\textup{x}}$ is the dominant pole of \eqref{eq:chain of integrators, inner closed loop dynamics}'s transfer function, i.e., the response behavior of $\xv(t)$ to $\vv(t)$ can be approximated by
\begin{equation}\label{eq:chain of integrators, inner closed loop dynamics approx}
	\dot{\hat{\thetav}} = - K_{\textup{x}} \hat{\thetav} + K_{\textup{x}} \vv.
\end{equation}
Given this first-order approximation of the linear inner loop~\eqref{eq:chain of integrators, inner closed loop dynamics}, we adopt \OTC~\eqref{eq:gryphonheart} as the outer controller; that is
\begin{equation}\label{eq:chain of integrators, outer controller}
	\vv 
	= \thetav 
	+ \Hm_{\uptheta}(t,\thetav)^{\#}_{\Rm_{\uptheta}(t,\thetav)} \big(
	\desired{\yv} - \hv_{\uptheta}(t,\thetav)
	\big) 
	-\big(
		\Imat - \Hm_{\uptheta}(t,\thetav)^{\#}_{\Rm_{\uptheta}(t,\thetav)} \Hm_{\uptheta}(t,\thetav)
	\big) \Rm_{\uptheta}(t,\thetav)^{-1} \rv_{\uptheta}(t,\thetav).
\end{equation}
Inserting \eqref{eq:chain of integrators, outer controller} into \eqref{eq:chain of integrators, inner controller} yields \eqref{eq:OTC for chain of integrators}.

\subsection{Numerical example and result}
For testing \eqref{eq:OTC for chain of integrators}, we adapt Example M1~\eqref{eq:example m1} as follows:
\begin{subequations}\label{eq:example m2}
	\begin{align}
		&\xv = \begin{bmatrix}
			\thetav^\transp &~ \dot{\thetav}^\transp &~ \cdots &~ \thetav^{(5)\transp} &~ \thetav^{(6)\transp}
		\end{bmatrix}^\transp, 
		\quad
		\thetav = \begin{bmatrix}
			\theta_{1} & \theta_{2}
		\end{bmatrix}^\transp,
		\\
		&\fv_{\scriptscriptstyle\textup{A},n}(\thetav) = \frac{1}{\thetav^\transp \thetav + 0.001},
		\quad
		\Bm_{n}(\thetav) = \big(\thetav^\transp \thetav + 1\big) \Imat,
		\quad \\
		&h_{\uptheta}(t,\thetav) = -6\,t+{\big(5\,\theta_{1}+25\,\theta_{2}^2-7\big)}^2+{\big(25\,\theta_{1}^2+5\,\theta_{2}-11\big)}^2 - 1,
		\quad
		\desired{y} = 0,
		\quad \\
		&\gv_{\uptheta}(t,\thetav) = \begin{bmatrix} 1 & 0\\[1ex] -1 & 0\\[1ex] 0 & 1\\[1ex] 0 & -1\\[1ex] \frac{209}{100\,\sin(\pi t/10)-143} & 1\\[1ex] \frac{\sin(\pi t/10)}{2}+\frac{73}{75} & 1\\[1ex] \frac{77}{8\,\left(5\,\sin(\pi t/10)+7\right)} & -1\\[1ex] \frac{15\,\sin(\pi t/10)}{38}-\frac{67}{152} & -1 \end{bmatrix} \thetav 
		+ \begin{bmatrix} -\frac{19}{20}\\[1ex] -\frac{19}{20}\\[1ex] -\frac{19}{20}\\[1ex] -\frac{19}{20}\\[1ex] -\frac{209\,\left(3\,\sin(\pi t/10)-10\right)}{10\,\left(100\,\sin(\pi t/10)-143\right)}\\[1ex] -\frac{3\,\sin(\pi t/10)}{10}-1\\[1ex] -\frac{77\,\left(\sin(\pi t/10)+5\right)}{40\,\left(5\,\sin(\pi t/10)+7\right)}\\[1ex] \frac{3\,\sin(\pi t/10)}{10}-1 \end{bmatrix}, 
		\quad \\
		&\sigma_{\uptheta}(t,\xv) = \sqrt{
			\big(\thetav - \pv(t)\big)^\transp \Pm(t) \big(\thetav - \pv(t)\big) + 0.001
		} + 0.001\, \big(\thetav - \pv(t)\big)^\transp \big(\thetav - \pv(t)\big),
		\quad \\
		&\quad\Pm(t) = \begin{bmatrix} \frac{9 t + 9 t \cos\big(4 \pi^{2}\,\sin(\pi t/100)\big)+40}{18 t + 40} & \frac{9 t \sin\big(4 \pi^{2}\,\sin(\pi t/100)\big)}{18 t + 40}\\[1ex] \frac{9 t \sin\big(4 \pi^{2}\,\sin(\pi t/100)\big)}{18 t + 40} & \frac{9 t {\sin\big(2\,\pi ^2\,\sin(\pi t/100)\big)}^2+20}{9 t + 20} \end{bmatrix},
		\quad \\%[1ex]
		&\quad\pv(t) = \begin{bmatrix} -\frac{2\,\sin(\pi t/5)}{3}-\frac{10\,\sin(3 \pi t/10)}{27}\\[1ex] \frac{2\,\cos(\pi t/5)}{3}-\frac{10\,\cos(3 \pi t/10)}{27} \end{bmatrix}.
	\end{align}
\end{subequations}
The resulting $\thetav(t)$ is illustrated in Figure~\ref{fig:chain of integrators example m1 result} and the expected exponential convergence and tracking behavior can be seen.

\afterpage{%
	\newgeometry{left=1.5cm,right=1.5cm,top=2cm,bottom=6cm}
\begin{figure}[!t]
	\centering
	\includegraphics[width=15cm]{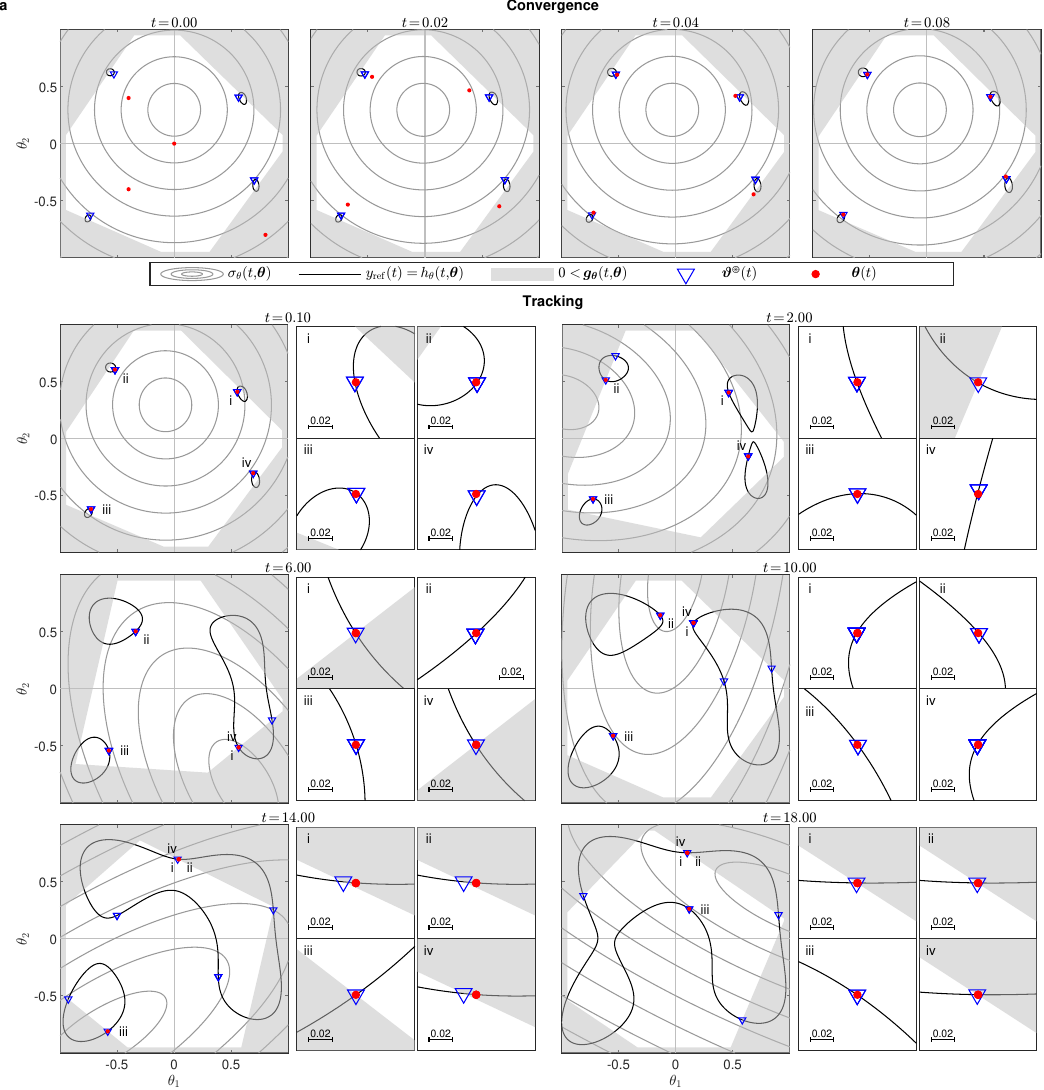}
	\includegraphics[width=15cm]{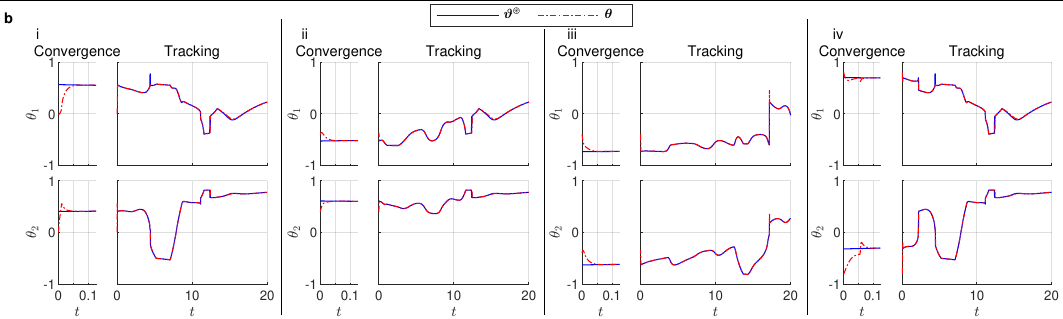}
	%\thisfloatpagestyle{empty}
	%	\vspace{-2ex}
	\caption{\label{fig:chain of integrators example m1 result}
		Numerical results of \OTC~($K_{\textup{x}} = 100$) controlling the exemplary system~\eqref{eq:chain of integrators}\eqref{eq:example m2} with $n=7$.
		The trajectories $\thetav(t)$ (marked as i, ii, iii, iv) starting from different initial conditions $\thetav(0)$ follow the corresponding $\varthetav^{\optimal}(t)$ over time $t$. 
		(a) Trajectories of $\thetav(t)$ (red) and $\varthetav^{\optimal}(t)$ (blue) in the two-dimensional space along with cost function and constraints. Each subplot corresponds to a snapshot at time instant $t$ as indicated.
		(b) Exponentially convergence and tracking of $\thetav(t)$ to the closest $\varthetav^{\optimal}(t)$. Each column represents the time evolution starting from a different initial condition.
		A zoomed-in perspective to the initial transient phase is provided on the left in each column.
	}
\end{figure}
%\restoregeometry
}

\subsection{Practical relevance}
One exemplary physical system that can be written as~\eqref{eq:chain of integrators} is a synchronous generator~\cite{BusawonDeAch2001}.
An other relevant example is a multibody dynamics with serial kinematics, i.e., a robot:
\begin{subequations}\label{eq:serial robot}
	\begin{align}
		\ddot{\thetav} &= -\Mm(\thetav)^{-1} \nv(\thetav,\dot{\thetav}) + \Mm(\thetav)^{-1} \uv, \\
		\yv &= \hv_{\scriptscriptstyle\textup{EE}}(\thetav),
	\end{align}
\end{subequations}
where $\thetav$ is the joint-space generalized coordinates, $\uv$ is the joint-space generalized force (e.g., torques), 
$\Mm(\thetav)$ accounts for inertia, $\nv(\thetav,\dot{\thetav})$ accounts for Coriolis force and gravity, and
$\yv$ is the \uls{e}nd-\uls{e}ffector pose.
System~\eqref{eq:serial robot} is explicitly in the form of \eqref{eq:chain of integrators} with $n=2$.
%%
%%	\begin{equation}
%	%		\uv
%	%		= \Mm(\thetav) \bigg(
%	%		\Mm(\thetav)^{-1} \nv(\thetav,\dot{\thetav})
%	%		+ K_{0} \Hm_{\scriptscriptstyle\textup{EE}}(t,\thetav)^{\#}_{\Rm_{\uptheta}(t,\thetav)} \big(
%	%		\desired{\yv} - \hv_{\scriptscriptstyle\textup{EE}}(t,\thetav)
%	%		\big)
%	%		- K_{0} \Big(
%	%		\Imat - \Hm_{\scriptscriptstyle\textup{EE}}(t,\thetav)^{\#}_{\Rm_{\uptheta}(t,\thetav)} \Hm_{\scriptscriptstyle\textup{EE}}(t,\thetav)
%	%		\Big) \Rm_{\uptheta}(t,\thetav)^{-1} \rv_{\uptheta}(t,\thetav)
%	%		%			
%	%		- \sum_{i=1}^{n-1} K_{i} \thetav^{(i)}
%	%		\bigg),
%	%	\end{equation}
%\begin{equation}
%	\begin{split}
%		\uv
%		&=
%		\nv(\thetav,\dot{\thetav})
%		+ \Mm(\thetav) K_{0} \Hm_{\scriptscriptstyle\textup{EE}}(t,\thetav)^{\#}_{\Rm_{\uptheta}(t,\thetav)} \big(
%		\desired{\yv} - \hv_{\scriptscriptstyle\textup{EE}}(t,\thetav)
%		\big) \\[-1ex]
%		&\qquad- \Mm(\thetav) K_{0} \Big(
%		\Imat - \Hm_{\scriptscriptstyle\textup{EE}}(t,\thetav)^{\#}_{\Rm_{\uptheta}(t,\thetav)} \Hm_{\scriptscriptstyle\textup{EE}}(t,\thetav)
%		\Big) \Rm_{\uptheta}(t,\thetav)^{-1} \rv_{\uptheta}(t,\thetav)
%		%			
%		- \Mm(\thetav) K_{1} \dot{\thetav},
%	\end{split}
%\end{equation}
%related to impedance control with null-space handling

\newpage
\section{Background material}
% !TEX root = ../supplementary.tex
%
\subsection{Newton-Raphson method}\label{sec:newton-raphson method}

The Newton-Raphson method is a root-finding algorithm which produces successively better approximations to a root of a system of nonlinear equations
\begin{equation}\label{eq:newton-raphson method, equation to be solved}
	\nullv = \rv(\chiv)
\end{equation}
with continuously differentiable $\rv: \realset^{m} \to \realset^{n}$.
Given a suitable initial guess $\chiv^{[0]}$, the iterative formula
\begin{equation}\label{eq:newton-raphson method}
	\chiv^{[k+1]} = \chiv^{[k]} - \Rm(\chiv^{[k]})^{\#} \rv(\chiv^{[k]})
\end{equation}
improves the guess at each iteration~\cite{Kelley1995}.
Herein, $\chiv^{[k]}$ denotes the estimate of a root at the $k$-th iteration, 
$\Rm(\chiv) = {\partial \rv(\chiv)}/{\partial \chiv}$ the Jacobian of $\rv$,
and $\Am^{\#}$ a suitable pseudoinverse of the matrix $\Am$.
Should $m = n$ and $\Rm(\chiv)$ be invertible for all $\chiv \in \realset^{m}$, the pseudoinverse $\Rm(\chiv^{[k]})^{\#}$ becomes the inverse $\Rm(\chiv^{[k]})^{-1}$.
The update
\begin{equation*}
	\chiv^{[k+1]} - \chiv^{[k]} = - \Rm(\chiv^{[k]})^{\#} \rv(\chiv^{[k]})
\end{equation*}
is known as a \textit{Newton-Raphson step}.
The convergence behavior of \eqref{eq:newton-raphson method} is described in Lemma~\ref{lemma:chord method}.

\paragraph{Application in optimization}
Optimization aims at finding critical points--which are local minima, local maxima, or saddle points (minimax points)--of a twice continuously differentiable scalar function $\rho: \realset^{m} \to \realset$~\cite{NocedalWri2006}. This is equivalent to finding the zeros of $\rho$'s gradient $\rv(\chiv) := \big[{\partial \rho(\chiv)}/{\partial \chiv}\big]^\transp$, i.e., finding the roots of \eqref{eq:newton-raphson method, equation to be solved}.
% !TEX root = ../supplementary.tex
%
\subsection{Stability of slowly varying systems}\label{sec:slowly varying systems}

\begin{theorem}[Stability theorem for slowly varying systems~\cite{Khalil2002}]\label{theorem:slowly varying system}
	Consider the system 
	\begin{equation}\label{eq:slowly varying system}
		\dot{\thetav} = \fv(\thetav,\vv)
	\end{equation}
	with state $\thetav \in \realset^{\dimtheta}$,
	input $\vv \in \mathcal{V} \subset \realset^{\dimv}$,
	and locally Lipschitz continuous $\fv: \realset^{\dimtheta}\times\mathcal{V} \to \realset^{\dimtheta}$.
	Suppose
	\begin{equation}\label{eq:slowly varying system, equilibrium}
		\psiv \in \big\{\psiv: \realset^{\dimv} \to \realset^{\dimtheta}~\big\lvert~\nullv = \fv\big(\psiv(\vv),\vv\big)~\forall \vv \in \mathcal{V}\big\}
	\end{equation}
	exists uniquely and is differentiable.
	Let $\Psim(\vv) := {\partial \psiv(\vv)}/{\partial \vv}$.
	Suppose $\|\Psim(\vv)\|_{2} \leqslant D_{\scriptscriptstyle\Psi}$ for all $\vv \in \realset^{\dimv}$.
	Suppose $\vv(t)$ is continuously differentiable and $\|\dot{\vv}(t)\|_{2} \leqslant D_{\textup{v}}$ for all $t \geqslant 0$.
	Suppose there exists for
	\begin{subequations}\label{eq:slowly varying system, error dynamics}
	\begin{align}
		\widetilde{\thetav}(t) &:= \thetav(t) - \psiv\big(\vv(t)\big), \\
		\dot{\widetilde{\thetav}} = \widetilde{\fv}(\widetilde{\thetav}, \vv) &:= \fv\big(\widetilde{\thetav} + \psiv(\vv), \vv\big)
	\end{align}
	\end{subequations}
	a Lyapunov function $U(\widetilde{\thetav},\vv)$ which satisfies
	\begin{subequations}\label{eq:slowly varying system, Lyapunov inequalities}
	\begin{align}
		b_{1} \big\|\widetilde{\thetav}\big\|_{2}^{2} \leqslant U(\widetilde{\thetav},\vv) \leqslant b_{2} \big\|\widetilde{\thetav}\big\|_{2}^{2}, \label{eq:slowly varying system, Lyapunov inequality 1}\\
		%
%		\frac{\partial U(\widetilde{\thetav},\vv)}{\partial \widetilde{\thetav}} \fv\big(\widetilde{\thetav} + \psiv(\vv), \vv\big) \leqslant -b_{3} \big\|\widetilde{\thetav}\big\|_{2}^{2}, \\
		\frac{\partial U(\widetilde{\thetav},\vv)}{\partial \widetilde{\thetav}} \widetilde{\fv}(\widetilde{\thetav}, \vv) \leqslant -b_{3} \big\|\widetilde{\thetav}\big\|_{2}^{2}, \label{eq:slowly varying system, Lyapunov inequality 2}\\
		\left\| \frac{\partial U(\widetilde{\thetav},\vv)}{\partial \widetilde{\thetav}} \right\|_{2} \leqslant b_{4} \big\|\widetilde{\thetav}\big\|_{2}, \label{eq:slowly varying system, Lyapunov inequality 3}\\
		\left\| \frac{\partial U(\widetilde{\thetav},\vv)}{\partial \vv} \right\|_{2} \leqslant b_{5} \big\|\widetilde{\thetav}\big\|_{2}^{2} \label{eq:slowly varying system, Lyapunov inequality 4}
	\end{align}
	\end{subequations}
	for all $\widetilde{\thetav} \in \mathcal{N}(\nullv,\delta)$ and $\vv \in \mathcal{V}$, where $b_{1} > 0$, $b_{2} > 0$, $b_{3} > 0$, $b_{4} > 0$, $b_{5} \geqslant 0$ are constants independent of $\vv$.
	Note that \eqref{eq:slowly varying system, error dynamics} is equivalent to \eqref{eq:slowly varying system}.
	
	If $\delta$ is such that 
	\begin{equation}\label{eq:slowly varying systems, inequality Dv}
		D_{\textup{v}} 
		< 
%		\frac{b_{1} b_{3}}{b_{2} b_{5}} \frac{\delta}{\delta + \frac{b_{4} D_{\scriptscriptstyle\Psi}}{b_{5}}}
%		= 
		\frac{b_{1} b_{3}}{b_{2}} \frac{\delta}{b_{5}\delta + b_{4} D_{\scriptscriptstyle\Psi}},
	\end{equation}
	then with initial condition %$\thetav(0) \in \mathcal{N}\big($\scalebox{0.89}{$\psiv\big(\vv(0)\big), \delta\,\sqrt{b_{1}/b_{2}}\,$}$\big)$, 
	\begin{equation}\label{eq:slowly varying systems, initial condition}
		\thetav(0) \in \mathcal{N}\hspace{-0.5ex}\left(\psiv\big(\vv(0)\big), \delta\,\sqrt{b_{1}/b_{2}}\,\right),
	\end{equation}
	the state $\thetav(t)$ of \eqref{eq:slowly varying system} satisfies
	\begin{subequations}
	\begin{align}
		\big\|\widetilde{\thetav}(t)\big\|_{2} 
		&\leqslant \sqrt{\frac{b_{2}}{b_{1}}}\,\big\|\widetilde{\thetav}(0)\big\|_{2} \exp(-\kappa t) + \frac{b_{4} D_{\scriptscriptstyle\Psi}}{2 b_{1}} \int_{0}^{t} \exp\big(-\kappa (t-\tau)\big) \big\|\dot{\vv}(\tau)\big\|_{2}\,\mathrm{d}\tau \\
		&\leqslant \sqrt{\frac{b_{2}}{b_{1}}}\,\big\|\widetilde{\thetav}(0)\big\|_{2} \exp(-\kappa t) + \frac{b_{4} D_{\scriptscriptstyle\Psi}}{2 b_{1}} \int_{0}^{t} \exp\big(-\kappa (t-\tau)\big) D_{\textup{v}}\,\mathrm{d}\tau \\
		&\leqslant \sqrt{\frac{b_{2}}{b_{1}}}\,\big\|\widetilde{\thetav}(0)\big\|_{2} \exp(-\kappa t) + \frac{b_{4} D_{\scriptscriptstyle\Psi} D_{\textup{v}}}{2 b_{1} \kappa},
	\end{align}
	\end{subequations}
	where %$\displaystyle\kappa = \frac{1}{2}\left(\frac{b_{3}}{b_{1}} - \frac{b_{5}}{b_{1}} D_{\textup{v}}\right)$.
	\begin{equation}\label{eq:slowly varying systems, kappa}
		\kappa = \frac{b_{3} - b_{5} D_{\textup{v}}}{2 b_{1}}.
	\end{equation}
\end{theorem}
\begin{proof}
	See Theorem 9.3 in \cite{Khalil2002}.
\end{proof}

\begin{remark}\textup{\cite{Khalil2002}}
	Inequalities~\eqref{eq:slowly varying system, Lyapunov inequality 1} and \eqref{eq:slowly varying system, Lyapunov inequality 2} state the usual requirements that $U$ be positive definite and has a negative definite derivative along the trajectories of system~\eqref{eq:slowly varying system, error dynamics}.
	On top of that, they show that $\widetilde{\thetav}^{\ast} = \nullv$ is exponentially stable.
\end{remark}

\begin{remark}\label{remark:application of theorm for slowly varying systems}%[Application in Lemma~\ref{lemma:x to x_ast, v3}]
	The following facts must be noted when applying Theorem~\ref{theorem:slowly varying system} to prove Lemma~\ref{lemma:x to x_ast, v3}:
	\begin{enumerate}
		\item As reasoned in the proof of Lemma~\ref{lemma:quasi Newton linear convergence}, the function $\etav$--defined in \eqref{eq:definition of eta}\eqref{eq:definition of eta_A}\eqref{eq:definition of eta_B}--is locally Lipschitz continuous.
		Note that $\etav$ in Lemma~\ref{lemma:x to x_ast, v3} corresponds to $\fv$ in Theorem~\ref{theorem:slowly varying system}.
		\item %By comparing \eqref{eq:slowly varying system, equilibrium} to Definition~\ref{definition:equilibrium x_ast}, one can say that  
		Analogously to Definition~\ref{definition:equilibrium x_ast}, one may define
		$\thetav^{\ast}(t) := \psiv\big(\vv(t)\big)$
		as an \emph{equilibrium point} of system~\eqref{eq:slowly varying system} \emph{at instant $t$}.
		This analogy makes \eqref{eq:slowly varying chi_optimal} in Lemma~\ref{lemma:x to x_ast, v3} be an instance of 
		\begin{equation*}
			%		\dot{\thetav}^{\ast} = \Psim(\vv) \dot{\vv} 
			%		~\Rightarrow~
			\big\|\dot{\thetav}^{\ast}\big\|_{2} = \big\|\Psim(\vv) \dot{\vv}\big\|_{2} \leqslant \big\|\Psim(\vv)\big\|_{2} \big\|\dot{\vv}\big\|_{2} \leqslant D_{\scriptscriptstyle\Psi} D_{\textup{v}},
		\end{equation*}
		i.e., $D_{\chi}$ in Lemma~\ref{lemma:x to x_ast, v3} corresponds to $D_{\scriptscriptstyle\Psi} D_{\textup{v}}$ in Theorem~\ref{theorem:slowly varying system}.
		\item Given $b_{i}$ in~\eqref{eq:x to x_ast, v3, Lyapunov inequalities} of Lemma~\ref{lemma:x to x_ast, v3},
		inequality~\eqref{eq:slowly varying systems, inequality Dv} in Theorem~\ref{theorem:slowly varying system} becomes
%		\begin{equation}
%			D_{\textup{v}} 
%			< \frac{b_{1} b_{3}}{b_{2}} \frac{\delta}{0\delta + b_{4} D_{\scriptscriptstyle\Psi}}
%			~\Leftrightarrow~
%			D_{\textup{v}} D_{\scriptscriptstyle\Psi}
%			< \frac{b_{1} b_{3}}{b_{2}} \frac{\delta}{b_{4}} 
%			= \frac{0.5 (1-M_{\textup{sup}}^{2})}{1} \frac{\delta}{K_{\textup{x}}^{-1}}
%			= K_{\textup{x}} 0.5 (1-M_{\textup{sup}}^{2}) \delta 
%			= \kappa \delta
%		\end{equation}
		\begin{equation*}
			D_{\textup{v}} D_{\scriptscriptstyle\Psi}
			< \frac{1}{2} (1-M_{\textup{sup}}^{2}) K_{\textup{x}} \delta 
			= \kappa \delta,
		\end{equation*}
		which is covered by \eqref{eq: x(t) is sufficiently close to chi_optimal(t)} in Lemma~\ref{lemma:x to x_ast, v3}.
		\item Given $b_{1}$ and $b_{2}$ in~\eqref{eq:x to x_ast, v3, Lyapunov inequalities} of Lemma~\ref{lemma:x to x_ast, v3}, the condition \eqref{eq:slowly varying systems, initial condition} in Theorem~\ref{theorem:slowly varying system} becomes \eqref{eq: x(0) is sufficiently close to chi_optimal(0)} in Lemma~\ref{lemma:x to x_ast, v3}.
		\item Given $b_{i}$ in~\eqref{eq:x to x_ast, v3, Lyapunov inequalities} of Lemma~\ref{lemma:x to x_ast, v3}, the constant $\kappa$~\eqref{eq:slowly varying systems, kappa} in Theorem~\ref{theorem:slowly varying system} becomes \eqref{eq:definition of kappa} in Lemma~\ref{lemma:x to x_ast, v3}.
	\end{enumerate}
\end{remark}

%\newpage
%\input{supplementary/emulated_rt_computer} % only the algorithm
%\input{supplementary/supplementary_figures}

%\backmatter

\clearpage
\bibliographystyle{naturemag} %\bibliographystyle{amsplain}
\bibliography{ref_hu}
% if required, the content of .bbl file can be included here once bbl is generated
%\input supplementary.bbl